\documentclass{article}

\usepackage[utf8]{inputenc}

\usepackage{wrapfig}

\usepackage{booktabs,caption}
\usepackage[flushleft]{threeparttable}

\usepackage[nointegrals]{wasysym}
\usepackage{enumitem}
\usepackage{adjustbox}

\usepackage{geometry}
\usepackage[english]{babel}

\usepackage{changepage}
\usepackage{epsfig}

\usepackage{multicol}
\usepackage{tikz}
\usepackage{rotating}

\usetikzlibrary{positioning}
\usetikzlibrary{decorations.markings}
\usetikzlibrary{lindenmayersystems}
\usetikzlibrary{patterns}
\usetikzlibrary[shadings]
\usetikzlibrary{fadings}

\tikzfading[
  name=fade out,
  inner color=transparent!0,
  outer color=transparent!100
]

\usepackage{mathtools}
\usepackage{amssymb}
\usepackage{array}
\usepackage{cancel}

\usepackage{xcolor}
\usepackage{graphicx}

\usepackage[numbers]{natbib}
\usepackage{hyperref}
\hypersetup{
    colorlinks=true,
    linkcolor=blue,
    citecolor=blue,
    urlcolor=red!75!black
}

\usepackage{url}

\definecolor{lejla}{RGB}{120,0,133}

\title{\large \textbf{Reconfiguration of Hamiltonian Cycles in Rectangular Grid Graphs}}

\author{Albi Kazazi\\
\small Department of Mathematics and Statistics, York University, Toronto, Canada}

\date{\today}

\begin{document}
\maketitle
\vspace*{-1.0\baselineskip}

\begingroup
\small
\noindent \textbf{ORCID:} \texttt{https://orcid.org/0009-0003-0120-4435}
\hfill
\textbf{Email:} \texttt{ak87@yorku.ca} \quad \textbf{Phone:} +1\,(416)\,878--7650

\medskip

\noindent
\begin{minipage}[t]{0.49\textwidth}\raggedright
\textbf{Mailing Address}\\
Albi Kazazi\\
24 Delfire Street\\
Maple, Ontario, Canada\\
L6A 2L9
\end{minipage}\hfill
\begin{minipage}[t]{0.49\textwidth}\raggedleft
\textbf{Affiliation Address}\\
Department of Mathematics and Statistics\\
York University\\
4700 Keele Street\\
Toronto, ON M3J 1P3, Canada
\end{minipage}
\par
\endgroup

\vspace*{-0.6\baselineskip}

\begin{abstract}
\noindent An \textit{\(m \times n\) grid graph} is the induced subgraph of the square lattice whose vertex set 
consists of all integer grid points \(\{(i,j) : 0 \leq i < m,\ 0 \leq j < n\}\). 
Let $H$ and $K$ be Hamiltonian cycles in an $m \times n$ grid graph $G$. 
We study the problem of reconfiguring $H$ into $K$, \textcolor{blue}{\textbullet} where the Hamiltonian cycles are viewed as vertices of a reconfiguration graph \textcolor{blue}{\textbullet}, using a sequence of local transformations 
called \textit{moves}. A \textit{box} of $G$ is a unit square face. A box with vertices $a, b, c, d$ 
is \textit{switchable} in $H$ if exactly two of its edges belong to $H$, and these edges are parallel. 
Given such a box with edges $ab$ and $cd$ in $H$, a \textit{switch move} removes $ab$ and $cd$, 
and adds $bc$ and $ad$. A \textit{double-switch move} consists of performing two consecutive 
switch moves. If, after a double-switch move, we obtain a Hamiltonian cycle, we say that the 
double-switch move is \textit{valid}. 

We prove that any Hamiltonian cycle $H$ can be transformed into any other Hamiltonian cycle $K$ 
via a sequence of valid double-switch moves, such that every intermediate graph remains a Hamiltonian cycle. \textcolor{blue}{\textbullet}Moreover, assuming $n \geq m$, the number of required moves is bounded by $mn^2$.\textcolor{blue}{\textbullet}
\end{abstract}

\vspace*{-0.6\baselineskip}
\section*{Acknowledgments}
\begingroup
\small
\noindent This paper is based on part of the author’s PhD dissertation at York University, 
written under the supervision of Professor Neal Madras.

\medskip

\noindent The author is grateful to his advisor, Neal Madras, for his patience and guidance through several rewrites, and for his thoughtful and insightful feedback on each of them. The author also thanks Nathan Clisby for creating an interactive online tool for visualizing and experimenting with Hamiltonian paths and cycles\cite{clisbyHamiltonian}.

\medskip

\noindent This research was supported in part by a Discovery Grant from NSERC Canada to the author’s advisor, Neal Madras.
\par
\endgroup

\section*{Introduction}

An \textit{\(m \times n\) grid graph} is the induced subgraph of the square lattice whose vertex set consists of all integer grid points 
\(\{(i,j) : 0 \leq i < m,\ 0 \leq j < n\}\) with edges between vertices at distance 1. We call the unit-square faces of the square lattice  \textit{boxes}. A Hamiltonian path (cycle) of a graph $G$ is a path (cycle) that visits each vertex of the graph exactly once.

\begingroup
\setlength{\intextsep}{0pt}
\setlength{\columnsep}{10pt}
\begin{wrapfigure}[]{r}{0cm}
\setlength{\intextsep}{0pt}
\setlength{\columnsep}{20pt}
\begin{adjustbox}{trim=0cm 0cm 0cm 0cm}
\begin{tikzpicture}[scale=1]

\begin{scope}[xshift=0cm]
   \draw[gray,very thin, step=0.5cm, opacity=0.4] (0,0) grid (3.5,2.5); 

\draw[gray,very thin, step=0.5cm, opacity=0.4] (0,0) grid (2.5,2); 

\draw[blue, line width=0.5mm] (0,2.5)--++(0,-2.5)--++(0.5,0)--++(0,0.5)--++(0.5,0)--++(0,-0.5)--++(0.5,0)--++(0,1)--++(0.5,0)--++(0,-1)--++(0.5,0)--++(0,1)--++(0.5,0)--++(0,-1)--++(0.5,0)--++(0,1.5)--++(-2,0)--++(0,0.5)--++(2,0)--++(0,0.5)--++(-2.5,0)--++(0,-1.5)--++(-0.5,0)--++(0,1.5)--++(-0.5,0); 


\node[below] at (1.75, 0) [scale=0.8]{\small{\begin{tabular}{c} Fig. I.1. A 1-complex Hamiltonian \\ cycle on an $8 \times 6$ grid graph. \end{tabular}}};;
\end{scope}

\end{tikzpicture}
\end{adjustbox}
\end{wrapfigure}

The question of whether an $m \times n$ grid graph has a Hamiltonian path was first studied by Itai et al. in \cite{itai1982hamilton}. They showed that for an $m \times n$ grid graph to have a Hamiltonian cycle, it is necessary and sufficient that at least one of $m$ and $n$ is even. Chen et al. gave an efficient algorithm to construct Hamiltonian paths in rectangular grid graphs \cite{chen2002efficient}. A \textit{solid} grid graph is a grid graph without holes, i.e. each bounded face of the graph is a box. Umans and Lenhart \cite{umans1997hamiltonian} gave a polynomial-time algorithm to find a Hamiltonian cycle in solid grid graphs, if one exists.

Given any two Hamiltonian cycles $H$ and $K$, the reconfiguration problem asks whether it is possible to transform $H$ into $K$ step-by-step, so that each intermediate step is also a Hamiltonian cycle of $G$. Nishat and Whitesides \cite{nishat2017bend} introduced the ``flip'' and ``transpose'' moves described  below, and a complexity measure called ``bend complexity" for Hamiltonian cycles in rectangular grid graphs. Roughly, a 1-complex Hamiltonian cycle is one in which every vertex of $G$ is connected to the boundary via a straight line.  They prove that using these two moves, it is possible to reconfigure any pair of 1-complex Hamiltonian cycles in $G$ into one another. \textcolor{blue}{\textbullet} Equivalently, the reconfiguration graph of 1-complex Hamiltonian cycles in rectangular grid graphs is connected.\textcolor{blue}{\textbullet}

\noindent We dispense with the need for bend complexity constraints, proving the connectivity of reconfiguration graph of all Hamiltonian cycles in rectangular grid graphs, by using a more general move, which we call a double-switch move.

\endgroup

\begingroup
\setlength{\intextsep}{0pt}
\setlength{\columnsep}{20pt}
\begin{wrapfigure}[]{l}{0cm}
\setlength{\intextsep}{0pt}
\setlength{\columnsep}{20pt}
\begin{adjustbox}{trim=0cm 0cm 0cm 0cm}
\begin{tikzpicture}[scale=1]

\begin{scope}[xshift=0cm]{
\draw[gray,very thin, step=0.5cm, opacity=0.5] (0,0) grid (1.5,2);

\draw[blue, line width=0.5mm] (0,0)--++(1.5,0)--++(0,2)--++(-1.5,0)--++(0,-1)--++(0.5,0)--++(0,0.5)--++(0.5,0)--++(0,-1)--++(-1,0)--++(0,-0.5);

\node at (0.25,0.75) [scale=1] {\small{Y}};
\node at (0.75,1.75) [scale=1] {\small{X}};

\draw [->,black, very thick] (1.25,2.1) to [out=10,in=170] (5.25,2.1);
\node[above] at  (3.25,2.3) [scale=0.8]{\small{$X \mapsto Y$}};

\node[above] at  (2,1) [scale=0.8]{\small{$\textrm{Sw}(X)$}};
\draw [->,black, very thick] (1.6,1) -- (2.4,1);

} \end{scope}

\begin{scope}[xshift=2.5cm]{
\draw[gray,very thin, step=0.5cm, opacity=0.5] (0,0) grid (1.5,2);

\draw[blue, line width=0.5mm] (0,0)--++(1.5,0)--++(0,2)--++(-0.5,0)--++(0,-1.5)--++(-1,0)--++(0,-0.5);
\draw[blue, line width=0.5mm] (0,1)--++(0.5,0)--++(0,1)--++(-0.5,0)--++(0,-1);

\node at (0.25,0.75) [scale=1] {\small{Y}};
\node at (0.75,1.75) [scale=1] {\small{X}};

\node[above] at  (2,1) [scale=0.8]{\small{$\textrm{Sw}(Y)$}};
\draw [->,black, very thick] (1.6,1) -- (2.4,1);

\node[below] at (0.75,-0.15) [scale=0.8]{\small{\begin{tabular}{c} Fig. I.2. Illustration of each switch of a general \\ double-switch move in a $4 \times 5$ grid graph. \end{tabular}}};;

} \end{scope}

\begin{scope}[xshift=5cm]{
\draw[gray,very thin, step=0.5cm, opacity=0.5] (0,0) grid (1.5,2);

\draw[blue, line width=0.5mm] (0,0)--++(1.5,0)--++(0,2)--++(-0.5,0)--++(0,-1.5)--++(-0.5,0)--++(0,1.5)--++(-0.5,0)--++(0,-2);

\node at (0.25,0.75) [scale=1] {\small{Y}};
\node at (0.75,1.75) [scale=1] {\small{X}};

} \end{scope}

\end{tikzpicture}
\end{adjustbox}
\end{wrapfigure}

\noindent Let $H$ be a Hamiltonian cycle of an $m\times n$ grid graph $G$. A box in $G$ with vertices $a, b, c, d$ is considered \textit{switchable} in $H$ if it has exactly two edges in $H$, and these edges are parallel. Let $abcd$ be a switchable box with edges $ab$ and $cd$ in $H$. We define a \textit{switch move} on the box $abcd$ in $H$ as follows: remove edges $ab$ and $cd$ and add edges $bc$ and $ad$. If $X$ is a switchable box of in $H$, we denote a switch move on $X$ by $\textrm{Sw}(X)$.

A \textit{double-switch move} (or simply a \textit{move}) is a pair of switch operations where we first switch $X$ and then find another switchable box $Y$ and switch it, and denote it by $X \mapsto Y$. If, after a double-switch move, we obtain a new Hamiltonian cycle, we call the move a \textit{valid move}.

\endgroup 

\begingroup
\setlength{\intextsep}{-10pt}
\setlength{\columnsep}{10pt}
\begin{wrapfigure}[]{r}{0cm}
\setlength{\intextsep}{0pt}
\setlength{\columnsep}{20pt}
\begin{adjustbox}{trim=0cm 0cm 0cm 0cm}
\begin{tikzpicture}[scale=1.5]

\begin{scope}[xshift=0cm]{
\draw[gray,very thin, step=0.5cm, opacity=0.5] (0,0) grid (0.5,1);

\draw[blue, line width=0.5mm] (0,1)--++(0,-0.5)--++(0.5,0)--++(0,0.5);
\draw[blue, line width=0.5mm] (0,0)--++(0.5,0);

\draw[fill=blue, opacity=1] (0,0) circle [radius=0.035];
\draw[fill=blue, opacity=1] (0,0.5) circle [radius=0.035];
\draw[fill=blue, opacity=1] (0,1) circle [radius=0.035];

\draw[fill=blue, opacity=1] (0.5,0) circle [radius=0.035];
\draw[fill=blue, opacity=1] (0.5,0.5) circle [radius=0.035];
\draw[fill=blue, opacity=1] (0.5,1) circle [radius=0.035];

\node[left] at (0,0) [scale=0.8] {\small{$a$}};
\node[left] at (0,0.5) [scale=0.8] {\small{$d$}};
\node[left] at (0,1) [scale=0.8] {\small{$f$}};

\node[right] at (0.5,0) [scale=0.8] {\small{$b$}};
\node[right] at (0.5,0.5) [scale=0.8] {\small{$c$}};
\node[right] at (0.5,1) [scale=0.8] {\small{$e$}};

\node at (0.25,0.25) [scale=1] {\small{X}};
\node at (0.25,0.75) [scale=1] {\small{Y}};

\node[above] at  (0.8,0.75) [scale=0.8]{\small{flip}};
\draw [->,black, very thick] (0.6,0.75) -- (1.1,0.75);

\node[below] at (0.85,0) [scale=0.8]{\small{\begin{tabular}{c} Fig. I.3. A flip move. \end{tabular}}};;

} \end{scope}

\begin{scope}[xshift=1.2cm]{
\draw[gray,very thin, step=0.5cm, opacity=0.5] (0,0) grid (0.5,1);

\draw[blue, line width=0.5mm] (0,0)--++(0,0.5)--++(0.5,0)--++(0,-0.5);
\draw[blue, line width=0.5mm] (0,1)--++(0.5,0);

\draw[fill=blue, opacity=1] (0,0) circle [radius=0.035];
\draw[fill=blue, opacity=1] (0,0.5) circle [radius=0.035];
\draw[fill=blue, opacity=1] (0,1) circle [radius=0.035];

\draw[fill=blue, opacity=1] (0.5,0) circle [radius=0.035];
\draw[fill=blue, opacity=1] (0.5,0.5) circle [radius=0.035];
\draw[fill=blue, opacity=1] (0.5,1) circle [radius=0.035];

\node[left] at (0,0) [scale=0.8] {\small{$a$}};
\node[left] at (0,0.5) [scale=0.8] {\small{$d$}};
\node[left] at (0,1) [scale=0.8] {\small{$f$}};

\node[right] at (0.5,0) [scale=0.8] {\small{$b$}};
\node[right] at (0.5,0.5) [scale=0.8] {\small{$c$}};
\node[right] at (0.5,1) [scale=0.8] {\small{$e$}};

\node at (0.25,0.25) [scale=1] {\small{X}};
\node at (0.25,0.75) [scale=1] {\small{Y}};

} \end{scope}

\end{tikzpicture}
\end{adjustbox}
\end{wrapfigure}

Let $X=abcd$ and $Y=dcef$ be boxes sharing the edge $cd$ of $G$. Assume that the edges $ab, fd, dc$ and $ce$ belong to $H$, and that the edges $fe,ad$ and $bc$ do not. A \textit{flip} move consists in removing the edges $fd, ce$ and $ab$, and adding the edges $ad, bc$ and $fe$. Effectively, this is the same as first switching $X$, and then switching $Y$. See Figure I.3.

\endgroup

\begingroup
\setlength{\intextsep}{0pt}
\setlength{\columnsep}{20pt}
\begin{wrapfigure}[]{l}{0cm}
\setlength{\intextsep}{0pt}
\setlength{\columnsep}{20pt}
\begin{adjustbox}{trim=0cm 0cm 0cm 0cm}
\begin{tikzpicture}[scale=1.5]

\begin{scope}[xshift=0cm]{
\draw[gray,very thin, step=0.5cm, opacity=0.5] (0,0) grid (1,1);

\draw[blue, line width=0.5mm] (0,1)--++(1,0)--++(0,-0.5)--++(-1,0);
\draw[blue, line width=0.5mm] (0.5,0)--++(0.5,0);

\draw[fill=blue, opacity=1] (0,0) circle [radius=0.035];
\draw[fill=blue, opacity=1] (0,0.5) circle [radius=0.035];
\draw[fill=blue, opacity=1] (0,1) circle [radius=0.035];

\draw[fill=blue, opacity=1] (0.5,0) circle [radius=0.035];
\draw[fill=blue, opacity=1] (0.5,0.5) circle [radius=0.035];
\draw[fill=blue, opacity=1] (0.5,1) circle [radius=0.035];

\draw[fill=blue, opacity=1] (1,0) circle [radius=0.035];
\draw[fill=blue, opacity=1] (1,0.5) circle [radius=0.035];
\draw[fill=blue, opacity=1] (1,1) circle [radius=0.035];

\node[left] at (0,1) [scale=0.8] {\small{$a$}};
\node[left] at (0,0.5) [scale=0.8] {\small{$d$}};
\node[left] at (0,0) [scale=0.8] {\small{$i$}};

\node[above] at (0.5,1) [scale=0.8] {\small{$b$}};
\node[right] at (0.5,0.6) [scale=0.8] {\small{$c$}};
\node[below] at (0.5,0) [scale=0.8] {\small{$h$}};

\node[above] at (1,1) [scale=0.8] {\small{$e$}};
\node[right] at (1,0.5) [scale=0.8] {\small{$f$}};
\node[right] at (1,0) [scale=0.8] {\small{$g$}};

\node at (0.25,0.25) [scale=1] {\small{W}};
\node at (0.25,0.75) [scale=1] {\small{X}};
\node at (0.75,0.75) [scale=1] {\small{Y}};
\node at (0.75,0.25) [scale=1] {\small{Z}};

\node[above] at  (1.5,0.75) [scale=0.8]{\small{transpose}};
\draw [->,black, very thick] (1.1,0.75) -- (1.9,0.75);

} \end{scope}

\node[below] at (1.5,-0.15) [scale=0.8]{\small{\begin{tabular}{c} Fig. I.4. A transpose move. \end{tabular}}};;

\begin{scope}[xshift=2cm]{
\draw[gray,very thin, step=0.5cm, opacity=0.5] (0,0) grid (1,1);

\draw[blue, line width=0.5mm] (0.5,0)--++(0,1)--++(0.5,0)--++(0,-1);
\draw[blue, line width=0.5mm] (0,0.5)--++(0,0.5);

\draw[fill=blue, opacity=1] (0,0) circle [radius=0.035];
\draw[fill=blue, opacity=1] (0,0.5) circle [radius=0.035];
\draw[fill=blue, opacity=1] (0,1) circle [radius=0.035];

\draw[fill=blue, opacity=1] (0.5,0) circle [radius=0.035];
\draw[fill=blue, opacity=1] (0.5,0.5) circle [radius=0.035];
\draw[fill=blue, opacity=1] (0.5,1) circle [radius=0.035];

\draw[fill=blue, opacity=1] (1,0) circle [radius=0.035];
\draw[fill=blue, opacity=1] (1,0.5) circle [radius=0.035];
\draw[fill=blue, opacity=1] (1,1) circle [radius=0.035];

\node[left] at (0,1) [scale=0.8] {\small{$a$}};
\node[left] at (0,0.5) [scale=0.8] {\small{$d$}};
\node[left] at (0,0) [scale=0.8] {\small{$i$}};

\node[above] at (0.5,1) [scale=0.8] {\small{$b$}};
\node[right] at (0.5,0.6) [scale=0.8] {\small{$c$}};
\node[below] at (0.5,0) [scale=0.8] {\small{$h$}};

\node[above] at (1,1) [scale=0.8] {\small{$e$}};
\node[right] at (1,0.5) [scale=0.8] {\small{$f$}};
\node[right] at (1,0) [scale=0.8] {\small{$g$}};

\node at (0.25,0.25) [scale=1] {\small{W}};
\node at (0.25,0.75) [scale=1] {\small{X}};
\node at (0.75,0.75) [scale=1] {\small{Y}};
\node at (0.75,0.25) [scale=1] {\small{Z}};

} \end{scope}

\end{tikzpicture}
\end{adjustbox}
\end{wrapfigure}

Consider the four boxes $X=abcd$, $Y=cbef$, $Z=cfgh$ and $W=dchi$ that are incident on the vertex $c$. Note that $X$ and $Y$ share the edge $cb$, $Y$ and $Z$ share $cf$, $Z$ and $W$ share $ch$, and $W$ and $X$ share $cd$. Assume that the edges $ab, be, ef, fc, cd$ and $hg$ belong to $H$ and that the edges $ad, bc, ch$ and $fg$ do not. A \textit{transpose} move consists in switching $X$ and then switching $Z$. See Figure I.4.

Nishat in \cite{nishat2020reconfiguration} showed that flip and transpose moves are always valid. The more general double-switch moves are sufficient for constructing algorithms that reconfigure arbitrary Hamiltonian cycles in grid graphs. This comes at the added cost of verifying the validity of each move. We provide such reconfiguration algorithms and prove the existence of all required moves.

\null


\noindent \textbf{Theorem.} Let $H$ and $K$ be any two Hamiltonian cycles in an $m \times n$ grid graph $G$ with $n \geq m$. Then there exists a sequence of at most $n^2m$ valid double-switch moves that reconfigures $H$ into $K$. 

\null 

\noindent See \cite{video} for an illustration. \textcolor{blue}{\textbullet} In particular, this yields an explicit $mn^2$ upper bound on the diameter of the reconfiguration graph of Hamiltonian cycles in $m \times n$ grid graphs\textcolor{blue}{\textbullet}. An analogous result for Hamiltonian paths is treated in \cite{kazazi_dissertation}. The extension makes use of two additional moves beyond the double-switch: the switch move and the \textit{backbite} move, the latter used to relocate the endpoints of a path and originally introduced by Mansfield \cite{mansfield1982monte}. For a description of the backbite move, see Appendix~A.2.

\begingroup
\setlength{\intextsep}{0pt}
\setlength{\columnsep}{20pt}
\begin{wrapfigure}[]{l}{0cm}
\begin{adjustbox}{trim=0cm 0cm 0cm 0cm}
\begin{tikzpicture}[scale=0.7]

\begin{scope}
 
{

\draw[gray,very thin, step=0.5cm, opacity=0.4] (0,0) grid (4.5,3.5); 

\filldraw[green!50!white, opacity=0.5] (2,0)--++(-0.5,0)--++(0,0.5)--++(-0.5,0)--++(0,0.5)--++(-0.5,0)--++(0,0.5)--++(-0.5,0)--++(0,0.5)--++(0.5,0)--++(0,0.5)--++(0.5,0)--++(0,0.5)--++(0.5,0)--++(0,0.5)--++(1,0)--++(0.5,0)--++(0,-0.5)--++(0.5,0)--++(0,-0.5)--++(0.5,0)--++(0,-0.5)--++(0.5,0)--++(0,-0.5)--++(-0.5,0)--++(0,-0.5)--++(-0.5,0)--++(0,-0.5)--++(-0.5,0)--++(0,-0.5)--++(-0.5,0);


\draw[blue, line width=0.5mm] 
(2,0)--++(-0.5,0)--++(0,0.5)--++(-0.5,0)--++(0,0.5)--++(-0.5,0)--++(0,0.5)--++(-0.5,0)--++(0,0.5)--++(0.5,0)--++(0,0.5)--++(0.5,0)--++(0,0.5)--++(0.5,0)--++(0,0.5)--++(0.5,0); 

\draw[blue, line width=0.5mm] 
(2,3.5)--++(0,-1)--++(-0.5,0)--++(0,-0.5)--++(-0.5,0)--++(0,-0.5)--++(0.5,0)--++(0,-0.5)--++(0.5,0)--++(0,1)--++(0.5,0)--++(0,-1)--++(0.5,0)--++(0,0.5)--++(0.5,0)--++(0,0.5)--++(-0.5,0)--++(0,0.5)--++(-0.5,0)--++(0,1);

\draw[blue, line width=0.5mm] 
(2.5,3.5)--++(0.5,0)--++(0,-0.5)--++(0.5,0)--++(0,-0.5)--++(0.5,0)--++(0,-0.5)--++(0.5,0)--++(0,-0.5)--++(-0.5,0)--++(0,-0.5)--++(-0.5,0)--++(0,-0.5)--++(-0.5,0)--++(0,-0.5)--++(-0.5,0);

\draw[blue, line width=0.5mm] 
(2.5,0)--++(0,0.5)--++(-0.5,0)--++(0,-0.5);




}   

\node[below] at (5,0) [scale=0.8]{\small{\begin{tabular}{c} Figure I.5. A solid grid graph shaded green, with two distinct \\ Hamiltonian cycles frozen under the double-switch move. \end{tabular}}};;

\end{scope}

\begin{scope}[xshift=5.5cm]

{
\draw[gray,very thin, step=0.5cm, opacity=0.4] (0,0) grid (4.5,3.5); 

\filldraw[green!50!white, opacity=0.5] (2,0)--++(-0.5,0)--++(0,0.5)--++(-0.5,0)--++(0,0.5)--++(-0.5,0)--++(0,0.5)--++(-0.5,0)--++(0,0.5)--++(0.5,0)--++(0,0.5)--++(0.5,0)--++(0,0.5)--++(0.5,0)--++(0,0.5)--++(1,0)--++(0.5,0)--++(0,-0.5)--++(0.5,0)--++(0,-0.5)--++(0.5,0)--++(0,-0.5)--++(0.5,0)--++(0,-0.5)--++(-0.5,0)--++(0,-0.5)--++(-0.5,0)--++(0,-0.5)--++(-0.5,0)--++(0,-0.5)--++(-0.5,0);

\draw[blue, line width=0.5mm] 
(2,0)--++(-0.5,0)--++(0,0.5)--++(-0.5,0)--++(0,0.5)--++(-0.5,0)--++(0,0.5)--++(-0.5,0)--++(0,0.5)--++(0.5,0)--++(0,0.5)--++(0.5,0)--++(0,0.5)--++(0.5,0)--++(0,0.5)--++(0.5,0); 

\draw[blue, line width=0.5mm] 
(2,3.5)--++(0,-0.5)--++(0.5,0)--++(0,0.5);

\draw[blue, line width=0.5mm] 
(2.5,3.5)--++(0.5,0)--++(0,-0.5)--++(0.5,0)--++(0,-0.5)--++(0.5,0)--++(0,-0.5)--++(0.5,0)--++(0,-0.5)--++(-0.5,0)--++(0,-0.5)--++(-0.5,0)--++(0,-0.5)--++(-0.5,0)--++(0,-0.5)--++(-0.5,0);

\draw[blue, line width=0.5mm] 
(2.5,0)--++(0,1)--++(0.5,0)--++(0,0.5)--++(0.5,0)--++(0,0.5)--++(-0.5,0)--++(0,0.5)--++(-0.5,0)--++(0,-1)--++(-0.5,0)--++(0,1)--++(-0.5,0)--++(0,-0.5)--++(-0.5,0)--++(0,-0.5)--++(0.5,0)--++(0,-0.5)--++(0.5,0)--++(0,-1);

}   

\end{scope}

\end{tikzpicture}

\end{adjustbox}
\end{wrapfigure}

\noindent \textbf{Scope and extensions.} In Section 4 we will use the condition that the boundary is rectangular to prove the theorem. This assumption cannot be relaxed to include all general solid grid graphs. For example, neither Hamiltonian cycle in Figure I.5. admits a valid double-switch move, so it is not possible to reconfigure one into the other through double-switch moves. While we believe there may be finer classes of graphs between rectangular grid graphs and solid grid graphs 
that can be reconfigured by the double-switch move, it seems likely that such classes would require imposing boundary conditions on general solid grid graphs.

We conjecture that the double-switch move should also suffice to reconfigure Hamiltonian cycles in three-dimensional rectangular grid graphs, although we do not yet have a proof. By a 
three-dimensional rectangular grid graph we mean the induced subgraph of the cubic lattice whose 
vertex set is \(\{(i,j,k) : 0 \leq i < m,\ 0 \leq j < n,\ 0 \leq k < p\}\) with edges between vertices at distance one. The arguments in this paper rely on Jordan’s curve theorem, which has no direct analogue in three dimensions. Thus, if this conjecture is true, the proof would seem to require different techniques.

\null 

\noindent \textbf{Applications and related work.}  A self-avoiding walk is a walk in a lattice where every vertex is unique. A Hamiltonian path in a grid graph is an example of a self-avoiding walk. Madras and Slade in \cite{madras2013self} present a comprehensive and rigorous study of self-avoiding walks. One application of the theorem is in chemical physics, drawing from the theory of self-avoiding walks. Researchers in \cite{oberdorf2006secondary}, \cite{jacobsen2008unbiased}, \cite{deutsch1997long}, and \cite{mansfield1982monte} use Monte Carlo methods to study statistical properties of polymer chains, which they abstracted as cycles and paths in the cubic lattice. 

\begingroup
\setlength{\intextsep}{0pt}
\setlength{\columnsep}{20pt}
\begin{wrapfigure}[]{r}{0cm}
\begin{adjustbox}{trim=0cm 0cm 0cm 0cm}
\begin{tikzpicture}[scale=0.7]

\begin{scope}[xshift=0cm, yshift=0cm]

\draw[gray,very thin, step=0.5cm, opacity=0.5] (1.49,1.49) grid (6,6); 

\draw[gray,very thin, step=0.5cm, opacity=0.5] (7.49,1.49) grid (12,6);


\begin{scope}[xshift=6cm, yshift=0cm]

{
\foreach \x in {1,3,5,6,7,9} \draw[blue, line width=0.5mm] (1+0.5*\x,1.5)--(1+0.5*\x+0.5,1.5); 
\foreach \x in {2,5,6} \draw[blue, line width=0.5mm] (1+0.5*\x,2)--(1+0.5*\x+0.5,2);
\foreach \x in {1,2,4,5,6,8} \draw[blue, line width=0.5mm] (1+0.5*\x,2.5)--(1+0.5*\x+0.5,2.5);
\foreach \x in {1,3,4,5,6,8,9} \draw[blue, line width=0.5mm] (1+0.5*\x,3)--(1+0.5*\x+0.5,3);
\foreach \x in {3,4,6,9} \draw[blue, line width=0.5mm] (1+0.5*\x,3.5)--(1+0.5*\x+0.5,3.5);
\foreach \x in {2,4,6} \draw[blue, line width=0.5mm] (1+0.5*\x,4)--(1+0.5*\x+0.5,4);
\foreach \x in {1,2,4,6,8} \draw[blue, line width=0.5mm] (1+0.5*\x,4.5)--(1+0.5*\x+0.5,4.5);
\foreach \x in {1,4,7,8} \draw[blue, line width=0.5mm]  (1+0.5*\x,5)--(1+0.5*\x+0.5,5);
\foreach \x in {2,4,5,8} \draw[blue, line width=0.5mm] (1+0.5*\x,5.5)--(1+0.5*\x+0.5,5.5);
\foreach \x in {1,2,3,4,5,6,8,9} \draw[blue, line width=0.5mm] (1+0.5*\x,6)--(1+0.5*\x+0.5,6);
}
{
\foreach \x in {1,2,4,5,6,8,9} \draw[blue, line width=0.5mm] 
(1.5, 1+0.5*\x)--(1.5, 1+0.5*\x+0.5); 
\foreach \x in {1,4,5,8} \draw[blue, line width=0.5mm] 
(2, 1+0.5*\x)--(2, 1+0.5*\x+0.5); 
\foreach \x in {1,3,5,7,8} \draw[blue, line width=0.5mm] 
(2.5, 1+0.5*\x)--(2.5, 1+0.5*\x+0.5);
\foreach \x in {1,2,6,8} \draw[blue, line width=0.5mm] 
(3, 1+0.5*\x)--(3, 1+0.5*\x+0.5);
\foreach \x in {1,5,7} \draw[blue, line width=0.5mm] 
(3.5, 1+0.5*\x)--(3.5, 1+0.5*\x+0.5);
\foreach \x in {5,7,8} \draw[blue, line width=0.5mm] 
(4, 1+0.5*\x)--(4, 1+0.5*\x+0.5);
\foreach \x in {2,4,6,8,9} \draw[blue, line width=0.5mm] 
(4.5, 1+0.5*\x)--(4.5, 1+0.5*\x+0.5);
\foreach \x in {1,2,4,5,6,9} \draw[blue, line width=0.5mm] 
(5, 1+0.5*\x)--(5, 1+0.5*\x+0.5);
\foreach \x in {1,2,5,6,8} \draw[blue, line width=0.5mm] 
(5.5, 1+0.5*\x)--(5.5, 1+0.5*\x+0.5);
\foreach \x in {1,2,3,5,6,7,8,9} \draw[blue, line width=0.5mm] 
(6, 1+0.5*\x)--(6, 1+0.5*\x+0.5);
}

\end{scope}


{
\foreach \x in {1,2,3,5,6,8,9} \draw[blue, line width=0.5mm] (1+0.5*\x,1.5)--(1+0.5*\x+0.5,1.5); 
\foreach \x in {2,6,8} \draw[blue, line width=0.5mm] (1+0.5*\x,2)--(1+0.5*\x+0.5,2);
\foreach \x in {1,4,6,7} \draw[blue, line width=0.5mm] (1+0.5*\x,2.5)--(1+0.5*\x+0.5,2.5);
\foreach \x in {1,4,6,7} \draw[blue, line width=0.5mm] (1+0.5*\x,3)--(1+0.5*\x+0.5,3);
\foreach \x in {3,6,7,8} \draw[blue, line width=0.5mm] (1+0.5*\x,3.5)--(1+0.5*\x+0.5,3.5);
\foreach \x in {2,3,5,6,7,9} \draw[blue, line width=0.5mm] (1+0.5*\x,4)--(1+0.5*\x+0.5,4);
\foreach \x in {2,5,6,7,9} \draw[blue, line width=0.5mm] (1+0.5*\x,4.5)--(1+0.5*\x+0.5,4.5);
\foreach \x in {3,6,7,8} \draw[blue, line width=0.5mm]  (1+0.5*\x,5)--(1+0.5*\x+0.5,5);
\foreach \x in {2,3,4,6,8} \draw[blue, line width=0.5mm] (1+0.5*\x,5.5)--(1+0.5*\x+0.5,5.5);
\foreach \x in {1,2,3,4,5,6,8,9} \draw[blue, line width=0.5mm] (1+0.5*\x,6)--(1+0.5*\x+0.5,6);
}
{
\foreach \x in {1,2,4,5,6,7,8,9} \draw[blue, line width=0.5mm] 
(1.5, 1+0.5*\x)--(1.5, 1+0.5*\x+0.5); 
\foreach \x in {2,4,5,7,8} \draw[blue, line width=0.5mm] 
(2, 1+0.5*\x)--(2, 1+0.5*\x+0.5); 
\foreach \x in {2,3,4,7} \draw[blue, line width=0.5mm] 
(2.5, 1+0.5*\x)--(2.5, 1+0.5*\x+0.5);
\foreach \x in {1,2,4,6,7} \draw[blue, line width=0.5mm] 
(3, 1+0.5*\x)--(3, 1+0.5*\x+0.5);
\foreach \x in {1,2,4,5,7,8} \draw[blue, line width=0.5mm] 
(3.5, 1+0.5*\x)--(3.5, 1+0.5*\x+0.5);
\foreach \x in {2,4,8} \draw[blue, line width=0.5mm] 
(4, 1+0.5*\x)--(4, 1+0.5*\x+0.5);
\foreach \x in {1,9} \draw[blue, line width=0.5mm] 
(4.5, 1+0.5*\x)--(4.5, 1+0.5*\x+0.5);
\foreach \x in {1,3,6,9} \draw[blue, line width=0.5mm] 
(5, 1+0.5*\x)--(5, 1+0.5*\x+0.5);
\foreach \x in {2,3,4,6,8} \draw[blue, line width=0.5mm] 
(5.5, 1+0.5*\x)--(5.5, 1+0.5*\x+0.5);
\foreach \x in {1,2,3,4,5,7,8,9} \draw[blue, line width=0.5mm] 
(6, 1+0.5*\x)--(6, 1+0.5*\x+0.5);
}

\end{scope}

\node[below] at (6.75,1.25) [scale=0.8]{\small{\begin{tabular}{c} Figure I.6. Two distinct Hamiltonian cycles on $10 \times 10$ grid graphs. \end{tabular}}};;
\end{tikzpicture}
\end{adjustbox}
\end{wrapfigure}

They use self avoiding-walks to model how a flexible polymer chain is arranged in a liquid solution. A polymer chain's concentration is the fraction of vertices of the lattice that are occupied by the vertices (monomers) of the polymer. The authors consider maximally concentrated polymers (high density polymers), where all the space is occupied by the polymer. These can be naturally represented as Hamiltonian paths or cycles. They view the set of Hamiltonian cycles in a rectangular grid graph as the state space of a Markov chain, with the double-switch move being the transition mechanism. Given a Hamiltonian cycle (a state in the state space), we choose two switchable boxes at random and perform a double-switch move. If the move is valid, then the new state is the resulting Hamiltonian cycle. Otherwise, we remain at the initial state and choose another pair of switchable boxes. The idea is that after a sufficiently large number of transitions, we obtain a set with many different states, which represents a reasonable random sample of the entire state space. The validity of these methods requires uniform random sampling of the entire state space, which in turn requires the Markov chain to be irreducible (i.e., any Hamiltonian cycle can be reconfigured into any other through a sequence of valid moves). The authors in \cite{oberdorf2006secondary, jacobsen2008unbiased, deutsch1997long, mansfield1982monte} assume irreducibility, but do not prove it. For a more detailed discussion on Monte Carlo methods and reconfiguration of self-avoiding walks, see Chapter 9 in \cite{madras2013self}.

\endgroup

Nishat, Whitesides, and Srinivasan extended the result of \cite{nishat2017bend} to 1-complex Hamiltonian paths in rectangular grid graphs \cite{nishat2024hamiltonian, nishat20231, nishat2023reconfiguration}, and to 1-complex Hamiltonian cycles in L-shaped grid graphs \cite{nishat2019reconfiguring}. The authors define a 1-complex $s$,$t$ Hamiltonian path to be a 1-complex Hamiltonian path that begins and ends at diagonally opposite corners $s$ and $t$ of a rectangular grid graph. We note that the results \cite{nishat2024hamiltonian}, \cite{nishat20231}, and  \cite{nishat2023reconfiguration} are extended to arbitrary $s$,$t$ Hamiltonian paths in \cite{kazazi_dissertation}. 

\null 

\noindent The rest of the paper is organized as follows. In Section 1 we introduce notation and definitions, prove some lemmas about the structure that a Hamiltonian cycle imparts on a grid graph, and some other lemmas characterizing the validity of double-switch moves. In Section 2, we state the algorithm required for the proof of the main result, and show that it depends on the existence of a further two algorithms, the MLC and the 1LC. In Section 3 we prove the MLC and 1LC algorithms. The 1LC proof depends on a lemma whose proof takes up all of Section 4.

\section{Preliminaries}

A \textit{grid graph} is a subgraph of the integer grid $\mathbb{Z}^2$. A \textit{lattice animal} is a finite connected subgraph of $\mathbb{Z}^2$. A Hamiltonian path (cycle) of a graph $G$ is a path (cycle) that visits each vertex of the graph exactly once. Assume that $G$ has a cut vertex $v$. Then $G$ cannot have a Hamiltonian cycle. Let $G_1, G_2$ be the components of $G\setminus v$. Let $H_1$ be a Hamiltonian path of $G\setminus G_1$ and let  $H_2$  be a Hamiltonian path of $G\setminus G_2$ such that $H_1$ and $H_2$ have $v$ as an end-vertex. Then a Hamiltonian path $H$ of $G$ can be obtained by concatenating $H_1$ and $H_2$. Since $H_1$ and $H_2$ are smaller than $H$ they are easier to find and to reconfigure. It follows that a graph that cannot be decomposed in this manner must be 2-connected. Thus, from here on, we will restrict our attention to 2-connected grid graphs.


\null 

\noindent  \textbf{Definitions.} Let $G$ be a grid graph and let $H$ be a Hamiltonian cycle of $G$. 
We denote the set of boxes of a grid graph $G$ by $\textrm{Boxes}(G)$. We will need some definitions to navigate $G$ and $H$. Position 
$G$ in the first quadrant so that its westernmost vertices have x-coordinate 0 and southernmost vertices have y-coordinate 0. We use the $x$ and $y$ coordinates to describe a rectangle in the graph and denote it by $R(k_1,k_2; l_1,l_2)$. This rectangle corresponds to the Cartesian product of the intervals $(k_1,k_2)$ and $(l_1,l_2)$. We will denote a box of $G$ by $R(k,l)$ where $k$ and $l$ are the coordinates of the corner of the box that is closest to the origin. That is, $R(k,l)=R(k,k+1;l,l+1)$.

We specify a vertex $v$ by $v(k,l)$, where $k$ and $l$ are the vertex coordinates. We denote a horizontal edge $e$ by $e(k_1,k_2; l_1)$, where $k_1, k_2$ are the x-coordinates of the vertices of $e$ and $l_1$ is the y-coordinate of the vertices of $e$. Similarly, we write $e(k_1; l_1, l_2)$ for vertical edges. It will be convenient to use the notation $\{u, v\}$ to describe edges of $G$ and the notation $(u, v)$ to describe directed edges of $G$. For a directed edge $e=(u,v)$, $u$ is said to be the \textit{tail} of $e$ and $v$ is said to be the \textit{head} of $e$.  

Let $G'$ be a subgraph of $G$. Then we write $G'+(x,y)$ to denote the translation of $G'$ by $(x,y)$.

\null

\noindent \textbf{Theorem 1.1.} Jordan's Curve Theorem for polygons (JCT). A polygon $Q$ divides the set of points of the plane not on $Q$ into two disjoint subsets Int (for ``Interior'') and Ext (for ``Exterior'') that have $Q$ as a common boundary and are such that any two points within a subset can be joined by a path that does not intersect $Q$ while any path joining a point of Int to a point of Ext must intersect $Q$. $\square$

\begingroup
\setlength{\intextsep}{0pt}
\setlength{\columnsep}{20pt}
\begin{wrapfigure}[]{r}{0cm}
\begin{adjustbox}{trim=0cm 0.75cm 0cm 0.25cm}
\begin{tikzpicture}[scale=1]

\begin{scope}[xshift=0cm]{
\draw[black,very thin, step=0.5cm, opacity=0.5] (0,0) grid (2.5,2);

\draw[black, thick] (0,0)--++(2.5,0)--++(0,2)--++(-2.5,0)--++(0,-2);

\node[left] at (0.1,-0.1) [scale=1]{\tiny{0}};
\node[left] at (0,2) [scale=1]
{\tiny{n-1}};
\node[below] at (2.5,0) [scale=1]{\tiny{m-1}};

\node[below] at (1.25, -0.15) [scale=0.8]{\small{\begin{tabular}{c} Fig. 1.1.  A $6 {\times} 5$ grid graph. \end{tabular}}};;

} \end{scope}

\end{tikzpicture}
\end{adjustbox}
\end{wrapfigure}

\null 

\noindent We record here a useful consequence of Jordan's curve theorem.

\null 

\noindent \textbf{Corollary 1.2.} Let $p_1$ and $p_2$ be points on the plane not on $Q$. If the segment $[p_1,p_2]$ intersects $Q$ exactly once at a point $q$, where $q$ is not a vertex of $Q$, then one of $p_1$ and $p_2$ is on Ext and the other is on Int. $\square$

\null

\noindent Recall that a \textit{solid} grid graph is a grid graph without holes. Note that an $m \times n$ grid graph is a solid grid graph such that the outer boundary is an $(m-1)\times(n-1)$ rectangle, with corners at $(0,0), (m-1,0), (m-1,n-1)$, and $(0,n-1)$. We call this rectangle the \textit{boundary of $G$}.

\endgroup 

\null

\noindent \textbf{Definitions.} Let $G$ be an $m \times n$ grid graph and let $X_1,X_2$ be two distinct boxes of $G$. If $X_1$ and $X_2$ share an edge of $G$, we say that $X_1$ and $X_2$ are adjacent. Define a \textit{walk of boxes in G} to be a sequence $X_1, ...,X_r$ of boxes in $G$, not necessarily distinct, such that for all $j \in \{1,2, ..., r-1 \}$, $X_j$ is adjacent to $X_{j+1}$ or $X_j = X_{j+1}$, We denote such a walk by $W(X_1,X_r)$. For each $j \in \{1,...,r-1\}$,  we call the edge of $G$ that $X_j$ and $X_{j+1}$ share a \textit{gluing edge of $W(X_1,X_r)$}, whenever $X_j$ and $X_{j+1}$ are distinct boxes. If for all $i,j \in \{1,2, ..., r \}$ with $i \neq j$, $X_i$ is distinct from $X_j$, we call the sequence a \textit{path of boxes in G} and denote it by $P(X_1,X_r)$. A cycle of boxes in $G$ is a walk $X_1, ..., X_r$ such that $X_1=X_r$ and for $i,j \in \{1, ..., r-1 \}$ $X_i \neq X_j$.

\null 

\noindent \textcolor{blue}{\textbullet}The rest of Section 1 contains definitions and technical results used in Sections 2-4. Our reconfiguration strategy relies on controlling which edges belong to the Hamiltonian
cycle by applying moves. To analyze when moves can add or remove specific edges while preserving the Hamiltonian property, we need to understand how the Hamiltonian cycle $H$ decomposes the grid graph into components that we call $H$-components. This decomposition is introduced in Section 1.1. In section 1.2 we define the Follow-the-wall construction and use it to build walks of boxes that respect the structure of $H$, called $H$-walks. $H$-walks are used extensively in many proofs throughout the rest of the paper. In Section 1.3 we prove some basic properties of $H$-components; In section 1.4 we give a more detailed description of double-switch moves and prove a lemma about their validity. \textcolor{blue}{\textbullet}

\subsection{The $H$-decomposition of $G$}

\begingroup
\setlength{\intextsep}{0pt}
\setlength{\columnsep}{20pt}
\begin{wrapfigure}[]{r}{0cm}
\begin{adjustbox}{trim=0cm 0cm 0cm 1cm}
\begin{tikzpicture}[scale=1.5]

 \begin{scope}[xshift=0cm]
{
\draw[gray,very thin, step=0.5cm, opacity=0.5] (0,0) grid (3,2.5);

\fill[blue!50!white,opacity=0.5](0.5, 0.5) rectangle  (2.5,2); 

\fill[green!50!white,opacity=0.5](0,0)--++(3,0)--++(0,0.5)--++(-3,0); 
\fill[green!50!white,opacity=0.5](0,2)--++(3,0)--++(0,0.5)--++(-3,0); 
\fill[green!50!white,opacity=0.5](0,0.5)--++(0.5,0)--++(0,1.5)--++(-0.5,0); 
\fill[green!50!white,opacity=0.5](2.5,0.5)--++(0.5,0)--++(0,1.5)--++(-0.5,0);

\draw[blue, line width=0.5mm] (1.5,1.5)--++(0,-0.5)--++(-0.5,0)--++(0,0.5)--++(1,0);
\draw[blue, line width=0.5mm] (1.5,2)--++(1,0)--++(0,-1.5);

\draw[red, line width=0.5mm, dotted](0.75,0.75)--++(0,1)--++(1.5,0)--++(0,-1)--++(-1.5,0);

\node[below] at (1.5,0) [scale=0.75] {\small{\begin{tabular}{c}  Fig. 1.2. An $m \times n$ grid graph $G$ shaded blue, \\ a subgraph $H$ of $G$ in blue, BBoxes($G$) shaded \\  green, an $H$-cycle in $G$ dotted red. \end{tabular}}};;

}
\end{scope}

\end{tikzpicture}
\end{adjustbox}
\end{wrapfigure}

\textbf{Definitions.} A \textit{walk} (of length $r$) in a graph is an alternating sequence $v_0e_1v_1e_2....e_rv_r$ of vertices and edges. Define a \textit{lazy walk} to be sequence of edges and vertices where every edge is in between two vertices that are its endpoints, and in between every two edges there is a vertex or multiple copies of a vertex. That is, a lazy walk is roughly a walk in which consecutive vertices can be the same, allowing the walk to remain at a vertex for one or more steps without traversing any edges.

Let $G$ be an $m \times n$ grid graph and let $H$ be any subgraph in $G$. Let $X_1,X_2$ be two adjacent boxes of $G$. If $E(X_1)\cap E(X_2) \cap E(H)=\emptyset$, we say that $X_1$ and $X_2$ are \textit{$H$-neighbours} or $X_1$ is \textit{$H$-adjacent} to $X_2$. Define an \textit{H-walk of boxes in G} ($H$-walk) to be a sequence $X_1, ...,X_r$ of boxes in $G$, not necessarily distinct, such that for all $j \in \{1,2, ..., r-1 \}$, $X_j$ is an $H$-neighbour of $X_{j+1}$ or $X_j = X_{j+1}$. 

If for all $i,j \in \{1,2, ..., r \}$ with $i \neq j$, $X_i$ is distinct from $X_j$, we call the sequence an \textit{H-path of boxes in G} and denote it by $P(X_1,X_r)$. Let $r \geq 4$. Define an \textit{H-cycle of boxes in G} ($H$-cycle) to be a set $X_1, X_2, ..., X_r=X_1$ of boxes in $G$ such that for each $j \in \{1,2, ..., r-1 \}$, $X_j$ is an $H$-neighbour of $X_{j+1}$ and the boxes $X_1, ..., X_{r-1}$ are distinct. Let $C$ be an $H$-cycle of boxes in $G$. We note that every box of $C$ has exactly two gluing edges. Proposition 1.3 will show that if $H$ is a Hamiltonian cycle of $G$, then there are no $H$-cycles of boxes in $G$.

\noindent Define a \textit{boundary box} of $G$ to be a box that is incident on the boundary of $G$ but that is not a box of $G$. denote the set of boundary boxes of $G$ by \textit{BBoxes(G)}. Define $G_{-1}$ to be the graph with vertex set $V(G_{-1})=V(G) \cup V(\text{BBoxes}(G))$ and edge set $E(G_{-1})=E(G) \cup E(\text{BBoxes}(G))$. We extend the definitions of $H$-walks, $H$-paths and $H$-cycles to $G_{-1}$. See Figure 1.2 for an illustration.

\endgroup 

\null 

\begingroup
\setlength{\intextsep}{0pt}
\setlength{\columnsep}{15pt}
\begin{wrapfigure}[]{l}{0cm}
\begin{adjustbox}{trim=0cm 0cm 0cm 0cm}
\begin{tikzpicture}[scale=1.25]

\draw[gray,very thin, step=0.5cm, opacity=0.5] (0,0) grid (4,3);

\fill[green!50!white,opacity=0.5](0.5, 0.5) rectangle  (1,2.5); 
\fill[green!50!white,opacity=0.5](3, 0.5) rectangle  (3.5,2.5); 
\fill[green!50!white,opacity=0.5](1, 0.5) rectangle  (3,1); 
\fill[green!50!white,opacity=0.5](1, 2) rectangle  (3,2.5); 

\draw[green!50!black, line width=0.5mm] 
(0.75,0.75)--++(0,1.5)--++(2.5,0)--++(0,-1.5)--++(-2.5,0);

\draw [blue, line width=1mm] (2.5,2)--++(0,0.5);

\foreach \x in {0,...,4}
\draw[orange, line width=0.5mm] 
(1+0.5*\x,0.5)--++(0,0.5);

\foreach \x in {0,...,4}
\draw[orange, line width=0.5mm] 
(1+0.5*\x,2)--++(0,0.5);

\foreach \x in {0,...,2}
\draw[orange, line width=0.5mm] 
(0.5,1+0.5*\x)--++(0.5,0);

\foreach \x in {0,...,2}
\draw[orange, line width=0.5mm] 
(3,1+0.5*\x)--++(0.5,0);

\draw [blue, line width=1mm] plot [smooth, tension=0.75] coordinates {(1.5,2)(2.25,1.75)(2.5,2)};
\draw [blue, line width=1mm] plot [smooth, tension=0.75] coordinates {(2.5,2.5)(2.25,2.7)(1.75,2.75)(1,2.5)};

\draw[fill=blue, opacity=1] (1.5,2) circle [radius=0.05];
\draw[fill=blue, opacity=1] (1,2.5) circle [radius=0.05];

\node[right] at (2.45,2.35) [scale=1] {\small{$e$}};
\node[below] at (1.5,2) [scale=1] {\small{$a$}};
\node[above] at (1,2.5) [scale=1] {\small{$b$}};
\node at (1.25,2.25) [scale=1] {\small{$X_i$}};

\node[below] at (2,0) [scale=0.75] {\small{\begin{tabular}{c}  Fig. 1.3. $P(a,b)$ in blue, gluing edges of $C$ in \\    orange, Q in dark green, $C$ shaded in light green.  \end{tabular}}};;

\end{tikzpicture}
\end{adjustbox}
\end{wrapfigure}

\noindent \textbf{Proposition 1.3.} Let $G$ be an $m \times n$ grid graph and let $H$ be Hamiltonian cycle of $G$. Then every $H$-cycle in $G_{-1}$ is contained in $\textrm{BBoxes}(G)$. 

\null

\noindent \textit{Proof.} 
By way of contradiction assume that there is an $H$-cycle $C$ in $G_{-1}$ with boxes $X_1, X_2, ..., X_r=X_1$ that has a box $X_i$ contained in $G$. Let $c_1, c_2, ..., c_r=c_1$ be the centers of the boxes of $C$. That is, if $X_j=R(k,l)$ then $c_j=(k+\frac{1}{2}, l+\frac{1}{2})$. Note that for each $j \in \{1,2,...,r-1\}$, $[c_j,c_{j+1}]$ intersects the gluing edge of $X_j$ and $X_{j+1}$ and $[c_j,c_{j+1}]$ intersects no other edge of $G_{-1}$.

We will first show that the set of segments $[c_j,c_{j+1}]$ is a non-self-intersecting polygon $Q$. Since $c_1=c_r$, $Q$ is a polygon. For a contradiction assume that $Q$ is self intersecting, so there are points $c_{j-1}$, $c_j$, $c_{j+1}$ and $c_i$ such that the segments $[c_{j-1},c_j]$, $[c_j,c_{j+1}]$ and $[c_j,c_i]$ are edges of $Q$. But then $X_j$ has three gluing edges, a contradiction.

By JCT, $Q$ divides the plane into two subsets Int and Ext such that any two points within a subset can be joined by a path that does not intersect $Q$ while a path joining a point of Int to a point of Ext must intersect $Q$. Let $V(\textrm{Int})$ be the vertices of $G_{-1}$ contained in Int and let $V(\textrm{Ext})$ be the vertices of $G_{-1}$  contained in Ext. Note that any box of $C$ contains at least one vertex in $V(\textrm{Int})$ and one vertex in $V(\textrm{Ext})$, so both sets are nonempty. In particular, this is true for $X_i$. For definiteness, assume that $a \in V(X_i)$ is contained in $V(\textrm{Int})$ and $b \in V(X_i)$ is contained in $V(\textrm{Ext})$. By JCT, $V(\textrm{Int}) \cup V(\textrm{Ext})=V(G)$ and $V(\textrm{Int}) \cap V(\textrm{Ext})=\emptyset$. 
Consider the subpath $P(a,b)$ of $H$. By JCT again, there is an edge $e \in P(a,b) \subset H$ intersecting $Q$ at some segment $[c_j, c_{j+1}]$. But then $e$ is a gluing edge of $C$, so $e$ cannot belong to $H$. $\square$

\noindent Let $G$ be an $m \times n$ grid graph and let $H$ be a Hamiltonian cycle of $G$. Let $X_1,..., X_r$ be a set of boxes in $G$ such that for any $i,j \in \{1,2, ..., r \}$, there is an $H$-path $P(X_i,X_j)$ between $X_i$ and $X_j$ contained in $G$. We call $X_i$ and $X_j$ the \textit{end-boxes} of $P(X_i,X_j)$. We say that $X_i$ and $X_j$ are $H$-path-connected in $G$ and the set of boxes $\{ X_1,..., X_r \}$ is an \textit{$H$-path connected set of boxes in $G$}. If the set $\{X_1,..., X_r\}$ contains no cycles of boxes we call it an \textit{$H$-tree}. A \textit{$H$-component} of $G$ is a maximal $H$-path connected set of boxes. 

\null 

\noindent  \textbf{Corollary 1.4.}  A Hamiltonian cycle $H$ partitions the boxes of $G$ into $H$-path-connected $H$-components which are maximal $H$-trees. In particular, if an $H$-path is contained in an $H$-component, that $H$-path is unique. $\square$

\subsection{The follow-the-wall construction}

\textbf{Definitions.} Let $G$ be a graph. A \textit{trail} is a walk in $G$ where all edges are distinct. A trail where the first and last vertices coincide is called a \textit{closed trail} or a \textit{circuit}. A \textit{directed walk} (of length $s$) is an alternating sequence $v_0e_1v_1e_2....e_sv_s$ of vertices and directed edges such that for $j \in \{1,...,s\}$, the directed edge $e_j$ has tail $v_{j-1}$ and head $v_j$. A \textit{directed trail} is a directed walk where all directed edges are distinct. We will use the notation $\overrightarrow{K}$ to denote directed trails. 

Let the box $X$ be incident on a directed edge $(u,v)$ of the integer grid. Then $X$ and its vertices not incident on the edge $e_j$ are either on the \textit{right} or on the \textit{left} side of $(u,v)$. See Figure 1.4 for an illustration and Appendix A.1 for a more precise definition of a box being on the right or left side of a directed edge. Note that if $X$ is on the right side of $(u,v)$, then the other box incident on $(u,v)$, say $X'$, is on the right side of $(v,u)$.

\begingroup
\setlength{\intextsep}{0pt}
\setlength{\columnsep}{20pt}
\begin{wrapfigure}[]{r}{0cm}
\begin{adjustbox}{trim=0cm 0cm 0cm 0cm}
\begin{tikzpicture}[scale=2]
\begin{scope}[xshift=0cm]{
\draw[gray,very thin, step=0.5cm, opacity=0.5] (0,0) grid (1,1);

\begin{scope}
[very thick,decoration={
    markings,
    mark=at position 0.6 with {\arrow{>}}}
    ]
    \draw[postaction={decorate}, blue, line width=0.5mm] (0.5,0)--++(0,0.5);
    \draw[postaction={decorate}, blue, line width=0.5mm] (0.5,0.5)--++(0.5,0);
\end{scope}

\draw[fill=blue, opacity=1] (1,0.5) circle [radius=0.035];
\draw[fill=blue, opacity=1] (0.5,0.5) circle [radius=0.035];
\draw[fill=blue, opacity=1] (0.5,0) circle [radius=0.035];

\node[left] at (0.5,0) [scale=1] {\small{$u$}};
\node[left] at (0.5,0.5) [scale=1] {\small{$v$}};
\node[above] at (1,0.5) [scale=1] {\small{$w$}};
\node at (0.65,0.15) [scale=1] {\small{$X_i$}};
\node at (0.8,0.35) [scale=1] {\small{$X_{i+1}$}};

\node[below] at (0.5, 0) [scale=0.8]{\small{\begin{tabular}{c} Fig. 1.4 (i).  $w$ is \\ right of $e_j$. \end{tabular}}};;

} \end{scope}
\begin{scope}[xshift=1.5cm]{
\draw[gray,very thin, step=0.5cm, opacity=0.5] (0,0) grid (1,1);

\begin{scope}
[very thick,decoration={
    markings,
    mark=at position 0.6 with {\arrow{>}}}
    ]
    \draw[postaction={decorate}, blue, line width=0.5mm] (0.5,0)--++(0,0.5);
    \draw[postaction={decorate}, blue, line width=0.5mm] (0.5,0.5)--++(0,0.5);
\end{scope}

\draw[fill=blue, opacity=1] (0.5,1) circle [radius=0.035];
\draw[fill=blue, opacity=1] (0.5,0.5) circle [radius=0.035];
\draw[fill=blue, opacity=1] (0.5,0) circle [radius=0.035];

\node[left] at (0.5,0) [scale=1] {\small{$u$}};
\node[left] at (0.5,0.5) [scale=1] {\small{$v$}};
\node[left] at (0.5,1) [scale=1] {\small{$w$}};
\node at (0.75,0.25) [scale=1] {\small{$X_i$}};
\node at (0.75,0.75) [scale=1] {\small{$X_{i+1}$}};

\node[below] at (0.55, 0) [scale=0.8]{\small{\begin{tabular}{c} Fig. 1.4 (ii).  $w\neq u$  \\ is collinear with $e_j$. \end{tabular}}};;

} \end{scope}
\begin{scope}[xshift=3cm]{
\draw[gray,very thin, step=0.5cm, opacity=0.5] (0,0) grid (1,1);

\begin{scope}
[very thick,decoration={
    markings,
    mark=at position 0.6 with {\arrow{>}}}
    ]
    \draw[postaction={decorate}, blue, line width=0.5mm] (0.5,0)--++(0,0.5);
    \draw[postaction={decorate}, blue, line width=0.5mm] (0.5,0.5)--++(-0.5,0);
\end{scope}

\draw[green!50!black, dotted, line width=0.5mm] (0.5,0.5)--++(0,0.5);

\draw[fill=blue, opacity=1] (0.5,0.5) circle [radius=0.035];
\draw[fill=blue, opacity=1] (0,0.5) circle [radius=0.035];
\draw[fill=blue, opacity=1] (0.5,0) circle [radius=0.035];

\node[right] at (0.5,0) [scale=1] {\small{$u$}};
\node[right] at (0.5,0.5) [scale=1] {\small{$v$}};
\node[above] at (0,0.5) [scale=1] {\small{$w$}};
\node[right] at (0.5,1) [scale=1] {\small{$v'$}};
\node at (0.75,0.25) [scale=1] {\small{$X_i$}};
\node at (0.75,0.75) [scale=1] {\small{$X_{i+1}$}};
\node at (0.25,0.75) [scale=1] {\small{$X_{i+2}$}};

\node[below] at (0.55, 0) [scale=0.8]{\small{\begin{tabular}{c} Fig. 1.4 (iii).  $w$ is \\ left of $e_j$. \end{tabular}}};;

} \end{scope}

\end{tikzpicture}
\end{adjustbox}
\end{wrapfigure}

\noindent Let $H$ be a cycle in an $m \times n$ grid graph $G$. Let $K=v_1, ..., v_{s+1}$ be a subpath of $H$. We orient $K$ to obtain a directed trail $\overrightarrow{K}=(v_1,v_2), ..., (v_s,v_{s+1})= e_1 ...., e_s$. We will use $\overrightarrow{K}$ to construct an $H$-walk of boxes in {$G_{-1}$} that we will call  \textit{the right $H$-walk induced by $\overrightarrow{K}$} and denote it by $W_{\overrightarrow{K}, \textrm{right}}(X_1,X_r)$, where $X_1$ and $X_r$ are the end-boxes of $W_{\overrightarrow{K}, \textrm{right} }(X_1,X_r)$. Roughly, $W_{\overrightarrow{K}, \textrm{right},}(X_1,X_r)$ will be the walk of boxes that a ``walker'' would encounter as they followed along the side of $\overrightarrow{K}$ when starting on the right side of the first edge $e_1$ of $\overrightarrow{K}$. We will call this construction the \textit{follow-the-wall} construction (FTW). This is very similar to the well-known hand-on-the-wall maze-solving algorithm. 

Let $X$ be the box of $G_{-1}$ on the right of the edge $e_1$ of $\overrightarrow{K}$ and let $e_j$ be the $j^{\textrm{th}}$ edge of $\overrightarrow{K}$. Then $X=X_1$ is the first box of the $H$-walk. Let $e_j=(u,v)$, $e_{j+1}=(v,w)$ and let $X_i$ be on the right of the edge $e_j$. There are three possibilities for the position of $w$ with respect to $(u,v)$: $w$ is right of $(u,v)$, $w \neq u$ is colinear with $(u,v)$ and $w$ is left of $(u,v)$. For the last case define $e_j'=(v, v')$ to be the edge in $G_{-1} \setminus H$ that is colinear with $e_j$ and set $e_j''=(v', v)$. See Figure 1.4.

\endgroup 

\null 

\begin{itemize}[topsep=0.2cm, partopsep=0.2pt]
    \item[(i)] \textit{$w$ is right of $e_j$.} Then $X_{i+1}$ is on the right of the edge $e_{j+1}$. Note that in this case, the walk has a repeated box since $X_i=X_{i+1}.$
    \item[(ii)] \textit{$w \neq u$ is collinear with $e_j$.} Then $X_{i+1}$ is on the right of the edge $e_{j+1}$. 
    \item[(iii)] \textit{$w$ is left of $e_j$.} Then $X_{i+1}$ is on the right of the edge $e_j'$ and $X_{i+2}$ is on the right of the edge $e_{j+1}$.
\end{itemize}

\begingroup
\setlength{\intextsep}{0pt}
\setlength{\columnsep}{20pt}
\begin{wrapfigure}[]{l}{0cm}
\begin{adjustbox}{trim=0cm 0.5cm 0cm 0cm}
\begin{tikzpicture}[scale=1.5]
\begin{scope}[xshift=0cm]
{
\node[inner sep=0pt] at (1.99, 1.5)[scale=0.27]
    {\includegraphics[width=.25\textwidth]{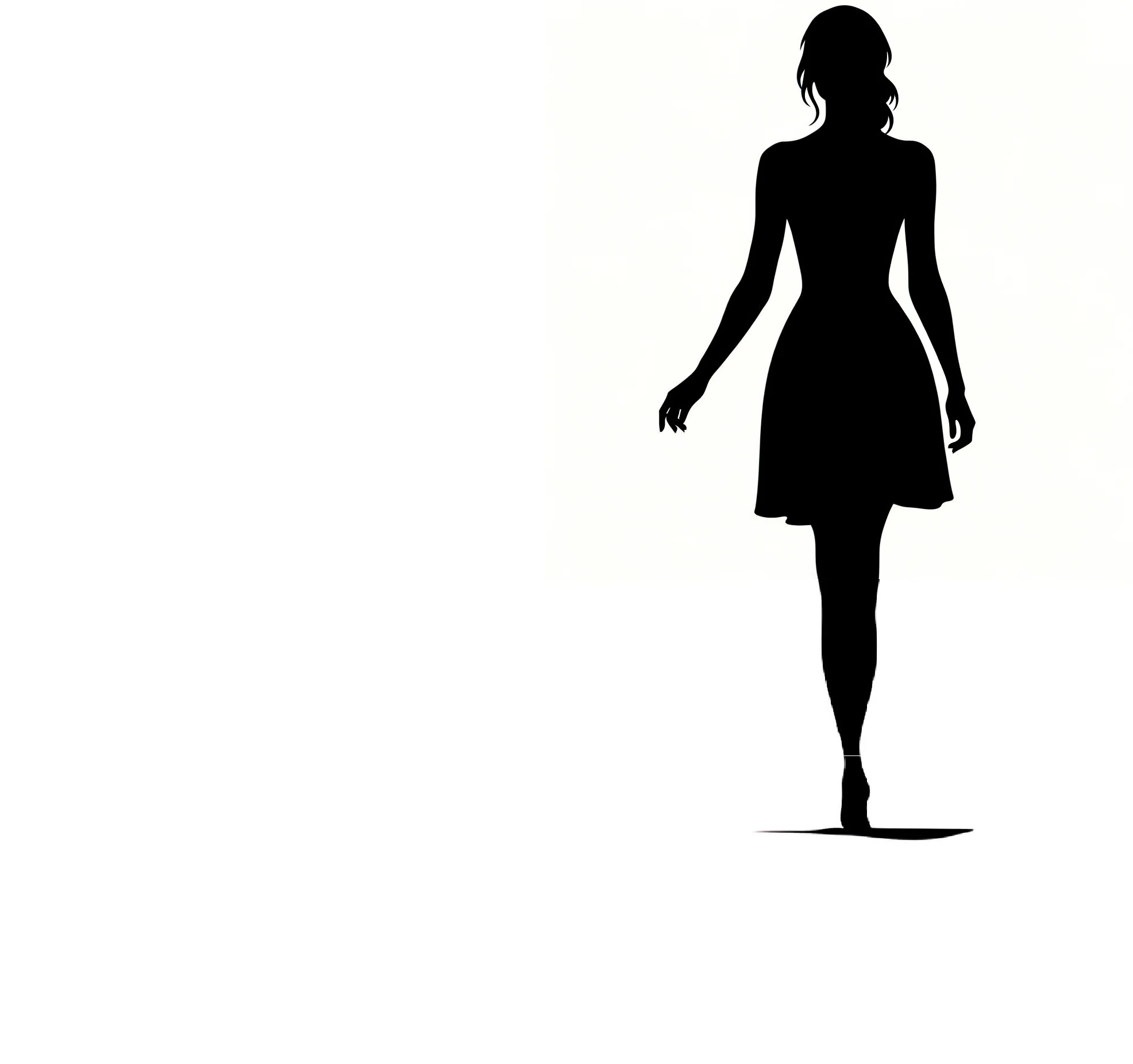}};

\draw[gray,very thin, step=1cm, opacity=0.5] (0,0) grid (3,3);

\fill[green!40!white,opacity=0.5](0,1) rectangle  (2,2);
\fill[orange!40!white,opacity=0.5](0,0)--++(3,0)--++(0,1)--++(-3,0); 
\fill[orange!40!white,opacity=0.5](2,1)--++(1,0)--++(0,1)--++(-1,0); 
\fill[orange!40!white,opacity=0.5](1,2)--++(2,0)--++(0,1)--++(-2,0); 

\begin{scope}
[very thick,decoration={
    markings,
    mark=at position 0.6 with {\arrow{>}}}
    ]
    \draw[postaction={decorate}, blue, line width=1mm] (0,1)--++(1,0);
    \draw[postaction={decorate}, blue, line width=1mm] (1,1)--++(1.03,0);
    \draw[postaction={decorate}, blue, line width=1mm] (2,1)--++(0,1);
    \draw[postaction={decorate}, blue, line width=1mm] (2.03,2)--++(-1.03,0);
    
\end{scope}


\foreach \x in {1,...,3}
\node at (\x-0.5, 0.5) [scale=1]
{\small{$Y_{\x}$}};

\foreach \x in {4,...,6}
\node at (\x-3.5, 1.5) [scale=1]
{\small{$Y_{\x}$}};

\foreach \x in {7,...,9}
\node at (\x-6.5, 2.5) [scale=1]
{\small{$Y_{\x}$}};

\node[above] at (0.6,1) [scale=1] {\small{$e_1$}};
\node[above] at (1.625,1) [scale=1] {\small{$e_2$}};
\node[left] at (2.05,1.65) [scale=1] {\small{$e_3$}};
\node[below] at (1.375,2) [scale=1] {\small{$e_4$}};

\node[below] at (1.5,0) [scale=0.75] {\small{\begin{tabular}{c} Fig. 1.5.  A subtrail $\overrightarrow{K}=\overrightarrow{K}(e_1,e_4)$ of a \\ 

Hamiltonian cycle of a grid graph $G$ in blue; \\

$\Phi (\overrightarrow{K}, \textrm{right} )=Y_1 Y_2Y_3Y_6Y_9Y_8$ shaded orange \\ 

and $\Phi (\overrightarrow{K}, \textrm{left} )=Y_4Y_5Y_5Y_5$ shaded green;
\\  silhouette in gray following the wall along $e_3$. 
\end{tabular}}};;

}
\end{scope}
\end{tikzpicture}
\end{adjustbox}
\end{wrapfigure}

\noindent We say that the edge $e_{j+1}$ \textit{adds} to the $H$-walk $W_{\overrightarrow{K},\textrm{right}}(X_1,X_r)$ the box $X_{i+1}$, in cases (i) and (ii), and boxes $X_{i+1}$ and $X_{i+2}$ in case (iii). If $e_s$ adds more than one box to $W_{\overrightarrow{K}, \textrm{right}}(X_1,X_r)$ then we will adopt the convention that $X_r$ is the last box added by the edge $e_s$ and that all the boxes added by each edge of $e_j$, $j\in \{1, ...,s \}$ are on the right side of $e_j$. Note that this convention is necessary for the box $X_{i+1}$ in Case (iii). We remark that the first edge $e_1$ can only add the single box $X_1$. The left $H$-walk $W_{\overrightarrow{K},\textrm{left}}(X_1,X_r)$ induced by $\overrightarrow{K}$ can be constructed analogously. 

Let $H$ denote either a path or a cycle in $G$. Let $\mathcal{K}(H)$ be the set of all directed subtrails of $H$. We can view the FTW construction as a function $\Phi$ that assigns an $H$-walk to elements of $\mathcal{K}(H) \times \{\text{right, left}\}$. We will take a closer look at $H$-trails in the case where $H$ is a Hamiltonian 
cycle of $G$. 

\null

\noindent Let $H=v_1, ..., v_r,v_1$ be a Hamiltonian cycle in $G$. Orient $H$ as a directed circuit $\overrightarrow{K}_H$. We observe that it is possible to choose a starting edge $e_j=(v_j, v_{j+1})$ of $\overrightarrow{K}_H$ so that the directed circuit $\overrightarrow{K}_H=e_j, e_{j+1}$, ...,  $e_{j-1}$ is such that $\textrm{Boxes} (\Phi (\overrightarrow{K}_H, \textrm{right})) \cup\textrm{Boxes} (\Phi (\overrightarrow{K}_H, \textrm{left})) = \textrm{Boxes}(G_{-1})$. We record an equivalent statement for reference as Observation 1.5 below. From here on, all circuits $\overrightarrow{K}_H$ will be assumed to satisfy Observation 1.5.

Let $H=v_1, v_2,   ... ,  v_r,v_1$ be a Hamiltonian cycle in $G$. Note that any subtrail of $\overrightarrow{K}_H$ is completely determined by its first and last edges. Therefore, it will be fitting to use the notation $\overrightarrow{K}(e_s,e_t)$ to denote the unique subtrail starting at edge $e_s$ and ending at edge $e_t$. 

\endgroup

\subsection{Properties of $H$-components}

\textbf{Observation 1.5.} Let $G$ be an $m \times n$ grid graph and let $H$ be a Hamiltonian cycle of $G$. Then for every box $X$ in $G_{-1}$  there is an edge $e$ of $\overrightarrow{K}_H$ and side $\in \{\text{right, left}\}$ such that $e$ adds $X$ to $W_{\overrightarrow{K}_H,\text{ side}}$.

\null 

\noindent \textbf{Lemma 1.6.} Let $G$ be an $m \times n$ grid graph and let $Q=v_1, ..., v_r,v_1$ be any cycle in $G$ and let $U$ be the region bounded by $Q$. Orient the edges of $G$ into the directed circuit $\overrightarrow{K}_Q=(v_1,v_2), ..., (v_r,v_1)$. Then $\textrm{Boxes} (\Phi (\overrightarrow{K}_Q, \textrm{right}))\subseteq U$ and $\textrm{Boxes} (\Phi (\overrightarrow{K}_Q, \textrm{left})) \subseteq G_{-1} \setminus U$ iff $\Phi((v_1,v_2), \textrm{right})$ is a box of $U$.

\null 

\noindent \textit{Proof.} This follows from Corollary 1.2, the definition of FTW, and induction on the edges of $\overrightarrow{K}_Q$. $\square$

\null 

\noindent Note that, by JCT, $H$ divides the boxes of $G_{-1}$ into the disjoint sets $\text{int}(H)$ and $\text{ext}(H)$, where $\text{int}(H)$ is the bounded region.

\null 

\noindent \textbf{Proposition 1.7.} Let $G$ be an $m \times n$ grid graph and let $H$ be a Hamiltonian cycle of $G$. Then $\text{int}(H)$ is an $H$-component of $G$.

\null 

\noindent \textit{Proof.} Let $H=v_1, ..., v_r,v_1$ and, for definiteness, assume that $\Phi((v_1,v_2), \text{right}) \in \text{int}(H)$. By Lemma 1.6,  $\textrm{Boxes} (\Phi (\overrightarrow{K}_H, \textrm{right}))$ $\subseteq \text{int}(H)$ and $\textrm{Boxes} (\Phi (\overrightarrow{K}_H, \textrm{left})) \subseteq G_{-1} \setminus \text{int}(H)=\text{ext}(H)$. Since $H$ is contained in $G$ and $H$ is the boundary of $\text{int}(H)$, $\text{int}(H)$ is contained in $G$. We check that $\text{int}(H)$ is $H$-path-connected and maximal.

\noindent Consider any box $Z \in \text{int}(H)$. By Observation 1.5, there is an edge $e \in \overrightarrow{K}_H$ such that $e$ adds $Z$ to $\Phi (\overrightarrow{K}_H, \textrm{left})$ or $e$ adds $Z$ to $\Phi (\overrightarrow{K}_H, \textrm{right})$. The former implies that $Z$ belongs to $\text{ext}(H)$, contradicting that $Z \in \text{int}(H) $. It must be the case that $Z \in \textrm{Boxes} (\Phi (\overrightarrow{K}_H, \textrm{right}))$, which is $H$-path connected. Note that this also shows that $ \text{int}(H)= \textrm{Boxes} (\Phi (\overrightarrow{K}_H, \textrm{right}))$.

\noindent To see that $\text{int}(H)$ is maximal, we note that $G \subset \text{int}(H) \cup \text{ext}(H)$ and that $\text{int}(H) \cap \text{ext}(H) = \emptyset$, so $\text{int}(H)$ cannot be extended. $\square$

\null 

\noindent \textbf{Corollary 1.8.} Let $G$ be an $m \times n$ grid graph, $H=v_1, ..., v_r,v_1$ be a Hamiltonian cycle of $G$, and assume that $\Phi((v_1,v_2), \textrm{right})$ is a box of $\text{int}(H)$. Then $\textrm{Boxes} (\Phi (\overrightarrow{K}_H, \textrm{right}))= \text{int}(H)$ and $\textrm{Boxes} (\Phi (\overrightarrow{K}_H, \textrm{left}))= \text{ext}(H)$. $\square$

\noindent Let $G$ be an $m \times n$ grid graph and let $H$ be a Hamiltonian cycle of $G$. Let $J_0, J_1, ..., J_s$ be the $H$-components of $G$, where $J_0=\text{int}(H)$. It follows from Proposition 1.7 that $\text{BBoxes}(G)=G_{-1}\setminus G$ is contained in $\text{ext}(H)$ and that all other components $J_1, ..., J_s$ of $G$ are contained in $\text{ext}(H) \setminus  \text{BBoxes}(G)$. We write this as Observation 1.10 below for reference. We call the $H$-components $J_1, ..., J_s$ \textit{cookies of $G$.}  If a cookie $J$ has more than one box, we call $J$ a \textit{large} cookie. Otherwise, we say that $J$ is a \textit{small} cookie.

Let $J$ be a cookie of $G$. We define a \textit{neck of $J$} to be a box $N_J$ of $J$ that is incident on a boundary edge $e_J$ of $G$ with $e_J \notin H$. We call $e_J$ a \textit{neck-edge} of $J$. Note that the other box incident on $e_J$ must be in $G_{-1}\setminus G$. With these definitions, Lemma 1.6 has the following corollary:

\null

\noindent \textbf{Corollary 1.9.} Let $G$ be an $m \times n$ grid graph and let $H$ be a Hamiltonian cycle of $G$,  and let $J_0, ..., J_s$ be the $H$-components of $G$. Then every edge of $H$ is incident on a box of $J_0$ and a box of $G_{-1} \setminus J_0$. $\square$

\null

\noindent \textbf{Observation 1.10.} Let $J$ be an $H$-component of $G$. Then $J \subseteq \text{int}(H)$ or $J \subset \text{ext}(H)$.

\null 

\noindent \textbf{Corollary 1.11.} Let $X$ and $Y$ be boxes of an $H$-component $J$. Then $X$ and $Y$ are on the same side of $\overrightarrow{K}_H$.

\null

\noindent \textbf{Lemma 1.12.} Let $G$ be an $m \times n$ grid graph, let $H$ be a Hamiltonian cycle of $G$, and let $J$ be a cookie of $G$. Then $J$ has a unique neck. 

\null 

\noindent \textit{Proof.} If $J$ only has one box, we are done, so assume that $J$ has more than one box. Let $Z$ be a box of $J$. By Observation 1.5, we may assume, WLOG, that $Z$ is on the left of $e_z \in \overrightarrow{K}_H$. If no edge is incident on $Z$, we choose one of the four neighbours beside it.

\begingroup 
\setlength{\intextsep}{0pt}
\setlength{\columnsep}{20pt}
\begin{wrapfigure}[]{r}{0cm}
\begin{adjustbox}{trim=0cm 0.5cm 0cm 0cm}
\begin{tikzpicture}[scale=1.5]
\usetikzlibrary{decorations.markings}
\begin{scope}[xshift=0cm] 

\draw[gray,very thin, step=0.5cm, opacity=0.5] (0,0) grid (2.5,2.5);

\fill[blue!40!white,opacity=0.5] plot [smooth, tension=0.75] coordinates {(2,1.5)(2,2) (1.85,2.25)(1.5,2.4)(1, 2.2)(0.75,1.5)(0.975,1)(1,0.5)  (1.15,0.65)(1.1,1.1)(1,1.25)(0.9,1.5)(0.95,1.75)(1,1.9)(1.25,2.2)(1.5,2.25)(1.6,2.25)(1.75,2.15)(1.85,2)(1.9,1.65)};

\draw[blue, line width=0.5mm] (0.5,0.5)--++(0.5,0);

\begin{scope}[very thick,decoration={
    markings,
    mark=at position 0.6 with {\arrow{>}}}
    ] 
    \draw[postaction={decorate}, blue, line width=0.5mm] (2,1.5)--++(0,0.5); 
    \draw[postaction={decorate}, blue, line width=0.5mm] (1,1)--++(0,-0.5); ; 

\end{scope}

\draw [blue, line width=0.5mm] plot [smooth, tension=0.75] coordinates {(2,2) (1.85,2.25)(1.5, 2.4)(1, 2.2)(0.75,1.5)(1,1)};

\node at  (0.75,0.75) [scale=0.8]{\small{$e_t$}};
\node at  (2.25,1.75) [scale=0.8]{\small{$e_z$}};
\node at  (1.3,0.75) [scale=0.8]{\small{$N_J$}};
\node at  (1.25,0.25) [scale=0.8]{\small{$Y$}};
\node at  (1.75,1.75) [scale=0.8]{\small{$Z$}};

\draw[black, line width=0.15mm] (1.2,0.45)--(1.2,0.55);
\draw[black, line width=0.15mm] (1.25,0.45)--(1.25,0.55);
\draw[black, line width=0.15mm] (1.3,0.45)--(1.3,0.55);

\node[below] at (1.25,0) [scale=0.8]{\small{\begin{tabular}{c} Fig. 1.6. $\Phi \big(\overrightarrow{K}(e_z, e_t) , \textrm{left}\big)$ \\  shaded in blue. \end{tabular}}};;

\end{scope}
\end{tikzpicture}
\end{adjustbox}
\end{wrapfigure}

\noindent We claim that there exists a subtrail $\overrightarrow{K}(e_z, e_{t+1})$ of $\overrightarrow{K}_H$ such that $\Phi \big(\overrightarrow{K}(e_z, e_t) , \textrm{left}\big)$ is contained in $J$ but $\Phi \big(\overrightarrow{K}(e_z, e_{t+1}) , \textrm{left}\big)$ is not. Assume for contradiction that for every subtrail of $\overrightarrow{K}_{H}$ starting at $e_z$, $\Phi \big(\overrightarrow{K}(e_z, e_j) , \textrm{left}\big)$ is contained in $J$, where $j \in \{z+1, ..., z\}$. But then $\Phi \big( \overrightarrow{K}_H, \text{left} \big)$ is contained in $J \subset \text{Boxes}(G)$, contradicting that 
$\Phi \big(\overrightarrow{K}_H,\textrm{left}\big)$ contains the boxes of $G_{-1}\setminus G$. It follows that $e_{t+1}$ adds the first box $Y$ of $\Phi \big(\overrightarrow{K}(e_z, e_{t+1}) , \textrm{left}\big)$ that is not contained in $J$. Note that, by definition of $H$-component, $Y$ must belong to $G_{-1} \setminus G$. (Since $Y$ is $H$-adjacent to the box $X$ preceding it, but $Y$ does not belong to $J$, it must be the case that $Y$ is not in $G$). Let $X$ be the box of $J$ preceding $Y$ in $\Phi \big(\overrightarrow{K}(e_z, e_{t+1}) , \textrm{left}\big)$. We have that $X$ and $Y$ are $H$-adjacent and share a boundary edge $e_J$ of $G$ that is not in $H$. By definition of neck of an $H$-component, $X=N_J$. 

To see that the neck of $J$ is unique, assume for contradiction that $J$ has at least two necks, $N_1$ and $N_2$. By Corollary 1.8, $\text{ext}(H)$ is $H$-path-connected. Let $N_1'$ and $N_2'$ be the boxes in BBoxes($G$) that are adjacent to $N_1$ and $N_2$, respectively. Then there is an $H$-cycle $P(N_1',N_2')$, $P(N_2,N_1)$ in $G_{-1}$ that is not contained in BBoxes($G$), which contradicts Proposition 1.3. $\square$

\endgroup 

\null 

\noindent \textbf{Lemma 1.13.} Let $G$ be an $m \times n$ grid graph, let $H$ be a Hamiltonian cycle of $G$, and let $J$ be a cookie of $G$ with neck $N_J$. Orient $H$ and let $(v_x, v_{x+1})$ and $(v_{y-1}, v_y)$ be edges of $N_J$ in $H$, where $\{v_x, v_y\}$ is a boundary edge of $G$. Then $V(J) = V(\overrightarrow{K} ((v_x, v_{x+1}),(v_{y-1}, v_y)))$. 

\null

\noindent \textit{Proof.} Let $\overrightarrow{K}_1=\overrightarrow{K} ((v_x, v_{x+1}),(v_{y-1},v_y))$ and $\overrightarrow{K}_2=\overrightarrow{K} ((v_y, v_{y+1}),(v_{x-1},v_x))$, and note that $\overrightarrow{K}_H=\overrightarrow{K}_1, \overrightarrow{K}_2$. For definiteness, assume that $N_J=\Phi ((v_x,v_{x+1}), \text{right})$. 

\noindent Let $v_z \in V(J)$. Assume for contradiction that $v_z \notin \overrightarrow{K}_1$. Then $v_z \in \overrightarrow{K}_2$. Let $Z$ be a box of $J$ on which $v_z$ is incident, and let $e_z$ be the edge of $\overrightarrow{K}_2$ containing $v_z$ that adds $Z$ to $\Phi (\overrightarrow{K}_H, \text{right})$. Note that we used Corollary 1.11 here. By proof of Lemma 1.12, there is a subtrail $\overrightarrow{K} (e_z,e_t)$ of $\overrightarrow{K}_H$ such that $\Phi(\overrightarrow{K} (e_z,e_t), \text{right})$ is contained in $J$ but $\Phi(\overrightarrow{K} (e_z,e_{t+1}), \text{right})$ is not contained in $J$. (If no edge of $H$ is incident on $Z$ we can just consider $\overrightarrow{K} (e_{z-1},e_t)$ instead.) Note that this implies that $e_t$ must add $N_J$ to $\Phi(\overrightarrow{K} (e_z,e_t), \text{right})$. Then $e_t=(v_x, v_{x+1})$ or $e_t=(v_{y-1}, v_y)$.

\begingroup 
\setlength{\intextsep}{0pt}
\setlength{\columnsep}{20pt}
\begin{wrapfigure}[]{l}{0cm}
\begin{adjustbox}{trim=0cm 0cm 0cm 0cm}
\begin{tikzpicture}[scale=1.5]

\begin{scope}[xshift=0cm]{
\draw[gray,very thin, step=0.5cm, opacity=0.5] (0,0) grid (1.5,1);

\fill[blue!50!white, opacity=0.5] (0.5,0.5) rectangle (1,1);

\fill[green!50!white, opacity=0.5] (0,0) rectangle (1.5,0.5);

\begin{scope}
[very thick,decoration={
    markings,
    mark=at position 0.6 with {\arrow{>}}}
    ]
    \draw[postaction={decorate}, orange, line width=0.5mm] (0,0.5)--++(0.5,0);
    \draw[postaction={decorate}, orange, line width=0.5mm] (1,0.5)--++(0.5,0);

    \draw[postaction={decorate}, blue, line width=0.5mm] (0.5,0.5)--++(0,0.5);
    \draw[postaction={decorate}, blue, line width=0.5mm] (1,1)--++(0,-0.5);

\end{scope}

\draw[fill=blue] (0.5,1) circle [radius=0.05];
\draw[fill=blue] (1,1) circle [radius=0.05];

\draw[fill=orange] (0,0.5) circle [radius=0.05];
\draw[fill=orange] (1.5,0.5) circle [radius=0.05];

\draw[fill=orange] (0.5,0.5) -- ++(135:0.05) arc[start angle=135, end angle=315, radius=0.05] -- cycle;
\draw[fill=blue] (0.5,0.5) -- ++(315:0.05) arc[start angle=315, end angle=495, radius=0.05] -- cycle;

\draw[fill=orange] (1,0.5) -- ++(225:0.05) arc[start angle=225, end angle=405, radius=0.05] -- cycle;
\draw[fill=blue] (1,0.5) -- ++(45:0.05) arc[start angle=45, end angle=225, radius=0.05] -- cycle;

{

\node at (0.75,0.75) [scale=0.8] {\small{$N_J$}};
\node at (0.75,0.25) [scale=0.8] {\small{$N_J'$}};

\node[below] at (0.5,0.5) [scale=0.8] {\small{$v_x$}};
\node[below] at (1,0.5) [scale=0.8] {\small{$v_y$}};
\node[above] at (1.0,1) [scale=0.8] {\small{$v_{y-1}$}};
\node[above] at (0.5,1) [scale=0.8] {\small{$v_{x+1}$}};
\node[below] at (0,0.5) [scale=0.8] {\small{$v_{x-1}$}};
\node[below] at (1.5,0.5) [scale=0.8] {\small{$v_{y+1}$}};

}

\node[below] at (0.75,0) [scale=0.8]{\small{\begin{tabular}{c} Fig. 1.7. $G_{-1} \setminus G$ \\  shaded in green. \end{tabular}}};;

} \end{scope}

\end{tikzpicture}
\end{adjustbox}
\end{wrapfigure}

Let $N_J'$ be the neighbour of $N_J$ in $G_{-1} \setminus G$. Note that $N_J' \in \Phi(\overrightarrow{K}$ $ (e_z,(v_x,v_{x+1})), \text{right})$ and so $\Phi(\overrightarrow{K} (e_z,(v_x,v_{x+1}))$ is not contained in $J$. Then it must be the case that $e_t=(v_{y-1}, v_y)$.

Since $e_z \in \overrightarrow{K}_2$ and $e_t=(v_{y-1}, v_y) \in \overrightarrow{K}_1$, there must be some $j_0 \in \{z, z+1, ..., y-2 \}$ such that for each $j \leq j_0$, $e_j \in \overrightarrow{K}_2$, but $e_{j_0+1} \in \overrightarrow{K}_1$. It follows that $e_{j_0+1}=(v_x, v_{x+1})$ or $e_{j_0+1}=(v_{y-1}, v_y)$.

\endgroup 

\noindent Note that if $e_{j_0+1}=(v_{y-1}, v_y)$ then $(v_{y-2},v_{y-1})=e_{j_0}$. But $(v_{y-2},v_{y-1}) \in \overrightarrow{K}_1$, contradicting $e_{j_0} \in \overrightarrow{K}_2$. Then it must be the case that $e_{j_0+1}=(v_x, v_{x+1})$. But then $\Phi(\overrightarrow{K} (e_z,e_{j_0}), \text{right})$ contains $\Phi(e_{j_0}, \text{right})$, which belongs to $G_{-1} \setminus G$, contradicting that $\Phi(\overrightarrow{K}$ $ (e_z,e_{j_0}), \text{right})$ is contained in $J$. Thus we must have that $v_z \in \overrightarrow{K}_1$. 

\null

\noindent It remains to check that $V(\overrightarrow{K}_1) \subseteq V(J)$. We will prove that $\overrightarrow{K}_1 \subseteq E(J)$. For $i \in  \{x, x+1, ..., y-1 \}$, let $e_i=(v_i, v_{i+1})$. Note that if $\overrightarrow{K}(e_x, e_i) \subseteq E(J)$ then $\Phi(\overrightarrow{K}(e_x, e_i), \text{right}) \subset J$. This follows by definition of FTW and induction on the edges of $\overrightarrow{K}(e_x, e_i)$.

We have that $e_x \in \overrightarrow{K}_1 \cap E(J)$. Assume for contradiction that there exists some $j_0 \in \{x+1, x+2, ..., y-3\}$\footnote{$\Phi(e_{y-1}, \text{right})=N_J$, so $j_0 \leq y-3$} such that $\overrightarrow{K} (e_x, e_{j_0})$ is contained in $E(J)$ but $\overrightarrow{K} (e_x, e_{j_0+1})$ is not. Since $e_{j_0+1} \notin E(J)$, we have that $\Phi (\overrightarrow{K}(e_x, e_{j_0}) , \text{right})$ is contained in $J$ but $\Phi (\overrightarrow{K}(e_x, e_{j_0+1}) , \text{right})$ is not. Let $Z$ be the first box of $\Phi (\overrightarrow{K}(e_x, e_{j_0+1}) , \text{right})$ that is not contained in $J$ and let $Z'$ be the box preceding $Z$ in $\Phi (\overrightarrow{K}(e_x, e_{j_0+1}) , \text{right})$. The fact that $Z \in J$ and $Z' \notin J$ are $H$-neighbours implies that $Z'=N_J$. Since $e_{j_0+1}$ adds $Z$ to $\Phi (\overrightarrow{K}(e_x, e_{j_0+1}) , \text{right})$, $e_{j_0+1}=e_x$ or $e_{j_0+1}=(v_y, v_{y+1})$. But both possibilities contradict that $j_0 \in \{x+1, x+2, ..., y-3\}$. $\square$

\null 

\noindent Let $H$ be a subgraph of $G$. A box of $G$ on vertices $a,b,c,d$ is \textit{switchable} in $H$ if it has exactly two edges in $H$ and the edges are parallel to each other. \textcolor{blue}{\textbullet}We call e box of $G$ with exactly three edges in $H$ a \textit{leaf}.\textcolor{blue}{\textbullet}

\null 

\noindent \textbf{Lemma 1.14.} Let $G$, $H$ and $J_0, ..., J_s$ be as in Corollary 1.9. Then a large cookie has exactly one box incident on a boundary edge of $G$, namely its neck. Furthermore, the neck of each large cookie is switchable.

\null 

\noindent \textit{Proof.} First we show that a large cookie has exactly one box in incident on a boundary edge of $G$, namely its neck. Let $J_i$ be a cookie. By Lemma 1.12, $J_i$ has a neck $N_{J_i}$ and $N_{J_i}$ is incident on the boundary of $G$. Suppose that there is another box $X$ of $J_i$ that is incident on $B_0$. Note that $X \in G_{-1} \setminus J_0$. Let $e$ be the boundary edge of $X$ and let $Y$ be the box in $G_{-1} \setminus G$ that is incident on $e$. Then $e \in H$ or $e \notin H$. Note that $e \notin H$ contradicts Lemma 1.12 so we only need to check the case where $e \in H$. Suppose that $e \in H$. By Corollary 1.9,  $Y \in J_0 \subset G$, contradicting that $Y \in G_{-1} \setminus G$.

Now we show that the neck of each large cookie is a switchable box. Let $X=R(k,l)$ be the neck of a cookie $J_i$. Let $v(k,l)=a$, $v(k+1,l)=b$, $v(k+1,l+1)=c$ and $v(k,l+1)=d$. For definiteness, assume that $\{a,b\}$ is the neck edge of $J_i$. Observe that we must have $0<k<m-1$. It follows that $\{a,d\} \in H$ and $ \{b,c\} \in H$. Since $J_i$ is not a small cookie $\{c,d\} \notin H$. Thus, $X$ is switchable. $\square$

\subsection{Moves}

\textbf{Definitions.} Let $G$ be an $m \times n$ grid graph and let $H$ be a Hamiltonian cycle of $G$. Let $abcd$ be a switchable box with edges $ab$ and $cd$ in $H$. A \textit{switch move} on the box $abcd$ in $H$  removes edges $ab$ and $cd$ and adds edges $bc$ and $ad$. Let $X \in G$ be a switchable box in $H$. We write $\textrm{Sw}(X)$ to denote a switch move.

\noindent A \textit{double-switch move} is pair of switch moves where we first switch $X$ and then find a switchable $Y$ and switch it. We denote a double-switch move by $X \mapsto Y$. See Figure 1.9. If after a double-switch, we get a new Hamiltonian cycle, then we call the move a \textit{valid move}. We call $X\mapsto X$ a \textit{trivial} move.

\noindent Orient $H$ as directed cycle $v_1,..., v_r,v_1$. Let $X$ be a switchable box in $H$ with edges $e_1=(v_s,v_{s+1})$ and $e_2=(v_t,v_{t+1})$. If $v_s$ is adjacent to $v_t$ in $G$, we say that $e_1$ and $e_2$ are \textit{parallel}; and if $v_s$ is adjacent to $v_{t+1}$ in $G$, we say that $e_1$ and $e_2$ are \textit{anti-parallel}. Similarly, we call the box $X$ a \textit{parallel (anti-parallel)} box if its edges are parallel (anti-parallel). 

\noindent Let $J$ be a large cookie with neck $N_J$. We call a valid flip move $N_J \mapsto N_J'$ a neck-shifting flip (NSF) move. Observe that after a NSF move, the large cookie $J$ necessarily becomes the large cookie $(J \cup N_J') \setminus N_J$, with new neck $N_J'$.

We define a \textit{cascade} to be a sequence of moves $\mu_1, ..., \mu_r$ such that for $0\leq j\leq r-1$:

1) $\mu_1$ is valid, 

2) if $\mu_1, ...,\mu_j$ have been applied then $\mu_{j+1}$ is valid, and

3) the sequence may contain NSF moves but does not otherwise create any new cookies.

\null 

\noindent Let $H$ be a Hamiltonian e-cycle of an $m \times n$ grid graph $G$ and let $J$ be a cookie of $H$ with
neck $N_J$ . Consider a cascade $\mu_1, . . . , \mu_r$ where $\mu_r$ is the nontrivial move $Z \mapsto N_J$.  We say
that the cascade $\mu_1, . . . , \mu_r$ \textit{collects} the cookie $J$. Note that all double-switch moves are invertible. For non-adjacent boxes $X$ and $Y$ , the moves $X \mapsto Y$ and $Y \mapsto X$. When $X$ and $Y$ are adjacent with $X$ switchable and $Y$ a leaf
(i.e. $X \mapsto Y$ is a flip move), $X$ must be switched first before $Y$ becomes switchable, so the order matters. 

\null 

\noindent \textbf{Lemma 1.15.} Let $G$ be an $m \times n$ grid graph and let $H$ be a Hamiltonian cycle of $G$. Orient $H$ as a directed cycle $v_1, ..., v_r, v_1$. Then every switchable box of $H$ is anti-parallel.

\null 

\noindent \textit{Proof.} Assume for contradiction that there is a box $X$ in $H$ with parallel edges $e_1$ and $e_2$. For definiteness assume that $X=R(k,l)$, $e_1=e(k,k+1;l)$, $e_2=e(k,k+1;l+1)$ and that $\textrm{Boxes} (\Phi (\overrightarrow{K}_H, \textrm{right}))= \text{int}(H)$. But then $X=\Phi (e_2, \textrm{right}) \in \text{int}(H)$ and $\Phi (e_1, \textrm{right}) \in \text{int}(H)$, contradicting Corollary 1.9. $\square$

\begingroup 
\setlength{\intextsep}{0pt}
\setlength{\columnsep}{10pt}
\begin{wrapfigure}[]{r}{0cm}
\begin{adjustbox}{trim=0cm 0.25cm 0cm 0.25cm}
\begin{tikzpicture}[scale=1.25]
\usetikzlibrary{decorations.markings}

\begin{scope}[xshift=0cm] 


\begin{scope}
[very thick,decoration={
    markings,
    mark=at position 0.6 with {\arrow{>}}}
    ]
    
    \draw[postaction={decorate}, blue, line width=0.5mm] (1,1)--++(0.5,0);
\end{scope}

\begin{scope}
[very thick,decoration={
    markings,
    mark=at position 0.6 with {\arrow{<}}}
    ]
    
    \draw[postaction={decorate}, blue, line width=0.5mm] (1,0.5)--++(0.5,0);
\end{scope}

\draw [blue, line width=0.5mm] plot [smooth, tension=0.75] coordinates {(1,0.5)(0.5,0.35)(0.25,0.75)(0.5,1.25)(1,1)};

\draw [blue, line width=0.5mm] plot [smooth, tension=0.75] coordinates {(1.5,1) (1.75,1.2) (2,1.25)(2.25,1) (2.3,0.7)(2.2,0.4) (1.85,0.3)(1.5,0.5)   };

\node[right] at  (0.3,0.2) [scale=0.8]{\small{$v_1$}};
\node[right] at  (1.1,0.2)  [scale=0.8]{\small{$v_r$}};

\node[above] at  (1,1) [scale=0.8]{\small{$v_s$}};
\draw[dotted, line width=0.2mm] (1.5,1)--++(0,0.6);
\node[above] at  (1.5,1.6) [scale=0.8]{\small{$v_{s+1}$}};
\node[above] at  (1.5,0.5)  [scale=0.8]{\small{$v_t$}};
\draw[dotted, line width=0.2mm] (1,0.5)--++(-0.25,-0.6);
\node[below] at  (0.75,0)  [scale=0.8]{\small{$v_{t+1}$}};

\node at  (1.25,0.75)  [scale=0.8]{\small{$X$}};

\draw [->,black, very thick] (2.1,1.65)--(3.1,1.65);
\node[above] at  (2.6,1.65) [scale=0.8]{\small{$\textrm{Sw}(X)$}};

\node[above] at  (0.25,1.15) [scale=1]{$H$};

\node[below] at (1.25,-0.25) [scale=0.8]{\small{\begin{tabular}{c} Fig. 1.8 (a).  A 
 directed \\ Hamiltonian cycle $H$,  \\  $X$ is switchable, $e_1$ and  \\ $e_2$ are anti-parallel. \end{tabular}}};;

\end{scope}

\begin{scope}[xshift=3cm] 


\draw[, blue, line width=0.5mm] (1,1)--++(0,-0.5);
\draw[ blue, line width=0.5mm] (1.5,1)--++(0,-0.5,0);

\draw [blue, line width=0.5mm] plot [smooth, tension=0.75] coordinates {(1,0.5)(0.5,0.35)(0.25,0.75)(0.5,1.25)(1,1)};

\draw [blue, line width=0.5mm] plot [smooth, tension=0.75] coordinates {(1.5,1) (1.75,1.2) (2,1.25)(2.25,1) (2.3,0.7)(2.2,0.4) (1.85,0.3)(1.5,0.5)   };

\node[right] at  (0.3,0.2) [scale=0.8]{\small{$v_1$}};
\node[right] at  (1.1,0.2)  [scale=0.8]{\small{$v_r$}};

\node[above] at  (1,1) [scale=0.8]{\small{$v_s$}};
\draw[dotted, line width=0.2mm] (1.5,1)--++(0,0.6);
\node[above] at  (1.5,1.6) [scale=0.8]{\small{$v_{s+1}$}};
\node[left] at  (1.5,0.5)  [scale=0.8]{\small{$v_t$}};
\draw[dotted, line width=0.2mm] (1,0.5)--++(-0.25,-0.6);
\node[below] at  (0.75,0)  [scale=0.8]{\small{$v_{t+1}$}};

\node at  (1.25,0.75)  [scale=0.8]{\small{$X$}};

\node[above] at  (2.25,1.1) [scale=1]{$H_1$};
\node[above] at  (0.165,1.01) [scale=1]{$H_2$};

\node[below] at (1.25,-0.25) [scale=0.8]{\small{\begin{tabular}{c} Fig. 1.8 (b).  Cycles \\ $H_1$ and $H_2$,  obtained \\   from $H$ after $Sw(X)$. \end{tabular}}};;

\end{scope}

\end{tikzpicture}
\end{adjustbox}
\end{wrapfigure}

\null 

\noindent We will show below that if we switch a switchable box of $H$ we get a cycle $H_1$ and a cycle $H_2$. We define a $(H_1,H_2)$\textit{-port} to be a switchable box of $H_1 \cup H_2$ that has one edge in $H_1$ and the other in $H_2$. 

\null 

\noindent \textbf{Lemma 1.16.} Let $G$ be  be an $m \times n$ grid graph and let $H$ be a Hamiltonian cycle of $G$. Orient $H$ as a directed cycle $v_1, ..., v_r, v_1$. Let $X$ be a switchable box of $H$ with edges $e_1=(v_s,v_{s+1})$ and $e_2=(v_t,v_{t+1})$, with $s+1<t$.

(i) \ $\textrm{Sw}(X)$ splits $H$ into two cycles, $H_1$ and $H_2$.

(ii) Suppose we apply $\textrm{Sw}(X)$. If $Y$ is an $(H_1,H_2)$-port then $X \mapsto Y$ is a valid double-switch move.

\noindent \textit{Proof.} Removing edges $e_1$ and $e_2$ splits $H$ into two disjoint paths $P_1=P(v_{s+1}, v_t)$ and $P_2=P(v_{t+1},v_r)$, $\{v_r,v_1\},P(v_1,v_s)=P(v_{t+1},v_s)$. By Lemma 1.15,  $e_1$ and $e_2$ are anti-parallel. Then we have that $v_s$ is adjacent to $v_{t+1}$ and $v_{s+1}$ is adjacent to $v_t$. Now ; adding $e_1'=(v_{s+1},v_t)$ gives a cycle $H_1=P_1, e_1'$; and adding $e_2'=(v_s,v_{t+1})$ gives a cycle $H_2=P_2,e_2'$. End of proof for (i).

The proof of (ii) is essentially the same as the proof for (i), so we omit it. $\square$ 

\endgroup

\begingroup 
\setlength{\intextsep}{0pt}
\setlength{\columnsep}{10pt}
\begin{wrapfigure}[]{l}{0cm}
\begin{adjustbox}{trim=0cm 0cm 0cm 0.25cm}
\begin{tikzpicture}[scale=1.3]
\begin{scope}[xshift=0cm]{
\draw[gray,very thin, step=0.5cm, opacity=0.5] (0,0) grid (1.5,1.5);

\draw[blue, line width=0.5mm] (0.5,1.5)--++(-0.5,0)--++(0,-1.5)--++(1.5,0)--++(0,0.5)--++(-1,0)--++(0,0.5)--++(1,0)--++(0,0.5)--++(-0.5,0);

\node at (0.25,0.75) [scale=1] {\small{X}};
\node at (0.75,0.75) [scale=1] {\small{Y}};

\draw [->,black, very thick] (1.25,1.65) to [out=45,in=135] (2.25,1.65);

\node[above] at  (1.75,1.8) [scale=0.8]{\small{$\textrm{Sw}(X)$}};

\node[below] at (0.75,-0.15) [scale=0.8]{\small{\begin{tabular}{c} Fig. 1.9 (a).  A  \\ Hamiltonian e-cycle  \\ $H$  on  a $4 \times 4$ grid; \\ $X$ is switchable. \end{tabular}}};;

} \end{scope}
\begin{scope}[xshift=2cm]
{
\draw[gray,very thin, step=0.5cm, opacity=0.5] (0,0) grid (1.5,1.5);

\draw[blue, line width=0.5mm] (0.5,1.5)--++(-0.5,0)--++(0,-0.5)--++(1.5,0)--++(0,0.5)--++(-0.5,0);

\draw[blue, line width=0.5mm] (0,0)--++(1.5,0)--++(0,0.5)--++(-1.5,0)--++(0,-0.5);

\node at (0.25,0.75) [scale=1] {\small{X}};
\node at (0.75,0.75) [scale=1] {\small{Y}};

\draw [->,black, very thick] (1.25,1.65) to [out=45,in=135] (2.25,1.65);
\node[above] at  (1.75,1.8) [scale=0.8]{\small{$\textrm{Sw}(Y)$}};

\node[below] at (0.75,-.15) [scale=0.8]{\small{\begin{tabular}{c} Fig. 1.9 (b). $H_{\textrm{cycle}}$  \\ and $H_{\textrm{path}}$ after  \\  switching $X$. Note    \\ that $Y$ is an 
 \\ ($H_{\textrm{cycle}}$, $H_{\textrm{path}}$)-port. \end{tabular}}};;

}
 \end{scope}
\begin{scope}[xshift=4cm]
{
\draw[gray,very thin, step=0.5cm, opacity=0.5] (0,0) grid (1.5,1.5);

\draw[blue, line width=0.5mm] (0.5,1.5)--++(-0.5,0)--++(0,-0.5)--++(0.5,0)--++(0,-0.5)--++(-0.5,0)--++(0,-0.5)--++(1.5,0)--++(0,0.5)--++(-0.5,0)--++(0,0.5)--++(0.5,0)--++(0,0.5)--++(-0.5,0);

\node at (0.25,0.75) [scale=1] {\small{X}};
\node at (0.75,0.75) [scale=1] {\small{Y}};

\node[below] at (0.75,-.15) [scale=0.8]{\small{\begin{tabular}{c} Fig. 1.9 (c).  A Hamil- \\ tonian e-cycle $H'$ \\ after switching $Y$.  \end{tabular}}};;

}
 \end{scope}

\end{tikzpicture}
\end{adjustbox}
\end{wrapfigure}

\null 

\noindent \textbf{Hamiltonian e-cycles.} Let $G$ be an $m \times n$ grid graph, let $H$ be a Hamiltonian cycle of $G$ and let $e$ be an edge of $H$ that lies in the boundary of $G$. We call the path $H'=H \setminus e$ a Hamiltonian \textit{e-cycle} of $G$. We remark that all the definitions and results about the case where $H$ is Hamiltonian cycles of $G$ translate immediately to the case where $H'$ is a Hamiltonian e-cycle of $G$. We may just add back the edge $e$ incident on the end-vertices of the $H'$ to obtain the cycle $H$. All relevant properties we have observed remain unchanged.

\endgroup 

\null

\noindent Section 2 contains algorithms we can use to reconfigure one Hamiltonian cycle (e-cycle) into another. Proofs of existence for the algorithms are in Section 3. Section 4 contains proofs of auxiliary results required in Section 3.

\subsection{Appendix}

\begingroup
\setlength{\intextsep}{0pt}
\setlength{\columnsep}{10pt}
\begin{wrapfigure}[]{r}{0cm}
\begin{adjustbox}{trim=0cm 1.5cm 0cm 1.25cm}
\begin{tikzpicture}[scale=1]
\usetikzlibrary{decorations.markings}

\begin{scope}[xshift=0cm] 

\begin{scope}
[very thick,decoration={
    markings,
    mark=at position 0.85 with {\arrow{>}}}
    ]
    \draw[postaction={decorate}] (1,-1)--(-0.5,2);
    \draw[postaction={decorate}] (1,-1)--(2.25,-0.5);
    \draw[postaction={decorate}] (2.25,-0.5)--(1.5,1);
    \draw[postaction={decorate}] (1,-1)--(1.5,1);
    \draw[postaction={decorate}] (0,1)--++(1.5,0.75);
\end{scope}

\draw[fill=black, opacity=1] (-0.5,2) circle [radius=0.05];
\draw[fill=black, opacity=1] (1,-1) circle [radius=0.05];
\draw[fill=black, opacity=1] (1.5,1) circle [radius=0.05];

\node[left] at  (-0.5,2) [scale=0.8]{\small${B}$};
\node[below] at  (1,-1) [scale=0.8]{\small${A}$};
\node[above] at  (1.5,1) [scale=0.8]{\small${P}$};

\node[right] at  (1.9,0.25) [scale=0.8]{\small{$\overrightarrow{d^{\parallel}}$}};

\node[below] at  (1.75,-0.7) [scale=0.8]{\small{$\overrightarrow{d^{\bot}}$}};

\node[left] at  (1.25,0) [scale=0.8]{\small{$\overrightarrow{d}$}};

\node[above] at  (0.75,1.4) [scale=0.8]{\small{$\overrightarrow{n}$}};


\node[below] at (1.25,-1.25) [scale=0.8]{\small{\begin{tabular}{c} Fig. 1.10. $P$ on the right of $\overrightarrow{AB}$. \end{tabular}}};;

\end{scope}

\end{tikzpicture}
\end{adjustbox}
\end{wrapfigure}

\textbf{A.1.} Let $A=(x_1,y_1)$, $B=(x_2,y_2)$ and $P=(x,y)$ be points in the plane that are not collinear. We define $(x_2-x_1, y_2-y_1)$ as the direction of the vector $\overrightarrow{AB}$. Then the direction of the normal $\overrightarrow{n}$ to $\overrightarrow{AB}$, obtained by rotating $\overrightarrow{AB}$ by $-\frac{\pi}{2}$, is $(y_2-y_1, x_1-x_2)$. We want to know whether the point $P$ is on the side of $\overrightarrow{AB}$ toward which $\overrightarrow{n}$ is pointing. Let $\overrightarrow{d}=\overrightarrow{AP}$. Let  $\overrightarrow{d^{\bot}}$ be the component of $\overrightarrow{d}$ that is perpendicular to $\overrightarrow{AB}$ and let $\overrightarrow{d^{\parallel}}$ be the component of $\overrightarrow{d}$ that is parallel to $\overrightarrow{AB}$. Note that:

\endgroup 

$\overrightarrow{d} \cdot \overrightarrow{n} = (\overrightarrow{d^{\parallel}} + \overrightarrow{d^{\bot}}) \cdot \overrightarrow{n}= \overrightarrow{d^{\bot}} \cdot \overrightarrow{n} 
= (x-x_1, y-y_1) \cdot (y_2-y_1, x_1-x_2)
= (x-x_1)(y_2-y_1)+(y-y_1)(x_1-x_2).$

\null

\noindent Point $P$ is on the \textit{right} of $\overrightarrow{AB}$ if $\overrightarrow{d} \cdot \overrightarrow{n} > 0$, and on the \textit{left} if $\overrightarrow{d} \cdot \overrightarrow{n} < 0$. Let $e=(u,v)$ be an edge of a grid graph $G$, where $u=v(k_1,l_1)$, $v=v(k_2,l_2)$. Let $X$ be a box of the square lattice that is incident on $(u,v)$. We say that $X$ is on the \textit{right} of the edge $(u,v)$ if there is a vertex $w=v(k,l)$ in $V(X)\setminus V(e)$ such that $(k-k_1, l_2-l_1) + (l-l_1, k_1-k_2)=1$ and we say that $X$ is on the \textit{left} of the edge $(u,v)$ if there is a vertex $w=v(k,l)$ in $V(X)\setminus V(e)$ such that $(k-k_1, l_2-l_1) + (l-l_1, k_1-k_2)=-1$.

\null 

\begingroup
\setlength{\intextsep}{0pt}
\setlength{\columnsep}{20pt}
\begin{wrapfigure}[]{r}{0cm}
\begin{adjustbox}{trim=0cm 0.25cm 0cm 0cm}
\begin{tikzpicture}[scale=1.25]

\begin{scope}[xshift=0cm] {
\draw[gray,very thin, step=0.5cm, opacity=0.5] (0,0) grid (1,1);

\draw[blue, line width=0.5mm] (0,0)--++(0,1)--++(0.5,0)--++(0,-1)--++(0.5,0)--++(0,1);

\draw[fill=blue, opacity=1] (0,0) circle [radius=0.035];
\draw[fill=blue, opacity=1] (1,1) circle [radius=0.035];
\draw[fill=blue, opacity=1] (0.5,0.5) circle [radius=0.035];
\draw[fill=blue, opacity=1] (0.5,0) circle [radius=0.035];

\node[right] at  (1,1) [scale=0.8]{\small{$v_r$}};
\node[left] at  (0,0) [scale=0.8]{\small{$v_1$}};
\node[right] at  (0.45,0.1) [scale=0.8]{\small{$v_s$}};
\node[right] at  (0.46,0.5) [scale=0.8]{\small{$v_{s{-}1}$}};


}\end{scope}

\begin{scope}[xshift=1.75cm] {
\draw[gray,very thin, step=0.5cm, opacity=0.5] (0,0) grid (1,1);

\draw[blue, line width=0.5mm] (0,0)--++(0,1)--++(0.5,0)--++(0,-1)--++(0.5,0)--++(0,1);

\draw[blue, line width=0.5mm] (0,0)--++(0.5,0);

\draw[fill=blue, opacity=1] (0,0) circle [radius=0.035];
\draw[fill=blue, opacity=1] (1,1) circle [radius=0.035];
\draw[fill=blue, opacity=1] (0.5,0.5) circle [radius=0.035];
\draw[fill=blue, opacity=1] (0.5,0) circle [radius=0.035];

\node[right] at  (1,1) [scale=0.8]{\small{$v_r$}};
\node[left] at  (0,0) [scale=0.8]{\small{$v_1$}};
\node[right] at  (0.45,0.1) [scale=0.8]{\small{$v_s$}};
\node[right] at  (0.46,0.5) [scale=0.8]{\small{$v_{s{-}1}$}};

\node[below] at (0.5,-.15) [scale=0.8]{\small{\begin{tabular}{c} Fig. I.11.  An illustration of a backbite move. \end{tabular}}};;

}\end{scope}

\begin{scope}[xshift=3.5cm] {
\draw[gray,very thin, step=0.5cm, opacity=0.5] (0,0) grid (1,1);

\draw[blue, line width=0.5mm] (0.5,0.5)--++(0,0.5)--++(-0.5,0)--++(0,-1)--++(1,0)--++(0,1);

\draw[fill=blue, opacity=1] (0,0) circle [radius=0.035];
\draw[fill=blue, opacity=1] (1,1) circle [radius=0.035];
\draw[fill=blue, opacity=1] (0.5,0.5) circle [radius=0.035];
\draw[fill=blue, opacity=1] (0.5,0) circle [radius=0.035];

\node[right] at  (1,1) [scale=0.8]{\small{$v_r$}};
\node[left] at  (0,0) [scale=0.8]{\small{$v_1$}};
\node[right] at  (0.45,0.1) [scale=0.8]{\small{$v_s$}};
\node[right] at  (0.46,0.5) [scale=0.8]{\small{$v_{s{-}1}$}};


}\end{scope}

\end{tikzpicture}
\end{adjustbox}
\end{wrapfigure}

\noindent \textbf{A.2.} Let $H$ be a Hamiltonian path $v_1, \ldots, v_r$ of an $m\times n$ grid graph $G$, and let $v_s$ be adjacent to $v_1$, $s\neq 2$. If we add the edge $\{v_1,v_s\}$, we obtain a cycle $v_1, \ldots, v_s, v_1$, and a path $v_s, \ldots, v_r$. Now, if we remove the edge \{$v_{s-1},v_s\}$, we obtain a new Hamiltonian path $H'=(H \setminus \{v_{s-1},v_s\}) \cup \{v_1,v_s\}$. This operation is called a backbite move. See Figure 1.11.

\section{Reconfiguration algorithm for cycles and canonical forms}

\textbf{Definitions.} Denote by $G_s$ the induced subgraph of $G$ on all the vertices with distance $s$ or greater from the boundary of $G$. Denote by $R_s$ the rectangular induced subgraph on vertices of $G$ with distance $s$ from the boundary. Then $R_s$ is the boundary of $G_s$ and the edges of $R_0$ are the boundary edges of $G$.

The main result stated in the Introduction is an immediate consequence of the slightly more general theorem below. Its proof takes up the remainder of the paper.

\null 

\noindent \textbf{Theorem 2.1.} Let $G$ be an $m\times n$ grid graph with $n \geq m$. Let $H$ and $K$ be two Hamiltonian cycles or Hamiltonian e-cycles of $G$ with the same edge $e$. Then there is a sequence of at most $n^2m$ valid double-switch moves that reconfigures $H$ into $K$.

\null

\noindent The $3 \times n$ and $4 \times n$ cases were done by Nishat in \cite{nishat2020reconfiguration}, so from here on, we will assume that $m,n \geq 5$, and that $m$ and $n$ are not both odd. First we will describe canonical forms for Hamiltonian cycles and e-cycles. Then we show that we can reconfigure any two canonical forms into one another. Then we show that any Hamiltonian cycle (e-cycle) can be reconfigured into a canonical form. Observing that double-switch moves are invertible completes the proof. That is $X \mapsto Y$ followed by $Y \mapsto X$ results in no net change. More specifically, suppose we want to reconfigure a Hamiltonian cycle (e-cycle) $H$ into a Hamiltonian cycle (e-cycle) $K$. Let $\mu_1, ..., \mu_k$ and $\nu_1,  ..., \nu_s$ be the sequences of moves that reconfigure $H$ and $K$ into the canonical forms $H_{can}$ and $K_{can}$, respectively. Then $\nu_s, \nu_{s-1}, ...,\nu_1$ reconfigures $K_{can}$ into $K$. Let $\eta_1, ..., \eta_t$ be the sequence of moves that reconfigures $H_{can}$ into $K_{can}$. Then the sequence of moves $\mu_1, ..., \mu_k, \eta_1, ..., \eta_t, \nu_s, \nu_{s-1}, ...,\nu_1 $ reconfigures $H$ into $K$.

\null

\noindent \textbf{Description of  canonical forms.} We shall write $\mathcal{H}_{\textrm{can}}(m,n)$ to denote the set of canonical forms of Hamiltonian cycles and e-cycles on an $m \times n$ grid graph. Then $H \in \mathcal{H}_{\textrm{can}}(m,n)$ if and only if $H$ can be constructed by the ``Canonical Form Builder"  algorithm described below.

Let $t=\Big\lfloor \frac{\min(m,n)-4}{2}  \Big\rfloor$. Let $k_1=|m-n|+2$ and $k_2=|m-n|+3$. If $\min(m,n)$ is even, let $D$ be the Hamiltonian cycle of the $2 \times k_1$ grid graph $G_{t+1}$. If $\min(m,n)$ is odd, let $D$ be any Hamiltonian cycle of the $3 \times k_2$ grid graph $G_{t+1}$. Let $U = D \cup  \bigcup_{i=0}^t R_i$.

\begin{center}
\begin{tikzpicture}[scale=0.7]

\draw[gray,very thin, step=0.5cm, opacity=0.5] (0,0) grid (11,5.5);

\draw[gray,very thin, step=0.5cm, opacity=0.5] (11.99,0) grid (20,6.5);

\begin{scope}[xshift=-1cm, yshift=-1cm]
{
\foreach \x in {} \fill[blue!40!white,opacity=0.5]
(1.5, 1+0.5*\x) rectangle  (2, 1+0.5*\x+0.5); 
\foreach \x in {} \fill[blue!40!white,opacity=0.5]
(2, 1+0.5*\x) rectangle  (2.5, 1+0.5*\x+0.5); 
\foreach \x in {} \fill[blue!40!white,opacity=0.5]
(2.5, 1+0.5*\x) rectangle  (3, 1+0.5*\x+0.5); 
\foreach \x in {1,3} \fill[red!40!white,opacity=0.5]
(3, 1+0.5*\x) rectangle  (3.5, 1+0.5*\x+0.5); 
\foreach \x in {8} \fill[red!40!white,opacity=0.5]
(3.5, 1+0.5*\x) rectangle  (4, 1+0.5*\x+0.5); 
\foreach \x in {5} \fill[red!40!white,opacity=0.5]
(4, 1+0.5*\x) rectangle  (4.5, 1+0.5*\x+0.5); 
\foreach \x in {} \fill[blue!40!white,opacity=0.5]
(4.5, 1+0.5*\x) rectangle  (5, 1+0.5*\x+0.5); 
\foreach \x in {} \fill[blue!40!white,opacity=0.5]
(5, 1+0.5*\x) rectangle  (5.5, 1+0.5*\x+0.5); 
\foreach \x in {} \fill[blue!40!white,opacity=0.5]
(5.5, 1+0.5*\x) rectangle  (6, 1+0.5*\x+0.5); 
}
{
\foreach \x in {1,2,3,4,5,6,7,8,9} \draw[blue, line width=0.5mm] (1+0.5*\x,1.5)--(1+0.5*\x+0.5,1.5); 
\foreach \x in {2,3,4,5,6,7,8} \draw[blue, line width=0.5mm] (1+0.5*\x,2)--(1+0.5*\x+0.5,2);
\foreach \x in {3,4,5,6,7} \draw[blue, line width=0.5mm] (1+0.5*\x,2.5)--(1+0.5*\x+0.5,2.5);
\foreach \x in {4,5,6} \draw[blue, line width=0.5mm] (1+0.5*\x,3)--(1+0.5*\x+0.5,3);
\foreach \x in {5} \draw[blue, line width=0.5mm] (1+0.5*\x,3.5)--(1+0.5*\x+0.5,3.5);
\foreach \x in {5} \draw[blue, line width=0.5mm] (1+0.5*\x,4)--(1+0.5*\x+0.5,4);
\foreach \x in {4,5,6} \draw[blue, line width=0.5mm] (1+0.5*\x,4.5)--(1+0.5*\x+0.5,4.5);
\foreach \x in {3,4,5,6,7} \draw[blue, line width=0.5mm]  (1+0.5*\x,5)--(1+0.5*\x+0.5,5);
\foreach \x in {2,3,4,5,6,7,8} \draw[blue, line width=0.5mm] (1+0.5*\x,5.5)--(1+0.5*\x+0.5,5.5);
\foreach \x in {1,2,3,4,5,6,7,8,9} \draw[blue, line width=0.5mm] (1+0.5*\x,6)--(1+0.5*\x+0.5,6);
}
{
\foreach \x in {1,2,3,4,5,6,7,8,9}  \draw[blue, line width=0.5mm] 
(1.5, 1+0.5*\x)--(1.5, 1+0.5*\x+0.5); 
\foreach \x in {2,3,4,5,6,7,8} \draw[blue, line width=0.5mm] 
(2, 1+0.5*\x)--(2, 1+0.5*\x+0.5); 
\foreach \x in {3,4,5,6,7} \draw[blue, line width=0.5mm] 
(2.5, 1+0.5*\x)--(2.5, 1+0.5*\x+0.5);
\foreach \x in {4,5,6} \draw[blue, line width=0.5mm] 
(3, 1+0.5*\x)--(3, 1+0.5*\x+0.5);
\foreach \x in {5} \draw[blue, line width=0.5mm] 
(3.5, 1+0.5*\x)--(3.5, 1+0.5*\x+0.5);
\foreach \x in {5} \draw[blue, line width=0.5mm] 
(4, 1+0.5*\x)--(4, 1+0.5*\x+0.5);
\foreach \x in {4,5,6} \draw[blue, line width=0.5mm] 
(4.5, 1+0.5*\x)--(4.5, 1+0.5*\x+0.5);
\foreach \x in {3,4,5,6,7} \draw[blue, line width=0.5mm] 
(5, 1+0.5*\x)--(5, 1+0.5*\x+0.5);
\foreach \x in {2,3,4,5,6,7,8} \draw[blue, line width=0.5mm] 
(5.5, 1+0.5*\x)--(5.5, 1+0.5*\x+0.5);
\foreach \x in {1,2,3,4,5,6,7,8,9} \draw[blue, line width=0.5mm] 
(6, 1+0.5*\x)--(6, 1+0.5*\x+0.5);
}

\node[below] at (3.75,0.75) [scale=0.75]
{\small{\begin{tabular}{c} Fig. 2.1(a). $U$  with $\min(m,n)$ even.  \\ $D_1$ in dashed yellow. \end{tabular}}};;

\node at (3.25,1.75) [scale=0.75]
{\small{$X_0$}};
\node at (3.75,5.25) [scale=0.75]
{\small{$X_1$}};
\node at (3.25,2.75) [scale=0.75]
{\small{$X_2$}};
\node at (4.25,3.75) [scale=0.75]
{\small{$X_3$}};

\end{scope}

\begin{scope}[xshift=4.5cm, yshift=-1cm]
{
\foreach \x in {} \fill[blue!40!white,opacity=0.5]
(1.5, 1+0.5*\x) rectangle  (2, 1+0.5*\x+0.5); 
\foreach \x in {} \fill[blue!40!white,opacity=0.5]
(2, 1+0.5*\x) rectangle  (2.5, 1+0.5*\x+0.5); 
\foreach \x in {} \fill[blue!40!white,opacity=0.5]
(2.5, 1+0.5*\x) rectangle  (3, 1+0.5*\x+0.5); 
\foreach \x in {1,3} \fill[red!40!white,opacity=0.5]
(3, 1+0.5*\x) rectangle  (3.5, 1+0.5*\x+0.5); 
\foreach \x in {8} \fill[red!40!white,opacity=0.5]
(3.5, 1+0.5*\x) rectangle  (4, 1+0.5*\x+0.5); 
\foreach \x in {5} \fill[red!40!white,opacity=0.5]
(4, 1+0.5*\x) rectangle  (4.5, 1+0.5*\x+0.5); 
\foreach \x in {} \fill[blue!40!white,opacity=0.5]
(4.5, 1+0.5*\x) rectangle  (5, 1+0.5*\x+0.5); 
\foreach \x in {} \fill[blue!40!white,opacity=0.5]
(5, 1+0.5*\x) rectangle  (5.5, 1+0.5*\x+0.5); 
\foreach \x in {} \fill[blue!40!white,opacity=0.5]
(5.5, 1+0.5*\x) rectangle  (6, 1+0.5*\x+0.5); 
}
{
\foreach \x in {1,2,3,5,6,7,8,9} \draw[blue, line width=0.5mm] (1+0.5*\x,1.5)--(1+0.5*\x+0.5,1.5); 
\foreach \x in {2,3,5,6,7,8} \draw[blue, line width=0.5mm] (1+0.5*\x,2)--(1+0.5*\x+0.5,2);
\foreach \x in {3,5,6,7} \draw[blue, line width=0.5mm] (1+0.5*\x,2.5)--(1+0.5*\x+0.5,2.5);
\foreach \x in {5,6} \draw[blue, line width=0.5mm] (1+0.5*\x,3)--(1+0.5*\x+0.5,3);
\foreach \x in {5,6} \draw[blue, line width=0.5mm] (1+0.5*\x,3.5)--(1+0.5*\x+0.5,3.5);
\foreach \x in {5,6} \draw[blue, line width=0.5mm] (1+0.5*\x,4)--(1+0.5*\x+0.5,4);
\foreach \x in {4,5,6} \draw[blue, line width=0.5mm] (1+0.5*\x,4.5)--(1+0.5*\x+0.5,4.5);
\foreach \x in {3,4,6,7} \draw[blue, line width=0.5mm]  (1+0.5*\x,5)--(1+0.5*\x+0.5,5);
\foreach \x in {2,3,4,6,7,8} \draw[blue, line width=0.5mm] (1+0.5*\x,5.5)--(1+0.5*\x+0.5,5.5);
\foreach \x in {1,2,3,4,5,6,7,8,9} \draw[blue, line width=0.5mm] (1+0.5*\x,6)--(1+0.5*\x+0.5,6);
}
{
\foreach \x in {1,2,3,4,5,6,7,8,9} \draw[blue, line width=0.5mm] 
(1.5, 1+0.5*\x)--(1.5, 1+0.5*\x+0.5); 
\foreach \x in {2,3,4,5,6,7,8} \draw[blue, line width=0.5mm] 
(2, 1+0.5*\x)--(2, 1+0.5*\x+0.5); 
\foreach \x in {3,4,5,6,7} \draw[blue, line width=0.5mm] 
(2.5, 1+0.5*\x)--(2.5, 1+0.5*\x+0.5);
\foreach \x in {1,3,4,5,6} \draw[blue, line width=0.5mm] 
(3, 1+0.5*\x)--(3, 1+0.5*\x+0.5);
\foreach \x in {1,3,5,8} \draw[blue, line width=0.5mm] 
(3.5, 1+0.5*\x)--(3.5, 1+0.5*\x+0.5);
\foreach \x in {8} \draw[blue, line width=0.5mm] 
(4, 1+0.5*\x)--(4, 1+0.5*\x+0.5);
\foreach \x in {4,6} \draw[blue, line width=0.5mm] 
(4.5, 1+0.5*\x)--(4.5, 1+0.5*\x+0.5);
\foreach \x in {3,4,5,6,7} \draw[blue, line width=0.5mm] 
(5, 1+0.5*\x)--(5, 1+0.5*\x+0.5);
\foreach \x in {2,3,4,5,6,7,8} \draw[blue, line width=0.5mm] 
(5.5, 1+0.5*\x)--(5.5, 1+0.5*\x+0.5);
\foreach \x in {1,2,3,4,5,6,7,8,9} \draw[blue, line width=0.5mm] 
(6, 1+0.5*\x)--(6, 1+0.5*\x+0.5);
}

\node[below] at (3.75,0.75) [scale=0.75]
{\small{\begin{tabular}{c} Fig. 2.1(b),  \\ A canonical form from U. \end{tabular}}};;

\node at (3.25,1.75) [scale=0.75]
{\small{$X_0$}};
\node at (3.75,5.25) [scale=0.75]
{\small{$X_1$}};
\node at (3.25,2.75) [scale=0.75]
{\small{$X_2$}};
\node at (4.25,3.75) [scale=0.75]
{\small{$X_3$}};

\draw[yellow, dashed, line width =0.3mm] (3.5,3.5)--++(0.5,0)--++(0,0.5)--++(-0.5,0)--++(0,-0.5);

\end{scope}

\begin{scope}[xshift=11cm, yshift=-1cm]

{\foreach \x in {6} \fill[red!40!white,opacity=0.5]
(1.5, 1+0.5*\x) rectangle  (2, 1+0.5*\x+0.5); 
\foreach \x in {10} \fill[red!40!white,opacity=0.5]
(2.5, 1+0.5*\x) rectangle  (3, 1+0.5*\x+0.5); 
}
{
\foreach \x in {1,2,3,4,5,6} \draw[blue, line width=0.5mm] (1+0.5*\x,1.5)--(1+0.5*\x+0.5,1.5); 
\foreach \x in {2,3,4,5} \draw[blue, line width=0.5mm] (1+0.5*\x,2)--(1+0.5*\x+0.5,2);
\foreach \x in {3,4} \draw[blue, line width=0.5mm] (1+0.5*\x,2.5)--(1+0.5*\x+0.5,2.5);

\foreach \x in {3,4} \draw[blue, line width=0.5mm] (1+0.5*\x,6)--(1+0.5*\x+0.5,6);
\foreach \x in {2,3,4,5} \draw[blue, line width=0.5mm] (1+0.5*\x,6.5)--(1+0.5*\x+0.5,6.5);
\foreach \x in {1,2,3,4,5,6} \draw[blue, line width=0.5mm] (1+0.5*\x,7)--(1+0.5*\x+0.5,7);
}
{
\foreach \x in {1,2,3,4,5,6,7,8,9,10,11}  \draw[blue, line width=0.5mm] 
(1.5, 1+0.5*\x)--(1.5, 1+0.5*\x+0.5); 
\foreach \x in {2,3,4,5,6,7,8,9,10} \draw[blue, line width=0.5mm] 
(2, 1+0.5*\x)--(2, 1+0.5*\x+0.5); 

\draw[blue, line width=0.5mm] (2.5,2.5)--++(0,0.5)--++(0.5,0)--++(0,0.5)--++(-0.5,0)--++(0,0.5)--++(0.5,0)--++(0,0.5)--++(-0.5,0)--++(0,1.5);

\draw[blue, line width=0.5mm] (3.5,2.5)--++(0,2.5)--++(-0.5,0)--++(0,0.5)--++(0.5,0)--++(0,0.5);

\foreach \x in {2,3,4,5,6,7,8,9,10} \draw[blue, line width=0.5mm] 
(4, 1+0.5*\x)--(4, 1+0.5*\x+0.5);
\foreach \x in {1,2,3,4,5,6,7,8,9,10,11}  \draw[blue, line width=0.5mm] 
(4.5, 1+0.5*\x)--(4.5, 1+0.5*\x+0.5);

}

\node[below] at (3,0.75) [scale=0.75]
{\small{\begin{tabular}{c} Fig. 2.2(a). $U$ with  \\$\min(m,n)$ odd. \end{tabular}}};;

\node at (1.75,4.25) [scale=0.75]
{\small{$X_0$}};
\node at (2.75,6.25) [scale=0.75]
{\small{$X_1$}};

\end{scope}

\begin{scope}[xshift=15cm, yshift=-1cm]

{\foreach \x in {6} \fill[red!40!white,opacity=0.5]
(1.5, 1+0.5*\x) rectangle  (2, 1+0.5*\x+0.5); 
\foreach \x in {10} \fill[red!40!white,opacity=0.5]
(2.5, 1+0.5*\x) rectangle  (3, 1+0.5*\x+0.5); 
}
{
\foreach \x in {1,2,3,4,5,6} \draw[blue, line width=0.5mm] (1+0.5*\x,1.5)--(1+0.5*\x+0.5,1.5); 
\foreach \x in {2,3,4,5} \draw[blue, line width=0.5mm] (1+0.5*\x,2)--(1+0.5*\x+0.5,2);
\foreach \x in {3,4} \draw[blue, line width=0.5mm] (1+0.5*\x,2.5)--(1+0.5*\x+0.5,2.5);

\foreach \x in {1} \draw[blue, line width=0.5mm] (1+0.5*\x,4)--(1+0.5*\x+0.5,4);

\foreach \x in {1} \draw[blue, line width=0.5mm] (1+0.5*\x,4.5)--(1+0.5*\x+0.5,4.5);

\foreach \x in {4} \draw[blue, line width=0.5mm] (1+0.5*\x,6)--(1+0.5*\x+0.5,6);
\foreach \x in {2,4,5} \draw[blue, line width=0.5mm] (1+0.5*\x,6.5)--(1+0.5*\x+0.5,6.5);
\foreach \x in {1,2,3,4,5,6} \draw[blue, line width=0.5mm] (1+0.5*\x,7)--(1+0.5*\x+0.5,7);
}
{
\foreach \x in {1,2,3,4,5,7,8,9,10,11}  \draw[blue, line width=0.5mm] 
(1.5, 1+0.5*\x)--(1.5, 1+0.5*\x+0.5); 
\foreach \x in {2,3,4,5,7,8,9,10} \draw[blue, line width=0.5mm] 
(2, 1+0.5*\x)--(2, 1+0.5*\x+0.5); 

\foreach \x in {10} \draw[blue, line width=0.5mm] 
(2.5, 1+0.5*\x)--(2.5, 1+0.5*\x+0.5); 

\foreach \x in {10} \draw[blue, line width=0.5mm] 
(3, 1+0.5*\x)--(3, 1+0.5*\x+0.5); 

\draw[blue, line width=0.5mm] (2.5,2.5)--++(0,0.5)--++(0.5,0)--++(0,0.5)--++(-0.5,0)--++(0,0.5)--++(0.5,0)--++(0,0.5)--++(-0.5,0)--++(0,1.5);

\draw[blue, line width=0.5mm] (3.5,2.5)--++(0,2.5)--++(-0.5,0)--++(0,0.5)--++(0.5,0)--++(0,0.5);

\foreach \x in {2,3,4,5,6,7,8,9,10} \draw[blue, line width=0.5mm] 
(4, 1+0.5*\x)--(4, 1+0.5*\x+0.5);
\foreach \x in {1,2,3,4,5,6,7,8,9,10,11}  \draw[blue, line width=0.5mm] 
(4.5, 1+0.5*\x)--(4.5, 1+0.5*\x+0.5);

}
\node[below] at (3,0.75) [scale=0.75]
{\small{\begin{tabular}{c} Fig. 2.2(b). A canonical \\  form from $U$. $D_2(H)$ in \\ dashed yellow. \end{tabular}}};;

\node at (1.75,4.25) [scale=0.75]
{\small{$X_0$}};
\node at (2.75,6.25) [scale=0.75]
{\small{$X_1$}};

\draw[yellow, dashed, line width =0.3mm] (2.5,2.5)--++(1,0)--++(0,3.5)--++(-1,0)--++(0,-3.5);
\end{scope}
\end{tikzpicture}
\end{center}

\noindent We define $(R_i,R_{i+1})$ to be the set of all the boxes of $G$ adjacent to both $R_i$ and $R_{i+1}$. Now we can state the Canonical Form Builder algorithm (CFB) that takes as inputs $m$ and $n$ and  outputs an element of $\mathcal{H}_{\textrm{can}}(m,n)$.

\begin{itemize}
    
    \item [Step 1.] Set $i=0$. Switch one of the $2(m-3)+2(n-3)$ switchable boxes of $(R_0,R_1)$ of the graph $U$. This switch removes some edge, say $e_1$, from $E(R_1)$. If $t=0$, stop. If $t>0$, go to Step 2.
    
    \item [Step 2.] Increase $i$ by $1$. Switch one of the switchable boxes of $(R_i,R_{i+1}) \setminus e_i$.

    \item [Step 3.] If $i \leq t$, go to Step 2. If $i=t+1$, then stop.

    We have arrived at a canonical form $H$. Record the switched boxes $X_0, ..., X_t$ in a list \textit{List(H)}. So $\textrm{List(H)}=(X_0, ..., X_t)$  consists of the faces of $G$ that were chosen to be switched to make $U$ into a canonical form, listed in order.

\end{itemize}

\noindent We observe that the CFB algorithm above works just as well for e-cycles if we remove $e$ from $U$.

\null 

\noindent \textbf{Reconfiguration between canonical forms.} Let $H,K \in \mathcal{H}_{\textrm{can}}(m,n)$. Let $\textrm{List(H)}=(X_0,...,X_t)$ and $\textrm{List(K)}=(Y_0,..., Y_t)$ be the switched boxes of $H$ and $K$ respectively. We will reconfigure $H$ into $K$, so the algorithm will run on $H$. 

Let $D(H)=H \cap G_{t+1}$ and note that $D(H)$ is an e-cycle of $G_{t+1}$. Using the result of Nishat in \cite{nishat2020reconfiguration}, $D(H)$ can be reconfigured into $D(K)$ by a sequence of valid moves. 

\begin{itemize}
    \item [Step 1.] Set $i=0$. If $t=0$, go to Step 3. If $t>0$, go to Step 2.

    \item [Step 2.] If $Y_i$ is switchable after switching $X_i$, switch both $X_i$ and $Y_i$. 
    
    If $Y_i$ is not switchable after switching $X_i$, switch $X_{i+1}$ and any other switchable box in $(R_{i+1}, R_{i+2})$, say $X_{i+1}'$, such that $X_{i+1} \mapsto X_{i+1}'$ is valid. We remark that the only the only case where $Y_i$ would not be switchable after switching $X_i$ occurs when $Y_i$ is adjacent to $X_{i+1}$. Note that there are many possible choices for $X_{i+1}'$. Now $Y_i$ is switchable. Switch both $X_i$ and $Y_i$. Update $\textrm{List(H)}$ by setting the $(i+1)^{st}$ slot to  $X_{i+1}'$. Increase $i$ by $1$. 

    \item [Step 3.] If $i<t$, go to Step 2. 
    
    If $i=t$ and $\min(m,n)$ is even, switch $X_i$ and $Y_i$, and then stop.

    If $i=t$ and $\min(m,n)$ is odd, go to Step 3.1.

     \item [Step 3.1.] Switch $X_t$ and any one of the four switchable boxes, say $X_t'$, located on the short sides of $D$. Run NRI's algorithm to reconfigure $D(H)$ into $D(K)$. Switch $X_t'$ and $Y_t$. Stop.

\end{itemize}

\noindent \textbf{Reconfiguration of a cycle into a canonical form (RtCF).} The RtCF algorithm takes as input a Hamiltonian e-cycle and outputs a canonical e-cycle. We will need the following proposition:

\null

\noindent \textbf{Proposition 2.2.}
Let $H \in \mathcal{H}$.

(a) If $H$ has more than one large cookie, then there is a cascade of length at most two that reduces

\hspace*{13pt}  the number of large cookies of $H$ by one. This is the ManyLargeCookies (MLC) algorithm.

(b) If $H$ has exactly one large cookie and at least one small cookie, then there is a cascade of length 

\hspace*{13pt} at most $\frac{1}{2}\max(m,n)+\min(m,n)+2$ that reduces the number of small cookies of $H$ by one and 

\hspace*{13pt} such that it does not increase the number of large cookies. This is the OneLargeCookie (1LC) 

\hspace*{13pt} algorithm.

\null

\noindent The proof for Proposition 2.2 will be given in Section 3.

\null 

\noindent Now we can describe the RtCF algorithm. Without loss of generality, assume $H \in \mathcal{H}$ is a Hamiltonian e-cycle of $G$. Suppose $H=H_0$ has $c_{1;\textrm{large}}$ large cookies and $c_{1;\textrm{small}}$ small cookies. We run MLC $c_{1;\textrm{large}}-1$ times and then run 1LC $c_{1;\textrm{small}}$ times to reconfigure $H_0$ into $H_0'$, where $H_0'$ has exactly one (necessarily large) cookie $C_1$. We define  $H_1=(G_1 \cap H_0')$ and observe that $H_1$ is a Hamiltonian e-cycle of $G_1$. This is the first iteration of (RtCF). Now we describe the $j^{\textrm{th}}$ iteration. We run MLC $c_{j;\textrm{large}}-1$ times and then run 1LC $c_{j;\textrm{small}}$ times to reconfigure $H_{j-1}$ into $H_{j-1}'$, where $H_{j-1}'$ has exactly one (necessarily large) cookie $C_j$. The RtCF algorithm stops when $j=t$. We give a summary of the algorithm below.

\begin{itemize}
    \item [Step 1.] Set $j=0$. Run MLC $c_{1;\textrm{large}}$ times and then 1LC $c_{1;\textrm{small}}$ times on $H_0$ to reconfigure $H_0$ into $H_0'$.
    
    \item [Step 2.] Increase $j$ by $1$. Set $H_j=G_j \cap H'_{j-1}$ and note that $H_j$ is a Hamiltonian e-cycle in $G_j$. Run MLC $c_{j+1;\textrm{large}}$ times and then 1LC $c_{j+1;\textrm{small}}$ times on $H_{j}$ to reconfigure into $H_{j}$ into $H_{j}'$. 

    \item [Step 3.] If $j<t$, go to Step 2. If $j=t$, stop. 

\end{itemize}

\noindent \textit{Proof of the RtCF algorithm}. Let $N_j$ be the neck of the only cookie $C_j$ of $H_{j-1}'$ in $G_{j-1}$. Define  $e_1(N_j)=N_j \cap R_{j-1}$, $e_2(N_j)=N_j \cap R_{j}$ and $\{ e_3(N_j),e_4(N_j)\} = N_j \cap H_{j-1}'$. We observe that when the RtCF algorithm stops, we have reconfigured $H$ into

$$H_c=D(H) \cup \bigcup_{j=0}^t \big(R_j \cup e_3(N_{j+1}) \cup e_4(N_{j+1})\big) \setminus (e_1(N_{j+1}) \cup e_2(N_{j+1}).$$

\noindent Now we can see that $H_c$ is an element of  $\mathcal{H}_{\textrm{can}}(m,n)$  by setting $X_{j-1}=N_j$ for $j=1,2,...,t+1$ and running CFB. See  Figure 2.3 on page 12 for an illustration of the RtCF algorithm on a $10 \times 10$ grid.

\null

\noindent \textbf{Bound for Theorem 2.1.} Recall that $n \geq m$. Note that it takes at most $2m$ moves to reconfigure canonical forms into one another. Now we count the moves required for RtCF to terminate. Observe that for each $j \in \{0, \ldots, t-1\}$, $H_j$ has at most $2 \big(\frac{n-2j}{2}+ \frac{m-2j}{2}\big)=n+m-4j$ cookies. This is the number of iterations of MLC or 1LC required for each $j$. It will follow from the proofs in Sections 3 and 4 that each application of MLC or 1LC in $H_j$ requires at most $\frac{1}{2}n+m-3j+3$ moves. So, RtCF requires at most:

\begin{align*}
 & \sum_{j=1}^{\lfloor (m-2)/2 \rfloor} (n+m-4j)\Big(\frac{n}{2}+m-3j+2\Big) \\
  &= \sum_{j=1}^{\lfloor (m-2)/2 \rfloor} \Big(12j^2 + (-5n - 7m - 8)j + \frac{n^2}{2} + \frac{3nm}{2} + 2n + m^2 + 2m\Big) \\
  &\leq \tfrac{1}{2}(m^3 - 3m^2 + 2m) 
      + (-5n - 7m - 8)\Big(\frac{m^2}{8} - \frac{m}{4}\Big) 
      + \tfrac{1}{2}\Big(\frac{n^2}{2} + \frac{3nm}{2} + 2n + m^2 + 2m\Big)(m-2) \\
  &= \frac{n^2m}{4} + \frac{nm^2}{8} + \frac{3nm}{4} + \frac{m^3}{8} - \frac{3m^2}{4} - \frac{n^2}{2} + m - 2n \\
  &= \frac{n^2m}{2} + \frac{nm}{4}\Big(3 - \frac{3m}{n} - \frac{2n}{m}\Big) + m - 2n.
\end{align*}

\noindent Let $x=\frac{m}{n}$. Then $\frac{3m}{n}+\frac{2n}{m}=3x+\frac{2}{x}$. Using calculus, we find that it attains a minimum of $2\sqrt{6}$ at $x=\frac{\sqrt{6}}{3}$. Then $\big( 3-\frac{3m}{n}-\frac{2n}{m} \big)$ can be at most $3-2\sqrt{6} \leq -1$. It follows that RtCF requires at most $\frac{n^2m}{2} -\frac{nm}{4}+m-2n$ to terminate. For a complete reconfiguration we need to run RtCF once for each e-cycle and reconfigure the resulting canonical forms. So, we need at most $2 \big( \frac{n^2m}{2}-\frac{nm}{4}+m-2n \big) +m=n^2m-\frac{nm}{2}-4n+3m < n^2m$ moves for a complete reconfiguration. We remark that this is a worst case scenario and conjecture that the typical number of moves required is of the order of $n^2$.

\null 

\begin{tikzpicture}[scale=0.7]

\begin{scope}[xshift=0cm]
{
\draw[gray,very thin, step=0.5cm, opacity=0.5] (0,0) grid (5.5,5.5);

\begin{scope}[xshift=-1cm,yshift=-1cm]

\node[below] at (3.75,1) [scale=1]
{\small{Initial configuration}};

{
\foreach \x in {} \fill[blue!40!white,opacity=0.5]
(1.5, 1+0.5*\x) rectangle  (2, 1+0.5*\x+0.5); 
\foreach \x in {2,3,4,5,6,7,8} \fill[blue!40!white,opacity=0.5]
(2, 1+0.5*\x) rectangle  (2.5, 1+0.5*\x+0.5); 
\foreach \x in {5,8} \fill[blue!40!white,opacity=0.5]
(2.5, 1+0.5*\x) rectangle  (3, 1+0.5*\x+0.5); 
\foreach \x in {1,2,3,4,5,6,8} \fill[blue!40!white,opacity=0.5]
(3, 1+0.5*\x) rectangle  (3.5, 1+0.5*\x+0.5); 
\foreach \x in {6} \fill[blue!40!white,opacity=0.5]
(3.5, 1+0.5*\x) rectangle  (4, 1+0.5*\x+0.5); 
\foreach \x in {2,4,6,8} \fill[blue!40!white,opacity=0.5]
(4, 1+0.5*\x) rectangle  (4.5, 1+0.5*\x+0.5); 
\foreach \x in {1,2,4,6,8,9} \fill[blue!40!white,opacity=0.5]
(4.5, 1+0.5*\x) rectangle  (5, 1+0.5*\x+0.5); 
\foreach \x in {2,3,4,8} \fill[blue!40!white,opacity=0.5]
(5, 1+0.5*\x) rectangle  (5.5, 1+0.5*\x+0.5); 
\foreach \x in {6} \fill[blue!40!white,opacity=0.5]
(5.5, 1+0.5*\x) rectangle  (6, 1+0.5*\x+0.5); 
}

{
\foreach \x in {1,2,3,5,6,8,9} \draw[blue, line width=0.5mm] (1+0.5*\x,1.5)--(1+0.5*\x+0.5,1.5); 
\foreach \x in {2,6,8} \draw[blue, line width=0.5mm] (1+0.5*\x,2)--(1+0.5*\x+0.5,2);
\foreach \x in {6,7} \draw[blue, line width=0.5mm] (1+0.5*\x,2.5)--(1+0.5*\x+0.5,2.5);
\foreach \x in {6,7} \draw[blue, line width=0.5mm] (1+0.5*\x,3)--(1+0.5*\x+0.5,3);
\foreach \x in {3,6,7,8} \draw[blue, line width=0.5mm] (1+0.5*\x,3.5)--(1+0.5*\x+0.5,3.5);
\foreach \x in {3,5,6,7,9} \draw[blue, line width=0.5mm] (1+0.5*\x,4)--(1+0.5*\x+0.5,4);
\foreach \x in {4,5,6,7,9} \draw[blue, line width=0.5mm] (1+0.5*\x,4.5)--(1+0.5*\x+0.5,4.5);
\foreach \x in {3,4,6,7,8} \draw[blue, line width=0.5mm]  (1+0.5*\x,5)--(1+0.5*\x+0.5,5);
\foreach \x in {2,3,4,6,8} \draw[blue, line width=0.5mm] (1+0.5*\x,5.5)--(1+0.5*\x+0.5,5.5);
\foreach \x in {1,2,3,4,5,6,8,9} \draw[blue, line width=0.5mm] (1+0.5*\x,6)--(1+0.5*\x+0.5,6);
}
{
\foreach \x in {1,2,3,4,5,6,7,8,9} \draw[blue, line width=0.5mm] 
(1.5, 1+0.5*\x)--(1.5, 1+0.5*\x+0.5); 
\foreach \x in {2,3,4,5,6,7,8} \draw[blue, line width=0.5mm] 
(2, 1+0.5*\x)--(2, 1+0.5*\x+0.5); 
\foreach \x in {2,3,4,6,7} \draw[blue, line width=0.5mm] 
(2.5, 1+0.5*\x)--(2.5, 1+0.5*\x+0.5);
\foreach \x in {1,2,3,4,6} \draw[blue, line width=0.5mm] 
(3, 1+0.5*\x)--(3, 1+0.5*\x+0.5);
\foreach \x in {1,2,3,4,5,8} \draw[blue, line width=0.5mm] 
(3.5, 1+0.5*\x)--(3.5, 1+0.5*\x+0.5);
\foreach \x in {2,4,8} \draw[blue, line width=0.5mm] 
(4, 1+0.5*\x)--(4, 1+0.5*\x+0.5);
\foreach \x in {1,9} \draw[blue, line width=0.5mm] 
(4.5, 1+0.5*\x)--(4.5, 1+0.5*\x+0.5);
\foreach \x in {1,3,6,9} \draw[blue, line width=0.5mm] 
(5, 1+0.5*\x)--(5, 1+0.5*\x+0.5);
\foreach \x in {2,3,4,6,8} \draw[blue, line width=0.5mm] 
(5.5, 1+0.5*\x)--(5.5, 1+0.5*\x+0.5);
\foreach \x in {1,2,3,4,5,7,8,9} \draw[blue, line width=0.5mm] 
(6, 1+0.5*\x)--(6, 1+0.5*\x+0.5);
}

{
\foreach \x in {7} \fill[red!60!white,opacity=0.5]
(4.5, 1+0.5*\x) rectangle  (5, 1+0.5*\x+0.5); 

\foreach \x in {9} \fill[red!60!white,opacity=0.5]
(4.5, 1+0.5*\x) rectangle  (5, 1+0.5*\x+0.5); 
}

\end{scope}
}
\end{scope}

\begin{scope}[xshift=7cm]
{
\draw[gray,very thin, step=0.5cm, opacity=0.5] (0,0) grid (5.5,5.5);

\begin{scope}[xshift=-1cm, yshift=-1cm]

\node[below] at (3.75,1) [scale=1]
{\small{First move}};

{
\foreach \x in {} \fill[blue!40!white,opacity=0.5]
(1.5, 1+0.5*\x) rectangle  (2, 1+0.5*\x+0.5); 
\foreach \x in {2,3,4,5,6,7,8} \fill[blue!40!white,opacity=0.5]
(2, 1+0.5*\x) rectangle  (2.5, 1+0.5*\x+0.5); 
\foreach \x in {5,8} \fill[blue!40!white,opacity=0.5]
(2.5, 1+0.5*\x) rectangle  (3, 1+0.5*\x+0.5); 
\foreach \x in {1,2,3,4,5,6,8} \fill[blue!40!white,opacity=0.5]
(3, 1+0.5*\x) rectangle  (3.5, 1+0.5*\x+0.5); 
\foreach \x in {6} \fill[blue!40!white,opacity=0.5]
(3.5, 1+0.5*\x) rectangle  (4, 1+0.5*\x+0.5); 
\foreach \x in {2,4,6,8} \fill[blue!40!white,opacity=0.5]
(4, 1+0.5*\x) rectangle  (4.5, 1+0.5*\x+0.5); 
\foreach \x in {1,2,4,6,7,8} \fill[blue!40!white,opacity=0.5]
(4.5, 1+0.5*\x) rectangle  (5, 1+0.5*\x+0.5); 
\foreach \x in {2,3,4,8} \fill[blue!40!white,opacity=0.5]
(5, 1+0.5*\x) rectangle  (5.5, 1+0.5*\x+0.5); 
\foreach \x in {6} \fill[blue!40!white,opacity=0.5]
(5.5, 1+0.5*\x) rectangle  (6, 1+0.5*\x+0.5); 
}
{
\foreach \x in {1,2,3,5,6,8,9} \draw[blue, line width=0.5mm] (1+0.5*\x,1.5)--(1+0.5*\x+0.5,1.5); 
\foreach \x in {2,6,8} \draw[blue, line width=0.5mm] (1+0.5*\x,2)--(1+0.5*\x+0.5,2);
\foreach \x in {6,7} \draw[blue, line width=0.5mm] (1+0.5*\x,2.5)--(1+0.5*\x+0.5,2.5);
\foreach \x in {6,7} \draw[blue, line width=0.5mm] (1+0.5*\x,3)--(1+0.5*\x+0.5,3);
\foreach \x in {3,6,7,8} \draw[blue, line width=0.5mm] (1+0.5*\x,3.5)--(1+0.5*\x+0.5,3.5);
\foreach \x in {3,5,6,7,9} \draw[blue, line width=0.5mm] (1+0.5*\x,4)--(1+0.5*\x+0.5,4);
\foreach \x in {4,5,6,9} \draw[blue, line width=0.5mm] (1+0.5*\x,4.5)--(1+0.5*\x+0.5,4.5);
\foreach \x in {3,4,6,8} \draw[blue, line width=0.5mm]  (1+0.5*\x,5)--(1+0.5*\x+0.5,5);
\foreach \x in {2,3,4,6,7,8} \draw[blue, line width=0.5mm] (1+0.5*\x,5.5)--(1+0.5*\x+0.5,5.5);
\foreach \x in {1,2,3,4,5,6,7,8,9} \draw[blue, line width=0.5mm] (1+0.5*\x,6)--(1+0.5*\x+0.5,6);
}
{
\foreach \x in {1,2,3,4,5,6,7,8,9} \draw[blue, line width=0.5mm] 
(1.5, 1+0.5*\x)--(1.5, 1+0.5*\x+0.5); 
\foreach \x in {2,3,4,5,6,7,8} \draw[blue, line width=0.5mm] 
(2, 1+0.5*\x)--(2, 1+0.5*\x+0.5); 
\foreach \x in {2,3,4,6,7} \draw[blue, line width=0.5mm] 
(2.5, 1+0.5*\x)--(2.5, 1+0.5*\x+0.5);
\foreach \x in {1,2,3,4,6} \draw[blue, line width=0.5mm] 
(3, 1+0.5*\x)--(3, 1+0.5*\x+0.5);
\foreach \x in {1,2,3,4,5,8} \draw[blue, line width=0.5mm] 
(3.5, 1+0.5*\x)--(3.5, 1+0.5*\x+0.5);
\foreach \x in {2,4,8} \draw[blue, line width=0.5mm] 
(4, 1+0.5*\x)--(4, 1+0.5*\x+0.5);
\foreach \x in {1,7} \draw[blue, line width=0.5mm] 
(4.5, 1+0.5*\x)--(4.5, 1+0.5*\x+0.5);
\foreach \x in {1,3,6,7} \draw[blue, line width=0.5mm] 
(5, 1+0.5*\x)--(5, 1+0.5*\x+0.5);
\foreach \x in {2,3,4,6,8} \draw[blue, line width=0.5mm] 
(5.5, 1+0.5*\x)--(5.5, 1+0.5*\x+0.5);
\foreach \x in {1,2,3,4,5,7,8,9} \draw[blue, line width=0.5mm] 
(6, 1+0.5*\x)--(6, 1+0.5*\x+0.5);
}

{
\foreach \x in {1} \fill[red!60!white,opacity=0.5]
(4.5, 1+0.5*\x) rectangle  (5, 1+0.5*\x+0.5); 

\foreach \x in {5} \fill[red!60!white,opacity=0.5]
(4.5, 1+0.5*\x) rectangle  (5, 1+0.5*\x+0.5); 
}

{
\foreach \x in {7} \fill[green!60!white,opacity=0.5]
(4.5, 1+0.5*\x) rectangle  (5, 1+0.5*\x+0.5); 

\foreach \x in {9} \fill[green!60!white,opacity=0.5]
(4.5, 1+0.5*\x) rectangle  (5, 1+0.5*\x+0.5);
 }
\end{scope}
}
\end{scope}

\begin{scope}[xshift=14cm]
{
\draw[gray,very thin, step=0.5cm, opacity=0.5] (0,0) grid (5.5,5.5);

\begin{scope}[xshift=-1cm,yshift=-1cm]

\node[below] at (3.75,1) [scale=1]
{\small{Second move}};

{
\foreach \x in {} \fill[blue!40!white,opacity=0.5]
(1.5, 1+0.5*\x) rectangle  (2, 1+0.5*\x+0.5); 
\foreach \x in {2,3,4,5,6,7,8} \fill[blue!40!white,opacity=0.5]
(2, 1+0.5*\x) rectangle  (2.5, 1+0.5*\x+0.5); 
\foreach \x in {5,8} \fill[blue!40!white,opacity=0.5]
(2.5, 1+0.5*\x) rectangle  (3, 1+0.5*\x+0.5); 
\foreach \x in {1,2,3,4,5,6,8} \fill[blue!40!white,opacity=0.5]
(3, 1+0.5*\x) rectangle  (3.5, 1+0.5*\x+0.5); 
\foreach \x in {6} \fill[blue!40!white,opacity=0.5]
(3.5, 1+0.5*\x) rectangle  (4, 1+0.5*\x+0.5); 
\foreach \x in {2,4,6,8} \fill[blue!40!white,opacity=0.5]
(4, 1+0.5*\x) rectangle  (4.5, 1+0.5*\x+0.5); 
\foreach \x in {2,4,5,6,7,8} \fill[blue!40!white,opacity=0.5]
(4.5, 1+0.5*\x) rectangle  (5, 1+0.5*\x+0.5); 
\foreach \x in {2,3,4,8} \fill[blue!40!white,opacity=0.5]
(5, 1+0.5*\x) rectangle  (5.5, 1+0.5*\x+0.5); 
\foreach \x in {6} \fill[blue!40!white,opacity=0.5]
(5.5, 1+0.5*\x) rectangle  (6, 1+0.5*\x+0.5); 
}

{
\foreach \x in {1,2,3,5,6,7,8,9} \draw[blue, line width=0.5mm] (1+0.5*\x,1.5)--(1+0.5*\x+0.5,1.5); 
\foreach \x in {2,6,7,8} \draw[blue, line width=0.5mm] (1+0.5*\x,2)--(1+0.5*\x+0.5,2);
\foreach \x in {6,7} \draw[blue, line width=0.5mm] (1+0.5*\x,2.5)--(1+0.5*\x+0.5,2.5);
\foreach \x in {6,7} \draw[blue, line width=0.5mm] (1+0.5*\x,3)--(1+0.5*\x+0.5,3);
\foreach \x in {3,6,8} \draw[blue, line width=0.5mm] (1+0.5*\x,3.5)--(1+0.5*\x+0.5,3.5);
\foreach \x in {3,5,6,9} \draw[blue, line width=0.5mm] (1+0.5*\x,4)--(1+0.5*\x+0.5,4);
\foreach \x in {4,5,6,9} \draw[blue, line width=0.5mm] (1+0.5*\x,4.5)--(1+0.5*\x+0.5,4.5);
\foreach \x in {3,4,6,8} \draw[blue, line width=0.5mm]  (1+0.5*\x,5)--(1+0.5*\x+0.5,5);
\foreach \x in {2,3,4,6,7,8} \draw[blue, line width=0.5mm] (1+0.5*\x,5.5)--(1+0.5*\x+0.5,5.5);
\foreach \x in {1,2,3,4,5,6,7,8,9} \draw[blue, line width=0.5mm] (1+0.5*\x,6)--(1+0.5*\x+0.5,6);
}
{
\foreach \x in {1,2,3,4,5,6,7,8,9} \draw[blue, line width=0.5mm] 
(1.5, 1+0.5*\x)--(1.5, 1+0.5*\x+0.5); 
\foreach \x in {2,3,4,5,6,7,8} \draw[blue, line width=0.5mm] 
(2, 1+0.5*\x)--(2, 1+0.5*\x+0.5); 
\foreach \x in {2,3,4,6,7} \draw[blue, line width=0.5mm] 
(2.5, 1+0.5*\x)--(2.5, 1+0.5*\x+0.5);
\foreach \x in {1,2,3,4,6} \draw[blue, line width=0.5mm] 
(3, 1+0.5*\x)--(3, 1+0.5*\x+0.5);
\foreach \x in {1,2,3,4,5,8} \draw[blue, line width=0.5mm] 
(3.5, 1+0.5*\x)--(3.5, 1+0.5*\x+0.5);
\foreach \x in {2,4,8} \draw[blue, line width=0.5mm] 
(4, 1+0.5*\x)--(4, 1+0.5*\x+0.5);
\foreach \x in {5,7} \draw[blue, line width=0.5mm] 
(4.5, 1+0.5*\x)--(4.5, 1+0.5*\x+0.5);
\foreach \x in {3,5,6,7} \draw[blue, line width=0.5mm] 
(5, 1+0.5*\x)--(5, 1+0.5*\x+0.5);
\foreach \x in {2,3,4,6,8} \draw[blue, line width=0.5mm] 
(5.5, 1+0.5*\x)--(5.5, 1+0.5*\x+0.5);
\foreach \x in {1,2,3,4,5,7,8,9} \draw[blue, line width=0.5mm] 
(6, 1+0.5*\x)--(6, 1+0.5*\x+0.5);
}

{
\foreach \x in {6} \fill[red!60!white,opacity=0.5]
(5, 1+0.5*\x) rectangle  (5.5, 1+0.5*\x+0.5); 

\foreach \x in {6} \fill[red!60!white,opacity=0.5]
(5.5, 1+0.5*\x) rectangle  (6, 1+0.5*\x+0.5);
 }
{
{
\foreach \x in {1} \fill[green!60!white,opacity=0.5]
(4.5, 1+0.5*\x) rectangle  (5, 1+0.5*\x+0.5); 

\foreach \x in {5} \fill[green!60!white,opacity=0.5]
(4.5, 1+0.5*\x) rectangle  (5, 1+0.5*\x+0.5); 
}
}

\end{scope}
}
\end{scope}

\begin{scope}[xshift=0cm, yshift=-7cm]  {
\draw[gray,very thin, step=0.5cm, opacity=0.5] (0,0) grid (5.5,5.5);

\begin{scope}[xshift=-1cm, yshift=-1cm]

\node[below] at (3.75,1) [scale=1]
{\small{Third move}};

{
\foreach \x in {} \fill[blue!40!white,opacity=0.5]
(1.5, 1+0.5*\x) rectangle  (2, 1+0.5*\x+0.5); 
\foreach \x in {2,3,4,5,6,7,8} \fill[blue!40!white,opacity=0.5]
(2, 1+0.5*\x) rectangle  (2.5, 1+0.5*\x+0.5); 
\foreach \x in {5,8} \fill[blue!40!white,opacity=0.5]
(2.5, 1+0.5*\x) rectangle  (3, 1+0.5*\x+0.5); 
\foreach \x in {1,2,3,4,5,6,8} \fill[blue!40!white,opacity=0.5]
(3, 1+0.5*\x) rectangle  (3.5, 1+0.5*\x+0.5); 
\foreach \x in {6} \fill[blue!40!white,opacity=0.5]
(3.5, 1+0.5*\x) rectangle  (4, 1+0.5*\x+0.5); 
\foreach \x in {2,4,6,8} \fill[blue!40!white,opacity=0.5]
(4, 1+0.5*\x) rectangle  (4.5, 1+0.5*\x+0.5); 
\foreach \x in {2,4,5,6,7,8} \fill[blue!40!white,opacity=0.5]
(4.5, 1+0.5*\x) rectangle  (5, 1+0.5*\x+0.5); 
\foreach \x in {2,3,4,6,8} \fill[blue!40!white,opacity=0.5]
(5, 1+0.5*\x) rectangle  (5.5, 1+0.5*\x+0.5); 
\foreach \x in {} \fill[blue!40!white,opacity=0.5]
(5.5, 1+0.5*\x) rectangle  (6, 1+0.5*\x+0.5); 
}
{
\foreach \x in {1,2,3,5,6,7,8,9} \draw[blue, line width=0.5mm] (1+0.5*\x,1.5)--(1+0.5*\x+0.5,1.5); 
\foreach \x in {2,6,7,8} \draw[blue, line width=0.5mm] (1+0.5*\x,2)--(1+0.5*\x+0.5,2);
\foreach \x in {6,7} \draw[blue, line width=0.5mm] (1+0.5*\x,2.5)--(1+0.5*\x+0.5,2.5);
\foreach \x in {6,7} \draw[blue, line width=0.5mm] (1+0.5*\x,3)--(1+0.5*\x+0.5,3);
\foreach \x in {3,6,8} \draw[blue, line width=0.5mm] (1+0.5*\x,3.5)--(1+0.5*\x+0.5,3.5);
\foreach \x in {3,5,6,8} \draw[blue, line width=0.5mm] (1+0.5*\x,4)--(1+0.5*\x+0.5,4);
\foreach \x in {4,5,6,8} \draw[blue, line width=0.5mm] (1+0.5*\x,4.5)--(1+0.5*\x+0.5,4.5);
\foreach \x in {3,4,6,8} \draw[blue, line width=0.5mm]  (1+0.5*\x,5)--(1+0.5*\x+0.5,5);
\foreach \x in {2,3,4,6,7,8} \draw[blue, line width=0.5mm] (1+0.5*\x,5.5)--(1+0.5*\x+0.5,5.5);
\foreach \x in {1,2,3,4,5,6,7,8,9} \draw[blue, line width=0.5mm] (1+0.5*\x,6)--(1+0.5*\x+0.5,6);
}
{
\foreach \x in {1,2,3,4,5,6,7,8,9} \draw[blue, line width=0.5mm] 
(1.5, 1+0.5*\x)--(1.5, 1+0.5*\x+0.5); 
\foreach \x in {2,3,4,5,6,7,8} \draw[blue, line width=0.5mm] 
(2, 1+0.5*\x)--(2, 1+0.5*\x+0.5); 
\foreach \x in {2,3,4,6,7} \draw[blue, line width=0.5mm] 
(2.5, 1+0.5*\x)--(2.5, 1+0.5*\x+0.5);
\foreach \x in {1,2,3,4,6} \draw[blue, line width=0.5mm] 
(3, 1+0.5*\x)--(3, 1+0.5*\x+0.5);
\foreach \x in {1,2,3,4,5,8} \draw[blue, line width=0.5mm] 
(3.5, 1+0.5*\x)--(3.5, 1+0.5*\x+0.5);
\foreach \x in {2,4,8} \draw[blue, line width=0.5mm] 
(4, 1+0.5*\x)--(4, 1+0.5*\x+0.5);
\foreach \x in {5,7} \draw[blue, line width=0.5mm] 
(4.5, 1+0.5*\x)--(4.5, 1+0.5*\x+0.5);
\foreach \x in {3,5,7} \draw[blue, line width=0.5mm] 
(5, 1+0.5*\x)--(5, 1+0.5*\x+0.5);
\foreach \x in {2,3,4,6,8} \draw[blue, line width=0.5mm] 
(5.5, 1+0.5*\x)--(5.5, 1+0.5*\x+0.5);
\foreach \x in {1,2,3,4,5,6,7,8,9} \draw[blue, line width=0.5mm] 
(6, 1+0.5*\x)--(6, 1+0.5*\x+0.5);
}

{
\draw[yellow, dashed, line width =0.3mm] (2,2) -- (5.5,2)--(5.5,5.5)--(2,5.5)--(2,2);

{

\foreach \x in {2} \fill[red!60!white,opacity=0.5]
(2.5, 1+0.5*\x) rectangle  (3, 1+0.5*\x+0.5); 

\foreach \x in {5} \fill[red!60!white,opacity=0.5]
(2.5, 1+0.5*\x) rectangle  (3, 1+0.5*\x+0.5); 

}

{
\foreach \x in {6} \fill[green!60!white,opacity=0.5]
(5, 1+0.5*\x) rectangle  (5.5, 1+0.5*\x+0.5); 

\foreach \x in {6} \fill[green!60!white,opacity=0.5]
(5.5, 1+0.5*\x) rectangle  (6, 1+0.5*\x+0.5); 
}

}

\end{scope}
}\end{scope}

\begin{scope}[xshift=7cm, yshift=-7cm]
{
\draw[gray,very thin, step=0.5cm, opacity=0.5] (0,0) grid (5.5,5.5);

\begin{scope}[xshift=-1cm, yshift=-1cm]

\node[below] at (3.75,1) [scale=1]
{\small{Fourth move}};

{
\foreach \x in {} \fill[blue!40!white,opacity=0.5]
(1.5, 1+0.5*\x) rectangle  (2, 1+0.5*\x+0.5); 
\foreach \x in {2,3,4,5,6,7,8} \fill[blue!40!white,opacity=0.5]
(2, 1+0.5*\x) rectangle  (2.5, 1+0.5*\x+0.5); 
\foreach \x in {2,8} \fill[blue!40!white,opacity=0.5]
(2.5, 1+0.5*\x) rectangle  (3, 1+0.5*\x+0.5); 
\foreach \x in {1,2,3,4,5,6,8} \fill[blue!40!white,opacity=0.5]
(3, 1+0.5*\x) rectangle  (3.5, 1+0.5*\x+0.5); 
\foreach \x in {6} \fill[blue!40!white,opacity=0.5]
(3.5, 1+0.5*\x) rectangle  (4, 1+0.5*\x+0.5); 
\foreach \x in {2,4,6,8} \fill[blue!40!white,opacity=0.5]
(4, 1+0.5*\x) rectangle  (4.5, 1+0.5*\x+0.5); 
\foreach \x in {2,4,5,6,7,8} \fill[blue!40!white,opacity=0.5]
(4.5, 1+0.5*\x) rectangle  (5, 1+0.5*\x+0.5); 
\foreach \x in {2,3,4,6,8} \fill[blue!40!white,opacity=0.5]
(5, 1+0.5*\x) rectangle  (5.5, 1+0.5*\x+0.5); 
\foreach \x in {} \fill[blue!40!white,opacity=0.5]
(5.5, 1+0.5*\x) rectangle  (6, 1+0.5*\x+0.5); 
}

{
\foreach \x in {1,2,3,5,6,7,8,9} \draw[blue, line width=0.5mm] (1+0.5*\x,1.5)--(1+0.5*\x+0.5,1.5); 
\foreach \x in {2,3,6,7,8} \draw[blue, line width=0.5mm] (1+0.5*\x,2)--(1+0.5*\x+0.5,2);
\foreach \x in {3,6,7} \draw[blue, line width=0.5mm] (1+0.5*\x,2.5)--(1+0.5*\x+0.5,2.5);
\foreach \x in {6,7} \draw[blue, line width=0.5mm] (1+0.5*\x,3)--(1+0.5*\x+0.5,3);
\foreach \x in {6,8} \draw[blue, line width=0.5mm] (1+0.5*\x,3.5)--(1+0.5*\x+0.5,3.5);
\foreach \x in {5,6,8} \draw[blue, line width=0.5mm] (1+0.5*\x,4)--(1+0.5*\x+0.5,4);
\foreach \x in {4,5,6,8} \draw[blue, line width=0.5mm] (1+0.5*\x,4.5)--(1+0.5*\x+0.5,4.5);
\foreach \x in {3,4,6,8} \draw[blue, line width=0.5mm]  (1+0.5*\x,5)--(1+0.5*\x+0.5,5);
\foreach \x in {2,3,4,6,7,8} \draw[blue, line width=0.5mm] (1+0.5*\x,5.5)--(1+0.5*\x+0.5,5.5);
\foreach \x in {1,2,3,4,5,6,7,8,9} \draw[blue, line width=0.5mm] (1+0.5*\x,6)--(1+0.5*\x+0.5,6);
}
{
\foreach \x in {1,2,3,4,5,6,7,8,9} \draw[blue, line width=0.5mm] 
(1.5, 1+0.5*\x)--(1.5, 1+0.5*\x+0.5); 
\foreach \x in {2,3,4,5,6,7,8} \draw[blue, line width=0.5mm] 
(2, 1+0.5*\x)--(2, 1+0.5*\x+0.5); 
\foreach \x in {3,4,5,6,7} \draw[blue, line width=0.5mm] 
(2.5, 1+0.5*\x)--(2.5, 1+0.5*\x+0.5);
\foreach \x in {1,3,4,5,6} \draw[blue, line width=0.5mm] 
(3, 1+0.5*\x)--(3, 1+0.5*\x+0.5);
\foreach \x in {1,2,3,4,5,8} \draw[blue, line width=0.5mm] 
(3.5, 1+0.5*\x)--(3.5, 1+0.5*\x+0.5);
\foreach \x in {2,4,8} \draw[blue, line width=0.5mm] 
(4, 1+0.5*\x)--(4, 1+0.5*\x+0.5);
\foreach \x in {5,7} \draw[blue, line width=0.5mm] 
(4.5, 1+0.5*\x)--(4.5, 1+0.5*\x+0.5);
\foreach \x in {3,5,7} \draw[blue, line width=0.5mm] 
(5, 1+0.5*\x)--(5, 1+0.5*\x+0.5);
\foreach \x in {2,3,4,6,8} \draw[blue, line width=0.5mm] 
(5.5, 1+0.5*\x)--(5.5, 1+0.5*\x+0.5);
\foreach \x in {1,2,3,4,5,6,7,8,9} \draw[blue, line width=0.5mm] 
(6, 1+0.5*\x)--(6, 1+0.5*\x+0.5);
}

{
\foreach \x in {2} \fill[red!60!white,opacity=0.5]
(3.5, 1+0.5*\x) rectangle  (4, 1+0.5*\x+0.5); 

\foreach \x in {6} \fill[red!60!white,opacity=0.5]
(3.5, 1+0.5*\x) rectangle  (4, 1+0.5*\x+0.5); 
}


{
\foreach \x in {2} \fill[green!60!white,opacity=0.5]
(2.5, 1+0.5*\x) rectangle  (3, 1+0.5*\x+0.5); 

\foreach \x in {5} \fill[green!60!white,opacity=0.5]
(2.5, 1+0.5*\x) rectangle  (3, 1+0.5*\x+0.5); 
}


\draw[yellow, dashed, line width =0.3mm] (2,2) -- (5.5,2)--(5.5,5.5)--(2,5.5)--(2,2);

\end{scope}
}
\end{scope}

\begin{scope}[xshift=14cm, yshift=-7cm]
 {
\draw[gray,very thin, step=0.5cm, opacity=0.5] (0,0) grid (5.5,5.5);

\begin{scope}[xshift=-1cm, yshift=-1cm]

\node[below] at (3.75,1) [scale=1]
{\small{Fifth move}};

{
\foreach \x in {} \fill[blue!40!white,opacity=0.5]
(1.5, 1+0.5*\x) rectangle  (2, 1+0.5*\x+0.5); 
\foreach \x in {2,3,4,5,6,7,8} \fill[blue!40!white,opacity=0.5]
(2, 1+0.5*\x) rectangle  (2.5, 1+0.5*\x+0.5); 
\foreach \x in {2,8} \fill[blue!40!white,opacity=0.5]
(2.5, 1+0.5*\x) rectangle  (3, 1+0.5*\x+0.5); 
\foreach \x in {1,2,3,4,5,6,8} \fill[blue!40!white,opacity=0.5]
(3, 1+0.5*\x) rectangle  (3.5, 1+0.5*\x+0.5); 
\foreach \x in {2} \fill[blue!40!white,opacity=0.5]
(3.5, 1+0.5*\x) rectangle  (4, 1+0.5*\x+0.5); 
\foreach \x in {2,4,6,8} \fill[blue!40!white,opacity=0.5]
(4, 1+0.5*\x) rectangle  (4.5, 1+0.5*\x+0.5); 
\foreach \x in {2,4,5,6,7,8} \fill[blue!40!white,opacity=0.5]
(4.5, 1+0.5*\x) rectangle  (5, 1+0.5*\x+0.5); 
\foreach \x in {2,3,4,6,8} \fill[blue!40!white,opacity=0.5]
(5, 1+0.5*\x) rectangle  (5.5, 1+0.5*\x+0.5); 
\foreach \x in {} \fill[blue!40!white,opacity=0.5]
(5.5, 1+0.5*\x) rectangle  (6, 1+0.5*\x+0.5); 
}

{
\foreach \x in {1,2,3,5,6,7,8,9} \draw[blue, line width=0.5mm] (1+0.5*\x,1.5)--(1+0.5*\x+0.5,1.5); 
\foreach \x in {2,3,5,6,7,8} \draw[blue, line width=0.5mm] (1+0.5*\x,2)--(1+0.5*\x+0.5,2);
\foreach \x in {3,5,6,7} \draw[blue, line width=0.5mm] (1+0.5*\x,2.5)--(1+0.5*\x+0.5,2.5);
\foreach \x in {6,7} \draw[blue, line width=0.5mm] (1+0.5*\x,3)--(1+0.5*\x+0.5,3);
\foreach \x in {6,8} \draw[blue, line width=0.5mm] (1+0.5*\x,3.5)--(1+0.5*\x+0.5,3.5);
\foreach \x in {6,8} \draw[blue, line width=0.5mm] (1+0.5*\x,4)--(1+0.5*\x+0.5,4);
\foreach \x in {4,6,8} \draw[blue, line width=0.5mm] (1+0.5*\x,4.5)--(1+0.5*\x+0.5,4.5);
\foreach \x in {3,4,6,8} \draw[blue, line width=0.5mm]  (1+0.5*\x,5)--(1+0.5*\x+0.5,5);
\foreach \x in {2,3,4,6,7,8} \draw[blue, line width=0.5mm] (1+0.5*\x,5.5)--(1+0.5*\x+0.5,5.5);
\foreach \x in {1,2,3,4,5,6,7,8,9} \draw[blue, line width=0.5mm] (1+0.5*\x,6)--(1+0.5*\x+0.5,6);
}
{
\foreach \x in {1,2,3,4,5,6,7,8,9} \draw[blue, line width=0.5mm] 
(1.5, 1+0.5*\x)--(1.5, 1+0.5*\x+0.5); 
\foreach \x in {2,3,4,5,6,7,8} \draw[blue, line width=0.5mm] 
(2, 1+0.5*\x)--(2, 1+0.5*\x+0.5); 
\foreach \x in {3,4,5,6,7} \draw[blue, line width=0.5mm] 
(2.5, 1+0.5*\x)--(2.5, 1+0.5*\x+0.5);
\foreach \x in {1,3,4,5,6} \draw[blue, line width=0.5mm] 
(3, 1+0.5*\x)--(3, 1+0.5*\x+0.5);
\foreach \x in {1,3,4,5,6,8} \draw[blue, line width=0.5mm] 
(3.5, 1+0.5*\x)--(3.5, 1+0.5*\x+0.5);
\foreach \x in {4,6,8} \draw[blue, line width=0.5mm] 
(4, 1+0.5*\x)--(4, 1+0.5*\x+0.5);
\foreach \x in {5,7} \draw[blue, line width=0.5mm] 
(4.5, 1+0.5*\x)--(4.5, 1+0.5*\x+0.5);
\foreach \x in {3,5,7} \draw[blue, line width=0.5mm] 
(5, 1+0.5*\x)--(5, 1+0.5*\x+0.5);
\foreach \x in {2,3,4,6,8} \draw[blue, line width=0.5mm] 
(5.5, 1+0.5*\x)--(5.5, 1+0.5*\x+0.5);
\foreach \x in {1,2,3,4,5,6,7,8,9} \draw[blue, line width=0.5mm] 
(6, 1+0.5*\x)--(6, 1+0.5*\x+0.5);
}

{
\foreach \x in {7} \fill[red!60!white,opacity=0.5]
(4.5, 1+0.5*\x) rectangle  (5, 1+0.5*\x+0.5); 

\foreach \x in {7} \fill[red!60!white,opacity=0.5]
(5, 1+0.5*\x) rectangle  (5.5, 1+0.5*\x+0.5); 

}

{
\foreach \x in {2} \fill[green!60!white,opacity=0.5]
(3.5, 1+0.5*\x) rectangle  (4, 1+0.5*\x+0.5); 

\foreach \x in {6} \fill[green!60!white,opacity=0.5]
(3.5, 1+0.5*\x) rectangle  (4, 1+0.5*\x+0.5); 

}


\draw[yellow, dashed, line width =0.3mm] (2,2) -- (5.5,2)--(5.5,5.5)--(2,5.5)--(2,2);

\end{scope}
}   
\end{scope}

\begin{scope}[xshift=0cm, yshift=-14cm]
 {
\draw[gray,very thin, step=0.5cm, opacity=0.5] (0,0) grid (5.5,5.5);

\begin{scope}[xshift=-1cm, yshift=-1cm]

\node[below] at (3.75,1) [scale=1]
{\small{Sixth move}};

{
\foreach \x in {} \fill[blue!40!white,opacity=0.5]
(1.5, 1+0.5*\x) rectangle  (2, 1+0.5*\x+0.5); 
\foreach \x in {2,3,4,5,6,7,8} \fill[blue!40!white,opacity=0.5]
(2, 1+0.5*\x) rectangle  (2.5, 1+0.5*\x+0.5); 
\foreach \x in {2,8} \fill[blue!40!white,opacity=0.5]
(2.5, 1+0.5*\x) rectangle  (3, 1+0.5*\x+0.5); 
\foreach \x in {1,2,3,4,5,6,8} \fill[blue!40!white,opacity=0.5]
(3, 1+0.5*\x) rectangle  (3.5, 1+0.5*\x+0.5); 
\foreach \x in {2} \fill[blue!40!white,opacity=0.5]
(3.5, 1+0.5*\x) rectangle  (4, 1+0.5*\x+0.5); 
\foreach \x in {2,4,6,8} \fill[blue!40!white,opacity=0.5]
(4, 1+0.5*\x) rectangle  (4.5, 1+0.5*\x+0.5); 
\foreach \x in {2,4,5,6,8} \fill[blue!40!white,opacity=0.5]
(4.5, 1+0.5*\x) rectangle  (5, 1+0.5*\x+0.5); 
\foreach \x in {2,3,4,6,7,8} \fill[blue!40!white,opacity=0.5]
(5, 1+0.5*\x) rectangle  (5.5, 1+0.5*\x+0.5); 
\foreach \x in {} \fill[blue!40!white,opacity=0.5]
(5.5, 1+0.5*\x) rectangle  (6, 1+0.5*\x+0.5); 
}

{
\foreach \x in {1,2,3,5,6,7,8,9} \draw[blue, line width=0.5mm] (1+0.5*\x,1.5)--(1+0.5*\x+0.5,1.5); 
\foreach \x in {2,3,5,6,7,8} \draw[blue, line width=0.5mm] (1+0.5*\x,2)--(1+0.5*\x+0.5,2);
\foreach \x in {3,5,6,7} \draw[blue, line width=0.5mm] (1+0.5*\x,2.5)--(1+0.5*\x+0.5,2.5);
\foreach \x in {6,7} \draw[blue, line width=0.5mm] (1+0.5*\x,3)--(1+0.5*\x+0.5,3);
\foreach \x in {6,8} \draw[blue, line width=0.5mm] (1+0.5*\x,3.5)--(1+0.5*\x+0.5,3.5);
\foreach \x in {6,8} \draw[blue, line width=0.5mm] (1+0.5*\x,4)--(1+0.5*\x+0.5,4);
\foreach \x in {4,6,7} \draw[blue, line width=0.5mm] (1+0.5*\x,4.5)--(1+0.5*\x+0.5,4.5);
\foreach \x in {3,4,6,7} \draw[blue, line width=0.5mm]  (1+0.5*\x,5)--(1+0.5*\x+0.5,5);
\foreach \x in {2,3,4,6,7,8} \draw[blue, line width=0.5mm] (1+0.5*\x,5.5)--(1+0.5*\x+0.5,5.5);
\foreach \x in {1,2,3,4,5,6,7,8,9} \draw[blue, line width=0.5mm] (1+0.5*\x,6)--(1+0.5*\x+0.5,6);
}
{
\foreach \x in {1,2,3,4,5,6,7,8,9} \draw[blue, line width=0.5mm] 
(1.5, 1+0.5*\x)--(1.5, 1+0.5*\x+0.5); 
\foreach \x in {2,3,4,5,6,7,8} \draw[blue, line width=0.5mm] 
(2, 1+0.5*\x)--(2, 1+0.5*\x+0.5); 
\foreach \x in {3,4,5,6,7} \draw[blue, line width=0.5mm] 
(2.5, 1+0.5*\x)--(2.5, 1+0.5*\x+0.5);
\foreach \x in {1,3,4,5,6} \draw[blue, line width=0.5mm] 
(3, 1+0.5*\x)--(3, 1+0.5*\x+0.5);
\foreach \x in {1,3,4,5,6,8} \draw[blue, line width=0.5mm] 
(3.5, 1+0.5*\x)--(3.5, 1+0.5*\x+0.5);
\foreach \x in {4,6,8} \draw[blue, line width=0.5mm] 
(4, 1+0.5*\x)--(4, 1+0.5*\x+0.5);
\foreach \x in {5} \draw[blue, line width=0.5mm] 
(4.5, 1+0.5*\x)--(4.5, 1+0.5*\x+0.5);
\foreach \x in {3,5,7} \draw[blue, line width=0.5mm] 
(5, 1+0.5*\x)--(5, 1+0.5*\x+0.5);
\foreach \x in {2,3,4,6,7,8} \draw[blue, line width=0.5mm] 
(5.5, 1+0.5*\x)--(5.5, 1+0.5*\x+0.5);
\foreach \x in {1,2,3,4,5,6,7,8,9} \draw[blue, line width=0.5mm] 
(6, 1+0.5*\x)--(6, 1+0.5*\x+0.5);
}

{
\foreach \x in {7} \fill[green!60!white,opacity=0.5]
(4.5, 1+0.5*\x) rectangle  (5, 1+0.5*\x+0.5); 

\foreach \x in {7} \fill[green!60!white,opacity=0.5]
(5, 1+0.5*\x) rectangle  (5.5, 1+0.5*\x+0.5); 

}

{
\foreach \x in {5} \fill[red!60!white,opacity=0.5]
(4.5, 1+0.5*\x) rectangle  (5, 1+0.5*\x+0.5); 

\foreach \x in {5} \fill[red!60!white,opacity=0.5]
(5, 1+0.5*\x) rectangle  (5.5, 1+0.5*\x+0.5); 

}


\draw[yellow, dashed, line width =0.3mm] (2,2) -- (5.5,2)--(5.5,5.5)--(2,5.5)--(2,2);

\end{scope}
}   
\end{scope}

\begin{scope}[xshift=7cm, yshift=-14cm]
{
\draw[gray,very thin, step=0.5cm, opacity=0.5]  (0,0) grid (5.5,5.5);

\begin{scope}[xshift=-1cm, yshift=-1cm]

\node[below] at (3.75,1) [scale=1]
{\small{Seventh move}};

{
\foreach \x in {} \fill[blue!40!white,opacity=0.5]
(1.5, 1+0.5*\x) rectangle  (2, 1+0.5*\x+0.5); 
\foreach \x in {2,3,4,5,6,7,8} \fill[blue!40!white,opacity=0.5]
(2, 1+0.5*\x) rectangle  (2.5, 1+0.5*\x+0.5); 
\foreach \x in {2,8} \fill[blue!40!white,opacity=0.5]
(2.5, 1+0.5*\x) rectangle  (3, 1+0.5*\x+0.5); 
\foreach \x in {1,2,3,4,5,6,8} \fill[blue!40!white,opacity=0.5]
(3, 1+0.5*\x) rectangle  (3.5, 1+0.5*\x+0.5); 
\foreach \x in {2} \fill[blue!40!white,opacity=0.5]
(3.5, 1+0.5*\x) rectangle  (4, 1+0.5*\x+0.5); 
\foreach \x in {2,4,6,8} \fill[blue!40!white,opacity=0.5]
(4, 1+0.5*\x) rectangle  (4.5, 1+0.5*\x+0.5); 
\foreach \x in {2,4,6,8} \fill[blue!40!white,opacity=0.5]
(4.5, 1+0.5*\x) rectangle  (5, 1+0.5*\x+0.5); 
\foreach \x in {2,3,4,5,6,7,8} \fill[blue!40!white,opacity=0.5]
(5, 1+0.5*\x) rectangle  (5.5, 1+0.5*\x+0.5); 
\foreach \x in {} \fill[blue!40!white,opacity=0.5]
(5.5, 1+0.5*\x) rectangle  (6, 1+0.5*\x+0.5); 
}
{
\foreach \x in {1,2,3,5,6,7,8,9} \draw[blue, line width=0.5mm] (1+0.5*\x,1.5)--(1+0.5*\x+0.5,1.5); 
\foreach \x in {2,3,5,6,7,8} \draw[blue, line width=0.5mm] (1+0.5*\x,2)--(1+0.5*\x+0.5,2);
\foreach \x in {3,5,6,7} \draw[blue, line width=0.5mm] (1+0.5*\x,2.5)--(1+0.5*\x+0.5,2.5);
\foreach \x in {6,7} \draw[blue, line width=0.5mm] (1+0.5*\x,3)--(1+0.5*\x+0.5,3);
\foreach \x in {6,7} \draw[blue, line width=0.5mm] (1+0.5*\x,3.5)--(1+0.5*\x+0.5,3.5);
\foreach \x in {6,7} \draw[blue, line width=0.5mm] (1+0.5*\x,4)--(1+0.5*\x+0.5,4);
\foreach \x in {4,6,7} \draw[blue, line width=0.5mm] (1+0.5*\x,4.5)--(1+0.5*\x+0.5,4.5);
\foreach \x in {3,4,6,7} \draw[blue, line width=0.5mm]  (1+0.5*\x,5)--(1+0.5*\x+0.5,5);
\foreach \x in {2,3,4,6,7,8} \draw[blue, line width=0.5mm] (1+0.5*\x,5.5)--(1+0.5*\x+0.5,5.5);
\foreach \x in {1,2,3,4,5,6,7,8,9} \draw[blue, line width=0.5mm] (1+0.5*\x,6)--(1+0.5*\x+0.5,6);
}
{
\foreach \x in {1,2,3,4,5,6,7,8,9} \draw[blue, line width=0.5mm] 
(1.5, 1+0.5*\x)--(1.5, 1+0.5*\x+0.5); 
\foreach \x in {2,3,4,5,6,7,8} \draw[blue, line width=0.5mm] 
(2, 1+0.5*\x)--(2, 1+0.5*\x+0.5); 
\foreach \x in {3,4,5,6,7} \draw[blue, line width=0.5mm] 
(2.5, 1+0.5*\x)--(2.5, 1+0.5*\x+0.5);
\foreach \x in {1,3,4,5,6} \draw[blue, line width=0.5mm] 
(3, 1+0.5*\x)--(3, 1+0.5*\x+0.5);
\foreach \x in {1,3,4,5,6,8} \draw[blue, line width=0.5mm] 
(3.5, 1+0.5*\x)--(3.5, 1+0.5*\x+0.5);
\foreach \x in {4,6,8} \draw[blue, line width=0.5mm] 
(4, 1+0.5*\x)--(4, 1+0.5*\x+0.5);
\foreach \x in {} \draw[blue, line width=0.5mm] 
(4.5, 1+0.5*\x)--(4.5, 1+0.5*\x+0.5);
\foreach \x in {3,5,7} \draw[blue, line width=0.5mm] 
(5, 1+0.5*\x)--(5, 1+0.5*\x+0.5);
\foreach \x in {2,3,4,5,6,7,8} \draw[blue, line width=0.5mm] 
(5.5, 1+0.5*\x)--(5.5, 1+0.5*\x+0.5);
\foreach \x in {1,2,3,4,5,6,7,8,9} \draw[blue, line width=0.5mm] 
(6, 1+0.5*\x)--(6, 1+0.5*\x+0.5);
}

{
\draw[yellow, dashed, line width =0.3mm] (2.5,2.5) -- (5,2.5)--(5,5)--(2.5,5)--(2.5,2.5);

}

{

\foreach \x in {4} \fill[red!60!white,opacity=0.5]
(3.5, 1+0.5*\x) rectangle  (4, 1+0.5*\x+0.5); 

\foreach \x in {4} \fill[red!60!white,opacity=0.5]
(4.5, 1+0.5*\x) rectangle  (5, 1+0.5*\x+0.5); 

}

{

\foreach \x in {5} \fill[green!60!white,opacity=0.5]
(4.5, 1+0.5*\x) rectangle  (5, 1+0.5*\x+0.5); 

\foreach \x in {5} \fill[green!60!white,opacity=0.5]
(5, 1+0.5*\x) rectangle  (5.5, 1+0.5*\x+0.5); 

}

\end{scope}
}    
\end{scope}

\begin{scope}[xshift=14cm, yshift=-14cm]
 {
\draw[gray,very thin, step=0.5cm, opacity=0.5] (0,0) grid (5.5,5.5);

\begin{scope}[xshift=-1cm, yshift=-1cm]

\node[below] at (3.75,1) [scale=1]
{\small{Eighth move}};

{
\foreach \x in {} \fill[blue!40!white,opacity=0.5]
(1.5, 1+0.5*\x) rectangle  (2, 1+0.5*\x+0.5); 
\foreach \x in {2,3,4,5,6,7,8} \fill[blue!40!white,opacity=0.5]
(2, 1+0.5*\x) rectangle  (2.5, 1+0.5*\x+0.5); 
\foreach \x in {2,8} \fill[blue!40!white,opacity=0.5]
(2.5, 1+0.5*\x) rectangle  (3, 1+0.5*\x+0.5); 
\foreach \x in {1,2,3,4,5,6,8} \fill[blue!40!white,opacity=0.5]
(3, 1+0.5*\x) rectangle  (3.5, 1+0.5*\x+0.5); 
\foreach \x in {2,4} \fill[blue!40!white,opacity=0.5]
(3.5, 1+0.5*\x) rectangle  (4, 1+0.5*\x+0.5); 
\foreach \x in {2,4,6,8} \fill[blue!40!white,opacity=0.5]
(4, 1+0.5*\x) rectangle  (4.5, 1+0.5*\x+0.5); 
\foreach \x in {2,6,8} \fill[blue!40!white,opacity=0.5]
(4.5, 1+0.5*\x) rectangle  (5, 1+0.5*\x+0.5); 
\foreach \x in {2,3,4,5,6,7,8} \fill[blue!40!white,opacity=0.5]
(5, 1+0.5*\x) rectangle  (5.5, 1+0.5*\x+0.5); 
\foreach \x in {} \fill[blue!40!white,opacity=0.5]
(5.5, 1+0.5*\x) rectangle  (6, 1+0.5*\x+0.5); 
}
{
\foreach \x in {1,2,3,5,6,7,8,9} \draw[blue, line width=0.5mm] (1+0.5*\x,1.5)--(1+0.5*\x+0.5,1.5); 
\foreach \x in {2,3,5,6,7,8} \draw[blue, line width=0.5mm] (1+0.5*\x,2)--(1+0.5*\x+0.5,2);
\foreach \x in {3,5,6,7} \draw[blue, line width=0.5mm] (1+0.5*\x,2.5)--(1+0.5*\x+0.5,2.5);
\foreach \x in {5,6} \draw[blue, line width=0.5mm] (1+0.5*\x,3)--(1+0.5*\x+0.5,3);
\foreach \x in {5,6} \draw[blue, line width=0.5mm] (1+0.5*\x,3.5)--(1+0.5*\x+0.5,3.5);
\foreach \x in {6,7} \draw[blue, line width=0.5mm] (1+0.5*\x,4)--(1+0.5*\x+0.5,4);
\foreach \x in {4,6,7} \draw[blue, line width=0.5mm] (1+0.5*\x,4.5)--(1+0.5*\x+0.5,4.5);
\foreach \x in {3,4,6,7} \draw[blue, line width=0.5mm]  (1+0.5*\x,5)--(1+0.5*\x+0.5,5);
\foreach \x in {2,3,4,6,7,8} \draw[blue, line width=0.5mm] (1+0.5*\x,5.5)--(1+0.5*\x+0.5,5.5);
\foreach \x in {1,2,3,4,5,6,7,8,9} \draw[blue, line width=0.5mm] (1+0.5*\x,6)--(1+0.5*\x+0.5,6);
}
{
\foreach \x in {1,2,3,4,5,6,7,8,9} \draw[blue, line width=0.5mm] 
(1.5, 1+0.5*\x)--(1.5, 1+0.5*\x+0.5); 
\foreach \x in {2,3,4,5,6,7,8} \draw[blue, line width=0.5mm] 
(2, 1+0.5*\x)--(2, 1+0.5*\x+0.5); 
\foreach \x in {3,4,5,6,7} \draw[blue, line width=0.5mm] 
(2.5, 1+0.5*\x)--(2.5, 1+0.5*\x+0.5);
\foreach \x in {1,3,4,5,6} \draw[blue, line width=0.5mm] 
(3, 1+0.5*\x)--(3, 1+0.5*\x+0.5);
\foreach \x in {1,3,5,6,8} \draw[blue, line width=0.5mm] 
(3.5, 1+0.5*\x)--(3.5, 1+0.5*\x+0.5);
\foreach \x in {6,8} \draw[blue, line width=0.5mm] 
(4, 1+0.5*\x)--(4, 1+0.5*\x+0.5);
\foreach \x in {4} \draw[blue, line width=0.5mm] 
(4.5, 1+0.5*\x)--(4.5, 1+0.5*\x+0.5);
\foreach \x in {3,4,5,7} \draw[blue, line width=0.5mm] 
(5, 1+0.5*\x)--(5, 1+0.5*\x+0.5);
\foreach \x in {2,3,4,5,6,7,8} \draw[blue, line width=0.5mm] 
(5.5, 1+0.5*\x)--(5.5, 1+0.5*\x+0.5);
\foreach \x in {1,2,3,4,5,6,7,8,9} \draw[blue, line width=0.5mm] 
(6, 1+0.5*\x)--(6, 1+0.5*\x+0.5);
}

{
\foreach \x in {6} \fill[red!60!white,opacity=0.5]
(3.5, 1+0.5*\x) rectangle  (4, 1+0.5*\x+0.5); 

\foreach \x in {6} \fill[red!60!white,opacity=0.5]
(4.5, 1+0.5*\x) rectangle  (5, 1+0.5*\x+0.5);
}

{

\foreach \x in {4} \fill[green!60!white,opacity=0.5]
(3.5, 1+0.5*\x) rectangle  (4, 1+0.5*\x+0.5); 

\foreach \x in {4} \fill[green!60!white,opacity=0.5]
(4.5, 1+0.5*\x) rectangle  (5, 1+0.5*\x+0.5); 

}

\draw[yellow, dashed, line width =0.3mm] (2.5,2.5) -- (5,2.5)--(5,5)--(2.5,5)--(2.5,2.5);

\end{scope}
}   
\end{scope}

\begin{scope}[xshift=0cm, yshift=-21cm]
    {
\draw[gray,very thin, step=0.5cm, opacity=0.5] (0,0) grid (5.5,5.5);

\begin{scope}[xshift=-1cm, yshift=-1cm]

\node[below] at (3.75,1) [scale=1]
{\small{Ninth move}};

{
\foreach \x in {} \fill[blue!40!white,opacity=0.5]
(1.5, 1+0.5*\x) rectangle  (2, 1+0.5*\x+0.5); 
\foreach \x in {2,3,4,5,6,7,8} \fill[blue!40!white,opacity=0.5]
(2, 1+0.5*\x) rectangle  (2.5, 1+0.5*\x+0.5); 
\foreach \x in {2,8} \fill[blue!40!white,opacity=0.5]
(2.5, 1+0.5*\x) rectangle  (3, 1+0.5*\x+0.5); 
\foreach \x in {1,2,3,4,5,6,8} \fill[blue!40!white,opacity=0.5]
(3, 1+0.5*\x) rectangle  (3.5, 1+0.5*\x+0.5); 
\foreach \x in {2,4,6} \fill[blue!40!white,opacity=0.5]
(3.5, 1+0.5*\x) rectangle  (4, 1+0.5*\x+0.5); 
\foreach \x in {2,4,6,8} \fill[blue!40!white,opacity=0.5]
(4, 1+0.5*\x) rectangle  (4.5, 1+0.5*\x+0.5); 
\foreach \x in {2,8} \fill[blue!40!white,opacity=0.5]
(4.5, 1+0.5*\x) rectangle  (5, 1+0.5*\x+0.5); 
\foreach \x in {2,3,4,5,6,7,8} \fill[blue!40!white,opacity=0.5]
(5, 1+0.5*\x) rectangle  (5.5, 1+0.5*\x+0.5); 
\foreach \x in {} \fill[blue!40!white,opacity=0.5]
(5.5, 1+0.5*\x) rectangle  (6, 1+0.5*\x+0.5); 
}
{
\foreach \x in {1,2,3,5,6,7,8,9} \draw[blue, line width=0.5mm] (1+0.5*\x,1.5)--(1+0.5*\x+0.5,1.5); 
\foreach \x in {2,3,5,6,7,8} \draw[blue, line width=0.5mm] (1+0.5*\x,2)--(1+0.5*\x+0.5,2);
\foreach \x in {3,5,6,7} \draw[blue, line width=0.5mm] (1+0.5*\x,2.5)--(1+0.5*\x+0.5,2.5);
\foreach \x in {5,6} \draw[blue, line width=0.5mm] (1+0.5*\x,3)--(1+0.5*\x+0.5,3);
\foreach \x in {5,6} \draw[blue, line width=0.5mm] (1+0.5*\x,3.5)--(1+0.5*\x+0.5,3.5);
\foreach \x in {5,6} \draw[blue, line width=0.5mm] (1+0.5*\x,4)--(1+0.5*\x+0.5,4);
\foreach \x in {4,5,6} \draw[blue, line width=0.5mm] (1+0.5*\x,4.5)--(1+0.5*\x+0.5,4.5);
\foreach \x in {3,4,6,7} \draw[blue, line width=0.5mm]  (1+0.5*\x,5)--(1+0.5*\x+0.5,5);
\foreach \x in {2,3,4,6,7,8} \draw[blue, line width=0.5mm] (1+0.5*\x,5.5)--(1+0.5*\x+0.5,5.5);
\foreach \x in {1,2,3,4,5,6,7,8,9} \draw[blue, line width=0.5mm] (1+0.5*\x,6)--(1+0.5*\x+0.5,6);
}
{
\foreach \x in {1,2,3,4,5,6,7,8,9} \draw[blue, line width=0.5mm] 
(1.5, 1+0.5*\x)--(1.5, 1+0.5*\x+0.5); 
\foreach \x in {2,3,4,5,6,7,8} \draw[blue, line width=0.5mm] 
(2, 1+0.5*\x)--(2, 1+0.5*\x+0.5); 
\foreach \x in {3,4,5,6,7} \draw[blue, line width=0.5mm] 
(2.5, 1+0.5*\x)--(2.5, 1+0.5*\x+0.5);
\foreach \x in {1,3,4,5,6} \draw[blue, line width=0.5mm] 
(3, 1+0.5*\x)--(3, 1+0.5*\x+0.5);
\foreach \x in {1,3,5,8} \draw[blue, line width=0.5mm] 
(3.5, 1+0.5*\x)--(3.5, 1+0.5*\x+0.5);
\foreach \x in {8} \draw[blue, line width=0.5mm] 
(4, 1+0.5*\x)--(4, 1+0.5*\x+0.5);
\foreach \x in {4,6} \draw[blue, line width=0.5mm] 
(4.5, 1+0.5*\x)--(4.5, 1+0.5*\x+0.5);
\foreach \x in {3,4,5,6,7} \draw[blue, line width=0.5mm] 
(5, 1+0.5*\x)--(5, 1+0.5*\x+0.5);
\foreach \x in {2,3,4,5,6,7,8} \draw[blue, line width=0.5mm] 
(5.5, 1+0.5*\x)--(5.5, 1+0.5*\x+0.5);
\foreach \x in {1,2,3,4,5,6,7,8,9} \draw[blue, line width=0.5mm] 
(6, 1+0.5*\x)--(6, 1+0.5*\x+0.5);
}

{
\foreach \x in {6} \fill[green!60!white,opacity=0.5]
(3.5, 1+0.5*\x) rectangle  (4, 1+0.5*\x+0.5); 

\foreach \x in {6} \fill[green!60!white,opacity=0.5]
(4.5, 1+0.5*\x) rectangle  (5, 1+0.5*\x+0.5); 
}

{
\draw[yellow, dashed, line width =0.3mm] (3,3) -- (4.5,3)--(4.5,4.5)--(3,4.5)--(3,3);


}

\end{scope}
}
\end{scope}

{
\node[below] at (9.5,-21.75) [scale=0.9]{\begin{tabular}{c}Figure 2.3. An illustration of RtCF on a $10 \times 10$ grid graph. Pink squares indicate boxes that  have just \\ been switched; green squares indicate boxes that we'll be switched next. The first three moves are in the \\ $j=0$ iteration of  RtCF, the next three moves are in the  $j=1$ iteration, and the last three moves are in \\   the $j=2$  iteration. No moves are required for the $j=3$ iteration, since $H_3$ has exactly one large cookie. \end{tabular}};;
}

\end{tikzpicture}

\newpage

\section{Existence of the MLC and 1LC Algorithms}

\noindent \textcolor{blue}{\textbullet} Recall that RtCF algorithm in Section 2 presupposes the existence of the moves required for its execution. The proofs of existence were deferred to the ManyLargeCookies (MLC)
and OneLargeCookie (1LC) algorithms, which we prove in Sections 3.1 and 3.2, respectively. These algorithms ensure that whenever RtCF requires a particular move, either the move is immediately available, or else there exists a cascade after which the required move becomes
available. Importantly, such cascades do not undo the progress already made: RtCF does not regress. The restrictions in the definition for cascades at the end of Section 1.4 were designed precisely for this reason.

Consider an iteration of RtCF on rectangle $G_i$ with Hamiltonian e-cycle $H_i$. At this stage, $H_i$ either has more than one large cookie, exactly one large cookie with at least one small
cookie, or exactly one large cookie with no small cookies. In the last case, $H_i$ is already in the desired form for this iteration, and RtCF proceeds to $G_{i+1}$. The MLC algorithm handles the
first case by finding the required cascades to collect large cookies when multiple large cookies are present. The 1LC algorithm handles the second case by finding the required cascades to collect small cookies when exactly one large cookie remains.

Why do we need two separate algorithms for what appears to be the same task? This is because small cookies can be harder to collect than large ones. A second large cookie $J$ always has a switchable neck $N_J$; to collect $J$ we need only find another switchable box $X$ such that $N_J \mapsto X$ is a valid move, or a cascade delivering such a switchable box. In Section 3.1, we show that it takes at most two moves to accomplish this (Proposition 3.8). Small cookies, by contrast, consist of a single non-switchable box. To collect a small cookie,
either the box $Y$ adjacent to it in $(R_i, R_{i+1})$ must be switchable, or we must find a cascade that makes $Y$ switchable. The latter task can require much longer cascades, and it is more
difficult to deal with. It requires Lemma 3.7, all of Section 3.2, most of Section 4, and several results from Chapter 1. Furthermore, the assumption that exactly one large cookie is present significantly shortens and simplifies the proofs of Proposition 3.10 and Lemmas
3.11–3.15 in Section 3.2, by precluding the possibility of several tedious cases. \textcolor{blue}{\textbullet}

\subsection{Existence of the MLC Algorithm}

\begingroup
\setlength{\intextsep}{0pt}
\setlength{\columnsep}{0pt}
\begin{wrapfigure}[]{r}{0cm}
\begin{adjustbox}{trim=0cm 0cm 0cm 0.75cm}
\begin{tikzpicture}[scale=1]
\begin{scope}[xshift=0cm, yshift=0cm]

\draw[gray,very thin, step=0.5cm, opacity=0.5] (0,0) grid (2.5,2.5);

\fill[orange!50!white, opacity=0.5] (1,0)--++(0.5,0)--++(0,1.5)--++(-0.5,0);
\fill[orange!50!white, opacity=0.5] (2,0)--++(0.5,0)--++(0,1.5)--++(-0.5,0);
\fill[orange!50!white, opacity=0.5] (1.5,0)--++(0.5,0)--++(0,0.5)--++(-0.5,0);

\draw[blue, line width=0.5mm] (0,0)--++(2.5,0)--++(0,2.5)--++(-1,0)--++(0,-0.5)--++(0.5,0)--++(0,-1.5)--++(-0.5,0)--++(0,1)--++(-0.5,0)--++(0,-1)--++(-0.5,0)--++(0,1.5)--++(0.5,0)--++(0,0.5)--++(-1,0)--++(0,-2.5);

\node at (1.75,1.25) [scale=0.8] {W};
\node at (1.25,1.25) [scale=0.8] {X};
\node at (2.25,1.25) [scale=0.8] {Y};

{
\node[below] at (1.25,0) [scale=0.8]{\small{\begin{tabular}{c} Fig. 3.1. The looping  \\ $H$-path of $W$ shaded orange. \end{tabular}}};;
}

\end{scope}
\end{tikzpicture}
\end{adjustbox}
\end{wrapfigure}

\textbf{Definitions.} Let $G$ be an $m \times n$ grid graph, let $H$ be a Hamiltonian cycle of $G$, and let $W$ be a switchable box in $H$. Let $X$ and $Y$ be the boxes adjacent to $W$ that are not its $H$-neighbours, and assume that $X$ and $Y$ belong to the same $H$-component. By Corollary 1.4, the $H$-path $P(X,Y)$ is unique. We call $P(X,Y)$ \textit{the looping $H$-path of $W$}. See Figure 3.1 for an illustration of the looping $H$-path of a switchable box $W$ in a Hamiltonian cycle of a $6 \times 6$ grid graph.

\endgroup

\null 

\noindent \textbf{Outline of the MLC algorithm.} Let $H$ be a Hamiltonian cycle of an $m \times n$ grid graph $G$ with multiple large
cookies. We first identify a large cookie $J$ with switchable
neck $N_J$. Consider what happens if $N_J$ is switched: this
would produce two cycles, $H_1$ and $H_2$. First we observe
that there must be some edge $\{v_1, v_2\}$ in $R_2$ (recall the nested rectangles from Chapter 2) with $v_1 \in H_1$ and $v_2 \in H_2$ (Lemma 3.7). The proximity of $\{v_1, v_2\}$  to the boundary constrains the possible configurations of edges in its vicinity. We analyze those configurations (Lemma 3.5) and show that either an $(H_1, H_2)$-port already exists near $\{v_1, v_2\}$, or a single-move cascade on the original $H$ yields a Hamiltonian cycle $H'$ where such a port exists after switching $N_J$.

\null 

\noindent \textbf{Proposition 3.1.} Let $G$ be an $m \times n$ grid graph, let $H$ be a Hamiltonian cycle of $G$, and let $P(X,Y)$ be the looping $H$-path of a switchable box $W$. Let $H'$ be the graph consisting of the cycles $H_1$ and $H_2$ obtained after switching $W$. Then a box $Z$ of $G$ belongs to the $H'$-cycle $P(X,Y), W,X$  if and only if $Z$ is incident on a vertex of $H_1$ and on a vertex of $H_2$.

\null

\begingroup
\setlength{\intextsep}{0pt}
\setlength{\columnsep}{10pt}
\begin{wrapfigure}[]{l}{0cm}
\begin{adjustbox}{trim=0cm 0.25cm 0cm 0cm}
\begin{tikzpicture}[scale=1.5]
\usetikzlibrary{decorations.markings}
\begin{scope}[xshift=0cm] 


\draw [orange, line width=0.5mm] (2.5,0)--++(0.5,0);
\begin{scope}
[very thick,decoration={
    markings,
    mark=at position 0.6 with {\arrow{>}}}
    ]
    
    \draw[postaction={decorate}, orange, line width=0.5mm] (0.5,1.5)--++(0,-0.5);
    \draw[postaction={decorate}, blue, line width=0.5mm] (1,1)--++(0,0.5);

    \draw[postaction={decorate}, blue, line width=0.5mm] (2,0)--++(0,0.5);
    \draw[postaction={decorate}, blue, line width=0.5mm] (2.5,0.5)--++(0,-0.5);
    
\end{scope}

\draw [orange, line width=0.5mm] plot [smooth, tension=0.75] coordinates {(3,0)(3.2,0.2)(3.25,1)(2.8,1.8)(1.75,2)(0.7,1.8)(0.5,1.5)};
\draw [orange, line width=0.5mm] plot [smooth, tension=0.75] coordinates {(0.5,1)(0.6,0.4)(0.8,0.15)(1.5,0)};

\draw [blue, line width=0.5mm] (2.5,0.5)--++(0,0.5);
\draw [blue, line width=0.5mm] (2,0.5)--++(0,0.5);
\draw [orange, line width=0.5mm, dashed] (2.5,0)--++(0,0.5);
\draw [orange, line width=0.5mm, dashed] (2,0)--++(0,0.5);

\draw [blue, line width=0.5mm] plot [smooth, tension=0.75] coordinates {(2,1)(1.9,1.25)(1.6,1.2)(1.4,0.8)(1.1,0.75)(1,1)};
\draw [blue, line width=0.5mm] plot [smooth, tension=0.75] coordinates {(1,1.5)(1.2,1.75)(2.25,1.75)(2.5,1)};

\draw [orange, line width=0.5mm] (2,0)--++(-0.5,0);



\fill[blue] (2.5,0) -- ++(0:0.05) arc [start angle=0, end angle=180, radius=0.05] -- cycle;
\fill[orange] (2.5,0) -- ++(180:0.05) arc [start angle=180, end angle=360, radius=0.05] -- cycle;

\fill[blue] (2,0.5) -- ++(0:0.05) arc [start angle=0, end angle=180, radius=0.05] -- cycle;
\fill[orange] (2,0.5) -- ++(180:0.05) arc [start angle=180, end angle=360, radius=0.05] -- cycle;

\fill[blue] (2.5,0.5) -- ++(0:0.05) arc [start angle=0, end angle=180, radius=0.05] -- cycle;
\fill[orange] (2.5,0.5) -- ++(180:0.05) arc [start angle=180, end angle=360, radius=0.05] -- cycle;

\fill[blue] (2,0) -- ++(0:0.05) arc [start angle=0, end angle=180, radius=0.05] -- cycle;
\fill[orange] (2,0) -- ++(180:0.05) arc [start angle=180, end angle=360, radius=0.05] -- cycle;

\draw[fill=orange, opacity=1] (0.5,1.5) circle [radius=0.05];
\draw[fill=orange, opacity=1] (0.5,1)  circle [radius=0.05];
\draw[fill=blue, opacity=1] (1,1.5) circle [radius=0.05];
\draw[fill=blue, opacity=1] (1,1) circle [radius=0.05];

\node[right] at  (2.5,0.5) [scale=0.8]{\small{$v_{y-1}$}};
\node[below] at  (2.5,0) [scale=0.8]{\small{$v_y$}};
\node[left] at  (2,0.5)  [scale=0.8]{\small{$v_{x+1}$}}; 
\node[below] at  (2,0) [scale=0.8]{\small{$v_x$}};

\node at  (0.75, 1.25) [scale=0.8]{Z};
\node at  (2.25, 0.25) [scale=0.8]{W};
\node at  (1.75, 0.25) [scale=0.8]{X};
\node at  (2.75, 0.25) [scale=0.8]{Y};

\node[below] at (1.75,-0.15) [scale=0.8]{\small{\begin{tabular}{c} Fig. 3.2. $\overrightarrow{K}_1$ in blue, $\overrightarrow{K}_2$ in orange. \end{tabular}}};;
\end{scope}

\end{tikzpicture}
\end{adjustbox}
\end{wrapfigure}

\noindent \textit{Proof.} Orient $H$ as a directed cycle $H=v_1, ..., v_r,v_1$. By Lemma 1.15, $W$ is anti-parallel. Let the edges of $W$ in $H$ be  $\{v_x, v_{x+1}\}$ and  $\{v_{y-1}, v_y\}$. For definiteness, assume that $X$ is adjacent to $\{v_x,v_{x+1}\}$, $Y$ is adjacent to $\{v_{y-1},v_y\}$ and that $W$ is on the right of $\{v_x,v_{x+1}\}$. Then we have that $\Phi((v_x,v_{x+1}),\text{left})=X$ and $\Phi((v_{y-1},v_y),\text{left})=Y$. Define $\overrightarrow{K}_1$ and $\overrightarrow{K}_2$ to be the subtrails $\overrightarrow{K}((v_x,v_{x+1})$, $(v_{y-1}, v_y))$ and $\overrightarrow{K}((v_{y-1},v_y), (v_x,v_{x+1}))$ of $\overrightarrow{K}_H$, respectively. By Lemma 1.16 (i), switching $W$ splits $H$ into two cycles $H_1$ and $H_2$, with  $V(H_1)=V(\overrightarrow{K}_1) \setminus \{v_x,v_y\}$ and $V(H_2)=V(\overrightarrow{K}_2) \setminus \{v_{x+1},v_{y-1}\}$.

\null

\noindent \textit{Proof of $(\implies)$.} Since  $P(X,Y)$ is unique, any $H$-walk of boxes between $X$ and $Y$ contains $P(X,Y)$. In particular, $\Phi(\overrightarrow{K}_1,\text{left})$ contains $P(X,Y)$ and $\Phi(\overrightarrow{K}_2,\text{left})$ contains $P(X,Y)$. Let $Z$ be a box of $P(X,Y)$. Then $Z$ is added to $\Phi(\overrightarrow{K}_1,\text{left})$ by an edge of $\overrightarrow{K}_1$ and $Z$ is added to $\Phi(\overrightarrow{K}_2,\text{left})$ by an edge of $\overrightarrow{K}_2$. By definition of FTW, $Z$ is incident on a vertex of $\overrightarrow{K}_1$ and a vertex of $\overrightarrow{K}_2$. Since for $i \in \{1,2\}$, $V(\overrightarrow{K}_i) \supset V(H_i)$, we have that $Z$ is incident on a vertex of $H_1$ and a vertex of $H_2$. Lastly, note that $W$ is incident on $v_{x+1} \in H_1$ and $v_x \in H_2$. End of proof for $(\implies)$.

\endgroup 

\null 

\noindent \textit{Proof of $(\impliedby)$.} Suppose we switch $W$ and obtain the graph $H'$ consisting of the cycles $H_1$ and $H_2$. Observe that $P(X,Y), W, X$ is the only $H'$-cycle in $G$. We will say that a box $Z$ of $G$ satisfies $(*)$ if $Z$ is incident on a vertex in $H_1$ and $H_2$. We will show that if a box $Z$ of $G$ satisfies $(*)$ then it must belong to an $H'$-cycle of boxes that satisfy $(*)$. Then, since there is only one $H'$-cycle in $G$, $Z=W$ or $Z \in P(X,Y)$.

Let $Z$ be a box in $G$ that satisfies $(*)$. For definiteness, assume that $Z=R(k,l)$. We will show that $Z$ has exactly two neighbours in $G$ that satisfy $(*)$ and that $Z$ is $H'$-adjacent to those two neighbours. Since $Z$ satisfies $(*)$, either $Z$ has two vertices in $H_1$ and two vertices in $H_2$ or $Z$ has one vertex in one of $H_1$ and $H_2$ and three vertices in the other.

\null 

\begingroup 
\setlength{\intextsep}{0pt}
\setlength{\columnsep}{10pt}
\begin{wrapfigure}[]{r}{0cm}
\begin{adjustbox}{trim=0cm 0cm 0cm 0.5cm}
\begin{tikzpicture}[scale=1.5]

\begin{scope}[xshift=0cm]{
\draw[gray,very thin, step=0.5cm, opacity=0.5] (0,0) grid (1.5,1);

\draw [blue, line width=0.5mm] plot [smooth, tension=0.75] coordinates {(0.5,0)(0.25,0.5)(0.5,0.75)(1,0.5)};

\draw [orange, line width=0.5mm] plot [smooth, tension=0.75] coordinates {(0.5,0.5)(1,0.75)(1.25,0.5)(1,0)};

\draw[blue, line width=0.5mm] (0.5,0)--++(0.5,0.5);

\draw[fill=orange] (0.5,0.5) circle [radius=0.05];
\draw[fill=orange] (1,0) circle [radius=0.05];
\draw[fill=blue] (0.5,0) circle [radius=0.05];
\draw[fill=blue] (1,0.5) circle [radius=0.05];

{

\node[left] at (0,0.5) [scale=1]{\tiny{+1}};
\node[left] at (0,0) [scale=1]
{\tiny{$\ell$}};

\node[above] at (0.5,1) [scale=1]
{\tiny{$k$}};
\node[above] at (1, 1) [scale=1]{\tiny{+1}};

\node at (0.85,0.15) [scale=0.8] {\small{$Z$}};
}

\node[below] at (0.75,0) [scale=0.8]{\small{\begin{tabular}{c} Fig. 3.3. \end{tabular}}};;

} \end{scope}

\end{tikzpicture}
\end{adjustbox}
\end{wrapfigure}

\noindent \textit{CASE 1: $Z$ has two vertices in $H_1$ and two vertices in $H_2$.} First we will check that the pair of vertices belonging to $H_i$, $i\in \{1,2\}$ must be adjacent in $Z$. Assume for contradiction that $v(k,l)$ and $v(k+1,l+1)$ belong to $H_1$ and $v(k+1,l)$ and $v(k,l+1)$ belong to $H_2$. See Figure 3.3. Let $Q$ be the closed polygon consisting of the subpath $P(v(k,l), v(k+1,l+1))$ of $H_1$ and the segment $[v(k,l), v(k+1,l+1)]$. Then, by Theorem 1.1, $v(k+1,l)$ and $v(k,l+1)$ are on different sides of $H_1$. It follows that the subpath $P(v(k+1,l),v(k,l+1))$ of $H_2$ intersects $Q$. Since $P(v(k+1,l),v(k,l+1))$ does not intersect $Q$ at the segment $[v(k,l), v(k+1,l+1)]$, it must intersect $Q$ at some vertex in $P(v(k,l), v(k+1,l+1))$. But this contradicts that $H_1$ and $H_2$ are disjoint. It follows that the pair of vertices belonging to $H_i$, $i\in \{1,2\}$ must be adjacent in $Z$. 

\begingroup 
\setlength{\intextsep}{0pt}
\setlength{\columnsep}{20pt}
\begin{wrapfigure}[]{l}{0cm}
\begin{adjustbox}{trim=0cm 0cm 0cm 0cm}
\begin{tikzpicture}[scale=1.5]

\begin{scope}[xshift=0cm]{
\draw[gray,very thin, step=0.5cm, opacity=0.5] (0,0) grid (1.5,1);

\fill[green!50!white, opacity=0.5] (0,0) rectangle (1.5,0.5);

\draw[orange, line width=0.5mm] (0,0.5)--++(1,0);

\draw[fill=blue] (1,0) circle [radius=0.05];
\draw[fill=blue] (0.5,0) circle [radius=0.05];
\draw[fill=blue] (0.5,1) circle [radius=0.05];

\draw[fill=orange] (0,0.5) circle [radius=0.05];
\draw[fill=orange] (0.5,0.5) circle [radius=0.05];
\draw[fill=orange] (1,0.5) circle [radius=0.05];

{

\draw[black, line width=0.15mm] (0.45,0.2)--++(0.1,0);
\draw[black, line width=0.15mm] (0.45,0.25)--++(0.1,0);
\draw[black, line width=0.15mm] (0.45,0.3)--++(0.1,0);

\draw[black, line width=0.15mm] (0.95,0.2)--++(0.1,0);
\draw[black, line width=0.15mm] (0.95,0.25)--++(0.1,0);
\draw[black, line width=0.15mm] (0.95,0.3)--++(0.1,0);

\draw[black, line width=0.15mm] (0.45,0.7)--++(0.1,0);
\draw[black, line width=0.15mm] (0.45,0.75)--++(0.1,0);
\draw[black, line width=0.15mm] (0.45,0.8)--++(0.1,0);

}

{

\node[left] at (0,0.5) [scale=1]{\tiny{+1}};
\node[left] at (0,0) [scale=1]
{\tiny{$\ell$}};

\node[above] at (0,1) [scale=1]{\tiny{-1}};
\node[above] at (0.5,1) [scale=1]
{\tiny{$k$}};
\node[above] at (1, 1) [scale=1]{\tiny{+1}};

\node at (0.75,0.25) [scale=0.8] {\small{$Z$}};
}

\node[below] at (0.75,0) [scale=0.8]{\small{\begin{tabular}{c} Fig. 3.4. \end{tabular}}};;

} \end{scope}

\end{tikzpicture}
\end{adjustbox}
\end{wrapfigure}

\noindent For definiteness assume that $v(k,l)$ and $v(k+1,l)$ belong to $H_1$ and that $v(k+1,l+1)$ and $v(k,l+1)$ belong to $H_2$. See Figure 3.4. Since $H_1$ and $H_2$ are disjoint $e(k;l,l+1) \notin H'$ and $e(k+1;l,l+1) \notin H'$ so $Z+(-1,0)$ and $Z+(1,0)$ are $H'$-adjacent to $Z$. Since $v(k,l) \in H_1 \cap V(Z+(-1,0))$ and $v(k,l+1) \in H_2 \cap V(Z+(-1,0))$, $Z+(-1,0)$ satisfies $(*)$. Similarly, $Z+(1,0)$ satisfies $(*)$. It remains to check that $Z+(0,1)$ and  $Z+(0,-1)$ do not satisfy $(*)$.

\noindent Assume for contradiction that one of $Z+(0,1)$ and $Z+(0,-1)$  satisfies $(*)$. By symmetry we may assume WLOG that $Z+(0,1)$ satisfies $(*)$. Then at least one of $v(k,l+2)$ and $v(k+1,l+2)$ belongs to $H_1$. By symmetry we may assume WLOG that $v(k,l+2) \in H_1$. It follows that $e(k;l+1,l+2) \notin H'$ and that $v(k-1,l+1) \in H_2$. Then we must have $e(k-1,k;l+1) \in H_2$ and $e(k,k+1;l+1) \in H_2$. But then, by Corollary 1.2, one of $v(k,l)$ and $v(k,l+2)$ belongs inside the region bounded by $H_2$ and the other belongs outside it. It follows that the subpath $P(v(k,l), v(k,l+2))$ of $H_1$ intersects $H_2$, contradicting that $H_1$ and $H_2$ are disjoint. End of Case 1.

\endgroup 

\null 

\noindent \textit{CASE 2: $Z$ has one vertex in one of $H_1$ and $H_2$ and three vertices in the other.} For definiteness assume that $v(k,l)$, $v(k,l+1)$ and $v(k+1,l+1)$ belong to $H_1$ and that $v(k+1,l)$ belongs to $H_2$. Then $e(k,k+1;l) \notin H'$ and $e(k+1;l,l+1) \notin H'$, so $Z+(1,0)$ and $Z+(0,-1)$ are $H'$-neighbours of $Z$. Since $v(k,l) \in H_1 \cap V(Z+(0,-1))$ and $v(k+1,l) \in H_2 \cap V(Z+(0,-1))$, $Z+(0,-1)$ satisfies $(*)$. Similarly, $Z+(1,0)$ satisfies $(*)$. It remains to check that $Z+(0,1)$ and $Z+(0,-1)$ do not satisfy $(*)$. 

Assume for contradiction that one of $Z+(-1,0)$ and $Z+(0,1)$  satisfies $(*)$. By symmetry we may assume WLOG that $Z+(0,1)$ satisfies $(*)$. Then one of $v(k,l+2)$ and $v(k+1,l+2)$ belongs to $H_2$. Note that if $v(k+1,l+2) \in H_2$ we run into the same contradiction as in Case 1, so we only need to check the case where $v(k,l+2) \in H_2$. Now, either $e(k,k+1;l+1) \in H_1$, or $e(k,k+1;l+1) \notin H'$. 

\null 

\begingroup 
\setlength{\intextsep}{0pt}
\setlength{\columnsep}{10pt}
\begin{wrapfigure}[]{r}{0cm}
\begin{adjustbox}{trim=0cm 0.5cm 0cm 0.5cm}
\begin{tikzpicture}[scale=1.5]

\begin{scope}[xshift=0cm]{
\draw[gray,very thin, step=0.5cm, opacity=0.5] (0,0) grid (1,1.5);

\fill[green!50!white, opacity=0.5] (0,0) rectangle (0.5,0.5);
\fill[green!50!white, opacity=0.5] (0,0.5) rectangle (1,1);

\draw[blue, line width=0.5mm] (0,1)--++(0.5,0);
\draw[red, line width=0.5mm] (0,1.5)--++(0.5,-1);

\draw[fill=blue] (0,0.5) circle [radius=0.05];
\draw[fill=blue] (0,1) circle [radius=0.05];
\draw[fill=blue] (0.5,1) circle [radius=0.05];

\draw[fill=orange] (0.5,0.5) circle [radius=0.05];
\draw[fill=orange] (0,1.5) circle [radius=0.05];

{

\draw[black, line width=0.15mm] (0.45,0.7)--++(0.1,0);
\draw[black, line width=0.15mm] (0.45,0.75)--++(0.1,0);
\draw[black, line width=0.15mm] (0.45,0.8)--++(0.1,0);

\draw[black, line width=0.15mm] (0.2,0.45)--++(0,0.1);
\draw[black, line width=0.15mm] (0.25,0.45)--++(0,0.1);
\draw[black, line width=0.15mm] (0.3,0.45)--++(0,0.1);

}

{

\node[left] at (0,1) [scale=1]{\tiny{+1}};
\node[left] at (0,0.5) [scale=1]
{\tiny{$\ell$}};
\node[left] at (0,0) [scale=1]{\tiny{-1}};

\node[above] at (0,1.5) [scale=1]
{\tiny{$k$}};
\node[above] at (0.5, 1.5) [scale=1]{\tiny{+1}};

\node at (0.25,0.75) [scale=0.8] {\small{$Z$}};
}

\node[below] at (0.5,0) [scale=0.8]{\small{\begin{tabular}{c} Fig. 3.5(a).\end{tabular}}};;

} \end{scope}

\begin{scope}[xshift=1.5cm]{
\draw[gray,very thin, step=0.5cm, opacity=0.5] (0,0) grid (1,1.5);

\fill[green!50!white, opacity=0.5] (0,0) rectangle (0.5,0.5);
\fill[green!50!white, opacity=0.5] (0,0.5) rectangle (1,1);

\draw[red, line width=0.5mm] (0,1.5)--++(0.5,-1);

\draw[fill=blue] (0,0.5) circle [radius=0.05];
\draw[fill=blue] (0,1) circle [radius=0.05];
\draw[fill=blue] (0.5,1) circle [radius=0.05];

\draw[fill=orange] (0.5,0.5) circle [radius=0.05];
\draw[fill=orange] (0,1.5) circle [radius=0.05];

{

\draw[black, line width=0.15mm] (0.45,0.7)--++(0.1,0);
\draw[black, line width=0.15mm] (0.45,0.75)--++(0.1,0);
\draw[black, line width=0.15mm] (0.45,0.8)--++(0.1,0);

\draw[black, line width=0.15mm] (0.2,0.45)--++(0,0.1);
\draw[black, line width=0.15mm] (0.25,0.45)--++(0,0.1);
\draw[black, line width=0.15mm] (0.3,0.45)--++(0,0.1);

\draw[black, line width=0.15mm] (0.2,0.95)--++(0,0.1);
\draw[black, line width=0.15mm] (0.25,0.95)--++(0,0.1);
\draw[black, line width=0.15mm] (0.3,0.95)--++(0,0.1);

}

{

\node[left] at (0,1) [scale=1]{\tiny{+1}};
\node[left] at (0,0.5) [scale=1]
{\tiny{$\ell$}};
\node[left] at (0,0) [scale=1]{\tiny{-1}};

\node[above] at (0,1.5) [scale=1]
{\tiny{$k$}};
\node[above] at (0.5, 1.5) [scale=1]{\tiny{+1}};

\node at (0.25,0.75) [scale=0.8] {\small{$Z$}};
}

\node[below] at (0.5,0) [scale=0.8]{\small{\begin{tabular}{c} Fig. 3.5(b). \end{tabular}}};;

} \end{scope}

\end{tikzpicture}
\end{adjustbox}
\end{wrapfigure}

\noindent \textit{CASE 2.1: $e(k,k+1;l+1) \in H_1$.} Then the segment $[v(k,l+2), v(k+1,l)]$ intersects $H_1$ at the point $(k \frac{1}{2}, l+1)$. by Corollary 1.1, the vertices $e(k,l+2)$ and $v(k+1,l)$ are on different sides of $H_1$, and we run into the same contradiction as in Case 1 again. End of Case 2.1. See Figure 3.5(a)

\null 

\noindent \textit{CASE 2.2: $e(k,k+1;l+1) \notin H_1$.} Consider the polygon $Q$ consisting of the segment $[v(k,l+2), v(k+1,l)]$ and the subpath $P(v(k,l+2), v(k+1,l))$ of $H_2$. By Corollary 1.1 the vertices $v(k,l+1)$ and $v(k+1,l+1)$ are on different sides of $Q$. BY JCT the subpath $P(v(k,l+1),v(k+1,l+1))$ of $H_1$ intersects $Q$. Since $P(v(k,l+1),v(k+1,l+1))$ does not intersect $Q$ at the segment $[v(k,l+2), v(k+1,l)]$, it must do so at some vertex of $P(v(k,l+2), v(k+1,l))$, contradicting that $H_1$ and $H_2$ are disjoint. See Figure 3.5(b). End of Case 2.2. End of Case 2. End of proof for $(\impliedby)$. $\square$

\endgroup

\null

\noindent \textbf{Corollary 3.2 (i).} Let $G$ be an $m \times n$ grid graph, let $H$ be a Hamiltonian cycle of $G$, and let $W$ be a switchable box in $H$. Let $H'$ be the graph consisting of the cycles $H_1$ and $H_2$ obtained after switching $W$. Let $a,b$ and $c$ be colinear vertices such that $b$ is adjacent to $a$ and $c$. Then, for $i \in \{1,2\}$, If $a$ and $c$ belong to $H_i$, so must $b$. See Figure 3.4. $\square$

\null

\noindent \textbf{Corollary 3.2 (ii).} Let $G$ be an $m \times n$ grid graph, let $H$ be a Hamiltonian cycle of $G$, and let $W$ be a switchable box in $H$. Let $H'$ be the graph consisting of the cycles $H_1$ and $H_2$ obtained after switching $W$. Let $Z$ be a box on vertices $a,b,c$, and $d$ such that $a$ and $b$ belong to $H_1$, and $c$ and $d$ belong to $H_2$. Then $a$ is adjacent to $b$, and $c$ is adjacent to $d$. See Figure 3.3. $\square$

\null 

\begingroup
\setlength{\intextsep}{0pt}
\setlength{\columnsep}{10pt}
\begin{wrapfigure}[]{l}{0cm}
\begin{adjustbox}{trim=0cm 0.5cm 0cm 0.25cm}
\begin{tikzpicture}[scale=1.5]
\usetikzlibrary{decorations.markings}
\begin{scope}[xshift=0.00cm] 


\draw [orange, line width=0.5mm] (2.5,0)--++(0.5,0);
\begin{scope}
[very thick,decoration={
    markings,
    mark=at position 0.6 with {\arrow{>}}}
    ]
    
    \draw[postaction={decorate}, orange, line width=0.5mm] (0.5,1.5)--++(0,-0.5);
    \draw[postaction={decorate}, blue, line width=0.5mm] (1,1)--++(0,0.5);
    
\end{scope}

\draw [orange, line width=0.5mm] plot [smooth, tension=0.75] coordinates {(3,0)(3.2,0.2)(3.25,1)(2.8,1.8)(1.75,2)(0.7,1.8)(0.5,1.5)};
\draw [orange, line width=0.5mm] plot [smooth, tension=0.75] coordinates {(0.5,1)(0.6,0.4)(0.8,0.15)(1.5,0)};

\draw [blue, line width=0.5mm] (2.5,0)--++(0,1);
\draw [blue, line width=0.5mm] (2,0)--++(0,1);
\draw [blue, line width=0.5mm] plot [smooth, tension=0.75] coordinates {(2,1)(1.9,1.25)(1.6,1.2)(1.4,0.8)(1.1,0.75)(1,1)};
\draw [blue, line width=0.5mm] plot [smooth, tension=0.75] coordinates {(1,1.5)(1.2,1.75)(2.25,1.75)(2.5,1)};

\draw [orange, line width=0.5mm] (2,0)--++(-0.5,0);

\draw[fill=orange, opacity=1] (0.8,0.15)  circle [radius=0.05];
\draw[fill=orange, opacity=1] (0.6,0.4) circle [radius=0.05];

\draw[fill=blue, opacity=1] (2.5,0.5) circle [radius=0.05];
\draw[fill=blue, opacity=1] (2,0.5) circle [radius=0.05];

\fill[orange] (2,0) -- ++(135:0.05) arc [start angle=135, end angle=315, radius=0.05] -- cycle;
\fill[blue] (2,0) -- ++(315:0.05) arc [start angle=315, end angle=495, radius=0.05] -- cycle;

\fill[orange] (2.5,0) -- ++(225:0.05) arc [start angle=225, end angle=405, radius=0.05] -- cycle;
    \fill[blue] (2.5,0) -- ++(45:0.05) arc [start angle=45, end angle=225, radius=0.05] -- cycle;

\draw[fill=orange, opacity=1] (0.5,1.5) circle [radius=0.05];
\draw[fill=orange, opacity=1] (0.5,1)  circle [radius=0.05];
\draw[fill=blue, opacity=1] (1,1.5) circle [radius=0.05];
\draw[fill=blue, opacity=1] (1,1) circle [radius=0.05];

\node[right] at  (2.5,0.5) [scale=0.8]{\small{$v_{y-1}$}};
\node[below] at  (2.5,0) [scale=0.8]{\small{$v_y$}};
\node[left] at  (2,0.5)  [scale=0.8]{\small{$v_{x+1}$}}; 
\node[below] at  (2,0) [scale=0.8]{\small{$v_x$}};

\node[right] at  (0.8,0.15) [scale=0.8]{\small{$v_1$}};
\node[right] at  (0.6,0.4) [scale=0.8]{\small{$v_r$}};

\node[right] at  (1,1) [scale=0.8]{\small{$v_s$}};
\node[right] at  (1,1.4) [scale=0.8]{\small{$v_{s+1}$}};
\node[left] at  (0.5,1.5) [scale=0.8]{\small{$v_t$}};
\node[left] at  (0.5,1) [scale=0.8]{\small{$v_{t+1}$}};

\node at  (0.75, 1.25) [scale=0.8]{Z};
\node at  (2.25, 0.25) [scale=0.8]{W};
\node at  (1.75, 0.25) [scale=0.8]{X};
\node at  (2.75, 0.25) [scale=0.8]{Y};

\node[below] at (1.75,-0.15) [scale=0.8]{\small{\begin{tabular}{c} Fig. 3.6. $P_1$ in blue, $P_2$ in orange. \end{tabular}}};;
\end{scope}

\end{tikzpicture}
\end{adjustbox}
\end{wrapfigure}

\noindent \textbf{Proposition 3.3.} Let $H$ be a Hamiltonian cycle of an $m \times n$ grid graph $G$, let $W$ be a switchable box in $H$ and let $P(X,Y)$ be the looping $H$-path of $W$. If $P(X,Y)$ has a switchable box $Z$, then $Z \mapsto W$ is a valid move.

\null 

\noindent \textit{Proof.} Let $H=v_1, v_2 ..., v_r, v_1$. Let $W$, $P(X,Y)$, $\{v_x, v_{x+1}\}$ and $\{v_{y-1}, v_y\}$ be as in Proposition 3.1. By Lemma 1.15, $W$ is anti-parallel. Let $P_1=P(v_x,v_y)$ and let $P_2=P(v_y,v_x)$. By Proposition 3.1, every box of $P(X,Y)$ is incident on a vertex of $P_1$ and a vertex of $P_2$. Let $Z$ be a switchable box of $P(X,Y)$. Let $(v_s,v_{s+1})$ and 

\endgroup

\noindent $(v_t, v_{t+1})$ be the edges of $Z$ in $H$. For definiteness, assume $s+1<t$. Proposition 3.1 implies that exactly one of $(v_s,v_{s+1})$ and $(v_t, v_{t+1})$ is in $P_1$ and the other is in $P_2$. WLOG assume that $(v_s,v_{s+1})$ is in $P_1$ and that $(v_t, v_{t+1})$ is in $P_2$. Then we can partition the edges of $H$ as follows: $P(v_1,v_s)$, $(v_s,v_{s+1})$, $P(v_{s+1},v_t)$, $(v_t,v_{t+1})$, $P(v_{t+1},v_r)$, $\{v_r,v_1\}$ where $1<x<s<y<t<r$. 

We we check that $Z\mapsto W$ is a valid move. After removing the edges $(v_s,v_{s+1})$ and $(v_t,v_{t+1})$ we are left with two paths: $P(v_{t+1},v_s)$ and $ P(v_{s+1},v_t)$. Note that adding the edge $\{v_s, v_{t+1}\}$ gives a cycle $H_1$ consisting of the path $P(v_s,v_{t+1})$ and the edge $\{v_s,v_{t+1}\}$. and adding the edge $\{v_{s+1}, v_t\}$ gives a cycle $H_2$ consisting of the path $P(v_{s+1},v_t)$ and the edge $\{v_{s+1}, v_t\}$. Now $1<x<s$ implies that $(v_x,v_{x+1}) \in H_1$ and $s<y<t$ implies that $(v_{y-1}, v_y) \in H_2$. It follows that $W$ is now an $(H_1,H_2)$-port. By Lemma 1.16 (ii), $Z \mapsto W$ is a valid move. $\square$ 

\null

\noindent \textbf{Observation 3.4.}  Let $X \mapsto Y$ be a valid move. If $X \in \text{ext}(H)$ and $Y \in \text{int}(H)$ then:

(i) \hspace{0.05cm} $X \mapsto Y$ increases the total number of cookies if and only if $X \in G_1$ and $Y\in G_0 \setminus G_1$.

(ii) \ $X \mapsto Y$ increases the total number of large cookies, leaving the total number of cookies unchanged, 

\hspace{0.55cm} if and only if $X \in G_1$, $Y \in G_1 \setminus G_2$ and $Y$ is adjacent to a small cookie. 

(iii) $X \mapsto Y$ decreases the total number of cookies if and only if $X \in G_0 \setminus G_1$ and $Y\in G_1$. 

\null

\noindent \textbf{Lemma 3.5.} Let $H$ be a Hamiltonian cycle of an $m \times n$ grid graph $G$ and let $Z$ be a switchable box in $Z \in \text{ext(H)} \cap ((G_0 \setminus G_1) \cup G_3)$. Assume that switching $Z$ splits $H$ into the cycles $H_1$ and $H_2$ that are such that there is $v_1 \in H_1 \cap R_2$ and $v_2 \in H_2 \cap R_2$ with $v_1$ adjacent to $v_2$. Then either $Z \mapsto Z'$ is a cascade, or there is a cascade $\mu, Z \mapsto Z'$ (of length two), with $Z \mapsto Z'$ nontrivial in either case.

\begingroup
\setlength{\intextsep}{0pt}
\setlength{\columnsep}{10pt}
\begin{wrapfigure}[]{r}{0cm}
\begin{adjustbox}{trim=0cm 0.5cm 0cm 0.5cm}
\begin{tikzpicture}[scale=1.5]

\begin{scope}[xshift=0cm]
{
\draw[gray,very thin, step=0.5cm, opacity=0.5] (0,0) grid (1.5,1.5);

\fill[green!50!white, opacity=0.5] (0.5,0.5)--++(1,0)--++(0,0.5)--++(-1,0);

\draw[yellow, line width =0.2mm] (0,0)--++(0,1.5);

\draw[fill=blue] (1,1) circle [radius=0.05];
\draw[fill=blue] (1,1.5) circle [radius=0.05];

\draw[fill=orange] (1,0.5) circle [radius=0.05];

{

\draw[black, line width=0.15mm] (0.95, 0.7)--++(0.1,0);
\draw[black, line width=0.15mm] (0.95, 0.75)--++(0.1,0);
\draw[black, line width=0.15mm] (0.95, 0.8)--++(0.1,0);
}

\node[left] at (0,0.5) [scale=1] {\tiny{-1}};
\node[left] at (-0.05,1) [scale=1] {\tiny{$\ell$}};
\node[left] at (0,1.5) [scale=1] {\tiny{+1}};

\node[above] at (0.5,1.55) [scale=1] {\tiny{1}};
\node[above] at (1,1.55) [scale=1] {\tiny{2}};
\node[above] at (1.5,1.55) [scale=1] {\tiny{3}};

\node[below] at (0.75,0) [scale=0.75] {\small{Fig. 3.7 (a).}};

}
\end{scope}
    
\begin{scope}[xshift=2.25cm]
{
\draw[gray,very thin, step=0.5cm, opacity=0.5] (0,0) grid (1.5,1.5);

\fill[green!50!white, opacity=0.5] (0.5,0.5)--++(1,0)--++(0,0.5)--++(-1,0);

\draw[yellow, line width =0.2mm] (0,0) -- (0,1.5);

\draw[fill=blue] (1,1) circle [radius=0.05];
\draw[fill=blue] (1,1.5) circle [radius=0.05];
\draw[fill=blue] (1.5,1) circle [radius=0.05];

\draw[fill=orange] (1,0.5) circle [radius=0.05];
\draw[fill=orange] (0.5,1) circle [radius=0.05];
\draw[fill=orange] (0.5,0.5) circle [radius=0.05];
\draw[fill=orange] (0,1) circle [radius=0.05];

\draw[blue, line width=0.5mm] (1,1.5)--++(0,-0.5)--++(0.5,0);
{
\draw[black, line width=0.15mm] (0.70,0.95)--++(0,0.1);
\draw[black, line width=0.15mm] (0.75,0.95)--++(0,0.1);
\draw[black, line width=0.15mm] (0.8,0.95)--++(0,0.1);

\draw[black, line width=0.15mm] (0.95, 0.7)--++(0.1,0);
\draw[black, line width=0.15mm] (0.95, 0.75)--++(0.1,0);
\draw[black, line width=0.15mm] (0.95, 0.8)--++(0.1,0);
}

\node[left] at (0,0.5) [scale=1] {\tiny{-1}};
\node[left] at (-0.05,1) [scale=1] {\tiny{$\ell$}};
\node[left] at (0,1.5) [scale=1] {\tiny{+1}};

\node[above] at (0.5,1.55) [scale=1] {\tiny{1}};
\node[above] at (1,1.55) [scale=1] {\tiny{2}};
\node[above] at (1.5,1.55) [scale=1] {\tiny{3}};

\node[below] at (0.75,0) [scale=0.75] {\small{Fig. 3.7 (b). Case 1.}};

}
\end{scope}
\end{tikzpicture}
\end{adjustbox}
\end{wrapfigure}

\null

\noindent In figures 3.7 through 3.11, vertices and edges of $H_1$ and $H_2$ are in blue and orange, respectively, and boxes of $\text{int}(H)$ are shaded in green.

\null 

\noindent \textit{Proof.} We will use the assumption that $Z \in  (G_0 \setminus G_1) \cup G_3$ repeatedly and implicitly throughout the proof. Switch $Z$ to obtain $H'$ consisting of the disjoint cycles $H_1$ and $H_2$. Note that now, if a vertex belongs to $H_i$ for $i\in \{1,2\}$, both edges incident on it must also belong to $H_i$. For definiteness, let $v_1 \in H_1 \cap R_2$ be the vertex $v(2,l)$ for some $l \in \{ 2, ..., n-3\}$ and let $v_2=v(2,l-1) \in H_2 \cap R_2$. Then $e(2;l-1,l) \notin H$, and by Corollary 3.2 (i), $v(2,l+1) \in H_1$ as well. By Proposition 3.1, $R(1,l-1)$ and $R(2,l-1)$ belong to $\text{int(H)}$. Now, either $v(1,l) \in H_1$, or $v(1,l) \in H_2$. See Figure 3.1.

\endgroup

\null 

\noindent \textit{CASE 1: $v(1,l) \in H_2$.} Then $e(1,2;l) \notin H$.  Corollary 3.2(i), $v(3,l) \in H_1$. It follows that $e(2;l,l+1) \in H'$ and $e(2,3;l) \in H'$. Now, by Corollary 3.2 (ii), $v(1,l-1) \in H_2$ and by Corollary 3.2 (i), $v(0,l) \in H_2$. At this point we must either have, $e(1;l,l+1) \in H'$ or $e(1;l,l+1) \notin H'$. 

\endgroup 

\null

\begingroup
\setlength{\intextsep}{0pt}
\setlength{\columnsep}{10pt}
\begin{wrapfigure}[]{r}{0cm}
\begin{adjustbox}{trim=0cm 0cm 0cm 0.5cm}
\begin{tikzpicture}[scale=1.5]
\begin{scope}[xshift=0cm]
{
\draw[gray,very thin, step=0.5cm, opacity=0.5] (0,0) grid (1.5,1.5);

\fill[green!50!white, opacity=0.5] (0,1)--++(1,0)--++(0,0.5)--++(-1,0);
\fill[green!50!white, opacity=0.5] (0.5,0.5)--++(1,0)--++(0,0.5)--++(-1,0);
\fill[green!50!white, opacity=0.5] (0,0)--++(1,0)--++(0,0.5)--++(-1,0);

\draw[yellow, line width =0.2mm] (0,0) -- (0,1.5);

\draw[blue, line width=0.5mm] (1,1.5)--++(0,-0.5)--++(0.5,0);
\draw[orange, line width=0.5mm] (0,1)--++(0.5,0)--++(0,-0.5);
\draw[orange, line width=0.5mm] (1,0.5)--++(0.5,0);
\draw[orange, line width=0.5mm] (0,0.5)--++(0.5,0);

\draw[fill=blue] (1,1) circle [radius=0.05];
\draw[fill=blue] (1,1.5) circle [radius=0.05];
\draw[fill=blue] (1.5,1) circle [radius=0.05];

\draw[fill=orange] (1,0.5) circle [radius=0.05];
\draw[fill=orange] (0.5,1) circle [radius=0.05];
\draw[fill=orange] (0.5,0.5) circle [radius=0.05];
\draw[fill=orange] (0,1) circle [radius=0.05];

\draw[fill=orange] (0,0.5) circle [radius=0.05];
\draw[fill=orange] (1.5,0.5) circle [radius=0.05];

{
\draw[black, line width=0.15mm] (0.70,0.95)--++(0,0.1);
\draw[black, line width=0.15mm] (0.75,0.95)--++(0,0.1);
\draw[black, line width=0.15mm] (0.8,0.95)--++(0,0.1);

\draw[black, line width=0.15mm] (0.70,0.45)--++(0,0.1);
\draw[black, line width=0.15mm] (0.75,0.45)--++(0,0.1);
\draw[black, line width=0.15mm] (0.8,0.45)--++(0,0.1);

\draw[black, line width=0.15mm] (0.95, 0.7)--++(0.1,0);
\draw[black, line width=0.15mm] (0.95, 0.75)--++(0.1,0);
\draw[black, line width=0.15mm] (0.95, 0.8)--++(0.1,0);

\draw[black, line width=0.15mm] (1.45, 0.7)--++(0.1,0);
\draw[black, line width=0.15mm] (1.45, 0.75)--++(0.1,0);
\draw[black, line width=0.15mm] (1.45, 0.8)--++(0.1,0);

\draw[black, line width=0.15mm] (0.45, 1.2)--++(0.1,0);
\draw[black, line width=0.15mm] (0.45, 1.25)--++(0.1,0);
\draw[black, line width=0.15mm] (0.45, 1.3)--++(0.1,0);
}

\node[left] at (0,0.5) [scale=1] {\tiny{-1}};
\node[left] at (-0.05,1)[scale=1] {\tiny{$\ell$}};
\node[left] at (0,1.5) [scale=1] {\tiny{+1}};

\node[above] at (0.5,1.5) [scale=1] {\tiny{1}};
\node[above] at (1,1.5) [scale=1] {\tiny{2}};
\node[above] at (1.5,1.5) [scale=1] {\tiny{3}};

\node[below] at (0.75,0) [scale=0.75] {\small{Fig. 3.8 (a). Case 1.1}};

}
\end{scope}

\begin{scope}[xshift=2.25cm]
{
\draw[gray,very thin, step=0.5cm, opacity=0.5] (0,0) grid (1.5,1.5);

\draw[yellow, line width =0.2mm] (0,0) -- (0,1.5);

\fill[green!50!white, opacity=0.5] (0,0.5)--++(1.5,0)--++(0,0.5)--++(-1.5,0);
\fill[green!50!white, opacity=0.5] (0.5,1)--++(0.5,0)--++(0,0.5)--++(-0.5,0);

\draw[blue, line width=0.5mm] (1,1.5)--++(0,-0.5)--++(0.5,0);
\draw[orange, line width=0.5mm] (0,1)--++(0.5,0)--++(0,0.5)--++(-0.5,0);
\draw[orange, line width=0.5mm] (0.5,0.5)--++(0.5,0);

\draw[fill=blue] (1,1) circle [radius=0.05];
\draw[fill=blue] (1,1.5) circle [radius=0.05];
\draw[fill=blue] (1.5,1) circle [radius=0.05];

\draw[fill=orange] (1,0.5) circle [radius=0.05];
\draw[fill=orange] (0.5,1) circle [radius=0.05];
\draw[fill=orange] (0.5,0.5) circle [radius=0.05];
\draw[fill=orange] (0,1) circle [radius=0.05];
\draw[fill=orange] (0,1.5) circle [radius=0.05];
\draw[fill=orange] (0.5,1.5) circle [radius=0.05];

{
\draw[black, line width=0.15mm] (0.70,0.95)--++(0,0.1);
\draw[black, line width=0.15mm] (0.75,0.95)--++(0,0.1);
\draw[black, line width=0.15mm] (0.8,0.95)--++(0,0.1);

\draw[black, line width=0.15mm] (0.95, 0.7)--++(0.1,0);
\draw[black, line width=0.15mm] (0.95, 0.75)--++(0.1,0);
\draw[black, line width=0.15mm] (0.95, 0.8)--++(0.1,0);

\draw[black, line width=0.15mm] (1.20,0.45)--++(0,0.1);
\draw[black, line width=0.15mm] (1.25,0.45)--++(0,0.1);
\draw[black, line width=0.15mm] (1.3,0.45)--++(0,0.1);
}

\node[left] at (0,0.5) [scale=1] {\tiny{-1}};
\node[left] at (-0.05,1) [scale=1] {\tiny{$\ell$}};
\node[left] at (0,1.5) [scale=1] {\tiny{+1}};

\node[above] at (0.5,1.55) [scale=1] {\tiny{1}};
\node[above] at (1,1.55) [scale=1] {\tiny{2}};
\node[above] at (1.5,1.55) [scale=1] {\tiny{3}};

\node[below] at (0.75,0) [scale=0.75] {\small{Fig. 3.8 (b). Case 1.2}};

}
\end{scope}

\end{tikzpicture}
\end{adjustbox}
\end{wrapfigure}

\noindent \textit{CASE 1.1: $e(1;l,l+1) \notin H'$.} Then $e(1;l-1,l) \in H_2$ and $e(0,1;l) \in H_2$. Since  $Z \in \text{ext(H)}$, by Proposition 3.1, $R(1,l-1) \in \text{int(H)}$. Then $R(0,l-1)$ must be a small cookie of $G$, so we must have that $e(0,1;l-1) \in H_2$ and that $e(1,2;l-1) \notin H'$. It follows that $e(2,3;l-1) \in H'$ and $e(3;l-1,l) \notin H'$. Now note that $R(2,l-1)$ is an $(H_1,H_2)$-port. Then, by Lemma 1.16 (ii) and Observation 3.4, $Z \mapsto R(2,l-1)$ is valid move that does not create new cookies. So, $Z \mapsto R(2,l-1)$ is the cascade we seek. End of Case 1.1

\null 

\noindent \textit{CASE 1.2: $e(1;l,l+1) \in H'$.} Proposition 3.1, and the assumption that $Z \in \text{ext}(H)$ imply that $R(1,l) \in \text{int}(H)$. Then, Lemma 1.14 implies that $R(0,l)$ is a small cookie of $G$. Note that if $e(2,3;l-1) \in H'$, then we're back to Case 1.1, so we may assume that $e(2,3;l-1) \notin H'$. It follows that $e(1,2; l-1) \in H_2$ and that $R(0,l-1)$ is not a small cookie of $H$. Then, by Observation 3.4 and Proposition 3.3, $R(0,l) \mapsto R(1,l)$, $Z \mapsto R(1,l-1)$ is the cascade we seek. End of Case 1.2. End of Case 1.

\endgroup

\begingroup
\setlength{\intextsep}{0pt}
\setlength{\columnsep}{10pt}
\begin{wrapfigure}[]{l}{0cm}
\begin{adjustbox}{trim=0cm 0cm 0cm 0.25cm}
\begin{tikzpicture}[scale=1.5]

\begin{scope}[xshift=0cm]
{
\draw[gray,very thin, step=0.5cm, opacity=0.5] (0,0) grid (1.5,1.5);

\fill[green!50!white, opacity=0.5] (0.5,0.5)--++(1,0)--++(0,0.5)--++(-1,0);

\draw[yellow, line width =0.2mm] (0,0)--++(0,1.5);

\draw[fill=blue] (1,1) circle [radius=0.05];
\draw[fill=blue] (1,1.5) circle [radius=0.05];
\draw[fill=blue] (0.5,1) circle [radius=0.05];

\draw[fill=orange] (1,0.5) circle [radius=0.05];

{

\draw[black, line width=0.15mm] (0.95, 0.7)--++(0.1,0);
\draw[black, line width=0.15mm] (0.95, 0.75)--++(0.1,0);
\draw[black, line width=0.15mm] (0.95, 0.8)--++(0.1,0);
}

\node[left] at (0,0.5) [scale=1] {\tiny{-1}};
\node[left] at (-0.05,1) [scale=1] {\tiny{$\ell$}};
\node[left] at (0,1.5) [scale=1] {\tiny{+1}};

\node[above] at (0.5,1.55) [scale=1] {\tiny{1}};
\node[above] at (1,1.55) [scale=1] {\tiny{2}};
\node[above] at (1.5,1.55) [scale=1] {\tiny{3}};

\node[below] at (0.75,0) [scale=0.75] {\small{Fig. 3.9 (a). Case 2.}};

}
\end{scope}
    
\begin{scope}[xshift=2.25cm]
{
\draw[gray,very thin, step=0.5cm, opacity=0.5] (0,0)grid (1.5,1.5);

\fill[green!50!white, opacity=0.5] (0,0.5)--++(1.5,0)--++(0,0.5)--++(-1.5,0);

\draw[blue, line width =0.5mm] (1,1)--++(0.5,0);
\draw[orange, line width =0.5mm] (0.5,0.5)--++(0.5,0);

\draw[yellow, line width =0.2mm] (0,0)--++(0,1.5);

\draw[fill=blue] (1,1) circle [radius=0.05];
\draw[fill=blue] (1,1.5) circle [radius=0.05];
\draw[fill=blue] (0.5,1) circle [radius=0.05];
\draw[fill=blue] (1.5,1) circle [radius=0.05];

\draw[fill=orange] (1,0.5) circle [radius=0.05];
\draw[fill=orange] (0.5,0.5) circle [radius=0.05];

{

\draw[black, line width=0.15mm] (0.95, 0.7)--++(0.1,0);
\draw[black, line width=0.15mm] (0.95, 0.75)--++(0.1,0);
\draw[black, line width=0.15mm] (0.95, 0.8)--++(0.1,0);

\draw[black, line width=0.15mm] (1.2, 0.45)--++(0,0.1);
\draw[black, line width=0.15mm] (1.25, 0.45)--++(0,0.1);
\draw[black, line width=0.15mm] (1.3, 0.45)--++(0,0.1);

\draw[black, line width=0.15mm] (0.45, 0.7)--++(0.1,0);
\draw[black, line width=0.15mm] (0.45, 0.75)--++(0.1,0);
\draw[black, line width=0.15mm] (0.45, 0.8)--++(0.1,0);

}

\node[left] at (0,0.5) [scale=1] {\tiny{-1}};
\node[left] at (-0.05,1) [scale=1] {\tiny{$\ell$}};
\node[left] at (0,1.5) [scale=1] {\tiny{+1}};

\node[above] at (0.5,1.55) [scale=1] {\tiny{1}};
\node[above] at (1,1.55) [scale=1] {\tiny{2}};
\node[above] at (1.5,1.55) [scale=1] {\tiny{3}};

\node[below] at (0.75,0) [scale=0.75] {\small{Fig. 3.9 (b). Case 2.1.}};

}
\end{scope}

\end{tikzpicture}
\end{adjustbox}
\end{wrapfigure}

\null 

\noindent \textit{CASE 2: $v(1,l) \in H_1$.} Then $e(2,3;l) \in H_1$ or $e(2,3;l) \notin H'$.

\null

\noindent \textit{CASE 2.1: $e(2,3;l) \in H_1$.} Note that if $e(2,3;l-1) \in H_2$, then we're back to essentially the same scenario as Case 1.1, so we may assume that $e(2,3;l-1) \notin H_2$. Then $e(1,2;l-1) \in H_2$. It follows that $e(1;l-1,l) \notin H'$ and so $R(0,l-1)$ is not a small cookie of $H$. Now, either $e(1,2;l) \in H_1$, or $e(1,2;l) \notin H'$.

\endgroup

\begingroup
\setlength{\intextsep}{0pt}
\setlength{\columnsep}{10pt}
\begin{wrapfigure}[]{r}{0cm}
\begin{adjustbox}{trim=0cm 0cm 0cm 0cm}
\begin{tikzpicture}[scale=1.5]

\begin{scope}[xshift=0cm]
{
\draw[gray,very thin, step=0.5cm, opacity=0.5] (0,0) grid (1.5,1.5);

\fill[green!50!white, opacity=0.5] (0,0.5)--++(1.5,0)--++(0,0.5)--++(-1.5,0);

\draw[yellow, line width =0.2mm] (0,0)--++(0,1.5);

\draw[blue, line width =0.5mm] (1,1)--++(0.5,0);
\draw[orange, line width =0.5mm] (0.5,0.5)--++(0.5,0);

\draw[blue, line width =0.5mm] (0.5,1)--++(0.5,0);

\draw[fill=blue] (1,1) circle [radius=0.05];
\draw[fill=blue] (1,1.5) circle [radius=0.05];
\draw[fill=blue] (0.5,1) circle [radius=0.05];
\draw[fill=blue] (1.5,1) circle [radius=0.05];

\draw[fill=orange] (1,0.5) circle [radius=0.05];
\draw[fill=orange] (0.5,0.5) circle [radius=0.05];

{

\draw[black, line width=0.15mm] (0.95, 0.7)--++(0.1,0);
\draw[black, line width=0.15mm] (0.95, 0.75)--++(0.1,0);
\draw[black, line width=0.15mm] (0.95, 0.8)--++(0.1,0);

\draw[black, line width=0.15mm] (1.2, 0.45)--++(0,0.1);
\draw[black, line width=0.15mm] (1.25, 0.45)--++(0,0.1);
\draw[black, line width=0.15mm] (1.3, 0.45)--++(0,0.1);

\draw[black, line width=0.15mm] (0.45, 0.7)--++(0.1,0);
\draw[black, line width=0.15mm] (0.45, 0.75)--++(0.1,0);
\draw[black, line width=0.15mm] (0.45, 0.8)--++(0.1,0);
}

\node[left] at (0,0.5) [scale=1] {\tiny{-1}};
\node[left] at (-0.05,1) [scale=1] {\tiny{$\ell$}};
\node[left] at (0,1.5) [scale=1] {\tiny{+1}};

\node[above] at (0.5,1.55) [scale=1] {\tiny{1}};
\node[above] at (1,1.55) [scale=1] {\tiny{2}};
\node[above] at (1.5,1.55) [scale=1] {\tiny{3}};

\node[below] at (0.75,0) [scale=0.75] {\small{Fig. 3.10 (a). Case 2.1(a).}};

}
\end{scope}

\begin{scope}[xshift=2.25cm]
{
\draw[gray,very thin, step=0.5cm, opacity=0.5] (0,0) grid (1.5,1.5);

\fill[green!50!white, opacity=0.5] (0,0.5)--++(1.5,0)--++(0,0.5)--++(-1.5,0);
\fill[green!50!white, opacity=0.5] (0.5,1)--++(0.5,0)--++(0,0.5)--++(-0.5,0);

\draw[yellow, line width =0.2mm] (0,0)--++(0,1.5);

\draw[blue, line width =0.5mm] (1,1)--++(0.5,0);
\draw[orange, line width =0.5mm] (0.5,0.5)--++(0.5,0);

\draw[blue, line width =0.5mm] (0,1)--++(0.5,0)--++(0,0.5)--++(-0.5,0);

\draw[blue, line width =0.5mm] (1,1)--++(0,0.5);

\draw[fill=blue] (1,1) circle [radius=0.05];
\draw[fill=blue] (1,1.5) circle [radius=0.05];
\draw[fill=blue] (0.5,1) circle [radius=0.05];
\draw[fill=blue] (1.5,1) circle [radius=0.05];

\draw[fill=blue] (0.5,1.5) circle [radius=0.05];
\draw[fill=blue] (0,1.5) circle [radius=0.05];
\draw[fill=blue] (0,1) circle [radius=0.05];

\draw[fill=orange] (1,0.5) circle [radius=0.05];
\draw[fill=orange] (0.5,0.5) circle [radius=0.05];

{

\draw[black, line width=0.15mm] (0.95, 0.7)--++(0.1,0);
\draw[black, line width=0.15mm] (0.95, 0.75)--++(0.1,0);
\draw[black, line width=0.15mm] (0.95, 0.8)--++(0.1,0);

\draw[black, line width=0.15mm] (1.2, 0.45)--++(0,0.1);
\draw[black, line width=0.15mm] (1.25, 0.45)--++(0,0.1);
\draw[black, line width=0.15mm] (1.3, 0.45)--++(0,0.1);

\draw[black, line width=0.15mm] (0.7, 0.95)--++(0,0.1);
\draw[black, line width=0.15mm] (0.75, 0.95)--++(0,0.1);
\draw[black, line width=0.15mm] (0.8, 0.95)--++(0,0.1);

\draw[black, line width=0.15mm] (0.45, 0.7)--++(0.1,0);
\draw[black, line width=0.15mm] (0.45, 0.75)--++(0.1,0);
\draw[black, line width=0.15mm] (0.45, 0.8)--++(0.1,0);
}

\node[left] at (0,0.5) [scale=1] {\tiny{-1}};
\node[left] at (-0.05,1) [scale=1] {\tiny{$\ell$}};
\node[left] at (0,1.5) [scale=1] {\tiny{+1}};

\node[above] at (0.5,1.55) [scale=1] {\tiny{1}};
\node[above] at (1,1.55) [scale=1] {\tiny{2}};
\node[above] at (1.5,1.55) [scale=1] {\tiny{3}};

\node[below] at (0.75,0) [scale=0.75] {\small{Fig. 3.10 (b). Case 2.1(b).}};

}
\end{scope}

\end{tikzpicture}
\end{adjustbox}
\end{wrapfigure}

\null

\noindent \textit{CASE 2.1 (a):  $e(1,2;l) \in H_1$.} By Observation 3.4 and Proposition 3.3,  $Z \mapsto R(1,l-1)$ is the cascade we seek. End of Case 2.1(a). 

\null

\noindent \textit{CASE 2.1 (b):  $e(1,2;l) \notin H'$.} Then we have that $e(2;l,l+1) \in H_1$, that $e(1;l,l+1) \in H_1$ and that $R(1,l) \in \text{int}(H)$. It follows that $R(0,l) \in \text{ext}(H)$, so $R(0,l)$ must be a small cookie. Then, after $R(0,l) \mapsto R(1,l)$, we are back to Case 2.1 (a). End of Case 2.1(b). End of Case 2.1.

\endgroup 

\begingroup
\setlength{\intextsep}{0pt}
\setlength{\columnsep}{10pt}
\begin{wrapfigure}[]{l}{0cm}
\begin{adjustbox}{trim=0cm 0cm 0cm 0.5cm}
\begin{tikzpicture}[scale=1.5]

\begin{scope}[xshift=0cm]
{
\draw[gray,very thin, step=0.5cm, opacity=0.5] (0,0) grid (1.5,1.5);

\fill[green!50!white, opacity=0.5] (0,0.5)--++(1.5,0)--++(0,0.5)--++(-1.5,0);
\fill[green!50!white, opacity=0.5] (0.5,0)--++(0.5,0)--++(0,0.5)--++(-0.5,0);

\draw[blue, line width =0.5mm] (0.5,1)--++(0.5,0);

\draw[orange, line width =0.5mm] (1,0)--++(0,0.5);
\draw[green!50!black, line width =0.5mm] (0,0)--++(0.5,0)--++(0,0.5)--++(-0.5,0);

\draw[yellow, line width =0.2mm] (0,0)--++(0,1.5);

\draw[fill=blue] (1,1) circle [radius=0.05];
\draw[fill=blue] (1,1.5) circle [radius=0.05];
\draw[fill=blue] (0.5,1) circle [radius=0.05];

\draw[fill=orange] (1,0.5) circle [radius=0.05];
\draw[fill=orange] (1,0) circle [radius=0.05];

{

\draw[black, line width=0.15mm] (0.95, 0.7)--++(0.1,0);
\draw[black, line width=0.15mm] (0.95, 0.75)--++(0.1,0);
\draw[black, line width=0.15mm] (0.95, 0.8)--++(0.1,0);

\draw[black, line width=0.15mm] (1.2, 0.95)--++(0,0.1);
\draw[black, line width=0.15mm] (1.25, 0.95)--++(0,0.1);
\draw[black, line width=0.15mm] (1.3, 0.95)--++(0,0.1);

\draw[black, line width=0.15mm] (0.7, 0.45)--++(0,0.1);
\draw[black, line width=0.15mm] (0.75, 0.45)--++(0,0.1);
\draw[black, line width=0.15mm] (0.8, 0.45)--++(0,0.1);

\draw[black, line width=0.15mm] (0.45, 0.7)--++(0.1,0);
\draw[black, line width=0.15mm] (0.45, 0.75)--++(0.1,0);
\draw[black, line width=0.15mm] (0.45, 0.8)--++(0.1,0);

}

\node[left] at (0,0.5) [scale=1] {\tiny{-1}};
\node[left] at (-0.05,1) [scale=1] {\tiny{$\ell$}};
\node[left] at (0,1.5) [scale=1] {\tiny{+1}};

\node[above] at (0.5,1.55) [scale=1] {\tiny{1}};
\node[above] at (1,1.55) [scale=1] {\tiny{2}};
\node[above] at (1.5,1.55) [scale=1] {\tiny{3}};

\node[below] at (0.75,0) [scale=0.75] {\small{Fig. 3.11. Case 2.2.}};

}
\end{scope}

\end{tikzpicture}
\end{adjustbox}
\end{wrapfigure}

\null 

\noindent \textit{CASE 2.2: $e(2,3;l) \notin H'$.} Then we have that $e(1,2;l) \in H_1$. Note that if $e(1,2;l-1) \in H_2$, then we're back to Case 2.1, so we may assume that $e(1,2;l-1) \notin H'$. It follows that $e(2;l-2,l-1) \in H_2$. Note that if $e(1; l-1,l) \in H_1$, then $R(0,l-1) \in \text{ext(H)}$ and $R(0,l) \in \text{ext(H)}$, contradicting Lemma 1.14, so we may assume that $e(1; l-1,l) \notin H'$. Then we must have that $e(1;l-2,l-1) \in H'$. This implies that $R(0,l-2)$ is a small cookie of $H$. Then after $R(0,l-2) \mapsto R(1,l-2)$, we're back to Case 2.1(a). End of Case 2.2. End of Case 2. $\square$.

\endgroup 

\null

\noindent \textbf{Observation 3.6.} Let $H$ be a Hamiltonian cycle of an $m \times n$ grid graph $G$, with $m,n \geq 5$ and let $J$ be a large cookie of $G$. Then $J \cap R_2 \neq \emptyset$.

\null

\noindent \textbf{Lemma 3.7.} Let $H$ be a Hamiltonian cycle of an $m \times n$ grid graph $G$, with $m,n \geq 5$ and assume that $G$ has at least two large cookies $J_1$ and $J_2$. Then switching the neck $N_{J_1}$ of $J_1$ splits $H$ into two cycles $H_1$ and $H_2$ such that there is $v_1 \in H_1 \cap R_2$ and $v_2 \in H_2 \cap R_2$ with $v_1$ adjacent to $v_2$.

\null

\noindent \textit{Proof.} Orient $H$. Let $\{v_x, v_y\}$ be the boundary edge of the neck $N_{J_1}$ of $J_1$. Define $\overrightarrow{K}_1$ and $\overrightarrow{K}_2$ to be the subtrails $\overrightarrow{K}((v_x,v_{x+1})$, $(v_{y-1}, v_y))$ and $\overrightarrow{K}((v_y,v_{y+1}), (v_{x-1},v_x))$ of $\overrightarrow{K}_H$, respectively. By Lemma 1.16 (i), switching $N_{J_1}$ gives two cycles $H_1$ and $H_2$, with  $V(H_1)=V(\overrightarrow{K}_1 \setminus \{(v_x, v_y)\})$ and $V(H_2)=V(\overrightarrow{K}_2)$. By Lemma 1.13, $V(J_1) = V(\overrightarrow{K}_1)$. By Observation 3.6, $V(J_1) \cap R_2 \neq \emptyset$. Then $V(H_1) \cap R_2 \neq \emptyset$.

Since $V(J_1) = V(\overrightarrow{K}_1)$, we have that $V(J_2) \subseteq V(\overrightarrow{K}_2)=V(H_2)$. By Observation 3.6, $V(J_2) \cap R_2 \neq \emptyset$. It follows that $ V(H_2) \cap R_2 \neq \emptyset$. 

We have shown that $V(H_1) = R_2 \neq \emptyset$ and that $ V(H_2) \cap R_2 \neq \emptyset$. It remains to check there is $v_1 \in H_1 \cap R_2$ and $v_2 \in H_2 \cap R_2$ with $v_1$ adjacent to $v_2$. $v_1\in H_1 \cap R_2$ and $R_2 = w_1, ..., w_s$, with $v_1=w_1$. Sweep $R_2$ starting at $w_1$. If there is $i \in \{1,...,s-2\} $ such that $v_{i} \in H_1$ and $v_{i+1} \in H_2$, we are done. If there is no such $i$ then $R_2 \cup H_1$, contradicting that $R_2 \cup H_2 \neq \emptyset$. $\square$

\null

\noindent \textbf{Proposition 3.8. (MLC Algorithm.)} Let $H$ be a Hamiltonian cycle of an $m \times n$ grid graph $G$. Assume that $G$ has more than one large cookie. Then there is a cascade of length at most two that reduces the number of large cookies of $G$ by one. 

\null 

\noindent \textit{Proof.} Let $J$ be a large cookie of $G$ with neck $N_J$. Switching $N_J$ splits $H$ into the cycles $H_1$ and $H_2$. By Lemma 3.7 there is $v_1 \in H_1 \cap R_2$ and $v_2 \in H_2\cap R_2$ with $v_1$ adjacent to $v_2$. By Lemma 3.5 there is a cascade of length at most two, whose last move is $N_j \mapsto N_J'$ with $N_J \neq N_J'$. By Observation 3.4, this cascade decreases the number of large cookies of $G$ by one. $\square$

\subsection{Existence of the 1LC Algorithm}

\begingroup
\setlength{\intextsep}{0pt}
\setlength{\columnsep}{10pt}
\begin{wrapfigure}[]{l}{0cm}
\begin{adjustbox}{trim=0cm -0.25cm 0cm 0cm}
\begin{tikzpicture}[scale=1.25]

\begin{scope}[xshift=0cm, yshift=0cm]
{
\draw[gray,very thin, step=0.5cm, opacity=0.5] (0,0) grid (1.5,1.5);

\draw[blue, line width=0.5mm](0.5,0.5)--++(0,0.5)--++(0.5,0)--++(0,-0.5);

\node[left] at (0,0.5) [scale=1]
{\tiny{$\ell$}};

\node[above] at (0.5,1.5) [scale=1]
{\tiny{k}};

\node[below] at (0.75, 0) [scale=0.75] {\small{}};

\node[below] at (0.75, 0) [scale=0.75] {\small{\begin{tabular}{c} Fig. 3.12(a). \\ A northern leaf. \end{tabular}}};;

}
\end{scope}

\begin{scope}[xshift=2.25cm, yshift=0cm]
{
\draw[gray,very thin, step=0.5cm, opacity=0.5] (0,0) grid (1.5,1.5);

\draw[blue, line width=0.5mm](0,0.5)--++(0.5,0)--++(0,0.5)--++(0.5,0)--++(0,-0.5)--++(0.5,0);

\node[left] at (0,0.5) [scale=1]
{\tiny{$\ell$}};

\node[above] at (0.5,1.5) [scale=1]
{\tiny{k}};

\node[below] at (0.75, 0) [scale=0.75] {\small{Fig. 3.12 (c). $A_0$-type.}};

}
\end{scope}

\begin{scope}[xshift=0cm, yshift=-2.5cm]
{
\draw[gray,very thin, step=0.5cm, opacity=0.5] (0,0) grid (1.5,1.5);

\draw[blue, line width=0.5mm](0,0.5)--++(0.5,0)--++(0,0.5);
\draw[blue, line width=0.5mm](1.5,0.5)--++(-0.5,0)--++(0,0.5);

\node[left] at (0,0.5) [scale=1]
{\tiny{$\ell$}};

\node[above] at (0.5,1.5) [scale=1]
{\tiny{k}};

\node[below] at (0.75, 0) [scale=0.75] {\small{Fig. 3.12 (b). $A$-type.}};

}
\end{scope}

\begin{scope}[xshift=2.25cm, yshift=-2.5cm]
{
\draw[gray,very thin, step=0.5cm, opacity=0.5] (0,0) grid (1.5,1.5);

\draw[blue, line width=0.5mm](0,0.5)--++(0.5,0)--++(0,1);
\draw[blue, line width=0.5mm](1.5,0.5)--++(-0.5,0)--++(0,1);

\node[left] at (0,0.5) [scale=1]
{\tiny{$\ell$}};

\node[above] at (0.5,1.5) [scale=1]
{\tiny{k}};

\node[below] at (0.75, 0) [scale=0.75] {\small{Fig. 3.12 (d). $A_1$-type.}};

}
\end{scope}

\begin{scope}[xshift=4.55cm, yshift=-1.75cm]
{
\draw[gray,very thin, step=0.5cm, opacity=0.5] (0,0) grid (1.5,3);

\draw[blue, line width=0.5mm] (0.5,0.5)--++(0,0.5)--++(0.5,0)--++(0,-0.5);

\draw[blue, line width=0.5mm] (0,1.5)--++(0.5,0)--++(0,1); 

\draw[blue, line width=0.5mm] (1.5,1.5)--++(-0.5,0)--++(0,1);



\begin{scope}
[very thick,decoration={
    markings,
    mark=at position 1 with {\arrow{>}}}
    ]   
\end{scope}

{
\node[left] at (0,0.5) [scale=1]
{\tiny{-1}};
\node[left] at (0,1) [scale=1]
{\tiny{$\ell$}};
\node[left] at (0,1.5) [scale=1]
{\tiny{+1}};
\node[left] at (0,2) [scale=1]
{\tiny{+2}};
\node[left] at (0,2.5) [scale=1]
{\tiny{+3}};

\node[above] at (0.5,3) [scale=1]
{\tiny{$k$}};
\node[above] at (1,3) [scale=1]
{\tiny{+1}};

}

\node[below] at (0.75, 0) [scale=0.75] {\small{\begin{tabular}{c} Fig. 3.12(e).  Northern:  \\leaf followed by a \\  northern $A_1$ type. \end{tabular}}};;

}
\end{scope}

\end{tikzpicture}
\end{adjustbox}
\end{wrapfigure}

\textbf{Definitions}.  Let $G$ be an $m \times n$ grid graph and $H$ be a Hamiltonian cycle of $G$. We call a subpath of $H$ on the edges $e(k;l,l+1), e(k,k+1;l+1)$ and $e(k+1;l,l+1)$ a \textit{northern leaf}. We will often say that $R(k,l)$ is a northern leaf to mean that $e(k;l,l+1), e(k,k+1;l+1)$ and $e(k+1;l,l+1)$ belong to $H$. Southern, eastern and western leaves are defined analogously. 

We call the subgraph of $H$ on the edges $e(k-1,k;l), e(k;l,l+1)$, $e(k+1,k+2;l)$, and $e(k+1;l,l+1)$ a \textit{northern $A$-type}. Suppose $H$ has a northern $A$-type. We call the subgraph $A \cup e(k,k+1;l+1)$ of $H$ a \textit{northern $A_0$-type}, and we call the subgraph $A \cup e(k;l+1,l+2) \cup e(k+1;l+1,l+2)$ of $H$ a \textit{northern $A_1$-type}. We make analogous definitions for eastern, southern and western $A$-types. See Figures 3.12.

Let $R(k,l-1)$ be a northern leaf. If $H$ has a northern $A$-type on $e(k-1,k;l+1), e(k;l+1,l+2)$ and $e(k+1,k+2;l+1), e(k+1;l+1,l+2)$ then we say that $A$-type \textit{follows} the northern leaf $R(k,l-1)$ \textit{northward}. We call the switchable box $R(k,l+1)$ the \textit{switchable middle-box} of the $A_1$-type. Analogous definitions apply for other compass directions.

\endgroup 

\noindent Let $A$ be a northern $A_0$-type in $H$ on the edges $e(k-1,k;0)$, $e(k;0,1)$, $e(k,k+1,1)$,  $e(k+1;0,1)$, $e(k+1,k+2;0)$ and let $j \in \big\{1,2,..\big\lfloor \frac{n}{2} \big\rfloor \big\}$. We define a \textit{northern j-stack of $A_0$'s starting at $A$} to be a subgraph \textit{$stack(j; A_0)$} of $H$, where $\textrm{stack}(j; A_0)= \bigcup_{i=0}^{j-1} \Big( A +(0,2i) \Big)$. If $j=\big\lfloor \frac{n}{2}   \big\rfloor$, we call the $j$-stack a \textit{full j-stack of $A_0$'s}. We note that $j$ is the number of $A_0$'s in $\textrm{stack}(j; A_0)$. Eastern, southern, and western $j$-stacks are defined analogously.

We denote the set of northern and southern small cookies by $SmallCookies\{N,S\}$ and the set of eastern and western small cookies by $SmallCookies\{E,W\}$. Assume that $C \in SmallCookies\{N,S\}$ is an easternmost or westernmost southern or northern small cookie. Then we call $C$ and \textit{outermost} small cookie in $SmallCookies\{N,S\}$. Outermost small cookies in $SmallCookies\{E,W\}$ are defined analogously.

Let $R(k,l)$ be a northern leaf. We say that the cascade $\mu_1, ..., \mu_r$ \textit{collects} $R(k,l)$, if $\mu_r$ is the move $Z \mapsto R(k,l)$. Note that, since $R(k,l)$ is not switchable, $Z$ must be a switchable box adjacent to $R(k,l)$.

\begingroup
\setlength{\intextsep}{0pt}
\setlength{\columnsep}{10pt}
\begin{wrapfigure}[]{r}{0cm}
\begin{adjustbox}{trim=0cm 0.25cm 0cm 0cm}
\begin{tikzpicture}[scale=1.25]

\begin{scope} 
{
\draw[gray,very thin, step=0.5cm, opacity=0.5] (0,0) grid (1.5,2);

\draw[yellow, line width =0.2mm] (0,0) -- (1.5,0);

{
\fill[blue!50!white, opacity=0.5] (0,1) rectangle (0.5,1.5);
\fill[blue!50!white, opacity=0.5] (1,1) rectangle (1.5,1.5);
\fill[blue!50!white, opacity=0.5] (0,1.5) rectangle (0.5,2);
\fill[blue!50!white, opacity=0.5] (1,1.5) rectangle (1.5,2);

\draw[blue, line width=0.5mm] (0,1)--++(0.5,0)--++(0,0.5); 
\draw[blue, line width=0.5mm] (1,1.5)--++(0,-0.5)--++(0.5,0); 

\fill[blue!40!white, opacity=0.4](0.5,0)--(0.5,0.5)--(1,0.5)--(1,0);
\draw[blue, line width=0.5mm] (0.5,0)--(0.5,0.5)--(1,0.5)--(1,0);
}
{
\foreach \x in {1, ...,4}
\node[left] at (0,0.5*\x) [scale=1]
{\tiny{\x}};

\node[above] at (0.5,2) [scale=1]
{\tiny{k}};
\node[above] at (1,2) [scale=1]
{\tiny{+1}};
\node[above] at (1.5,2) [scale=1]
{\tiny{+2}};

}
{
\draw[black, line width=0.15mm] (1.2,1.45)--(1.2,1.55);
\draw[black, line width=0.15mm] (1.25,1.45)--(1.25,1.55);
\draw[black, line width=0.15mm] (1.3,1.45)--(1.3,1.55);

\draw[black, line width=0.15mm] (0.2,1.45)--(0.2,1.55);
\draw[black, line width=0.15mm] (0.25,1.45)--(0.25,1.55);
\draw[black, line width=0.15mm] (0.3,1.45)--(0.3,1.55);

\draw[black, line width=0.15mm] (0.7,0.95)--(0.7,1.05);
\draw[black, line width=0.15mm] (0.75,0.95)--(0.75,1.05);
\draw[black, line width=0.15mm] (0.8,0.95)--(0.8,1.05);
}

\node[below] at (0.75, 0) [scale=0.75] {\small{\begin{tabular}{c} Fig. 3.13. \end{tabular}}};;

}
\end{scope}

\end{tikzpicture} 
\end{adjustbox}
\end{wrapfigure}

\null 

\noindent Given a small cookie $C$, we want to show that there is a cascade that collects $C$. For definiteness, assume that $C$ is the northern small cookie $R(k,0)$. If $e(k,k+1; 2) \in H$, then $C+(0,1) \mapsto C$ is the cascade we seek, so we only need to consider the case where $e(k,k+1; 2) \notin H$. Then we must have $e(k-1,k;2) \in H$, $e(k+1,k+2; 2) \in H$, $e(k;2,3) \in H$ and $e(k+1; 2,3) \in H$. Note that if $e(k-1,k;3) \in H$, then $C +(0,2) \mapsto C+(-1,2)$ followed by $C+(0,1) \mapsto C$ is the cascade we seek, so we consider the case where $e(k-1,k;3) \notin H$ and, by symmetry, where $e(k+1,k+2; 3) \notin H$. 
Now, we either have $e(k;3,4) \in H$  and $e(k+1;3,4) \in H$ or  $e(k,k+1;3) \in H$. That is, $C$ is followed northward by an $A_0$-type or by an $A_1$-type. See Figure 3.13. From this point onward, we will omit the compass direction when it does not introduce ambiguity. We coalesce this paragraph into the following lemma:

\endgroup

\null

\noindent \textbf{Proposition 3.10. (1LC Algorithm).} Let $G$ be an $m \times n$ grid graph, let $H$ be a Hamiltonian cycle of $G$. If $H$ has exactly one large cookie and at least one small cookie, then there is a cascade of length at most $\frac{1}{2}\max(m,n)+\min(m,n)+2$ moves that reduces the number of small cookies of $H$ by one and such that it does not increase the number of large cookies.

\null 

\noindent The proof of Proposition 3.10 requires the following two Lemmas.

\null 

\noindent \textbf{Lemma 3.11.} Let $G$ be an $m \times n$ grid graph, let $H$ be a Hamiltonian cycle of $G$, and let $C \in SmallCookies\{N,S\}$ be an easternmost small cookie. Assume that $G$ has only one large cookie. Then there cannot be a full $j$-stack of $A_0's$ starting at the $A_0$-type that contains $C$. 

\null

\noindent \textbf{Lemma 3.12.} Let $G$ be an $m \times n$ grid graph, let $H$ be a Hamiltonian cycle of $G$, and let $C \in SmallCookies\{N,S\}$ be an easternmost small cookie. Assume that there is a $j$-stack of $A_0$ starting at the $A_0$-type containing $C$. Let $L$ be the leaf in the top ($j^{\text{th}}$) $A_0$ of the stack. Assume that $L$ is followed by an $A_1$-type and that $G$ has only one large cookie. Then there is a cascade of at most $\min(m,n)+3$ moves that collects $L$.

\null

\noindent \textit{Proof of Proposition 3.10.} Since there is at least one small cookie, at least one of $SmallCookies\{N,S\}$ and $SmallCookies\{E,W\}$ is nonempty. WLOG assume that $SmallCookies\{N,S\}$ is nonempty. Let $C \in SmallCookies\{N,S\}$ be an easternmost small cookie. We will require the following two lemmas.

\setlength{\intextsep}{0pt}
\setlength{\columnsep}{20pt}
\begin{center}
\begin{adjustbox}{trim=0cm 0.25cm 0cm 0cm}
\begin{tikzpicture}[scale=1.35]

\begin{scope}[xshift=0cm]
 \draw[gray,very thin, step=0.5cm, opacity=0.5] (0,0) grid (1.5,3);

\draw[blue, line width=0.5mm] (0,0)--++(0.5,0)--++(0,0.5)--++(0.5,0)--++(0,-0.5)--++(0.5,0);

\draw[blue, line width=0.5mm] (0,1)--++(0.5,0)--++(0,0.5)--++(0.5,0)--++(0,-0.5)--++(0.5,0);

\draw[blue, line width=0.5mm] (0,2)--++(0.5,0)--++(0,0.5)--++(0.5,0)--++(0,-0.5)--++(0.5,0);

\draw[blue, line width=0.5mm] (0.5,3)--++(0.5,0);

\draw [->,black, very thick] (1.6,1.75)--(2.4,1.75);
\node[above] at  (2,1.75) [scale=0.8]{\small{$\mu_s$}};

\node at  (0.75, 2.25) [scale=0.7]{\small{$L_{j_0-1}$}};
\node at  (0.75, 1.25) [scale=0.7]{\small{$L_{j_0-2}$}};
\node at  (0.75, 0.25) [scale=0.7]{\small{$L_1$}};

\node[below] at (0.75,0) [scale=0.8]{\small{\begin{tabular}{c} Fig. 3.14 (a). \end{tabular}}};;   
\end{scope}
\begin{scope}[xshift=2.5cm]
 \draw[gray,very thin, step=0.5cm, opacity=0.5] (0,0) grid (1.5,3);

\draw[blue, line width=0.5mm] (0,0)--++(0.5,0)--++(0,0.5)--++(0.5,0)--++(0,-0.5)--++(0.5,0);

\draw[blue, line width=0.5mm] (0,1)--++(0.5,0)--++(0,0.5)--++(0.5,0)--++(0,-0.5)--++(0.5,0);

\draw[blue, line width=0.5mm] (0,2)--++(1.5,0);
\draw[blue, line width=0.5mm] (0.5,3)--++(0,-0.5)--++(0.5,0)--++(0,0.5);

\draw [->,black, very thick] (1.6,1.75)--(2.4,1.75);
\node[above] at  (2,1.75) [scale=0.8]{\small{$\mu_{s+1}$}};

\node at  (0.75, 1.25) [scale=0.7]{\small{$L_{j_0-2}$}};
\node at  (0.75, 0.25) [scale=0.7]{\small{$L_1$}};

\node[below] at (0.75,0) [scale=0.8]{\small{\begin{tabular}{c} Fig. 3.14 (a). \end{tabular}}};;   
\end{scope}
\begin{scope}[xshift=5cm]
 \draw[gray,very thin, step=0.5cm, opacity=0.5] (0,0) grid (1.5,3);

\draw[blue, line width=0.5mm] (0,0)--++(0.5,0)--++(0,0.5)--++(0.5,0)--++(0,-0.5)--++(0.5,0);

\draw[blue, line width=0.5mm] (0,1)--++(1.5,0);

\draw[blue, line width=0.5mm] (0,2)--++(0.5,0)--++(0,-0.5)--++(0.5,0)--++(0,0.5)--++(0.5,0);

\draw[blue, line width=0.5mm] (0.5,3)--++(0,-0.5)--++(0.5,0)--++(0,0.5);

\draw [->,black, very thick] (1.6,1.75)--(2.4,1.75);
\node[above] at  (2,1.75) [scale=0.8]{\small{$\mu_{s+2}$}};

\node at  (0.75, 0.25) [scale=0.7]{\small{$L_1$}};

\node[below] at (0.75,0) [scale=0.8]{\small{\begin{tabular}{c} Fig. 3.14 (c). \end{tabular}}};;   
\end{scope}
\begin{scope}[xshift=7.5cm]
 \draw[gray,very thin, step=0.5cm, opacity=0.5] (0,0) grid (1.5,3);

\draw[blue, line width=0.5mm] (0,0)--++(1.5,0);

\draw[blue, line width=0.5mm] (0,1)--++(0.5,0)--++(0,-0.5)--++(0.5,0)--++(0,0.5)--++(0.5,0);

\draw[blue, line width=0.5mm] (0,2)--++(0.5,0)--++(0,-0.5)--++(0.5,0)--++(0,0.5)--++(0.5,0);

\draw[blue, line width=0.5mm] (0.5,3)--++(0,-0.5)--++(0.5,0)--++(0,0.5);

\node[below] at (0.75,0) [scale=0.8]{\small{\begin{tabular}{c} Fig. 3.14 (d). \end{tabular}}};;   
\end{scope}

\end{tikzpicture}
\end{adjustbox}
\end{center}

\noindent For definiteness assume that $C$ is a small northern cookie on the southern boundary.  Let $Q(j)$ be the statement ``There is a $j$-stack of $A_0$'s starting at the $A_0$-type containing $C$''. Note that $C$ is contained in an $A_0$-type, so $Q(1)$ is true. By Lemma 3.11,  there is a $j_0 \in \Big\{2,3, ...\Big\lfloor \frac{n}{2} \Big\rfloor \Big\}$ such that for each $j \in \{2.., j_0-1\}$, $Q(j)$ is true for each $j < j_0$ but $Q(j_0)$ is not true.

{
}

\noindent For $j \in \{1.., j_0-1\}$ let $L_j$ be the northern leaf of the $j^{\text{th}}$ $A_0$-type in the stack. Note that $L_1=C$. Lemma 3.9 implies that $L_{j_0-1}$ is either followed by an $A_1$-type or there is a cascade that collects $L_{j_0-1}$. If $L_{j_0-1}$ is followed by an $A_1$-type, then by Lemma 3.12, we can find a cascade that collects it, so we only need to check the case in which there is a cascade $\mu_1, ..., \mu_s$ that collects $L_{j_0-1}$. Note that $\mu_s$ must be the move $L_{j_0-1}+(1,0) \mapsto L_{j_0-1}$.
Then $\mu_1, ..., \mu_s, L_{j_0-2}+(0,1) \mapsto L_{j_0-2} , ..., L_1+(0,1) \mapsto L_1$ is a cascade that collects $C$. See Figures 3.14 for an illustration with $j_0-1=3$. Note that $j \leq \frac{n}{2}$, and that by Lemma 3.12, there are at most $\min(m,n)+3$ moves required to collect $L$. After that, we need at most another $j-1$ flips to collect $C$, so $C$ can be collected after at most $\frac{1}{2}\max(m,n) + \min(m,n)+2$ moves. See Figure 3.14 for an illustration with $j_0-1=3$. $\square$

\null

\begingroup
\setlength{\intextsep}{0pt}
\setlength{\columnsep}{10pt}
\begin{wrapfigure}[]{r}{0cm}
\begin{adjustbox}{trim= 0cm 0.5cm 0cm 0.5cm}
\begin{tikzpicture}[scale=1.25]

\begin{scope}[xshift=0cm, yshift=0cm ]

\draw[gray,very thin, step=0.5cm, opacity=0.5] (0,0) grid (1.5,4);

\fill[blue!40!white, opacity=0.5] (0,1) rectangle (0.5,1.5);
\fill[blue!40!white, opacity=0.5] (1,1) rectangle (1.5,1.5);
\fill[blue!40!white, opacity=0.5] (0.5,2) rectangle (1,2.5);
\fill[blue!40!white, opacity=0.5] (0.5,0) rectangle (1,0.5);

\fill[blue!40!white, opacity=0.5] (0,3.5) rectangle (1.5,4);

\fill[blue!40!white, opacity=0.5] (0,3) rectangle (0.5,3.5);
\fill[blue!40!white, opacity=0.5] (1,3) rectangle (1.5,3.5);

\draw[green!50!black, dotted, line width=0.5mm] (1.5,0)--++(0,1);

\draw[blue, line width=0.5mm] (0,1)--++(0.5,0)--++(0,0.5)--++(0.5,0)--++(0,-0.5)--++(0.5,0);

\draw[blue, line width=0.5mm] (0,2)--++(0.5,0)--++(0,0.5)--++(0.5,0)--++(0,-0.5)--++(0.5,0);

\draw[blue, line width=0.5mm] (0,3)--++(0.5,0)--++(0,0.5)--++(0.5,0)--++(0,-0.5)--++(0.5,0)--++(0,0.5);

\draw[blue, line width=0.5mm] (0.5,0)--(0.5,0.5)--(1,0.5)--(1,0)--++(0.5,0);

\draw[green!50!black, thick] (1.5,1.5) circle [radius=0.05];

{
\tikzset
  {
    myCircle/.style=    {orange}
  }
 
\foreach \x in {0,...,2}
\fill[myCircle, blue,] (0.75,1.62575+0.125*\x) circle (0.035);
}

{
\node[left] at (0,4) [scale=1]
{\tiny{n-1}};
\node[left] at (0,3.5) [scale=1]
{\tiny{-2}};

\node[above] at (0.5,4) [scale=1]
{\tiny{k}};
\node[above] at (1,4) [scale=1]
{\tiny{+1}};
\node[above] at (1.5,4) [scale=1]
{\tiny{+2}};

}

\node[below] at (0.75,0) [scale=0.75] {\small Fig 3.15(a). Case 1.1.};

\end{scope}

\begin{scope}[xshift=2.25cm, yshift=0cm ]

\draw[gray,very thin, step=0.5cm, opacity=0.5] (0,0) grid (1.5,3.5);

{
\draw[blue, line width=0.5mm] (0.5,0)--(0.5,0.5)--(1,0.5)--(1,0)--++(0.5,0);

\draw[blue, line width=0.5mm] (0,1)--++(0.5,0)--++(0,0.5)--++(0.5,0)--++(0,-0.5)--++(0.5,0);

\draw[blue, line width=0.5mm] (0,2)--++(0.5,0)--++(0,0.5)--++(0.5,0)--++(0,-0.5)--++(0.5,0);

\draw[blue, line width=0.5mm] (0,3)--++(0.5,0)--++(0,0.5)--++(0.5,0)--++(0,-0.5)--++(0.5,0)--++(0,0.5);

\fill[blue!40!white, opacity=0.5] (0,1) rectangle (0.5,1.5);
\fill[blue!40!white, opacity=0.5] (1,1) rectangle (1.5,1.5);
\fill[blue!40!white, opacity=0.5] (0.5,2) rectangle (1,2.5);
\fill[blue!40!white, opacity=0.5] (0.5,0) rectangle (1,0.5);

\fill[blue!40!white, opacity=0.5] (0,3) rectangle (0.5,3.5);
\fill[blue!40!white, opacity=0.5] (1,3) rectangle (1.5,3.5);

}

{
\tikzset
  {
    myCircle/.style=    {orange}
  }
 
\foreach \x in {0,...,2}
\fill[myCircle, blue,] (0.75,1.62575+0.125*\x) circle (0.035);
}

{
\node[left] at (0,3.5) [scale=1]
{\tiny{n-1}};
\node[left] at (0,3) [scale=1]
{\tiny{-2}};

\node[above] at (0.5,3.5) [scale=1]
{\tiny{k}};
\node[above] at (1,3.5) [scale=1]
{\tiny{+1}};

}

\node[below] at (0.75,0) [scale=0.75] {\small Fig 3.15(b). Case 1.2.};

\end{scope}

\end{tikzpicture}
\end{adjustbox}
\end{wrapfigure}

\noindent It remains to prove Lemmas 3.11 and 3.12. 

\null 

\noindent \textit{Proof of Lemma 3.11.} Assume for contradiction that $C$ is in a full stack of $A_0$'s starting at the $A_0$ that contains $C$. For definiteness, assume that $C=R(k,0)$ is a small northern cookie on the southern boundary. First we check that $m-1 >k+2$. If $m-1=k+2$, then we must have $e(k+2;0,1) \in H$ and $e(k+2;1,2) \in H$. But then $H$ misses $v(k+2,3)$ (in green in Figure 3.15 (a)). The number $j$ of $A_0$'s in the full stack is even or an odd so there are two cases to check. Note that for each odd $i \in \{1, 2, ..., j \}$, the leaf of 
the i$^{\text{th}}$ $A_0$ belongs to $\text{ext}(H)$ and for each even $i \in \{1, 2, ..., j\}$, the leaf of the i$^{\text{th}}$ $A_0$ belongs to $\text{int}(H)$.

\endgroup

\null 

\noindent \textit{CASE 1: $j$ is even.} Note that the top leaf of the stack is in $\text{int}(H)$. Now, $n-1$ is either even or odd. 

\null 

\noindent \textit{CASE 1.1: $n-1$ is even.} We have that $R(k,n-3) \in \text{int}(H)$. But then we must have $R(k,n-2) \in \text{ext}(H)$ and $R(k+1,n-2) \in \text{ext}(H)$, contradicting Lemma 1.14. End of Case 1.1. See Figure 3.15(a).

\null 

\noindent \textit{CASE 1.2: $n-1$ is odd.} Then we must have that $R(k+1,n-2)$ is a small southern cookie. But this contradicts our assumption that $C$
is the easternmost small cookie in $SmallCookies\{N,S\}$.  End of Case 1.2. End of Case 1. See Figure 1.15 (b).

\null 

\begingroup
\setlength{\intextsep}{0pt}
\setlength{\columnsep}{20pt}
\begin{wrapfigure}[]{l}{0cm}
\begin{adjustbox}{trim= 0cm 0.5cm 0cm 0.5cm}
\begin{tikzpicture}[scale=1.25]

\begin{scope}[xshift=0cm, yshift=0cm ]

\draw[gray,very thin, step=0.5cm, opacity=0.5] (0,0) grid (1.5,2.5);

{
\draw[blue, line width=0.5mm] (0.5,0)--(0.5,0.5)--(1,0.5)--(1,0)--++(0.5,0);

\draw[blue, line width=0.5mm] (0,1)--++(0.5,0)--++(0,0.5)--++(0.5,0)--++(0,-0.5)--++(0.5,0);

\draw[blue, line width=0.5mm] (0,2)--++(0.5,0)--++(0,0.5)--++(0.5,0)--++(0,-0.5)--++(0.5,0);


\fill[blue!40!white, opacity=0.5] (0,1) rectangle (0.5,1.5);
\fill[blue!40!white, opacity=0.5] (1,1) rectangle (1.5,1.5);
\fill[blue!40!white, opacity=0.5] (0.5,2) rectangle (1,2.5);
\fill[blue!40!white, opacity=0.5] (0.5,0) rectangle (1,0.5);

}

{
\tikzset
  {
    myCircle/.style=    {orange}
  }
 
\foreach \x in {0,...,2}
\fill[myCircle, blue,] (0.75,1.62575+0.125*\x) circle (0.035);
}

{
\node[left] at (0,2.5) [scale=1]
{\tiny{n-1}};
\node[left] at (0,2) [scale=1]
{\tiny{-2}};

\node[above] at (0.5,2.5) [scale=1]
{\tiny{k}};
\node[above] at (1,2.5) [scale=1]
{\tiny{+1}};

}

\node[below] at (0.75,0) [scale=0.75] {\small Fig 3.16(a). Case 2.1.};

\end{scope}

\begin{scope}[xshift=2.25cm, yshift=0cm ]

\draw[gray,very thin, step=0.5cm, opacity=0.5] (0,0) grid (2,3);

{
\draw[blue, line width=0.5mm] (0.5,0)--(0.5,0.5)--(1,0.5)--(1,0)--++(0.5,0);

\draw[blue, line width=0.5mm] (0,1)--++(0.5,0)--++(0,0.5)--++(0.5,0)--++(0,-0.5)--++(0.5,0);

\draw[blue, line width=0.5mm] (0,2)--++(0.5,0)--++(0,0.5)--++(0.5,0)--++(0,-0.5)--++(0.5,0);


\fill[blue!40!white, opacity=0.5] (0,1) rectangle (0.5,1.5);
\fill[blue!40!white, opacity=0.5] (1,1) rectangle (1.5,1.5);
\fill[blue!40!white, opacity=0.5] (0.5,2) rectangle (1,2.5);
\fill[blue!40!white, opacity=0.5] (0.5,0) rectangle (1,0.5);

\fill[orange!50!white, opacity=0.5] (1.5,0.5) rectangle (2,1);
\fill[orange!50!white, opacity=0.5] (1.5,2.5) rectangle (2,3);

\draw[orange!90!black, line width=0.5mm] (1.5,0)--++(0.5,0);

\draw[orange!90!black, line width=0.5mm] (1.5,0)--++(0.5,0);
\draw[orange!90!black, line width=0.5mm] (1.5,1)--++(0,-0.5)--++(0.5,0);

\draw[orange!90!black, line width=0.5mm] (1.5,2)--++(0,-0.5)--++(0.5,0);

\draw[orange!90!black, line width=0.5mm] (0,3)--++(1.5,0)--++(0,-0.5)--++(0.5,0)--++(0,0.5);

}

{
\tikzset
  {
    myCircle/.style=    {orange}
  }
 
\foreach \x in {0,...,2}
\fill[myCircle, blue,] (0.75,1.62575+0.125*\x) circle (0.035);

\foreach \x in {0,...,2}
\fill[myCircle, orange,] (1.75,1.12575+0.125*\x) circle (0.035);
}

{
\node[left] at (0,3) [scale=1]
{\tiny{n-1}};
\node[left] at (0,2.5) [scale=1]
{\tiny{-2}};

\node[left] at (0,0) [scale=1] {\tiny{0}};
\node[left] at (0,0.5) [scale=1] {\tiny{1}};
\node[left] at (0,1) [scale=1] {\tiny{2}};

\node[above] at (0.5,3) [scale=1]
{\tiny{k}};
\node[above] at (1,3) [scale=1]
{\tiny{+1}};
\node[above] at (1.5,3) [scale=1]
{\tiny{+2}};
\node[above] at (2,3) [scale=1]
{\tiny{+3}};

}

\node[below] at (1,0) [scale=0.75] {\small Fig 3.16(b). Case 2.2(a).};

\end{scope}

\begin{scope}[xshift=5cm, yshift=0cm ]

\draw[gray,very thin, step=0.5cm, opacity=0.5] (0,0) grid (2,3);

{
\draw[blue, line width=0.5mm] (0.5,0)--(0.5,0.5)--(1,0.5)--(1,0)--++(0.5,0);

\draw[blue, line width=0.5mm] (0,1)--++(0.5,0)--++(0,0.5)--++(0.5,0)--++(0,-0.5)--++(0.5,0);

\draw[blue, line width=0.5mm] (0,2)--++(0.5,0)--++(0,0.5)--++(0.5,0)--++(0,-0.5)--++(0.5,0);


\fill[blue!40!white, opacity=0.5] (0,1) rectangle (0.5,1.5);
\fill[blue!40!white, opacity=0.5] (1,1) rectangle (1.5,1.5);
\fill[blue!40!white, opacity=0.5] (0.5,2) rectangle (1,2.5);
\fill[blue!40!white, opacity=0.5] (0.5,0) rectangle (1,0.5);

\fill[orange!50!white, opacity=0.5] (1.5,0) rectangle (2,1);
\fill[orange!50!white, opacity=0.5] (1.5,2.5) rectangle (2,3);


\draw[orange!90!black, line width=0.5mm] (1.5,1)--++(0,-1);

\draw[orange!90!black, line width=0.5mm] (1.5,2)--++(0,-0.5)--++(0.5,0);

\draw[orange!90!black, line width=0.5mm] (0,3)--++(1.5,0)--++(0,-0.5)--++(0.5,0)--++(0,0.5);

}

{
\tikzset
  {
    myCircle/.style=    {orange}
  }
 
\foreach \x in {0,...,2}
\fill[myCircle, blue,] (0.75,1.62575+0.125*\x) circle (0.035);

\foreach \x in {0,...,2}
\fill[myCircle, orange,] (1.75,1.12575+0.125*\x) circle (0.035);
}

{
\node[left] at (0,3) [scale=1]
{\tiny{n-1}};
\node[left] at (0,2.5) [scale=1]
{\tiny{-2}};

\node[left] at (0,0) [scale=1] {\tiny{0}};
\node[left] at (0,0.5) [scale=1] {\tiny{1}};
\node[left] at (0,1) [scale=1] {\tiny{2}};

\node[above] at (0.5,3) [scale=1]
{\tiny{k}};
\node[above] at (1,3) [scale=1]
{\tiny{+1}};
\node[above] at (1.5,3) [scale=1]
{\tiny{+2}};
\node[above] at (2,3) [scale=1]
{\tiny{+3}};

}

\node[below] at (1,0) [scale=0.75] {\small Fig 3.16(c). Case 2.2(b).};

\end{scope}

\end{tikzpicture}
\end{adjustbox}
\end{wrapfigure}

\noindent \textit{CASE 2: $j$ is odd.} Note that the top leaf of the stack is in $\text{ext}(H)$. Again, $n-1$ is either even or odd. 

\null

\noindent \textit{CASE 2.1: $n-1$ is odd.} We have that  $R(k,n-2) \in \text{ext}(H)$. But then, the fact that $e(k,k+1;n-1) \in H$ implies that $R(k,n-2)$ is not a cookie neck, contradicting Lemma 1.14. End of Case 2.1. See Figure 3.16 (a).

\null

\noindent \textit{CASE 2.2: $n-1$ is even.} We have that $e(k+1,k+2;0) \in H$. Then, either $e(k+2, k+3;0) \in H$, or $e(k+2;0,1) \in H$.

\null

\noindent \textit{CASE 2.2(a): $e(k+2, k+3;0) \in H$.} Then we must have $e(k+2;1,2) \in H$. Note that for $i \in \{ 1,3, ..., n-2\}$, $e(k+2;i,i+1)\in H$ implies $e(k+2;i+2,i+3)\in H$. Then, for $i \in \{ 1,3, ..., n-2\}$, we have that $e(k+2;i,i+1)\in H$. Note that we must also have $e(k+2,k+3;n-2) \in H$. Then $R(k+2,n-2)$ must be a southern small cookie, contradicting the easternmost assumption. End of Case 2.2 (a).

\null

\noindent \textit{CASE 2.2(b): $e(k{+}2;0,1) \in H$.} Note that if $e(k+2,k+3;1) \in H$, then $R(k+2,0)$ must be a small cookie, contradicting the easternmost assumption. Then $e(k+2,k+3;1) \notin H$. But then we have $e(k+2;1,2) \in H$, and we are back to Case 2.2(a). End of Case 2.2(b). End of Case 2.2. End of Case 2. $\square$

\endgroup

\null

\noindent We will need Lemmas 3.13-3.16 to prove Lemma 3.12.

\null 

\noindent \textbf{Lemma 3.13.} Let $G$ be an $m \times n$ grid graph, and let $H$ be a Hamiltonian cycle of $G$. Let $C$ be a small cookie of $G$.  Assume that $G$ has only one large cookie, and that there is a $j$-stack of $A_0$ starting at the $A_0$-type containing $C$. Let $L$ be the leaf in the top ($j^{\text{th}}$) $A_0$ of the stack, and assume that $L$ is followed by an $A_1$-type. Let $X$ and $Y$ be the boxes adjacent to the middle-box of the $A_1$-type that are not its $H$-neighbours. If $P(X,Y)$ has no switchable boxes, then either:

(i) there is a cascade of length at most $\min(m,n)$, which avoids the stack of $A_0$'s, and after 

\hspace{0.4 cm} which $P(X,Y)$, gains a switchable box, or 

(ii) there is a cascade of length at most $\min(m,n)+1$, that collects $L$ and avoids the stack of $A_0$'s.

\null 

\noindent We postpone the proof of Lemma 3.13 until Section 4. It takes up all of the section.

\null

\noindent \textbf{Lemma 3.14.} Let $G$ be an $m \times n$ grid graph, let $H$ be a Hamiltonian cycle of $G$, and let $C \in SmallCookies\{N,S\}$ be an easternmost small cookie. Assume that $G$ has only one large cookie, and that there is a $j$-stack of $A_0$ starting at the $A_0$-type containing $C$. Let $L$ be the leaf contained in the top ($j^{\text{th}}$) $A_0$ of the stack. Assume that $L$ is followed by an $A_1$-type with looping $H$-path $P(X,Y)$. Let $X'$ be the box of $G$ that shares edges with $X$ and $Y$. Then $X' \in G_2$.

\null 

\setlength{\intextsep}{0pt}
\setlength{\columnsep}{20pt}
\begin{wrapfigure}[]{l}{0cm}
\begin{adjustbox}{trim=0cm 0cm 0cm 0cm} 
\begin{tikzpicture}[scale=1.25]
\begin{scope}[xshift=0cm] 
{
\draw[gray,very thin, step=0.5cm, opacity=0.5] (0,0) grid (1.5,2);

{
\fill[blue!40!white, opacity=0.5] (0,1)--++(0.5,0)--++(0,1)--++(-0.5,0);
\fill[blue!40!white, opacity=0.5] (1,1)--++(0.5,0)--++(0,1)--++(-0.5,0);

\draw[blue, line width=0.5mm] (0,1)--++(0.5,0)--++(0,1); 
\draw[blue, line width=0.5mm] (1.5,1)--++(-0.5,0)--++(0,1); 

\fill[blue!40!white, opacity=0.4](0.5,0)--++(0,0.5)--++(0.5,0)--++(0,-0.5);
\draw[blue, line width=0.5mm] (0.5,0)--++(0,0.5)--++(0.5,0)--++(0,-0.5);

}


{

\draw[black, line width=0.15mm] (0.2,1.95)-++(0,0.1);
\draw[black, line width=0.15mm] (0.25,1.95)-++(0,0.1);
\draw[black, line width=0.15mm] (0.3,1.95)-++(0,0.1);

\draw[black, line width=0.15mm] (1.2,1.95)-++(0,0.1);
\draw[black, line width=0.15mm] (1.25,1.95)-++(0,0.1);
\draw[black, line width=0.15mm] (1.3,1.95)-++(0,0.1);

}

{

\node[left] at (0,1.5) [scale=1]{\tiny{+1}};
\node[left] at (0,1) [scale=1]{\tiny{$\ell$}};
\node[left] at (0,0.5) [scale=1]{\tiny{-1}};

\node[above] at (0,2) [scale=1]
{\tiny{-1}};
\node[above] at (0.5,2) [scale=1]
{\tiny{k}};
\node[above] at (1,2) [scale=1]
{\tiny{+1}};
}

\node at (1.25, 1.25) [scale=0.8]{Y};
\node at (0.25, 1.25) [scale=0.8]{X};
\node at (0.75, 0.25) [scale=0.8]{L};

\node at (0.75, 1.25) [scale=0.8]{X$'$};

\node[below] at (0.75, 0) [scale=0.75] {\small{Fig. 3.17 (a). Case 1.}};

}
\end{scope}

\begin{scope}[xshift=2.25cm] 
{
\draw[gray,very thin, step=0.5cm, opacity=0.5] (0,0) grid (2,2);

{

\fill[blue!40!white, opacity=0.5](0,0)--++(0.5,0)--++(0,1)--++(-0.5,0);
\fill[blue!40!white, opacity=0.5](0.5,0.5)--++(0.5,0)--++(0,1.5)--++(-0.5,0);
\fill[blue!40!white, opacity=0.5](1,0)--++(0.5,0)--++(0,1)--++(-0.5,0);

\fill[blue!40!white, opacity=0.5](1.5,2)--++(0,-0.5)--++(0.5,0)--++(0,0.5);

\draw[blue, line width=0.5mm] (0,1)--++(0.5,0)--++(0,1); 
\draw[blue, line width=0.5mm] (1.5,1)--++(-0.5,0)--++(0,1);

\draw[blue, line width=0.5mm] (0.5,0)--++(0,0.5)--++(0.5,0)--++(0,-0.5);

\draw[blue, line width=0.5mm] (1,2)--++(0.5,0)--++(0,-0.5)--++(0.5,0)--++(0,0.5);

}


{

\draw[black, line width=0.15mm] (0.7,1.95)-++(0,0.1);
\draw[black, line width=0.15mm] (0.75,1.95)-++(0,0.1);
\draw[black, line width=0.15mm] (0.8,1.95)-++(0,0.1);

}

{

\node[left] at (0,1.5) [scale=1]{\tiny{+1}};
\node[left] at (0,1) [scale=1]{\tiny{$\ell$}};
\node[left] at (0,0.5) [scale=1]{\tiny{-1}};

\node[above] at (0,2) [scale=1]
{\tiny{-1}};
\node[above] at (0.5,2) [scale=1]
{\tiny{k}};
\node[above] at (1,2) [scale=1]
{\tiny{+1}};
\node[above] at (1.5,2) [scale=1]
{\tiny{+2}};

}

\node at (1.25, 1.25) [scale=0.8]{Y};
\node at (0.25, 1.25) [scale=0.8]{X};
\node at (0.75, 0.25) [scale=0.8]{L};
\node at (0.75, 1.25) [scale=0.8]{X$'$};

\node[below] at (1, 0) [scale=0.75] {\small{Fig. 3.17 (b). Case 2.}};

}
\end{scope}
\end{tikzpicture}
\end{adjustbox}
\end{wrapfigure}

\noindent \textit{Proof.} For definiteness, assume that $L$ is the northern leaf $R(k,l-2)$, and that $X=R(k-1,l)$. Then $X'=R(k,l)$ and $Y=R(k+1,l)$. Note that $l-2 \geq 0$ and $l+2 \leq n-1$. Either $P(X,Y)$ is contained in $\text{ext}(H)$, or $P(X,Y)$ is contained in $\text{int}(H)$, so there are two cases to check.

\null

\noindent \textit{CASE 1: $P(X,Y) \subset \text{ext}(H)$.} By Lemma 1.14, we must have that $m-1 >k+2$ and $k-1>0$. To see that $n-1 > l+2$, assume for contradiction that $n-1 = l+2$. By Lemma 1.14, $X+(1,0)$ and $Y+(1,0)$ are cookie necks. But this contradicts the assumption that there is only large cookie in $G$. See Figure 3.17 (a). End of Case 1. 

\null

\noindent \textit{CASE 2: $P(X,Y) \subset \text{int}(H)$.} By Lemma 1.14, we must have that $m-1 >k+2$ and $k-1>0$. To see that $n-1 > l+2$, assume for contradiction that $n-1 = l+2$. Lemma 1.14 implies that $X'+(0,1)$
is the neck of the large cookie of $G$. But now $X'+(2,1)$ must be a small cookie of $G$, contradicting the easternmost assumption. See Figure 3.17 (b).

\null

\endgroup 

\noindent \textbf{Lemma 3.15.} Let $G$ be an $m \times n$ grid graph, let $H$ be a Hamiltonian cycle of $G$, and let $C \in SmallCookies\{N,S\}$ be an easternmost small cookie. Assume that $G$ has only one large cookie, and that there is a $j$-stack of $A_0$ starting at the $A_0$-type containing $C$. Let $L$ be the leaf in the top ($j^{\text{th}}$) $A_0$ of the stack. Assume that $L \in \text{int}(H)$ and that $L$ is followed by an $A_1$-type with looping $H$-path $P(X,Y)$. Let $X'$ be the box of $G$ that shares edges with $X$ and $Y$. If $X'$ is not in $G_3$ then, either $X' \mapsto W$ is a cascade, or there is a cascade $\mu, X' \mapsto W$, of length two, with $X' \mapsto W$ nontrivial in either case.

\null

\noindent \textit{Proof.} Suppose that $X'$ is not in $G_3$. By Lemma 3.14, $X' \in G_2 \setminus G_3$. For definiteness assume that $L$ is a northern leaf, and let $X'=R(k,l)$. The assumption that $L \in \text{int}(H)$ implies that $l-2>0$.

\null 

\begingroup
\setlength{\intextsep}{0pt}
\setlength{\columnsep}{10pt}
\begin{wrapfigure}[]{r}{0cm}
\begin{adjustbox}{trim= 0cm 0cm 0cm 0.5cm}
\begin{tikzpicture}[scale=1.25]

\begin{scope}[xshift=0cm, yshift=0cm ]

\draw[gray,very thin, step=0.5cm, opacity=0.5] (0,0) grid (2,1.5);

{
\draw[blue, line width=0.5mm] (0.5,0)--(0.5,0.5)--(1,0.5)--(1,0)--++(0.5,0);

\draw[blue, line width=0.5mm] (0,1)--++(0.5,0)--++(0,0.5)--++(0.5,0)--++(0,-0.5)--++(0.5,0);

\draw[blue, thick] (2,1) circle [radius=0.05];

\draw[blue, line width=0.5mm] (2,0)--++(0,0.5);

\draw[blue, line width=0.5mm] (1.5,0)--++(0.5,0);
\draw[blue, line width=0.5mm] (1.5,1)--++(0,-0.5)--++(0.5,0);

}

{

\node[above] at (0.5,1.5) [scale=1]
{\tiny{k}};
\node[above] at (1,1.5) [scale=1]
{\tiny{+1}};
\node[above] at (1.5,1.5) [scale=1]
{\tiny{+2}};
\node[above] at (2,1.5) [scale=1]
{\tiny{+3}};

\node at (0.75,0.25) [scale=0.8]
{\small{$C$}};

}

\node[below] at (1,0) [scale=0.75] {\small Fig 3.17$\frac{1}{2}$. $m-1=k+3$};

\end{scope}

\end{tikzpicture}
\end{adjustbox}
\end{wrapfigure}

\noindent Now we check that $m-1 > k+3$. By Lemma 1.14, $m-1 > k+2$. Assume for contradiction that $m-1=k+3$. Note that we must have $e(k+1,k+2;0) \in H$, $e(k+2,k+3;0) \in H$, and $e(k+3;0,1) \in H$. This implies that we must have $e(k+2; 1,2)\in H$ and $e(k+2,k+3; 1)\in H$. But now $H$ misses $v(k+3;2)$. It follows that we must have $m-1 > k+3$. See Figure 3.17$\frac{1}{2}$(b). By symmetry, $0 < k-2$. It follows that $l+3=n-1$. 

The same argument used in Case 2.2 of Lemma 3.11 (see Figures 3.16 (b) and (c)) shows that we have $e(k+2;l+1,l+2)\in H$ and $e(k+2,k+3;l+1)\in H$. Now either $e(k,k+1;l+2) \in H$ or $e(k,k+1;l+2) \notin H$. See figures 3.18 (a) and (b).

\setlength{\intextsep}{0pt}
\setlength{\columnsep}{20pt}
\begin{wrapfigure}[]{l}{0cm}
\begin{adjustbox}{trim=0cm 0.5cm 0cm 0.25cm} 
\begin{tikzpicture}[scale=1.25]

\begin{scope}[xshift=0cm] 
{
\draw[gray,very thin, step=0.5cm, opacity=0.5] (0,0) grid (2,2.5);

{

\fill[blue!40!white, opacity=0.5](0,0)--++(0.5,0)--++(0,1)--++(-0.5,0);
\fill[blue!40!white, opacity=0.5](0.5,0.5)--++(0.5,0)--++(0,1.5)--++(-0.5,0);
\fill[blue!40!white, opacity=0.5](1,0)--++(0.5,0)--++(0,1)--++(-0.5,0);

\fill[blue!40!white, opacity=0.5](1.5,2)--++(0,-0.5)--++(0.5,0)--++(0,0.5);

\draw[blue, line width=0.5mm] (0,1)--++(0.5,0)--++(0,1); 
\draw[blue, line width=0.5mm] (1.5,1)--++(-0.5,0)--++(0,1);

\draw[blue, line width=0.5mm] (0.5,0)--++(0,0.5)--++(0.5,0)--++(0,-0.5);

\draw[blue, line width=0.5mm] (1.5,2)--++(0,-0.5)--++(0.5,0);

\draw[blue, line width=0.5mm] (0.5,2)--++(0.5,0);

\draw[blue, line width=0.5mm] (1.5,1)--++(0,-0.5); 

}

{

\node[left] at (0,2) [scale=1]{\tiny{+2}};
\node[left] at (0,1.5) [scale=1]{\tiny{+1}};
\node[left] at (0,1) [scale=1]{\tiny{$\ell$}};
\node[left] at (0,0.5) [scale=1]{\tiny{-1}};

=
\node[above] at (0.5,2.5) [scale=1]
{\tiny{k}};
\node[above] at (1,2.5) [scale=1]
{\tiny{+1}};
\node[above] at (1.5,2.5) [scale=1]
{\tiny{+2}};

}

\node at (0.75, 0.25) [scale=0.8]{L};
\node at (0.75, 1.25) [scale=0.8]{X$'$};

\node[below] at (1, 0) [scale=0.75] {\small{Fig. 3.18 (a). Case 1.}};

}
\end{scope}

\begin{scope}[xshift=2.75cm] 
{
\draw[gray,very thin, step=0.5cm, opacity=0.5] (0,0) grid (2,2.5);

{

\fill[blue!40!white, opacity=0.5](0,0)--++(0.5,0)--++(0,1)--++(-0.5,0);
\fill[blue!40!white, opacity=0.5](0.5,0.5)--++(0.5,0)--++(0,1.5)--++(-0.5,0);
\fill[blue!40!white, opacity=0.5](1,0)--++(0.5,0)--++(0,1)--++(-0.5,0);
\fill[blue!40!white, opacity=0.5](0.5,2)--++(0.5,0)--++(0,0.5)--++(-0.5,0);

\fill[blue!40!white, opacity=0.5](1.5,2)--++(0,-0.5)--++(0.5,0)--++(0,0.5);

\draw[blue, line width=0.5mm] (0,1)--++(0.5,0)--++(0,1); 
\draw[blue, line width=0.5mm] (1.5,1)--++(-0.5,0)--++(0,1);

\draw[blue, line width=0.5mm] (0.5,0)--++(0,0.5)--++(0.5,0)--++(0,-0.5);

\draw[blue, line width=0.5mm] (1.5,2)--++(0,-0.5)--++(0.5,0);

\draw[blue, line width=0.5mm] (1.5,2)--++(0.5,0);

\draw[blue, line width=0.5mm] (0.5,2)--++(0,0.5)--++(-0.5,0);
\draw[blue, line width=0.5mm] (1,2)--++(0,0.5)--++(0.5,0);

\draw[blue, line width=0.5mm] (1.5,1)--++(0,-0.5); 
}


{

\draw[black, line width=0.15mm] (0.7,1.95)-++(0,0.1);
\draw[black, line width=0.15mm] (0.75,1.95)-++(0,0.1);
\draw[black, line width=0.15mm] (0.8,1.95)-++(0,0.1);

\draw[black, line width=0.15mm] (1.2,1.95)-++(0,0.1);
\draw[black, line width=0.15mm] (1.25,1.95)-++(0,0.1);
\draw[black, line width=0.15mm] (1.3,1.95)-++(0,0.1);

}

{

\node[left] at (0,2) [scale=1]{\tiny{+2}};
\node[left] at (0,1.5) [scale=1]{\tiny{+1}};
\node[left] at (0,1) [scale=1]{\tiny{$\ell$}};
\node[left] at (0,0.5) [scale=1]{\tiny{-1}};

=
\node[above] at (0.5,2.5) [scale=1]
{\tiny{k}};
\node[above] at (1,2.5) [scale=1]
{\tiny{+1}};
\node[above] at (1.5,2.5) [scale=1]
{\tiny{+2}};

}

\node at (0.75, 0.25) [scale=0.8]{L};
\node at (0.75, 1.25) [scale=0.8]{X$'$};

\node[below] at (1, 0) [scale=0.75] {\small{Fig. 3.18 (b). Case 2.}};

}
\end{scope}

\end{tikzpicture}
\end{adjustbox}
\end{wrapfigure}

\null 

\noindent \textit{CASE 1: $e(k,k+1;l+2) \in H$. } By Lemma 1.16, $X' \mapsto X' +(1,1)$ is a valid move, and by Observation 3.4, $X' \mapsto X' +(1,1)$ creates no new cookies. End of Case 1.

\null

\noindent \textit{CASE 2: $e(k,k+1;l+2) \notin H$. } Lemma 1.14 implies that $e(k+1,k+2;l+2)$ cannot be in $H$ either. It follows that $X'+(0,2)$ must be the neck of the large cookie. The assumption that there is only one large cookie implies that $e(k+2;l+2,l+3) \notin H$. Then we must have $e(k+2, k+3;l+2) \in H$. Then, by Lemma 1.16, $X'+(1,1) \mapsto X'+(2,1)$, $X' \mapsto X'+(1,0)$ is the cascade we seek. $\square$

\endgroup 

\null 

\noindent \textbf{Lemma 3.16.}  Let $H$ be a Hamiltonian cycle of an $m \times n$ grid graph $G$. Let $X' \in G_3 \cap \text{ext}(H)$ be a switchable box, and let $P(X,Y)$ be the looping $H$-path of $X'$. Assume that $P(X,Y)$ has a switchable box in $G_0 \setminus G_2$. Then switching $X'$ splits $H$ into two cycles $H_1$ and $H_2$ such that there is $v_1 \in H_1 \cap R_2$ and $v_2 \in H_2 \cap R_2$ with $v_1$ adjacent to $v_2$.

\null

\noindent \textit{Proof.} Let $Z$ be a switchable box of $P(X,Y)$ in $G_0 \setminus G_2$. Orient $H$. Let $(v_x, v_{x+1})$ and $(v_{y-1}, v_y)$ be the edges of $X'$ in $H$. Define $\overrightarrow{K}_1$ and $\overrightarrow{K}_2$ to be the subtrails $\overrightarrow{K}((v_x,v_{x+1})$, $(v_{y-1}, v_y))$ and $\overrightarrow{K}((v_y,v_{y+1}),$ $(v_{x-1},v_x))$ of $\overrightarrow{K}_H$, respectively. By Lemma 1.16 (i), switching $X'$ gives two cycles $H_1$ and $H_2$, with  $V(H_1)=V(\overrightarrow{K}_1 \setminus \{v_x, v_y\})$ and $V(H_2)=V(\overrightarrow{K}_2)$. By Proposition 3.1, $Z$ has a vertex in $H_1$ and another in $H_2$, and the same holds for $X'$. Since $X' \in G_3$ and $Z \in G_0 \setminus G_2$, by JCT, $H_1 \cap R_2 \neq \emptyset$. Similarly, $H_2 \cap R_2 \neq \emptyset$. Now the argument in the last paragraph of Lemma 3.7 shows that there must be $v_1 \in H_1 \cap R_2$ and $v_2 \in H_2 \cap R_2$ with $v_1$ adjacent to $v_2$. $\square$

\null 

\noindent \textit{Proof of Lemma 3.12.} For definiteness assume that $L$ is the northern leaf $R(k,l)$. Let $P(X,Y)$ be the looping $H$-path following $L$, with $X=R(k-1,l+2)$ and $Y=R(k+1,l+2)$. By Lemma 3.13, either there is a cascade that collects $L$, or a cascade after which $P(X,Y)$ gains a switchable box, with both cascades having length at most $\min(m,n)+1$, and both avoiding the $j$-stack of $A_0$'s starting at $C$. If the former, we are done, so may assume that $P(X,Y)$ has a switchable box $Z$. Let $J$ be the large cookie of $G$ and let $N_J$ be the neck of $J$. Note that $N_J$ is not a box of $P(X,Y)$. Now, $P(X,Y)$ is either contained in \text{ext}(H) or \text{int}(H). 

\null

\noindent \textit{CASE 1: $P(X,Y) \subseteq $} \text{ext}(H). Then $X' \subset$ \text{int}(H).  By Lemma 3.14, $X' \in G_2$. By Proposition 3.3, $Z \mapsto X'$ is a valid move. By Observation 3.4, $Z \mapsto X'$ does not create additional cookies. Then $Z \mapsto X'$, $L+(0,1) \mapsto L$ is a cascade that collects $L$. 

\null

\noindent \textit{CASE 2: $P(X,Y) \subseteq $} $\text{int}(H)$ .Then $X' \subset$ \text{ext}(H).  If $Z \subset G_2$, by Proposition 3.3 $Z \mapsto X'$ is a valid move, and by Observation 3.4, $Z \mapsto X'$ does not create additional cookies. Then $Z \mapsto X'$, $L \mapsto L+(0,1)$ is a cascade that collects $L$.

Suppose then that $Z \subset G_0 \setminus G_2$. By Lemma 3.15, we only need to check the case where $X' \in G_3$. Note that switching $X'$ splits $H$ into two cycles $H_1$ and $H_2$. By Lemma 3.16 there is $v_1 \in H_1 \cap R_2$ and $v_2 \in H_2 \cap R_2$ with $v_1$ adjacent to $v_2$.  By Lemma 3.5 there is a cascade $X' \mapsto W$, or $\mu, X' \mapsto W$, with $X' \mapsto W$ nontrivial. Note that here, $X'$ plays the role that $Z$ played in Lemma 3.5. Then $\mu, X' \mapsto W, L+(0,1) \mapsto L$ or $X' \mapsto W, L+(0,1) \mapsto L$ is a cascade that collects $L$. 

We have just shown that if $P(X,Y)$ has a switchable box, then the cascade required to collect $L$ has length at most three. By Lemma 3.13, the cascade after which  $P(X,Y)$ to gains a switchable box has length at most $\min(m,n)$. Thus, at most $\min(m,n)+3$ moves are required to collect $L$. $\square$

\null 

\subsection{Summary }

\textcolor{blue}{\textbullet} In Section 3 we proved the MLC and 1LC algorithms. The proof of the MLC algorithm is fully contained here, while the proof of the 1LC algorithm depends on Lemma 3.13, whose proof is given in Section 4. 

Proposition 3.3 characterizes when double-switch moves are valid and serves as the primary tool for both algorithms. 

The MLC algorithm handles the case where H has multiple large cookies. To collect a large cookie $J$ with switchable neck $N_J$, we look for a switchable box $Z$ in the looping $H$-path of $N_J$. Proposition 3.8 shows that either such a $Z$ already exists, or a there is single
preparatory move that produces one. 

The 1LC algorithm handles the case where $H$ has exactly one large cookie and at least one small cookie. It collects outermost small cookies. Suppose that $C$ is an outermost small cookie. Either $C$ can be collected immediately by a single move, or $C$ is followed by a $j$-stack of $A_0$-types and an $A_1$-type with switchable middle-box $X'$. If the latter, let $P(X, Y)$ be the $H$-path determined by $X'$. If $P(X, Y)$ contains a switchable box $Z$, then $C$ can be collected
by either switching $X'$ directly (if $Z \in G_2$) or by using Lemma 3.5 to find a cascade of length at most two that enables switching $X$
(if $Z \in G_0 \setminus G_2$). In both cases, a cascade of flips then
collects $C$. The existence of such a switchable box $Z$ is guaranteed by Lemma 3.13, whose proof takes up Section 4. \textcolor{blue}{\textbullet}

\null

\newpage

\section{Looping fat paths, turns and weakenings}

\begingroup
\setlength{\intextsep}{0pt}
\setlength{\columnsep}{20pt}
\begin{wrapfigure}[]{r}{0cm}
\begin{adjustbox}{trim=0cm 0cm 0cm 1.25cm}
\begin{tikzpicture}[scale=1.25]

\begin{scope}[xshift=0]
\draw[gray,very thin, step=0.5cm, opacity=0.5] (0,0) grid (2.5,3);

\fill[blue!50!white, opacity=0.5] (0,1.5)--++(0.5,0)--++(0,0.5)--++(0.5,0)--++(0,-1)--++(0.5,0)--++(0,1)--++(0.5,0)--++(0,-0.5)--++(0.5,0)--++(0,-0.5)--++(-0.5,0)--++(0,-0.5)--++(-0.5,0)--++(0,-0.5)--++(-0.5,0)--++(0,0.5)--++(-0.5,0)--++(0,0.5)--++(-0.5,0);

\draw[blue, line width=0.5mm] (1,3)--++(0,-0.5)--++(0.5,0)--++(0,0.5);
\draw[blue, line width=0.5mm] (0,1.5)--++(0.5,0)--++(0,0.5)--++(0.5,0)--++(0,-1)--++(0.5,0)--++(0,1)--++(0.5,0)--++(0,-0.5)--++(0.5,0)--++(0,-0.5)--++(-0.5,0)--++(0,-0.5)--++(-0.5,0)--++(0,-0.5);

\draw[blue, line width=0.5mm] (0,1)--++(0.5,0)--++(0,-0.5)--++(0.5,0)--++(0,-0.5);

\draw[red, line width=0.25mm, dotted](0.75,1.75)--++(0,-1)--++(1,0)--++(0,1);

{
\node[below] at (1.25,0) [scale=0.75] {\small{\begin{tabular}{c} Fig. 4.1. A southern looping fat  \\ path $G\langle N[P(X,Y] \rangle$. $N[P(X,Y]$  \\  shaded in light blue. $P(X,Y)$ \\  traced in red. \end{tabular}}};;

\node at  (1.75, 1.75) [scale=0.8]{X};
\node at  (0.75, 1.75) [scale=0.8]{Y};

\node[right] at (2.5,2.5) [scale=1]{\tiny{$\ell$}};
\node[above] at (1,3) [scale=1]{\tiny{$k$}};

}    
\end{scope}

\end{tikzpicture}
\end{adjustbox}
\end{wrapfigure}

\textbf{Definitions}. Let $J=\{X_1, X_2, ..., X_r\}$ be a collection of boxes in an $m \times n$ grid graph $G$, and let $H$ be a Hamiltonian cycle of $G$. We will use the notation $G\langle J \rangle$ to denote the subgraph of $G$ with vertex set $V(G\langle J \rangle)=V(J)$ and edge set $E(G\langle J \rangle)=E(J)\cap E(H)$. The boxes of $G\langle J \rangle$ are the boxes of $J$. We call $G\langle J \rangle$ the \textit{subgraph of $G$ induced by $J$}. 

Suppose that the southern leaf $R(k,l)$ is followed by an $A_1$-type. Let $X=R(k+1, l-2)$ and $Y=R(k-1, l-2)$. Let $P(X,Y)$ be a southern looping $H$-path following the southern leaf $R(k,l)$. The set of all boxes in $P(X,Y)$, along with their $H$-neighbours, is called the \textit{$H$-neighbourhood of $P(X,Y)$}, and is denoted by $N[P(X,Y)]$. Consider the subgraph $F=G\langle N[P(X,Y)] \rangle$ of $G$ induced by $N[P(X,Y)]$. We define a \textit{short weakening of $F$} to be a cascade of length three or less after  which, the edge $\{v(k,l-1),v(k+1,l-1)\}$ is in the resulting Hamiltonian cycle of $G$. We say that $F$ is a \textit{southern looping fat path} if $F$ has no short weakening. We define western, northern and eastern looping fat paths analogously.

\null 

\noindent Assume that $G$ has only one large cookie, and that there is a $j$-stack of $A_0$ starting at the $A_0$-type containing an outermost southern small cookie $C$. Let $L$ be the leaf in the top ($j^{\text{th}}$) $A_0$ of the stack, and assume that $L$ is followed by an $A_1$-type with looping $H$ path $P(X,Y)$. Let $F=G\langle N[P(X,Y)] \rangle$. If $F$ has no short weakening, we say that $F$ a southern looping fat path \textit{anchored} at the outermost small southern cookie $C$. Analogous definitions apply for northern, eastern, and western looping fat paths. 

We remark that if $P(X,Y)$ has a switchable box then, by Proposition 3.3 and the proof of Lemma 3.12, $F$ has a short weakening, and thus it cannot be a looping fat path. 

Throughout the remainder of this section, we assume that $G$ has exactly one large cookie, and that all looping fat paths considered are anchored at some outermost small cookie.

\null

\endgroup 

\null 

\begingroup
\setlength{\intextsep}{0pt}
\setlength{\columnsep}{10pt}
\begin{wrapfigure}[]{l}{0cm}
\begin{adjustbox}{trim=0cm 0cm 0cm 0.5cm}
\begin{tikzpicture}[scale=1.25]
\begin{scope}[xshift=0cm, yshift=0cm]

\draw[gray,very thin, step=0.5cm, opacity=0.5] (0.5,0.5) grid (2,2);

{
\node[right] at (2,2) [scale=1]{\tiny{$\ell$}};

\node[above] at (0.5,2) [scale=1]{\tiny{$k$}};
}
{
\draw[blue, line width=0.5mm] (0.5,2)--++(0.5,0)--++(0,-0.5)--++(0.5,0)--++(0,-0.5)--++(0.5,0)--++(0,-0.5);

\fill[blue!40!white, opacity=0.5] (0.5,2)--++(0.5,0)--++(0,-0.5)--++(0.5,0)--++(0,-0.5)--++(0.5,0)--++(0,-0.5)--++(-1,0)--++(0,0.5)--++(-0.5,0);
}

{
\node[below] at (1.25,0.5) [scale=0.8]{\small{\begin{tabular}{c} Fig. 4.2. $S_{\rightarrow}(k,l;$ \\ $k+3,l-3)$. \end{tabular}}};;
}

\end{scope}
\end{tikzpicture}
\end{adjustbox}
\end{wrapfigure}

\noindent We define below a subgraph of $G$ consisting of the union of translations of two adjacent and perpendicular edges of $G$. Let $r \in \mathbb{N}$, and let \textit{the stairs from (k,l) to (k+r,l-r) east} be denoted by $S_{\rightarrow}(k,l;k+r,l-r)$ and be defined as:

\setlength{\abovedisplayskip}{0pt}
\setlength{\belowdisplayskip}{10pt}

\[S_{\rightarrow}(k,l;k+r,l-r) = \bigcup_{j=0}^{r-1} \Big( e(k,k+1;l)+(j,-j) \Big )\cup \Big( e(k+1;l-1,l)+(j,-j) \Big).\]

\noindent We define $d(S)=r$ to be the \textit{length} of $S_{\rightarrow}(k,l;k+r,l-r)$. We say that $S_{\rightarrow}(k,l;k+r,l-r)$ starts at $v(k,l)$ and ends at $v(k+r,l-r)$. The subscripted arrow indicates the direction from $v(k,l)$ of the first edge of the subgraph. By choosing an ``up", ``down", ``left" or ``right" arrow for direction and a sign for the third and fourth arguments of $S_{\square}(k,l;k\pm r,l \pm r)$ we may describe any of the eight possible steps subgraphs starting at the vertex $v(k,l)$.  See Figure 4.2.

\endgroup

\null

\noindent \textbf{Definitions.} \textcolor{blue}{\textbullet} Let $H$ be a Hamiltonian cycle of an $m \times n$ grid graph $G$. Let $T$ be the subgraph of $H$ on the edges $S_{\downarrow}(k+1,l;k',l'+1)$, $e(k;l-1,l)$ and $e(k'-1,k';l')$, where $k'=k+d(T)$, $l'=l-d(T)$, where $d(T)=d(S)+1$ is \textit{the length of T} and $d(T) \geq 2 $. We call $T$ a \textit{north-east turn}. If both $e(k,k+1;l)$ and $e(k';l',l'+1)$ belong to $G \setminus H$, call $T$ an \textit{open north-east turn}. If exactly 

\begingroup
\setlength{\intextsep}{0pt}
\setlength{\columnsep}{20pt}
\begin{wrapfigure}[]{r}{0cm}
\begin{adjustbox}{trim=0cm 0cm 0cm 0cm}
\begin{tikzpicture}[scale=1.25]
\begin{scope}[xshift=0cm, yshift=0cm]

\draw[gray,very thin, step=0.5cm, opacity=0.5] (0,0) grid (2.5,2.5);

{
\node[left] at (0,2) [scale=1]{\tiny{$\ell$}};
\node[above] at (0.5,2.5) [scale=1]{\tiny{$k$}};

\node[left] at (0,0.5) [scale=1]{\tiny{$\ell'$}};
\node[above] at (2,2.5) [scale=1]{\tiny{$k'$}};
}
{

\draw[blue, line width=0.5mm] (1,2)--++(0,-0.5)--++(0.5,0)--++(0,-0.5)--++(0.5,0)--++(0,-0.5)--++(-0.5,0);
\draw[blue, line width=0.5mm] (0.5,1.5)--++(0,0.5);

\fill[blue!40!white, opacity=0.5] (0.5,2)--++(0.5,0)--++(0,-0.5)--++(0.5,0)--++(0,-0.5)--++(0.5,0)--++(0,-0.5)--++(-1,0)--++(0,0.5)--++(-0.5,0);

}

{
\node[below] at (1.25,-0) [scale=0.75]{\small{\begin{tabular}{c} Fig. 4.3. A half-open \\ northeast turn. \end{tabular}}};
}

\end{scope}
\end{tikzpicture}
\end{adjustbox}
\end{wrapfigure}

\noindent one of $e(k,k+1;l)$ and $e(k';l',l'+1)$ is in $H$, then $T$ is a \textit{half-open north-east turn}. If both $e(k,k+1;l)\in H$ and $e(k';l',l'+1) \in H$, then $T$ is a \textit{closed north-east turn}. See Figure 4.3 For any north-east turn $T$, we say that $R(k,l-1)$ is the \textit{northern leaf} of $T$ and $R(k'-1,l')$ is the \textit{eastern leaf} of $T$. If $e(k,k+1;l)\notin H$ we call $R(k,l-1)$ an \textit{open northern leaf of T} and if $e(k,k+1;l)\in H$ we call $R(k,l-1)$ a \textit{closed northern leaf of T}. We note that the two leaves of a turn will determine its ``leaf prefix": If a turn has a north leaf and an east leaf then the turn is a north-east turn. 

We say that a looping fat path $F$ \textit{has a turn} (open, half-open or closed) to mean that there exists some turn $T$ of $H$ such that $E(F) \supset E(T)$.

\null

\noindent \textbf{Sketch of proof of Lemma 3.13.} Let $H$ be a Hamiltonian cycle of an $m \times n$ grid graph $G$. Assume that $P(X, Y)$ is a looping $H$-path with no switchable boxes, following a leaf $L$. It follows that $P(X, Y)$ is contained in a looping fat path $F$. In Section 4.2, we show that every looping fat path must have a turn. In Sections 4.3 we show that given a turn, we can find a cascade we call a weakening (precise definition in Section 4.3) that collects at least one of the leaves of the turn. Then we show that after such a cascade, either $P(X, Y)$ gains a switchable box, or we can extend the cascade by a single move to collect $L$. The rest of the Section is organized as follows. Section 4.1 proves structural properties of fat paths, on which the later sections build. Section 4.2 shows that every fat path contains a turn (Proposition 4.7). In Section 4.3 we define weakenings, prove that turns have weakenings, and give a proof of Lemma 3.13. \textcolor{blue}{\textbullet}

\null

\subsection{Properties of looping fat paths}

\textbf{Lemma 4.1.} Let $F=G\langle N[P(X,Y)] \rangle$ be a looping fat path. Let $W=R(k,l)$ \footnote{ $W=R(k,l)$ is not related to the southern leaf $R(k,l)$ in the definitions in page 18.} be a box of $P=P(X,Y)$ with the $H$-neighbour $Z=W+(0,-1)$ southward in $N[P] \setminus P$. Then:

(a) $Z$ has exactly one $H$-neighbour in $P$ and $W$ has no other $H$-neighbour in $N[P] \setminus P$.

(b) If $W$ is not an end-box (i.e. $X$ or $Y$) of $P$, then the $H$-neighbours of $W$ in $P$ are $W+(-1,0)$ 

\hspace{0.5cm} and $W+(1,0)$. Furthermore, $S_{\rightarrow}(k-1,l;k,l-1) \in H$, $S_{\uparrow}(k+1,l-1;k+2,l) \in H$, 

\hspace{0.5cm} and $(k,k+1;l+1)\in H$. 

(c) If $W$ is an end-box of $P$, then $e(k,k+1;l+1) \in H$ and exactly one of $e(k;l,l+1)$ and 

\hspace{0.5cm}  $e(k+1;l,l+1)$ belong to $H$.

(d) $Z$ is a leaf or $Z$ is a switchable box in $H$.

\noindent Analogous statements apply when $Z$ is west, north or east of $W$.

\null 

\noindent \textit{Proof of (a).} Note that if $Z$ has more than one $H$-neighbour in $P$ then we can make an $H$-cycle. To see that $W$ has no other $H$-neighbour in $N[P] \setminus P$, assume for contradiction that $W$ has at least two $H$-neighbours in $N[P] \setminus P$. 

If $W$ is an end-box of $P$, then, by definition of $A_1$, $W$ has at most two $H$-neighbours, and at least one of them must belong to $P$, contradicting our assumption that $W$ has at least two $H$-neighbours in $N[P] \setminus P$.  

\noindent If $W$ is not an end-box, then $W$ must have four $H$ neighbours: two in $N[P] \setminus P$ and two in $P$. By definition of a looping fat path, at least one of the neighbours of $W$ in $P$, say $W'$, is not an end-box. But then $W'$ must be switchable, contradicting that $F$ is a looping fat path. $\square$

\null 

\begingroup 
\setlength{\intextsep}{0pt}
\setlength{\columnsep}{10pt}
\begin{wrapfigure}[]{l}{0cm}
\begin{adjustbox}{trim=0cm 0cm 0cm 0.25cm}
\begin{tikzpicture}[scale=1.5]
\begin{scope}[xshift=0cm]{
\draw[gray,very thin, step=0.5cm, opacity=0.5] (0,0) grid (1.5,1.5);

\draw[blue, line width=0.5mm] (0,0.5)--++(0.5,0)--++(0,-0.5);
\draw[blue, line width=0.5mm] (0,1)--++(0.5,0)--++(0,0.5);

\draw[blue, line width=0.5mm] (1,0.5)--++(0,0.5);

\node[right] at (1.5,0) [scale=1]
{\tiny{-1}};
\node[right] at (1.5,0.5) [scale=1]
{\tiny{$\ell$}};
\node[right] at (1.5,1) [scale=1]
{\tiny{+1}};
\node[right] at (1.5,1.5) [scale=1]
{\tiny{+2}};

\node[above] at (0, 1.5) [scale=1]
{\tiny{-1}};
\node[above] at (0.5, 1.5) [scale=1]
{\tiny{$k$}};
\node[above] at (1, 1.5) [scale=1]
{\tiny{+1}};

\node at (0.75,0.25) [scale=0.8] {\small{Z}};
\node at (0.75,0.75) [scale=0.8] {\small{W}};

{
\draw[black, line width=0.15mm] (0.70,0.45)--(0.7,0.55);
\draw[black, line width=0.15mm] (0.75,0.45)--(0.75,0.55);
\draw[black, line width=0.15mm] (0.8,0.45)--(0.8,0.55);

\draw[black, line width=0.15mm] (0.70,0.95)--(0.7,1.05);
\draw[black, line width=0.15mm] (0.75,0.95)--(0.75,1.05);
\draw[black, line width=0.15mm] (0.8,0.95)--(0.8,1.05);

\draw[black, line width=0.15mm] (0.45,0.70)--(0.55,0.70);
\draw[black, line width=0.15mm] (0.45,0.75)--(0.55,0.75);
\draw[black, line width=0.15mm] (0.45,0.80)--(0.55,0.80);

}

\node[below] at (0.75,0) [scale=0.8]{\small{\begin{tabular}{c} Fig. 4.4 (a).
 \end{tabular}}};;

} \end{scope}

\begin{scope}[xshift=2.25cm]{
\draw[gray,very thin, step=0.5cm, opacity=0.5] (0,0) grid (1.5,1.5);

\draw[blue, line width=0.5mm] (0,0.5)--++(0.5,0)--++(0,-0.5);
\draw[blue, line width=0.5mm] (1.5,0.5)--++(-0.5,0)--++(0,-0.5);

\draw[blue, line width=0.5mm] (0.5,1)--++(0.5,0);

\node[right] at (1.5,0) [scale=1]
{\tiny{-1}};
\node[right] at (1.5,0.5) [scale=1]
{\tiny{$\ell$}};
\node[right] at (1.5,1) [scale=1]
{\tiny{+1}};

\node[above] at (0, 1.5) [scale=1]
{\tiny{-1}};
\node[above] at (0.5, 1.5) [scale=1]
{\tiny{$k$}};
\node[above] at (1, 1.5) [scale=1]
{\tiny{+1}};
\node[above] at (1.5, 1.5) [scale=1]
{\tiny{+2}};

\node at (0.75,0.25) [scale=0.8] {\small{Z}};
\node at (0.75,0.75) [scale=0.8] {\small{W}};

{
\draw[black, line width=0.15mm] (0.70,0.45)--(0.7,0.55);
\draw[black, line width=0.15mm] (0.75,0.45)--(0.75,0.55);
\draw[black, line width=0.15mm] (0.8,0.45)--(0.8,0.55);


\draw[black, line width=0.15mm] (0.95,0.70)--(1.05,0.70);
\draw[black, line width=0.15mm] (0.95,0.75)--(1.05,0.75);
\draw[black, line width=0.15mm] (0.95,0.80)--(1.05,0.80);

\draw[black, line width=0.15mm] (0.45,0.70)--(0.55,0.70);
\draw[black, line width=0.15mm] (0.45,0.75)--(0.55,0.75);
\draw[black, line width=0.15mm] (0.45,0.80)--(0.55,0.80);

}

\node[below] at (0.75,0) [scale=0.8]{\small{\begin{tabular}{c} Fig. 4.4 (b).
 \end{tabular}}};;

} \end{scope}

\end{tikzpicture}
\end{adjustbox}
\end{wrapfigure}

\noindent \textit{Proof of (b).} First we show that the $H$-neighbours of $W$ are $W+(1,0)$ and $W+(-1,0)$. Assume for contradiction that $W+(1,0)$ is not an $H$-neighbour of $W$. Then the $H$-neighbours of $W$ in $P$ must be $W+(-1,0)$ and $W+(0,1)$. It follows that $S_{\rightarrow}(k-1,l;k,l-1)\in H$ and $S_{\rightarrow}(k-1,l+1;k,l+2)\in H$. Note that, by the definition of $A_1$ and looping fat paths, $W+(-1,0)$ is not an end-box of $P$. But then $W+(-1,0)$ is a switchable box of $P$, contradicting that $F$ is a looping fat path. Therefore, the $H$-neighbours of $W$ in $P$ are $W+(-1,0)$ and $W+(1,0)$. It follows that $S_{\rightarrow}(k-1,l;k,l-1) \in H$ and $S_{\uparrow}(k+1,l-1;k+2,l) \in H$, and by part (a), $(k,k+1;l+1)\in H$. See Figure 4.4 (a) and (b). End of proof for (b).

\endgroup

\begingroup 
\setlength{\intextsep}{0pt}
\setlength{\columnsep}{20pt}
\begin{center}

\begin{adjustbox}{trim=0cm 0cm 0cm 0.25cm}
\begin{tikzpicture}[scale=1.5]
\begin{scope}[xshift=0cm]{
\draw[gray,very thin, step=0.5cm, opacity=0.5] (0,0) grid (1.5,1.5);

\draw[blue, line width=0.5mm] (0.5,0.5)--++(0,0.5)--++(0.5,0);

\node[right] at (1.5,0) [scale=1]
{\tiny{-1}};
\node[right] at (1.5,0.5) [scale=1]
{\tiny{$\ell$}};
\node[right] at (1.5,1) [scale=1]
{\tiny{+1}};

\node[above] at (0.5, 1.5) [scale=1]
{\tiny{$k$}};
\node[above] at (1, 1.5) [scale=1]
{\tiny{+1}};

\node at (0.75,0.25) [scale=0.8] {\small{Z}};
\node at (0.75,0.75) [scale=0.8] {\small{W}};

{
\draw[black, line width=0.15mm] (0.70,0.45)--++(0,0.1);
\draw[black, line width=0.15mm] (0.75,0.45)--++(0,0.1);
\draw[black, line width=0.15mm] (0.8,0.45)--++(0,0.1);

\draw[black, line width=0.15mm] (0.95,0.70)--++(0.1,0);
\draw[black, line width=0.15mm] (0.95,0.75)--++(0.1,0);
\draw[black, line width=0.15mm] (0.95,0.80)--++(0.1,0);

}

\node[below] at (0.75,0) [scale=0.8]{\small{\begin{tabular}{c} Fig. 4.5 (a).
 \end{tabular}}};;

} \end{scope}

\begin{scope}[xshift=2.5cm]{
\draw[gray,very thin, step=0.5cm, opacity=0.5] (0,0) grid (1.5,1.5);

\draw[blue, line width=0.5mm] (0.5,0.5)--++(0,0.5)--++(1,0);

\draw[blue, line width=0.5mm] (1,0)--++(0,0.5)--++(0.5,0);

\node[right] at (1.5,0) [scale=1]
{\tiny{-1}};
\node[right] at (1.5,0.5) [scale=1]
{\tiny{$\ell$}};
\node[right] at (1.5,1) [scale=1]
{\tiny{+1}};

\node[above] at (0.5, 1.5) [scale=1]
{\tiny{$k$}};
\node[above] at (1, 1.5) [scale=1]
{\tiny{+1}};
\node[above] at (1.5, 1.5) [scale=1]
{\tiny{+2}};

\node at (0.75,0.25) [scale=0.8] {\small{Z}};
\node at (0.75,0.75) [scale=0.8] {\small{W}};

{
\draw[black, line width=0.15mm] (0.70,0.45)--(0.7,0.55);
\draw[black, line width=0.15mm] (0.75,0.45)--(0.75,0.55);
\draw[black, line width=0.15mm] (0.8,0.45)--(0.8,0.55);

\draw[black, line width=0.15mm] (0.95,0.70)--(1.05,0.70);
\draw[black, line width=0.15mm] (0.95,0.75)--(1.05,0.75);
\draw[black, line width=0.15mm] (0.95,0.80)--(1.05,0.80);

}

\node[below] at (0.75,0) [scale=0.8]{\small{\begin{tabular}{c} Fig. 4.5 (b). Case 1: \\ $F$ is eastern. \end{tabular}}};;

} \end{scope}

\begin{scope}[xshift=5cm]{
\draw[gray,very thin, step=0.5cm, opacity=0.5] (0,0) grid (1.5,1.5);

\draw[blue, line width=0.5mm] (0.5,0)--++(0,1)--++(0.5,0);
\draw[blue, line width=0.5mm] (1,0)--++(0,0.5);

\node[right] at (1.5,0) [scale=1]
{\tiny{-1}};
\node[right] at (1.5,0.5) [scale=1]
{\tiny{$\ell$}};
\node[right] at (1.5,1) [scale=1]
{\tiny{+1}};

\node[above] at (0.5, 1.5) [scale=1]
{\tiny{$k$}};
\node[above] at (1, 1.5) [scale=1]
{\tiny{+1}};

\node at (0.75,0.25) [scale=0.8] {\small{Z}};
\node at (0.75,0.75) [scale=0.8] {\small{W}};

{
\draw[black, line width=0.15mm] (0.70,0.45)--(0.7,0.55);
\draw[black, line width=0.15mm] (0.75,0.45)--(0.75,0.55);
\draw[black, line width=0.15mm] (0.8,0.45)--(0.8,0.55);

\draw[black, line width=0.15mm] (0.95,0.70)--(1.05,0.70);
\draw[black, line width=0.15mm] (0.95,0.75)--(1.05,0.75);
\draw[black, line width=0.15mm] (0.95,0.80)--(1.05,0.80);

}

\node[below] at (0.75,0) [scale=0.8]{\small{\begin{tabular}{c}Fig. 4.5 (c). Case 1: \\ $F$  is southern. \end{tabular}}};;

} \end{scope}

\begin{scope}[xshift=7.5cm]{
\draw[gray,very thin, step=0.5cm, opacity=0.5] (0,0) grid (1.5,1.5);

\draw[blue, line width=0.5mm] (0,0.5)--++(0.5,0)--++(0,-0.5);
\draw[blue, line width=0.5mm] (1.5,0.5)--++(-0.5,0)--++(0,-0.5);

\draw[blue, line width=0.5mm] (0.5,1)--++(0.5,0);

\node[right] at (1.5,0) [scale=1]
{\tiny{-1}};
\node[right] at (1.5,0.5) [scale=1]
{\tiny{$\ell$}};
\node[right] at (1.5,1) [scale=1]
{\tiny{+1}};

\node[above] at (0, 1.5) [scale=1]
{\tiny{-1}};
\node[above] at (0.5, 1.5) [scale=1]
{\tiny{$k$}};
\node[above] at (1, 1.5) [scale=1]
{\tiny{+1}};
\node[above] at (1.5, 1.5) [scale=1]
{\tiny{+2}};

\node at (0.75,0.25) [scale=0.8] {\small{Z}};
\node at (0.75,0.75) [scale=0.8] {\small{W}};

{
\draw[black, line width=0.15mm] (0.70,0.45)--(0.7,0.55);
\draw[black, line width=0.15mm] (0.75,0.45)--(0.75,0.55);
\draw[black, line width=0.15mm] (0.8,0.45)--(0.8,0.55);


\draw[black, line width=0.15mm] (0.95,0.70)--(1.05,0.70);
\draw[black, line width=0.15mm] (0.95,0.75)--(1.05,0.75);
\draw[black, line width=0.15mm] (0.95,0.80)--(1.05,0.80);

\draw[black, line width=0.15mm] (0.45,0.70)--(0.55,0.70);
\draw[black, line width=0.15mm] (0.45,0.75)--(0.55,0.75);
\draw[black, line width=0.15mm] (0.45,0.80)--(0.55,0.80);

}

\node[below] at (0.75,0) [scale=0.8]{\small{\begin{tabular}{c} Fig. 4.5 (d). Case 2.
 \end{tabular}}};;

} \end{scope}

\end{tikzpicture}
\end{adjustbox}

\end{center}

\noindent \textit{Proof of (c).} By part (a) and the assumption that $W$ is an end-box of $P$, $W$ has exactly one $H$-neighbour in $P$ and no other $H$-neighbours in $N[P]\setminus P$. It follows that $W$ has exactly two edges in $H$ and two edges not in $H$. Assume for contradiction that the other edge of $W$ not in $H$ is $e(k,k+1;l+1)$. But then $e(k;l-1,l) \in H$ and  $e(k+1;l-1,l) \in H$ and the $W$ is switchable, contradicting that $F$ is a looping fat path. It follows that $e(k,k+1;l+1)\in H$. Since $W$ has exactly two edges in $H$, we have that exactly one of $e(k;l,l+1)$ and $e(k+1;l,l+1)$ belong to $H$. See Figure 4.5 (a). End of proof for (c).

\endgroup

\null 

\noindent \textit{Proof of (d).} $W$ is either an end-box of $P$ or it is not.

\null

\noindent \textit{CASE 1: W is an end-box of P.} By part (c), we may assume WLOG that $e(k;l,l+1)\in H$ and $e(k+1;l,l+1)\notin H$. Then $F$ is eastern or southern. Suppose that $F$ is eastern. Then $e(k+1,k+2;l+1) \in H$. It follows that $S_{\uparrow}(k+1,l-1;k+2,l)\in H$. But then $W+(1,0) \in P$ is switchable, contradicting that $F$ is a looping fat path. So $F$ must be southern. Then $e(k;l-1,l) \in H$. It follows that $S_{\uparrow}(k+1,l-1;k+2,l)\in H$. Now, either $e(k,k+1;l-1) \in H$, or  $e(k,k+1;l-1)\notin H$. Either way, (d) is satisfied. See figures 4.5 (b) and (c). End of Case 1. 

\endgroup 

\null

\noindent \textit{CASE 2: W is not end-box of P.} By part (b), the $H$-neighbours of $W$ in $P$ are $W+(-1,0)$ and $W+(1,0)$ and we have that $S_{\rightarrow}(k-1,l;k,l-1) \in H$, $S_{\uparrow}(k+1,l-1;k+2,l) \in H$. Then, either $e(k,k+1;l-1) \in H$ or  $e(k,k+1;l-1) \notin H$. Either way, (d) is satisfied. See Figure 4.5(d). $\square$

\null 

\noindent \textbf{Definitions.} Let $G$ be an $m \times n$ grid graph, let $H$ be a Hamiltonian cycle of $G$ and let
$J$ be an $H$-subtree of an $H$-component of $G$. We say that a box $Z$ of $J$ is a \textit{border box} of $J$ if $Z$ is an $H$-neighbour of a box $Z'\in G_{-1} \setminus J$. We will call the edge that $Z$ and $Z'$ share a \textit{shadow edge} of $J$. We call the set of all shadow edges of $J$ the \textit{shadow border} of $J$ and denote it by $hb(J)$. We define the \textit{shadow of $J$} to be the graph $h(J)$ with vertex set $V(h(J))=V(J)$ and edge set $E(h(J))=(E(J) \cap E(H)) \cup hb(J)$. The boxes of the shadow of $J$ are the same as the boxes of $J$. We note that shadow edges cannot be incident on boxes of $P$.

\null 

\noindent \textbf{Observation 4.2.} Let $F=G\langle N[P(X,Y)] \rangle$ be a looping fat path. Then the shadow edges of the $H$-subtree $N[P]$ can only be incident on boxes of $N[P] \setminus P$. Moreover, exactly one of the two boxes incident on a shadow edge of $N[P]$ belongs to $N[P] \setminus P$, and the other belongs to $G_{-1} \setminus N[P]$. 

\null

\noindent \textbf{Lemma 4.3.} Let $G$ be an $m \times n$ grid graph, let $H$ be a Hamiltonian cycle of $G$ and let $F=G\langle N[P(X,Y)] \rangle $ be a looping fat path in $G$. Then $E(h(F))$ is a Hamiltonian cycle of $h(F)$. 

\null 

\noindent \textit{Proof.} Since every vertex of $h(F)$ is incident on some edge of $E(h(F))$, it is sufficient to show that $E(h(F))$ is a cycle. We will prove:

(i) Every vertex of $h(F)$ has degree two in $h(F)$.

(ii) $E(h(F))$ is connected.

\noindent \textit{Proof of (i).} Let $v=v(k,l) \in V(F)$ and let $R(k-1,l-1)=Z$. Either $v$ has a shadow edge incident on it or it does not.

\begingroup 
\setlength{\intextsep}{0pt}
\setlength{\columnsep}{5pt}
\begin{wrapfigure}[]{l}{0cm}
\begin{adjustbox}{trim=0cm 0.25cm 0cm 0.25cm}
\begin{tikzpicture}[scale=1.5]

\begin{scope}[xshift=0cm]{
\draw[gray,very thin, step=0.5cm, opacity=0.5] (0,0) grid (1,1);

\fill[blue!50!white, opacity=0.5] (0.5,0) rectangle (1,1);

\draw[blue, line width=0.5mm] (0.5,0)--++(0,1);

\draw[fill=blue] (0.5,0) circle [radius=0.05];
\draw[fill=blue] (0.5,0.5) circle [radius=0.05];
\draw[fill=blue] (0.5,1) circle [radius=0.05];
{

\node[left] at (0,0) [scale=1]{\tiny{-1}};
\node[left] at (0,0.5) [scale=1]
{\tiny{$\ell$}};

\node[above] at (0.5,1) [scale=1]
{\tiny{$k$}};
\node[above] at (1, 1) [scale=1]{\tiny{+1}};

\node[right] at (0.5,0.25) [scale=0.8] {\small{$e_1$}};

\node[right] at (0.5,0.75) [scale=0.8] {\small{$e_2$}};

\node at (0.25,0.25) [scale=0.8] {\small{$Z$}};
}

\node[below] at (0.5,0) [scale=0.8]{\small{\begin{tabular}{c} Fig. 4.6. Case 1.1. \end{tabular}}};;

} \end{scope}

\end{tikzpicture}
\end{adjustbox}
\end{wrapfigure}

\noindent \textit{CASE 1. $v$ has no shadow edge incident on it.} Then there are two edges $e_1$ and $e_2$ of $H$ incident on $v$. These edges are either colinear or not colinear, so there are two cases to check. 

\null 

\noindent \textit{CASE 1.1. $e_1$ and $e_2$ are colinear.} For definiteness assume that $e_1=e(k;l-1,l)$, $e_2=e(k;l,l+1)$ and that $Z+(1,0) \in F$. Since $e(k,k+1;l)$ is not a shadow edge, $Z+(1,1) \in F$ as well, $v(k;l-1)\in F$ and $v(k;l+1)\in F$. Since $e(k-1,k;l)$ is not a shadow edge, we have that $\text{deg}_{h(F)}(v)=2$. See Figure 4.6. End of Case 1.1.

\null 

\begingroup 
\setlength{\intextsep}{0pt}
\setlength{\columnsep}{10pt}
\begin{wrapfigure}[]{r}{0cm}
\begin{adjustbox}{trim=0cm 0cm 0cm 0.5cm}
\begin{tikzpicture}[scale=1.5]

\begin{scope}[xshift=0cm]{
\draw[gray,very thin, step=0.5cm, opacity=0.5] (0,0) grid (1,1);

\fill[blue!50!white, opacity=0.5] (0.5,0) rectangle (1,0.5);


    


\draw[blue, line width=0.5mm] (0.5,0)--++(0,0.5)--++(0.5,0);

\draw[fill=blue] (1,0.5) circle [radius=0.05];
\draw[fill=blue] (0.5,0.5) circle [radius=0.05];
\draw[fill=blue] (0.5,0) circle [radius=0.05];

{



}

{

\node[left] at (0,0) [scale=1]{\tiny{-1}};
\node[left] at (0,0.5) [scale=1]
{\tiny{$\ell$}};

\node[above] at (0.5,1) [scale=1]
{\tiny{$k$}};
\node[above] at (1, 1) [scale=1]{\tiny{+1}};

\node[left] at (0.5,0.25) [scale=0.8] {\small{$e_1$}};

\node[below] at (0.75,0.5) [scale=0.8] {\small{$e_2$}};

\node at (0.75,0.75) [scale=0.8] {\small{$Z$}};
}

\node[below] at (0.5,0) [scale=0.8]{\small{\begin{tabular}{c} Fig. 4.7 (a). Case  \\ 1.2. $Z+(0,-1) \in F$. \end{tabular}}};;

} \end{scope}

\begin{scope}[xshift=1.75cm]{
\draw[gray,very thin, step=0.5cm, opacity=0.5] (0,0) grid (1,1);

\fill[blue!50!white, opacity=0.5] (0,0) rectangle (0.5,1);
\fill[blue!50!white, opacity=0.5] (0.5,0.5) rectangle (1,1);


    


\draw[blue, line width=0.5mm] (0.5,0)--++(0,0.5)--++(0.5,0);

\draw[fill=blue] (1,0.5) circle [radius=0.05];
\draw[fill=blue] (0.5,0.5) circle [radius=0.05];
\draw[fill=blue] (0.5,0) circle [radius=0.05];

{



}

{

\node[left] at (0,0) [scale=1]{\tiny{-1}};
\node[left] at (0,0.5) [scale=1]
{\tiny{$\ell$}};

\node[above] at (0.5,1) [scale=1]
{\tiny{$k$}};
\node[above] at (1, 1) [scale=1]{\tiny{+1}};

\node[left] at (0.5,0.25) [scale=0.8] {\small{$e_1$}};

\node[below] at (0.75,0.5) [scale=0.8] {\small{$e_2$}};

\node at (0.75,0.75) [scale=0.8] {\small{$Z$}};
}

\node[below] at (0.5,0) [scale=0.8]{\small{\begin{tabular}{c} Fig. 4.7 (b). Case  \\ 1.2. $Z+(0,-1) \notin F$. \end{tabular}}};;

} \end{scope}

\end{tikzpicture}
\end{adjustbox}
\end{wrapfigure}

\noindent \textit{CASE 1.2. $e_1$ and $e_2$ are not colinear.} For definiteness assume that $e_1=e(k;l-1,l)$, $e_2=e(k,k+1;l)$. Then $Z+(1,0) \in F$ or $Z+(1,0)  \notin F$.

If $Z+(1,0)  \in F$, then $v(k;l-1)\in F$ and $v(k+1;l)\in F$. As $v$ has no shadow edges incident on it, it follows that $\text{deg}_{(h(F))}(v)=2$. See Figure 4.7 (a).

If $Z+(1,0)  \notin F$, then at least one of $Z$, $Z+(0,1)$ and $Z+(1,1)$ belong to $F$. If $Z \in F$, since $e(k-1,k;l)$, is not a shadow edge, we have that $Z+(0,1) \in F$. Similarly, $Z+(1,1) \in F$. Then $v(k;l-1) \in F$ and $v(k+1;l) \in F$.  As $v$ has no shadow edges incident on it, it follows that $\text{deg}_{(h(F))}(v)=2$. The cases where $Z+(0,1) \in F$ and $Z+(1,1) \in F$ are similar and we omit them. See Figure 4.7(b). End of Case 1.2.

\endgroup 

\null 

\noindent \textit{CASE 2. $v$ has a shadow edge incident on it.} For definiteness, let $e(k;l-1,l)$ be the shadow edge on which $v$ is incident. There are three possibilities: $e(k-1,k;l) \in H$ and $e(k;l,l+1) \in H$, $e(k-1,k;l) \in H$ and $e(k,k+1;l) \in H$, and, $e(k;l,l+1) \in H$ and $e(k,k+1;l) \in H$. By symmetry, we only need to check the first two.

\null 

\noindent \textit{CASE 2.1: $e(k-1,k;l) \in H$ and $e(k;l,l+1) \in H$.} Then $e(k,k+1;l) \notin H$. By Observation 4.2, exactly one of $Z$ and $Z+(1,0)$ belongs to $N[P(X,Y)] \setminus P(X,Y)=N[P] \setminus P$ and the other belongs to $G_{-1} \setminus N[P]$. Note that by Lemma 4.1 (d), $Z+(1,0) \notin N[P] \setminus P$. Then we must have $Z+(1,0) \in G_{-1} \setminus N[P]$ and $Z \in N[P] \setminus P$. So, we have that $v(k-1,l) \in F$ and $v(k,l-1) \in F$. By Corollary 1.9, $Z+(0,1) \in G_{-1} \setminus N[P]$. In order to have $\text{deg}_{(h(F))}(v) =2$ we need to check that $e(k,k+1;l)$ and $e(k;l,l+1)$ do not belong to $h(F)$. See Figure 4.8 (a).

\begingroup 
\setlength{\intextsep}{0pt}
\setlength{\columnsep}{10pt}
\begin{wrapfigure}[]{l}{0cm}
\begin{adjustbox}{trim=0cm 0cm 0cm 0cm}
\begin{tikzpicture}[scale=1.5]

\begin{scope}[xshift=0cm]{
\draw[gray,very thin, step=0.5cm, opacity=0.5] (0,0) grid (1,1);

\fill[green!50!white, opacity=0.5] (0,0) rectangle (0.5,0.5);

\fill[orange!50!white, opacity=0.5] (0.5,0) rectangle (1,0.5);

\fill[orange!50!white, opacity=0.5] (0,0.5) rectangle (0.5,1);

\draw[blue, line width=0.5mm] (0,0.5)--++(0.5,0)--++(0,0.5);

\draw[green!50!black, dotted, line width=0.5mm] (0.5,0)--++(0,0.5);

\draw[fill=blue] (0,0.5) circle [radius=0.05];
\draw[fill=blue] (0.5,0.5) circle [radius=0.05];
\draw[fill=blue] (0.5,0) circle [radius=0.05];

{


\draw[black, line width=0.15mm] (0.7,0.45)--++(0,0.1);
\draw[black, line width=0.15mm] (0.75,0.45)--++(0,0.1);
\draw[black, line width=0.15mm] (0.8,0.45)--++(0,0.1);

}

{

\node[left] at (0,0) [scale=1]{\tiny{-1}};
\node[left] at (0,0.5) [scale=1]
{\tiny{$\ell$}};

\node[above] at (0.5,1) [scale=1]
{\tiny{$k$}};
\node[above] at (1, 1) [scale=1]{\tiny{+1}};

\node at (0.25,0.25) [scale=0.8] {\small{$Z$}};
}

\node[below] at (0.5,0) [scale=0.8]{\small{\begin{tabular}{c} Fig. 4.8 (a). Case 2.1. \end{tabular}}};;

} \end{scope}

\begin{scope}[xshift=2cm]{
\draw[gray,very thin, step=0.5cm, opacity=0.5] (0,0) grid (1,1);

{
\draw [blue!40!white, line width=2mm, dotted, opacity=0.5] plot [smooth, tension=0.75] coordinates {(0.25,0.25) (-0.25,0.35)(-0.35,1)(0.25,1.35)(0.65,1.15)(0.75,0.75)};

}

\fill[green!50!white, opacity=0.5] (0,0) rectangle (0.5,0.5);
\fill[green!50!white, opacity=0.5] (0.5,0.5) rectangle (1,1);
\fill[orange!50!white, opacity=0.5] (0.5,0) rectangle (1,0.5);
\fill[orange!50!white, opacity=0.5] (0,0.5) rectangle (0.5,1);

    


\draw[blue, line width=0.5mm] (0,0.5)--++(0.5,0)--++(0,0.5);

\draw[green!50!black, dotted, line width=0.5mm] (0.5,0)--++(0,0.5);

\draw[green!50!black, dotted, line width=0.5mm] (0.5,0.5)--++(0.5,0);

\draw[fill=blue] (0,0.5) circle [radius=0.05];
\draw[fill=blue] (0.5,0.5) circle [radius=0.05];
\draw[fill=blue] (0.5,0) circle [radius=0.05];

{



}

{

\node[left] at (0,0) [scale=1]{\tiny{-1}};
\node[left] at (0,0.5) [scale=1]
{\tiny{$\ell$}};

\node[above] at (0.5,1) [scale=1]
{\tiny{$k$}};
\node[above] at (1, 1) [scale=1]{\tiny{+1}};

\node at (0.25,0.25) [scale=0.8] {\small{$Z$}};
}

\node[below] at (0.5,0) [scale=0.8]{\small{\begin{tabular}{c} Fig. 4.8 (b). Case 2.1. \\  $e(k,k+1;l) \in hb(F)$. \end{tabular}}};;

} \end{scope}

\end{tikzpicture}
\end{adjustbox}
\end{wrapfigure}

Assume for contradiction that $e(k,k+1;l)$ is a shadow edge of $F$. By Observation 4.2, $Z+(1,1) \in N[P] \setminus P$. Since $F$ is an $H$-subtree, there is an $H$-path $P(Z, Z+(1,1))$ contained in $F$. But then $P(Z, Z+(1,1)), Z+(1,0),Z$ is an $H$-cycle, which contradicts Proposition 1.3. Thus $e(k,k+1;l)$ is not a shadow edge of $F$. It follows that $Z+(1,1)$ belongs to $G_{-1} \setminus F$. See Figure 4.8 (b).

Similarly, if $e(k;l,l+1) \in h(F)$, then by Corollary 1.9, $Z+(1,1) \in N[P]$ again, and we obtain the same contradiction as above. End of Case 2.1.

\endgroup

\null 

\begingroup 
\setlength{\intextsep}{0pt}
\setlength{\columnsep}{10pt}
\begin{wrapfigure}[]{r}{0cm}
\begin{adjustbox}{trim=0cm 0.25cm 0cm 0.25cm}
\begin{tikzpicture}[scale=1.5]

\begin{scope}[xshift=0cm]{
\draw[gray,very thin, step=0.5cm, opacity=0.5] (0,0) grid (1,1);

\fill[green!50!white, opacity=0.5] (0,0) rectangle (0.5,0.5);

\fill[orange!50!white, opacity=0.5] (0.5,0) rectangle (1,0.5);

\fill[orange!50!white, opacity=0.5] (0,0.5) rectangle (0.5,1);

\draw[blue, line width=0.5mm] (0,0.5)--++(1,0);

\draw[green!50!black, dotted, line width=0.5mm] (0.5,0)--++(0,0.5);

\draw[fill=blue] (0,0.5) circle [radius=0.05];
\draw[fill=blue] (0.5,0.5) circle [radius=0.05];
\draw[fill=blue] (0.5,0) circle [radius=0.05];

{

\draw[black, line width=0.15mm] (0.45,0.7)--++(0.1,0);
\draw[black, line width=0.15mm] (0.45,0.75)--++(0.1,0);
\draw[black, line width=0.15mm] (0.45,0.8)--++(0.1,0);


}

{

\node[left] at (0,0) [scale=1]{\tiny{-1}};
\node[left] at (0,0.5) [scale=1]
{\tiny{$\ell$}};

\node[above] at (0.5,1) [scale=1]
{\tiny{$k$}};
\node[above] at (1, 1) [scale=1]{\tiny{+1}};

\node at (0.25,0.25) [scale=0.8] {\small{$Z$}};
}

\node[below] at (0.5,0) [scale=0.8]{\small{\begin{tabular}{c} Fig. 4.9 (a). Case 2.2. \end{tabular}}};;

} \end{scope}

\begin{scope}[xshift=2cm]{
\draw[gray,very thin, step=0.5cm, opacity=0.5] (0,0) grid (1,1);

{
}

\fill[green!50!white, opacity=0.5] (0,0) rectangle (0.5,0.5);
\fill[green!50!white, opacity=0.5] (0.5,0.5) rectangle (1,1);
\fill[orange!50!white, opacity=0.5] (0.5,0) rectangle (1,0.5);
\fill[orange!50!white, opacity=0.5] (0,0.5) rectangle (0.5,1);

\begin{scope}
[very thick,decoration={
    markings,
    mark=at position 0.6 with {\arrow{>}}}
    ]
    
    \draw[postaction={decorate}, blue, line width=0.5mm] (0,0.5)--++(0.5,0);
    
    \draw[postaction={decorate}, blue, line width=0.5mm] (0.5,0.5)--++(0.5,0);
    
\end{scope}

\draw[green!50!black, dotted, line width=0.5mm] (0.5,0)--++(0,0.5);

\draw[green!50!black, dotted, line width=0.5mm] (0.5,0.5)--++(0,0.5);

\draw[fill=blue] (0,0.5) circle [radius=0.05];
\draw[fill=blue] (0.5,0.5) circle [radius=0.05];
\draw[fill=blue] (0.5,0) circle [radius=0.05];

{



}

{

\node[left] at (0,0) [scale=1]{\tiny{-1}};
\node[left] at (0,0.5) [scale=1]
{\tiny{$\ell$}};

\node[above] at (0.5,1) [scale=1]
{\tiny{$k$}};
\node[above] at (1, 1) [scale=1]{\tiny{+1}};

\node at (0.25,0.25) [scale=0.8] {\small{$Z$}};
}

\node[below] at (0.5,0) [scale=0.8]{\small{\begin{tabular}{c} Fig. 4.9 (b). Case 2.2. \\  $e(k;l,l+1) \in hb(F)$. \end{tabular}}};;

} \end{scope}

\end{tikzpicture}
\end{adjustbox}
\end{wrapfigure}

\noindent \textit{CASE 2.2: $e(k-1,k;l) \in H$ and $e(k,k+1;l) \in H$.} Then $e(k;l,l+1) \notin H$. By Observation 4.2, exactly one of $Z$ and $Z+(1,0)$ belongs to $N[P] \setminus P$  and the other belongs to $G_{-1} \setminus N[P]$. By symmetry, we may assume WLOG that $Z \in N[P] \setminus P$ and $Z+(1,0) \in G_{-1} \setminus N[P]$. By Corollary 1.9, $Z+(0,1) \in G_{-1} \setminus N[P]$. In order to have $\text{deg}_{h(F)}(v)=2$ we need to check that $e(k;l,l+1)$ and $e(k,k+1;l)$ do not belong to $h(F)$.

Assume for contradiction that $e(k;l,l+1)$ is a shadow edge of $F$. Then $Z+(1,1) \in F$. But now, if we orient $H$ as directed cycle $\overrightarrow{K}_H$, $Z$ and $Z+(1,1)$ are on different sides of $\overrightarrow{K}_H$, contradicting  Corollary 1.11. Therefore, $e(k;l,l+1)$ cannot be a shadow edge of $F$. See Figure 4.9 (b). It follows that $Z+(1,1) \in G_{-1} \setminus N[P]$.

Similarly, if $e(k,k+1;l) \in h(F)$, then by Corollary 1.9, $Z+(1,1) \in N[P]$ again, and we obtain the same contradiction as above. End of Case 2.2. End of Case 2.  End of proof for (i).

\null

\noindent \textit{Proof of (ii).} Let $u,v$ be vertices in $F$. Orient the subpath $P=P(u,v)$ of $H$ from $u$ to $v$, labelled $u=u_0, u_1, ...$. If $E(P) \subset E(h(F))$, then we're done. Otherwise, let $s$ be the number of shadow edges of $h(F)$ that are incident on $P$. Let $u_{i_1}$ be the first vertex of $P$ after $u$ such that $u_{i_1} \in F$ but $u_{i_1+1} \notin F$. Let $u_{i'_1}$ be the first vertex of $P$ after $u_{i_1}$ that is in $F$. For $j \in \{ 2,...,s \}$ let $u_{i_j}$ be the first vertex of $P$ after $u_{i'_{j-1}}$ such that $u_{i_j} \in F$ but $u_{i_j+1} \notin F$, and let $u_{i'_j}$ be the first vertex of $P$ after $u_{i_j}$ that is in $F$. Note that $P(u_{i'_s},v) \subset P(u,v)$ and that $P(u_{i'_s},v)$ is contained in $F$. 

\noindent We claim that for $1 \leq j \leq s$, $\{u_{i_j},u_{i'_j} \}$ is a shadow edge of $F$. It follows from this claim that $P(u, u_{i_1}), (u_{i_1}, u_{i'_1}),P( u_{i'_1},u_{i_2})$, $(u_{i_2}, u_{i'_2}), ..., P( u_{i'_{s-1}},u_{i_s}), (u_{i_s}, u_{i'_s}), P(u_{i'_s,v}) $ is contained in $F$. It remains to check that the claim is true.

\noindent \textit{Proof of Claim.} For definiteness let $u_{i_j}=v(k,l)$, $u_{i_j+1}=v(k+1,l)$. Let $R(k-1, l-1)=Z$. Since $u_{i_j+1}$ is not in $F$, $Z+(1,0)$ and $Z+(1,1)$ belong to $G_{-1} \setminus F$. Since $u_{i_j}$ is in $F$, at least one of $Z$ and $Z+(0,1)$ is a box of $F$.  

We will first show that $u_{i_j-1}=v(k-1,l)$.
Assume for contradiction that $u_{i_j-1}=v(k,l-1)$ or $u_{i_j-1}=v(k,l+1)$. By symmetry we only need to check one of the two. For definiteness assume that $u_{i_j-1}=v(k,l-1)$.

\begingroup 
\setlength{\intextsep}{0pt}
\setlength{\columnsep}{20pt}
\begin{wrapfigure}[]{l}{0cm}
\begin{adjustbox}{trim=0cm 0cm 0cm 0cm}
\begin{tikzpicture}[scale=1.5]

\begin{scope}[xshift=0cm]{
\draw[gray,very thin, step=0.5cm, opacity=0.5] (0,0) grid (1,1.5);

\fill[blue!50!white, opacity=0.5] (0,0) rectangle (0.5,0.5);
\fill[green!50!white, opacity=0.5] (0,0.5) rectangle (0.5,1);

\fill[orange!50!white, opacity=0.5] (0.5,0.5) rectangle (1,1.5);
\fill[orange!50!white, opacity=0.5] (0,1) rectangle (0.5,1.5);

\begin{scope}
[very thick,decoration={
    markings,
    mark=at position 0.6 with {\arrow{>}}}
    ]
    
    \draw[postaction={decorate}, blue, line width=0.5mm] (0.5,1)--++(0.5,0);
    \draw[postaction={decorate}, blue, line width=0.5mm] (0.5,0.5)--++(0,0.5);
    
\end{scope}

\draw[blue, line width=0.5mm] (0,0.5)--++(0,0.5);
\draw[blue, line width=0.5mm] (0,0)--++(0.5,0);

\draw[fill=blue] (1,1) circle [radius=0.05];
\draw[fill=blue] (0.5,1) circle [radius=0.05];

{


\draw[black, line width=0.15mm] (0.2,0.45)--++(0,0.1);
\draw[black, line width=0.15mm] (0.25,0.45)--++(0,0.1);
\draw[black, line width=0.15mm] (0.3,0.45)--++(0,0.1);

}

{
\node[right] at (1,0) [scale=1]{\tiny{-2}};
\node[right] at (1,0.5) [scale=1]{\tiny{-1}};
\node[right] at (1,1) [scale=1]
{\tiny{$\ell$}};

\node[above] at (0,1.5) [scale=1]{\tiny{-1}};
\node[above] at (0.5,1.5) [scale=1]
{\tiny{$k$}};
\node[above] at (1, 1.5) [scale=1]{\tiny{+1}};

\node[above] at (0.5,1) [scale=0.8] {\small{$u_{_{i_j}}$}};
\node[above] at (1,1) [scale=0.8] {\small{$u_{_{i_j+1}}$}};

\node at (0.25,0.75) [scale=0.8] {\small{$Z$}};
}

\node[below] at (0.5,0) [scale=0.8]{\small{\begin{tabular}{c} Fig. 4.10 (a). \end{tabular}}};;

} \end{scope}

\begin{scope}[xshift=1.75cm]{
\draw[gray,very thin, step=0.5cm, opacity=0.5] (0,0) grid (1,1.5);

\fill[blue!50!white, opacity=0.5] (0,0) rectangle (0.5,0.5);
\fill[green!50!white, opacity=0.5] (0,0.5) rectangle (0.5,1);

\fill[orange!50!white, opacity=0.5] (0.5,0.5) rectangle (1,1.5);
\fill[orange!50!white, opacity=0.5] (0,1) rectangle (0.5,1.5);

\begin{scope}
[very thick,decoration={
    markings,
    mark=at position 0.6 with {\arrow{>}}}
    ]
    
    \draw[postaction={decorate}, blue, line width=0.5mm] (0.5,1)--++(0.5,0);
    \draw[postaction={decorate}, blue, line width=0.5mm] (0.5,0.5)--++(0,0.5);
    
\end{scope}

\draw[blue, line width=0.5mm] (0,0.5)--++(0,0.5);
\draw[blue, line width=0.5mm] (0,0)--++(0.5,0);

\draw[blue, line width=0.5mm] (0.5,0.5)--++(0.5,0);

\draw[fill=blue] (1,1) circle [radius=0.05];
\draw[fill=blue] (0.5,1) circle [radius=0.05];

{

\draw[black, line width=0.15mm] (-0.05,0.2)--++(0.1,0);
\draw[black, line width=0.15mm] (-0.05,0.25)--++(0.1,0);
\draw[black, line width=0.15mm] (-0.05,0.3)--++(0.1,0);

\draw[black, line width=0.15mm] (0.45,0.2)--++(0.1,0);
\draw[black, line width=0.15mm] (0.45,0.25)--++(0.1,0);
\draw[black, line width=0.15mm] (0.45,0.3)--++(0.1,0);

\draw[black, line width=0.15mm] (0.2,0.45)--++(0,0.1);
\draw[black, line width=0.15mm] (0.25,0.45)--++(0,0.1);
\draw[black, line width=0.15mm] (0.3,0.45)--++(0,0.1);

}

{
\node[right] at (1,0) [scale=1]{\tiny{-2}};
\node[right] at (1,0.5) [scale=1]{\tiny{-1}};
\node[right] at (1,1) [scale=1]
{\tiny{$\ell$}};

\node[above] at (0,1.5) [scale=1]{\tiny{-1}};
\node[above] at (0.5,1.5) [scale=1]
{\tiny{$k$}};
\node[above] at (1, 1.5) [scale=1]{\tiny{+1}};

\node[above] at (0.5,1) [scale=0.8] {\small{$u_{_{i_j}}$}};
\node[above] at (1,1) [scale=0.8] {\small{$u_{_{i_j+1}}$}};

\node at (0.25,0.75) [scale=0.8] {\small{$Z$}};
}

\node[below] at (0.5,0) [scale=0.8]{\small{\begin{tabular}{c} Fig. 4.10 (b).  \end{tabular}}};;

} \end{scope}

\end{tikzpicture}
\end{adjustbox}
\end{wrapfigure}

\noindent Note that if $Z+(0,1) \in F$, by Observation 4.2, $Z+(0,1) \in N[P(X,Y)] \setminus P(X,Y)$. But then $Z+(0,1)$ is neither a leaf nor a switchable box, contradicting Lemma 4.1 (d). It remains to check the case where $Z+(0,1) \in G_{-1} \setminus F$ and $Z \in N[P(X,Y)] \setminus P(X,Y)$. See Figure 4.10 (a). Using Lemma 4.1 (d) again, we have that $Z+(0,-1) \in P(X,Y)$, $e(k-1;l-1,l) \in H$ and $e(k-1,k;l-1) \notin H$. By Lemma 4.1 (b) and (c), we have that $e(k-1,k;l-2)\in H$. Note that if $Z+(0,-1)$ is not an end-box of $P(X,Y)$, by Lemma 4.1 (b), $e(k,k+1;l-1) \in H$, $e(k-1;l-2,l-1) \notin H$ and $e(k;l-2,l-1) \notin H$. But then after $Z \mapsto Z+(1,0)$, $Z+(0,-1) \in P(X,Y)$ is switchable, contradicting that $F$ is a looping fat path. See Figure 4.10(b) It remains to check the case where $Z+(0,-1)$ is an end-box of $Z$. WLOG assume that $Z+(0,-1)=X$. There are three possibilities: $X=Y+(2,0)$, $X=Y+(0,-2)$, and $X=Y+(-2,0)$. 

\endgroup 

\null 

\noindent \textit{CASE 1: $Y=X+(2,0)$.} Then $F$ must be northern, so $e(k+1;l-2,l-1)\in H$ and $e(k+1;l-1,l) \in H$. But then there is an $H$-cycle $P(Y,X), X+(0,1), X+(0,2), X+(1,2), X+(2,2), X+(2,1), Y$, which contradicts Proposition 1.3.  See Figure 4.11 (a). End of Case 1.

\null 

\begingroup 
\setlength{\intextsep}{0pt}
\setlength{\columnsep}{10pt}
\begin{wrapfigure}[]{r}{0cm}
\begin{adjustbox}{trim=0cm 0.5cm 0cm 0.25cm}
\begin{tikzpicture}[scale=1.5]

\begin{scope}[xshift=0cm]{
\draw[gray,very thin, step=0.5cm, opacity=0.5] (0,0) grid (1.5,1.5);

\fill[blue!50!white, opacity=0.5] (0,0) rectangle (0.5,0.5);
\fill[green!50!white, opacity=0.5] (0,0.5) rectangle (0.5,1);

\fill[orange!50!white, opacity=0.5] (0.5,0.5) rectangle (1,1.5);
\fill[orange!50!white, opacity=0.5] (0,1) rectangle (0.5,1.5);

\begin{scope}
[very thick,decoration={
    markings,
    mark=at position 0.6 with {\arrow{>}}}
    ]
    
    \draw[postaction={decorate}, blue, line width=0.5mm] (0.5,1)--++(0.5,0);
    \draw[postaction={decorate}, blue, line width=0.5mm] (0.5,0.5)--++(0,0.5);
    
\end{scope}

\draw[blue, line width=0.5mm] (0,0.5)--++(0,0.5);
\draw[blue, line width=0.5mm] (0,0)--++(0.5,0);

\draw[blue, line width=0.5mm] (1,0)--++(0,1);

\draw[fill=blue] (1,1) circle [radius=0.05];
\draw[fill=blue] (0.5,1) circle [radius=0.05];

{


\draw[black, line width=0.15mm] (0.2,0.45)--++(0,0.1);
\draw[black, line width=0.15mm] (0.25,0.45)--++(0,0.1);
\draw[black, line width=0.15mm] (0.3,0.45)--++(0,0.1);

}

{
\node[left] at (0,0) [scale=1]{\tiny{-2}};
\node[left] at (0,0.5) [scale=1]{\tiny{-1}};
\node[left] at (0,1) [scale=1]
{\tiny{$\ell$}};

\node[above] at (0,1.5) [scale=1]{\tiny{-1}};
\node[above] at (0.5,1.5) [scale=1]
{\tiny{$k$}};
\node[above] at (1, 1.5) [scale=1]{\tiny{+1}};

\node[above] at (0.5,1) [scale=0.8] {\small{$u_{_{i_j}}$}};
\node[above] at (1,1) [scale=0.8] {\small{$u_{_{i_j+1}}$}};

\node at (0.25,0.75) [scale=0.8] {\small{$Z$}};
\node at (0.25,0.25) [scale=0.8] {\small{$X$}};
\node at (1.25,0.25) [scale=0.8] {\small{$Y$}};
}

\node[below] at (0.75,0) [scale=0.8]{\small{\begin{tabular}{c} Fig. 4.11 (a). Case 1. \end{tabular}}};;

} \end{scope}

\begin{scope}[xshift=2.25cm]{
\draw[gray,very thin, step=0.5cm, opacity=0.5] (0,0) grid (1,1.5);

\fill[blue!50!white, opacity=0.5] (0.5,0) rectangle (0.5,0.5);
\fill[green!50!white, opacity=0.5] (0.5,0.5) rectangle (1,1);

\fill[orange!50!white, opacity=0.5] (0.5,1) rectangle (1,1.5);

\fill[blue!50!white, opacity=0.5] (0,0) rectangle (1,0.5);

\begin{scope}
[very thick,decoration={
    markings,
    mark=at position 0.6 with {\arrow{>}}}
    ]
    
    \draw[postaction={decorate}, blue, line width=0.5mm] (1,0.5)--++(0,0.5);
    
\end{scope}

\draw[blue, line width=0.5mm] (0.5,0.5)--++(0,0.5);
\draw[blue, line width=0.5mm] (0.5,0)--++(0.5,0);

\draw[blue, line width=0.5mm] (0,0)--++(0.5,0);
\draw[blue, line width=0.5mm] (0,0.5)--++(0.5,0);
\draw[blue, line width=0.5mm] (1,0)--++(0,0.5);

\draw[fill=blue] (1,1) circle [radius=0.05];

{

\draw[black, line width=0.15mm] (0.45,0.2)--++(0.1,0);
\draw[black, line width=0.15mm] (0.45,0.25)--++(0.1,0);
\draw[black, line width=0.15mm] (0.45,0.3)--++(0.1,0);

\draw[black, line width=0.15mm] (-0.05,0.2)--++(0.1,0);
\draw[black, line width=0.15mm] (-0.05,0.25)--++(0.1,0);
\draw[black, line width=0.15mm] (-0.05,0.3)--++(0.1,0);


}

{
\node[left] at (0,0) [scale=1]{\tiny{-2}};
\node[left] at (0,0.5) [scale=1]{\tiny{-1}};
\node[left] at (0,1) [scale=1]
{\tiny{$\ell$}};

\node[above] at (0,1.5) [scale=1]{\tiny{-2}};
\node[above] at (0.5,1.5) [scale=1]
{\tiny{-1}};
\node[above] at (1, 1.5) [scale=1]{\tiny{$k$}};

\node[above] at (1,1) [scale=0.8] {\small{$u_{_{i_j}}$}};

\node at (0.75,0.75) [scale=0.8] {\small{$Z$}};
\node at (0.75,0.25) [scale=0.8] {\small{$X$}};
}

\node[below] at (0.5,0) [scale=0.8]{\small{\begin{tabular}{c} Fig. 4.11 (b). \\ Case 2.1. \end{tabular}}};;

} \end{scope}

\begin{scope}[xshift=4cm]{
\draw[gray,very thin, step=0.5cm, opacity=0.5] (0,0) grid (1,1.5);

\fill[orange!50!white, opacity=0.5] (0,1) rectangle (1,1.5);
\fill[green!50!white, opacity=0.5] (0,0.5) rectangle (0.5,1);

\fill[blue!50!white, opacity=0.5] (0,0) rectangle (1,0.5);

\begin{scope}
[very thick,decoration={
    markings,
    mark=at position 0.6 with {\arrow{>}}}
    ]
    
    \draw[postaction={decorate}, blue, line width=0.5mm] (0.5,0.5)--++(0,0.5);
     \draw[postaction={decorate}, blue, line width=0.5mm] (0.5,1)--++(0.5,0);
    
\end{scope}

\draw[blue, line width=0.5mm] (0,0)--++(0,1);
\draw[blue, line width=0.5mm] (0,0)--++(1,0);
\draw[blue, line width=0.5mm] (0.5,0.5)--++(0.5,0);

\draw[fill=blue] (0.5,1) circle [radius=0.05];
\draw[fill=blue] (1,1) circle [radius=0.05];
{

\draw[black, line width=0.15mm] (0.55,0.2)--++(-0.1,0);
\draw[black, line width=0.15mm] (0.55,0.25)--++(-0.1,0);
\draw[black, line width=0.15mm] (0.55,0.3)--++(-0.1,0);

\draw[black, line width=0.15mm] (1.05,0.2)--++(-0.1,0);
\draw[black, line width=0.15mm] (1.05,0.25)--++(-0.1,0);
\draw[black, line width=0.15mm] (1.05,0.3)--++(-0.1,0);

}

{
\node[left] at (0,0) [scale=1]{\tiny{-2}};
\node[left] at (0,0.5) [scale=1]{\tiny{-1}};
\node[left] at (0,1) [scale=1]
{\tiny{$\ell$}};

\node[above] at (0,1.5) [scale=1]{\tiny{-1}};
\node[above] at (1,1.5) [scale=1]{\tiny{+1}};
\node[above] at (0.5,1.5) [scale=1]{\tiny{$k$}};

\node[above] at (0.5,1) [scale=0.8]{\small{$u_{_{i_j}}$}};

\node[above] at (1,1) [scale=0.8]{\small{$u_{_{i_j+1}}$}};

\node at (0.25,0.75) [scale=0.8]{\small{$Z$}};
\node at (0.25,0.25) [scale=0.8]{\small{$X$}};
}

\node[below] at (0.5,0) [scale=0.8]{\small{\begin{tabular}{c} Fig. 4.11 (c). \\ Case 2.2. \end{tabular}}};;

} \end{scope}

\end{tikzpicture}
\end{adjustbox}
\end{wrapfigure}

\noindent \textit{CASE 2: $Y=X+(0,-2)$.} Then $F$ can be western or eastern.

\null 

\noindent  \textit{CASE 2.1: $F$ is western.} Then $e(k-2,k-1;l-2) \in H$ and $e(k;l-2,l-1) \in H$. It follows that $X+(-1,0) \in P(X,Y)$ and that $e(k-2;l-2,l-1) \notin H$. But then $X+(-1,0)$ is switchable, contradicting that $F$ is a looping fat path. See Figure 4.11 (b). End of Case 2.1.

\null 

\noindent \textit{CASE 2.2: $F$ is eastern.} Then $e(k-1;l-2,l-1) \in H$ and $e(k,k+1;l-2) \in H$. It follows that $e(k,k+1;l-1) \in H$. But then $X+(1,0) \in P(X,Y)$ is switchable, contradicting  that $F$ is a looping fat path. See Figure 4.11 (c). End of Case 2.2. End of Case 2.

\endgroup 

\null 

\begingroup 
\setlength{\intextsep}{0pt}
\setlength{\columnsep}{10pt}
\begin{wrapfigure}[]{l}{0cm}
\begin{adjustbox}{trim=0cm 0.25cm 0cm 0.5cm}
\begin{tikzpicture}[scale=1.5]

\begin{scope}[xshift=0cm]{
\draw[gray,very thin, step=0.5cm, opacity=0.5] (0,0) grid (2,1.5);

\fill[blue!50!white, opacity=0.5] (1,0.5) rectangle (1.5,1);
\fill[blue!50!white, opacity=0.5] (0,0.5) rectangle (0.5,1.5);

\fill[green!50!white, opacity=0.5] (1,1) rectangle (1.5,1.5);

\fill[orange!50!white, opacity=0.5] (1.5,1) rectangle (2,1.5);

\begin{scope}
[very thick,decoration={
    markings,
    mark=at position 0.6 with {\arrow{>}}}
    ]
    
    \draw[postaction={decorate}, blue, line width=0.5mm] (1.5,1.5)--++(0.5,0);
    \draw[postaction={decorate}, blue, line width=0.5mm] (1.5,1)--++(0,0.5);
    
\end{scope}

\draw[blue, line width=0.5mm] (1,1)--++(0,0.5);
\draw[blue, line width=0.5mm] (1,0.5)--++(0.5,0);

\draw[blue, line width=0.5mm] (0,0.5)--++(0.5,0)--++(0,1);
\draw[blue, line width=0.5mm] (0.5,0)--++(0.5,0);

\draw[blue, line width=0.5mm] (1,0.5)--++(0,0.5);
\draw[blue, line width=0.5mm] (1.5,1)--++(0.5,0);

\draw[fill=blue] (2,1.5) circle [radius=0.05];
\draw[fill=blue] (1.5,1.5) circle [radius=0.05];

{

\draw[black, line width=0.15mm] (1.2,0.95)--++(0,0.1);
\draw[black, line width=0.15mm] (1.25,0.95)--++(0,0.1);
\draw[black, line width=0.15mm] (1.3,0.95)--++(0,0.1);

}

{
\node[left] at (0,0.5) [scale=1]{\tiny{-2}};
\node[left] at (0,1) [scale=1]{\tiny{-1}};
\node[left] at (0,1.5) [scale=1]{\tiny{$\ell$}};

\node[above] at (1,1.5) [scale=1]{\tiny{-1}};
\node[above] at (1.5,1.5) [scale=1]{\tiny{$k$}};
\node[above] at (2,1.5) [scale=1]{\tiny{+1}};

\node[below] at (1.65,1.5) [scale=0.8]{\small{$u_{_{i_j}}$}};
\node[below] at (2.15,1.45) [scale=0.8]{\small{$u_{_{i_j+1}}$}};

\draw[fill=blue] (2,1.5) circle [radius=0.05];
\draw[fill=blue] (1.5,1.5) circle [radius=0.05];

\node at (1.25,1.25) [scale=0.8]{\small{$Z$}};
\node at (1.25,0.75) [scale=0.8]{\small{$X$}};

\node at (0.25,0.75) [scale=0.8]{\small{$Y$}};
}

\node[below] at (1,0) [scale=0.8]{\small{\begin{tabular}{c} Fig. 4.12. Case 3. \end{tabular}}};;

} \end{scope}

\end{tikzpicture}
\end{adjustbox}
\end{wrapfigure}

\null 

\noindent \textit{CASE 3: $Y=X+(-2,0)$.} Then $F$ must be northern, so $e(k-1;l-2,l-1) \in H$. It follows that $e(k,k+1;l-1) \in H$. But then $X+(0,1) \mapsto X+(1,1)$, $X+(-1,0) \mapsto X$ is a short weakening, contradicting that $F$ is a looping fat path. See Figure 4.12. End of Case 3. This concludes the proof that $u_{i_j-1}=v(k-1,l)$.

\begingroup 
\setlength{\intextsep}{0pt}
\setlength{\columnsep}{20pt}
\begin{wrapfigure}[]{r}{0cm}
\begin{adjustbox}{trim=0cm 0.75cm 0cm 0.25cm}
\begin{tikzpicture}[scale=1.5]

{
\draw [orange!50!white, line width=2mm, dotted, opacity=0.5] plot [smooth, tension=0.75] coordinates {(1.25,0.25) (1.75,0.1)(1.9,-0.25)(1.75,-0.6)(1.25,-0.75)};

\draw [blue!40!white, line width=2mm, dotted, opacity=0.5] plot [smooth, tension=0.75] coordinates {(0.25,0.25) (-0.25,0.1)(-0.4,-0.25)(-0.25,-0.6)(0.25,-0.75)};
}

\begin{scope}[xshift=0cm]{
\draw[gray,very thin, step=0.5cm, opacity=0.5] (0,0) grid (1.5,1);

\fill[blue!40!white, opacity=0.5] (0,0) rectangle (0.5,0.5);
\fill[green!40!white, opacity=0.5] (0.5,0) rectangle (1,0.5);
\fill[orange!50!white, opacity=0.5] (1,0) rectangle (1.5,1);

\begin{scope}
[very thick,decoration={
    markings,
    mark=at position 0.6 with {\arrow{>}}}
    ]
    
    \draw[postaction={decorate}, blue, line width=0.5mm] (0.5,0.5)--++(0.5,0);
    
    \draw[postaction={decorate}, blue, line width=0.5mm] (1,0.5)--++(0.5,0);
    
\end{scope}

\draw[blue, line width=0.5mm] (0.5,0)--++(0.5,0);


{

\draw[black, line width=0.15mm] (0.45,0.2)--++(0.1,0);
\draw[black, line width=0.15mm] (0.45,0.25)--++(0.1,0);
\draw[black, line width=0.15mm] (0.45,0.3)--++(0.1,0);
}

\node[above] at (1,0.5) [scale=0.8]{\small{$u_{_{i_j}}$}};
\node[above] at (1.5,0.5) [scale=0.8]{\small{$u_{_{i_j+1}}$}};

\draw[fill=blue] (1,0.5) circle [radius=0.05];
\draw[fill=blue] (1.5,0.5) circle [radius=0.05];

{
\node[left] at (0,0) [scale=1]
{\tiny{-1}};
\node[left] at (0,0.5) [scale=1]
{\tiny{$\ell$}};
\node[left] at (0,1) [scale=1]
{\tiny{+1}};

\node[above] at (0.5,1) [scale=1]
{\tiny{-1}};
\node[above] at (1,1) [scale=1]
{\tiny{$k$}};
\node[above] at (1.5, 1) [scale=1]
{\tiny{+1}};


\node at (0.75,0.25) [scale=0.8] {\small{$Z$}};

}


} \end{scope}

\begin{scope}[yshift=-1.5cm]{
\draw[gray,very thin, step=0.5cm, opacity=0.5] (0,0) grid (1.5,1);

\fill[blue!40!white, opacity=0.5] (0,0.5) rectangle (0.5,1);
\fill[green!40!white, opacity=0.5] (0.5,0.5) rectangle (1,1);
\fill[orange!50!white, opacity=0.5] (1,0) rectangle (1.5,1);

\begin{scope}
[very thick,decoration={
    markings,
    mark=at position 0.6 with {\arrow{<}}}
    ]
    
    \draw[postaction={decorate}, blue, line width=0.5mm] (0.5,0.5)--++(0.5,0);
    
    \draw[postaction={decorate}, blue, line width=0.5mm] (1,0.5)--++(0.5,0);
    
\end{scope}

\draw[blue, line width=0.5mm] (0.5,1)--++(0.5,0);


{

\draw[black, line width=0.15mm] (0.45,0.7)--++(0.1,0);
\draw[black, line width=0.15mm] (0.45,0.75)--++(0.1,0);
\draw[black, line width=0.15mm] (0.45,0.8)--++(0.1,0);
}

\node[above] at (1,0.5) [scale=0.8]{\small{$u_{_{i'_j}}$}};
\node[above] at (1.5,0.5) [scale=0.8]{\small{$u_{_{i'_j-1}}$}};

\draw[fill=blue] (1,0.5) circle [radius=0.05];
\draw[fill=blue] (1.5,0.5) circle [radius=0.05];

{
\node[left] at (0,0) [scale=1]
{\tiny{-1}};
\node[left] at (0,0.5) [scale=1]
{\tiny{$b$}};
\node[left] at (0,1) [scale=1]
{\tiny{+1}};

\node[above] at (0.5,1) [scale=1]
{\tiny{-1}};
\node[above] at (1,1) [scale=1]
{\tiny{$a$}};
\node[above] at (1.5, 1) [scale=1]
{\tiny{+1}};


\node at (0.75,0.75) [scale=0.8] {\small{$Z'$}};

}

\node[below] at (0.75,0) [scale=0.8]{\small{\begin{tabular}{c} Fig. 4.13. 
 \end{tabular}}};;

} \end{scope}

\end{tikzpicture}
\end{adjustbox}
\end{wrapfigure}

\noindent Now, by Corollary 1.9, exactly one of $Z$ and $Z+(0,1)$ belongs to $F$. WLOG, assume that $Z \in F$. By Observation 4.2, since $Z+(1,0)$ is not in $F$, $Z \in N[P(X,Y)] \setminus P(X,Y)$. Note that this means that $e(k;l-1,l)$ is a shadow edge of $F$. By Lemma 4.1 (d), $e(k-1,k;l-1) \in H$, $e(k-1;l-1,l) \notin H$ and $Z+(-1,0) \in P(X,Y)$.
 
\noindent It remains to check that $u_{i'_j}=v(k,l-1)$. Assume for contradiction that $u_{i'_j}=v(a,b) \neq v(k,l-1)$. For definiteness assume that $u_{i'_j-1}=v(a+1,b)$. Note that, by the proof that $u_{i_j-1}=v(k-1,l)$, we have that $u_{i'_j+1}=v(a-1,b)$. Let $Z'=R(a-1,b)$. By Corollary 1.9, exactly one of $Z'$ and $Z'+(0,-1)$ belongs to $F$. Note that if $Z' \notin F$, then $Z'+(0,-1)= \Phi ((u_{i'_j},u_{i'_j+1}), \text{left}) \in F$. But then, $Z \in F$ and $ Z=\Phi ((u_{i_j-1},u_{i_j}), \text{right})$ contradicting Corollary 1.11. Thus we must have $Z' \in F$.

By Lemma 4.1 (d), we must have $e(a-1,a;b+1)\in H$, $e(a-1;b,b+1) \notin H$ and $Z'+(-1,0) \in P(X,Y)$.  Let $P(Z+(1,0), Z'+(1,0))$ be the $H$-path from $Z+(1,0)$ to $Z'+(1,0)$ contained in the $H$-walk $\Phi(\overrightarrow{K}(((u_{i_j},u_{i_j+1}), (u_{i'_j-1},u_{i'_j})), \text{right})$. Observe that $P(Z+(1,0), Z'+(1,0)) \subset G_{-1}\setminus F$. Since $F$ is an $H$-subtree, $F$ is $H$-path-connected so there is an $H$-path $P(Z',Z)$ contained in $F$. But then, then $P(Z',Z), P(Z+(1,0), Z'+(1,0))$ is an $H$-cycle, which contradicts Proposition 1.3. See Figure 4.13. Thus we must have $u_{i'_j}=v(k,l-1)$. $\square$

\endgroup 

\null 

\noindent \textbf{Proposition 4.4.} Let $G$ be an $m\times n$ grid graph, let $H$ be a Hamiltonian cycle of $G$ and let $F=G\langle N[P(X,Y)] \rangle $ be a looping fat path in $G$ following a leaf $L$. Then $h(F)$ does not have consecutive colinear edges other than the left and right collinear edges in the $A_1$-type of $F$ following $L$.

\null 

\noindent \textit{Proof.}  We need to check that $h(F)$ does not have consecutive colinear edges in the case where one of those colinear edges is a left or right colinear edge of the $A_1$-type of $F$ and in the case where neither of those colinear edges is a left or right colinear edge of the $A_1$-type of $F$. We divide the proof into Lemmas 4.5 and 4.6.

\null 

\begingroup 
\setlength{\intextsep}{0pt}
\setlength{\columnsep}{20pt}
\begin{wrapfigure}[]{r}{0cm}
\begin{adjustbox}{trim=0cm 0.5cm 0cm 0.25cm}
\begin{tikzpicture}[scale=1.5]
\begin{scope}[xshift=0cm]{
\draw[gray,very thin, step=0.5cm, opacity=0.5] (0,0) grid (1.5,1.5);

\draw[green!50!black, dotted, line width=0.5mm] (1,0.5)--++(0,0.5);


\draw[blue, line width=0.5mm] (0,0.5)--++(0.5,0)--++(0,-0.5);
\draw[blue, line width=0.5mm] (1.5,0.5)--++(-0.5,0)--++(0,-0.5);
\draw[blue, line width=0.5mm] (0.5,1.5)--++(0,-0.5)--++(0.5,0)--++(0,0.5);

\draw[fill=blue] (1,0.5) circle [radius=0.05];

\node[right] at (1.5,0.5) [scale=1]
{\tiny{$\ell{+2}$}};

\node[above] at (0.5, 1.5) [scale=1]
{\tiny{$k{-}1$}};

\node at (0.75,1.25) [scale=0.8] {\small{L}};
\node at (1.25,0.25) [scale=0.8] {\small{X}};
\node at (0.25,0.25) [scale=0.8] {\small{Y}};

\node[below] at (0.75,0) [scale=0.8]{\small{\begin{tabular}{c} Fig. 4.14. 
 \end{tabular}}};;

} \end{scope}

\end{tikzpicture}
\end{adjustbox}
\end{wrapfigure}

\noindent \textbf{Lemma 4.5.} The shadow of $F$ does not have a pair of consecutive colinear edges in the case where one edge of the pair is one of the left or right colinear edges of the $A_1$-type that follows $L$.

\null 

\noindent \textit{Proof.} For definiteness, assume that $F$ is a southern looping fat path following the leaf $R(k-1,l+3)$. See Figure 4.14. Assume for contradiction that $h(F)$ does have consecutive colinear edges and one of those edges is one of the right or left colinear edges of the $A_1$-type that follows $e(k-1,k;l+3)$. For definiteness, assume that one of those consecutive colinear edges is one of the right colinear edges of the $A_1$-type that follows $L$. Then, either $e(k;l+1,l+2)$, $e(k;l+2,l+3)$ is a pair of consecutive colinear edges or $e(k;l-1,l)$, $e(k;l,l+1)$ is a pair of consecutive colinear edges. If the former, then $\text{deg}_{(h(F))}(v(k,l+2)) =3$ contradicting Lemma 4.3, so we only need to check the latter. Suppose then, that $e(k;l-1,l)$, $e(k;l,l+1)$ is a pair of consecutive colinear edges. Note that $X$ and $X+(0,-1)$ belong to $F$. If the edge $e(k;l-1,l)$ is in $h(F)$ then it is either a shadow edge or it belongs to $H$. We show that both cases lead to contradictions.  

\endgroup 

\null 

\begingroup 
\setlength{\intextsep}{0pt}
\setlength{\columnsep}{10pt}
\begin{wrapfigure}[]{r}{0cm}
\begin{adjustbox}{trim=0cm 0.25cm 0cm 0.25cm}
\begin{tikzpicture}[scale=1.5]
\begin{scope}[xshift=0cm]{
\draw[gray,very thin, step=0.5cm, opacity=0.5] (0,0) grid (2,1.5);

\fill[blue!50!white, opacity=0.5](1,0.5) rectangle (1.5,1.5);
\fill[blue!50!white, opacity=0.5](0,1) rectangle (0.5,1.5);

\draw[blue, line width=0.5mm] (0,1.5)--++(0.5,0)--++(0,-1);
\draw[blue, line width=0.5mm] (1.5,1.5)--++(-0.5,0)--++(0,-1);
\draw[blue, line width=0.5mm] (1,0)--++(0,0.5);
\draw[blue, line width=0.5mm] (1,0)--++(0,0.5);

\draw[blue, line width=0.5mm] (1.5,0)--++(0,0.5)--++(0.5,0);

\draw[blue, line width=0.5mm] (1.5,1.5)--++(0,-0.5)--++(0.5,0);

\draw[black, line width=0.15mm] (1.45,0.70)--++(0.1,0);
\draw[black, line width=0.15mm] (1.45,0.75)--++(0.1,0);
\draw[black, line width=0.15mm] (1.45,0.80)--++(0.1,0);

\node[left] at (0,0) [scale=1]
{\tiny{-1}};
\node[left] at (0,0.5) [scale=1]
{\tiny{$\ell$}};
\node[left] at (0,1) [scale=1]
{\tiny{+1}};
\node[left] at (0,1.5) [scale=1]
{\tiny{+2}};

\node[above] at (1, 1.5) [scale=1]
{\tiny{$k$}};
\node[above] at (1.5, 1.5) [scale=1]
{\tiny{+1}};
\node[above] at (2, 1.5) [scale=1]
{\tiny{+2}};

\node at (1.25,1.25) [scale=0.8] {\small{X}};
\node at (0.25,1.25) [scale=0.8] {\small{Y}};

\node[below] at (1,0) [scale=0.8]{\small{\begin{tabular}{c} Fig. 4.15 (a). Case 1.1.
 \end{tabular}}};;

} \end{scope}

\begin{scope}[xshift=2.75 cm]{
\draw[gray,very thin, step=0.5cm, opacity=0.5] (0,0) grid (2,1.5);

\fill[blue!50!white, opacity=0.5](1,1) rectangle (2,1.5);


\fill[green!50!white, opacity=0.5] (1,0.5) rectangle (1.5,1);

\draw[blue, line width=0.5mm] (0,1.5)--++(0.5,0)--++(0,-1);

\draw[blue, line width=0.5mm] (1.5,1.5)--++(-0.5,0)--++(0,-1);

\draw[blue, line width=0.5mm] (1,0.5)--++(0,-0.5);

\draw[green!50!black, dotted, line width=0.5mm] (1,0.5)--++(0.5,0);

\draw[black, line width=0.15mm] (1.45,1.2)--++(0.1,0);
\draw[black, line width=0.15mm] (1.45,1.25)--++(0.1,0);
\draw[black, line width=0.15mm] (1.45,1.3)--++(0.1,0);


\draw[fill=blue] (1,0.5) circle [radius=0.05];

\node[left] at (0,0.5) [scale=1]
{\tiny{$\ell$}};
\node[left] at (0,1) [scale=1]
{\tiny{+1}};
\node[left] at (0,1.5) [scale=1]
{\tiny{+2}};

\node[above] at (1, 1.5) [scale=1]
{\tiny{$k$}};
\node[above] at (1.5, 1.5) [scale=1]
{\tiny{+1}};
\node[above] at (2, 1.5) [scale=1]
{\tiny{+2}};

\node at (1.25,1.25) [scale=0.8] {\small{X}};
\node at (0.25,1.25) [scale=0.8] {\small{Y}};

\node[below] at (1,0) [scale=0.8]{\small{\begin{tabular}{c} Fig. 4.15 (b). Case 1.2. 
 \end{tabular}}};;

} \end{scope}

\end{tikzpicture}
\end{adjustbox}
\end{wrapfigure}

\noindent \textit{CASE 1: $e(k;l-1,l) \in H$.} Then exactly one of $X+(0,-1)$ and $X+(1,0)$ must belong to $P$.

\null

\noindent \textit{CASE 1.1: $X+(0,-1)\in P$.} Note that if $e(k+1;l,l+1) \in H$, then $X+(0,-1)\in P$ is switchable, contradicting the definition of a fat path, so we only need to check the case where $e(k+1;l,l+1) \notin H$. Then $S_{\uparrow}(k+1,l-1;k+2,l) \in H$ and $S_{\downarrow}(k+1,l+2;k+2,l+1) \in H$. Now exactly one of $X+(0,-2)$ and $X+(1,-1)$ belong to $P$. Suppose that $X+(0,-2) \in P$ (Figure 4.15 (a)). Then $e(k,k+1;l-1) \in H$ or $e(k,k+1;l-1) \notin H$. If the former, then $X+(-2,0)$ must be an end-box of $P$, but this contradicts the fact that the other end-box of $P$ is $Y=X+(-2,0)$; and if the latter then $X+(0,-2)$ is switchable, contradicting that $F$ is a looping fat path. The case where $X+(1,-1)$ belongs to $P$ is very similar so we omit the proof. End of Case 1.1

\endgroup 

\null 

\noindent \textit{CASE 1.2: $X+(1,0)\in P$.} Then $X+(0,-1) \in N[P] \setminus P$ and $X+(0,-2) \notin N[P]$. It follows that $e(k,k+1;l)$ is a shadow edge of $F$. But then $\text{deg}_{(h(F))}(v(k,l)) =3$ contradicting Lemma 4.3.  See Figure 4.15 (b). End of Case 1.2.

\null 

\noindent \textit{CASE 2: $e(k;l-1,l) \in hb(F)$.} Let $Z=R(k-1,l-1)$. Then exactly one of $Z$ and $Z+(1,0)$ belongs to $N[P] \setminus P$.

\begingroup 
\setlength{\intextsep}{0pt}
\setlength{\columnsep}{20pt}
\begin{wrapfigure}[]{l}{0cm}
\begin{adjustbox}{trim=0cm 0.5cm 0cm 0cm}
\begin{tikzpicture}[scale=1.5]

\begin{scope}[xshift=0cm]{
\draw[gray,very thin, step=0.5cm, opacity=0.5] (0,0) grid (1.5,1.5);

\fill[green!50!white, opacity=0.5](0.5,0) rectangle (1,0.5);

\fill[blue!50!white, opacity=0.5](1,1) rectangle (1.5,1.5);

\draw[green!50!black, dotted, line width=0.5mm] (1,0)--++(0,0.5);

{
\draw[black, line width=0.15mm] (0.95,0.2)--++(0.1,0);
\draw[black, line width=0.15mm] (0.95,0.25)--++(0.1,0);
\draw[black, line width=0.15mm] (0.95,0.3)--++(0.1,0);

\draw[black, line width=0.15mm] (0.7,0.45)--++(0,0.1);
\draw[black, line width=0.15mm] (0.75,0.45)--++(0,0.1);
\draw[black, line width=0.15mm] (0.8,0.45)--++(0,0.1);
}

\draw[blue, line width=0.5mm] (0,1.5)--++(0.5,0)--++(0,-1);
\draw[blue, line width=0.5mm] (1.5,1.5)--++(-0.5,0)--++(0,-1);

\node[left] at (0,0.5) [scale=1]
{\tiny{$\ell$}};
\node[left] at (0,1) [scale=1]
{\tiny{+1}};

\node[above] at (0.5, 1.5) [scale=1]
{\tiny{-1}};
\node[above] at (1, 1.5) [scale=1]
{\tiny{$k$}};

\node at (1.25,1.25) [scale=0.8] {\small{X}};
\node at (0.25,1.25) [scale=0.8] {\small{Y}};
\node at (0.75,0.25) [scale=0.8] {\small{Z}};

\node[below] at (0.75,0) [scale=0.8]{\small{\begin{tabular}{c} Fig. 4.16 (a). Case 2. \\ $ Z \in N[P] \setminus P$.
 \end{tabular}}};;

} \end{scope}

\begin{scope}[xshift=2.25cm]{
\draw[gray,very thin, step=0.5cm, opacity=0.5] (0,0) grid (1.5,1.5);

\fill[green!50!white, opacity=0.5](1,0) rectangle (1.5,0.5);

\fill[blue!50!white, opacity=0.5](1,1) rectangle (1.5,1.5);

\draw[green!50!black, dotted, line width=0.5mm] (1,0)--++(0,0.5);

{
\draw[black, line width=0.15mm] (0.95,0.2)--++(0.1,0);
\draw[black, line width=0.15mm] (0.95,0.25)--++(0.1,0);
\draw[black, line width=0.15mm] (0.95,0.3)--++(0.1,0);

\draw[black, line width=0.15mm] (1.2,0.45)--++(0,0.1);
\draw[black, line width=0.15mm] (1.25,0.45)--++(0,0.1);
\draw[black, line width=0.15mm] (1.3,0.45)--++(0,0.1);
}

\draw[blue, line width=0.5mm] (0,1.5)--++(0.5,0)--++(0,-1);
\draw[blue, line width=0.5mm] (1.5,1.5)--++(-0.5,0)--++(0,-1);

\node[left] at (0,0.5) [scale=1]
{\tiny{$\ell$}};
\node[left] at (0,1) [scale=1]
{\tiny{+1}};

\node[above] at (0.5, 1.5) [scale=1]
{\tiny{-1}};
\node[above] at (1, 1.5) [scale=1]
{\tiny{$k$}};
\node[above] at (1.5, 1.5) [scale=1]
{\tiny{1}};

\node at (1.25,1.25) [scale=0.8] {\small{X}};
\node at (0.25,1.25) [scale=0.8] {\small{Y}};
\node at (1.25,0.25) [scale=0.8] {\small{Z}};

\node[below] at (0.75,0) [scale=0.8]{\small{\begin{tabular}{c} Fig. 4.16 (b). Case 2. \\ $ Z+(1,0) \in N[P] \setminus P$.
 \end{tabular}}};;

} \end{scope}

\end{tikzpicture}
\end{adjustbox}
\end{wrapfigure}

If $Z \in N[P] \setminus P$ then Lemma 4.3 implies that $e(k-1,k;l) \notin E(h(F))$, but  this contradicts Lemma 4.1 (d) (Figure 4.16 (a)); and if $Z+(1,0) \in N[P] \setminus P$, then Lemma 4.3 implies that $e(k,k+1,l)\notin H$, but  this also contradicts to Lemma 4.1 (d) (Figure 4.16 (b)). End of Case 2. $\square$

\null

\noindent \textbf{Lemma 4.6.} The shadow of $F$ does not have a pair of consecutive colinear edges in the case where neither edge of the pair is a left or right colinear edge of the $A_1$-type that follows $L$.

\null 

\noindent \textit{Proof.}  It is a fact that none of the boxes of $P$ discussed in this lemma can be end-boxes of $P$. The justifications are straightforward but distracting so we will omit them and use this fact repeatedly and implicitly throughout the proof. 

Assume for contradiction that there is a pair of consecutive collinear edges in $h(F)$ where neither edge of the pair is a left or right colinear edge of the $A_1$-type that follows $L$. For definiteness, we may assume that these edges are the horizontal edges $e(k'-1,k';l')$ and $e(k',k'+1;l')$. Let $R(k'-1,l')=Z$. Corollary 1.9 implies that exactly one of $Z$ and $Z+(0,-1)$ is in $F$. For definiteness, assume that $Z \in F$. Then, by Lemma 4.3, $e(k';l',l'+1)$ is neither in $H$ nor in $hb(F)$. Now, either both $e(k'-1,k';l')$ and $e(k',k'+1;l')$ belong to $H$, or at least one of them belongs to $hb(F)$.

\null

\noindent \textit{CASE 1: Both $e(k'-1,k';l')$ and $e(k',k'+1;l')$ belong to $H$.} By Lemma 4.3, $e(k';l',l'+1)$ is not a shadow edge of $F$. Then, Lemma 4.1 (a) implies that either exactly one of $Z$ and $Z+(1,0)$ belongs to $N[P] \setminus P$ or neither does.

\null 

\begingroup 
\setlength{\intextsep}{0pt}
\setlength{\columnsep}{20pt}
\begin{wrapfigure}[]{r}{0cm}
\begin{adjustbox}{trim=0cm 0cm 0cm 0.5cm}
\begin{tikzpicture}[scale=1.5]
\begin{scope}[xshift=0cm]{
\draw[gray,very thin, step=0.5cm, opacity=0.5] (0,0) grid (1,1);

\fill[blue!50!white, opacity=0.5](0,0.5) rectangle (0.5,1);

\draw[blue, line width=0.5mm] (0,0.5)--++(1,0);
\draw[blue, line width=0.5mm] (0,1)--++(0.5,0);

\node[right] at (1,0.5) [scale=1]
{\tiny{$\ell'$}};
\node[right] at (1,1) [scale=1]
{\tiny{+1}};

\node[above] at (0, 1) [scale=1]
{\tiny{$-1$}};
\node[above] at (0.5, 1) [scale=1]
{\tiny{$k'$}};

\node at (0.25,0.75) [scale=0.8] {\small{Z}};

\node[below] at (0.5,-0.1) [scale=0.8]{\small{\begin{tabular}{c} Fig. 4.17 (a).  \\  Case 1.1.
 \end{tabular}}};;

} \end{scope}

\begin{scope}[xshift=1.75cm]{
\draw[gray,very thin, step=0.5cm, opacity=0.5] (0,0) grid (1,1);

\fill[blue!50!white, opacity=0.5](0,0.5) rectangle (0.5,1);

\fill[green!50!white, opacity=0.5](0.5,0.5) rectangle (1,1);

\draw[blue, line width=0.5mm] (0,0.5)--++(1,0);

\draw[black, line width=0.15mm] (0.7,0.45)--++(0,0.1);
\draw[black, line width=0.15mm] (0.75,0.45)--++(0,0.1);
\draw[black, line width=0.15mm] (0.8,0.45)--++(0,0.1);

\node[right] at (1,0.5) [scale=1]
{\tiny{$\ell'$}};
\node[right] at (1,1) [scale=1]
{\tiny{+1}};

\node[above] at (0, 1) [scale=1]
{\tiny{$-1$}};
\node[above] at (0.5, 1) [scale=1]
{\tiny{$k'$}};

\node at (0.25,0.75) [scale=0.8] {\small{Z}};

\node[below] at (0.5,-0.1) [scale=0.8]{\small{\begin{tabular}{c} Fig. 4.17 (b).   \\ Case 1.2.
 \end{tabular}}};;

} \end{scope}
\end{tikzpicture}
\end{adjustbox}
\end{wrapfigure}

\noindent  \textit{CASE 1.1: Neither $Z$ nor $Z+(1,0)$ belongs to $N[P] \setminus P$.} This implies that $Z$ and $Z+(1,0)$ are in $P$ (Figure 4.17 (a)). Now, at least one of $e(k'-1, k'; l'+1)$ and $e(k', k'+1; l'+1)$ belongs to $H$. For definiteness assume that $e(k'-1, k'; l'+1) \in H$. But then $Z$ is a switchable box in $P$, contradicting the definition of a fat path. End of Case 1.1

\null 

\noindent  \textit{CASE 1.2: Exactly one of $Z$ and $Z+(1,0)$ belongs to $N[P] \setminus P$}. For definiteness, assume that $Z+(1,0) \in N[P] \setminus P$. Then $Z \in P$. But now, the fact that $Z$ is not an end-box and Lemma 4.1 (b) imply that $e(k',k'+1;l') \notin H$, contradicting the assumption of Case 1. End of Case 1.2. End of Case 1. 

\endgroup 

\null 

\noindent \textit{CASE 2: At least one of $e(k'-1,k';l')$ and $e(k',k'+1;l')$ belongs to $hb(F)$.} Either $e(k'-1,k';l') \in hb(F)$, or $e(k'-1,k';l') \notin hb(F)$.

\null

\begingroup 
\setlength{\intextsep}{0pt}
\setlength{\columnsep}{10pt}
\begin{wrapfigure}[]{l}{0cm}
\begin{adjustbox}{trim=0cm 0cm 0cm 0.5cm}
\begin{tikzpicture}[scale=1.5]

\begin{scope}[xshift=0cm]{
\draw[gray,very thin, step=0.5cm, opacity=0.5] (0,0) grid (1,1);

\fill[green!50!white, opacity=0.5](0,0.5) rectangle (0.5,1);

\draw[blue, line width=0.5mm] (0.5,0.5)--++(0.5,0);
\draw[green!50!black, line width=0.5mm] (0,0.5)--++(0.5,0);

{
\draw[black, line width=0.15mm] (0.45,0.7)--++(0.1,0);
\draw[black, line width=0.15mm] (0.45,0.75)--++(0.1,0);
\draw[black, line width=0.15mm] (0.45,0.8)--++(0.1,0);

\draw[black, line width=0.15mm] (0.2,0.45)--++(0,0.1);
\draw[black, line width=0.15mm] (0.25,0.45)--++(0,0.1);
\draw[black, line width=0.15mm] (0.3,0.45)--++(0,0.1);
}

\node[right] at (1,0.5) [scale=1]
{\tiny{$\ell'$}};
\node[right] at (1,1) [scale=1]
{\tiny{+1}};

\node[above] at (0, 1) [scale=1]
{\tiny{$-1$}};
\node[above] at (0.5, 1) [scale=1]
{\tiny{$k'$}};
;

\node at (0.25,0.75) [scale=0.8] {\small{Z}};

\node[below] at (0.5,0) [scale=0.8]{\small{\begin{tabular}{c} Fig. 4.18 (a)  Case 2.1
 \end{tabular}}};;

} \end{scope}

\begin{scope}[xshift=2cm]{
\draw[gray,very thin, step=0.5cm, opacity=0.5] (0,0) grid (1,1);

\fill[blue!50!white, opacity=0.5](0,0.5) rectangle (0.5,1);
\fill[green!50!white, opacity=0.5](0.5,0.5) rectangle (1,1);

\draw[green!50!black, line width=0.5mm] (0.5,0.5)--++(0.5,0);
\draw[blue, line width=0.5mm] (0,0.5)--++(0.5,0);

{
\draw[black, line width=0.15mm] (0.45,0.7)--++(0.1,0);
\draw[black, line width=0.15mm] (0.45,0.75)--++(0.1,0);
\draw[black, line width=0.15mm] (0.45,0.8)--++(0.1,0);

\draw[black, line width=0.15mm] (0.7,0.45)--++(0,0.1);
\draw[black, line width=0.15mm] (0.75,0.45)--++(0,0.1);
\draw[black, line width=0.15mm] (0.8,0.45)--++(0,0.1);
}

\node[right] at (1,0.5) [scale=1]
{\tiny{$\ell'$}};
\node[right] at (1,1) [scale=1]
{\tiny{+1}};

\node[above] at (0, 1) [scale=1]
{\tiny{$-1$}};
\node[above] at (0.5, 1) [scale=1]
{\tiny{$k'$}};
;

\node at (0.25,0.75) [scale=0.8] {\small{Z}};

\node[below] at (0.5,0) [scale=0.8]{\small{\begin{tabular}{c} Fig. 4.18 (b).  \\ Case 2.2(a).
 \end{tabular}}};;

} \end{scope}

\begin{scope}[xshift=4cm]{
\draw[gray,very thin, step=0.5cm, opacity=0.5] (0,0) grid (1,1);

\fill[green!50!white, opacity=0.5](0,0.5) rectangle (0.5,1);

\draw[green!50!black, line width=0.5mm] (0.5,0.5)--++(0.5,0);
\draw[green!50!black, line width=0.5mm] (0.5,0.5)--++(0,0.5);
\draw[blue, line width=0.5mm] (0,0.5)--++(0.5,0);

{
\draw[black, line width=0.15mm] (0.45,0.7)--++(0.1,0);
\draw[black, line width=0.15mm] (0.45,0.75)--++(0.1,0);
\draw[black, line width=0.15mm] (0.45,0.8)--++(0.1,0);

\draw[black, line width=0.15mm] (0.7,0.45)--++(0,0.1);
\draw[black, line width=0.15mm] (0.75,0.45)--++(0,0.1);
\draw[black, line width=0.15mm] (0.8,0.45)--++(0,0.1);
}

\node[right] at (1,0.5) [scale=1]
{\tiny{$\ell'$}};
\node[right] at (1,1) [scale=1]
{\tiny{+1}};

\node[above] at (0, 1) [scale=1]
{\tiny{$-1$}};
\node[above] at (0.5, 1) [scale=1]
{\tiny{$k'$}};
;

\node at (0.25,0.75) [scale=0.8] {\small{Z}};

\node[below] at (0.5,0) [scale=0.8]{\small{\begin{tabular}{c} Fig. 4.18 (c). \\  Case 2.2(b).
 \end{tabular}}};;

} \end{scope}

\end{tikzpicture}
\end{adjustbox}
\end{wrapfigure}

\noindent \textit{CASE 2.1: $e(k'-1,k';l') \in hb(F)$.} It follows that $Z \in N[P] \setminus P$. But then $Z$ can be neither switchable nor a leaf, contradicting Lemma 4.1(d). See Figure 4.18 (a). End of Case 2.1.

\null

\noindent \textit{CASE 2.2: $e(k'-1,k';l') \notin hb(F)$.} Then we must have that $e(k'-1,k';l') \in H$ and $e(k',k'+1;l') \in hb(F)$. Now, either $Z+(1,0) \in F$ or $Z+(1,0) \notin F$. 

\null 

\noindent \textit{CASE 2.2(a):  $Z+(1,0) \in F$.} Then, by Observation 4.2, $Z+(1,0) \in N[P] \setminus P$. But $Z+(1,0)$ can be neither switchable nor a leaf, contradicting Lemma 4.1(d). See Figure 4.18 (b). End of Case 2.2(a).

\null 

\noindent \textit{CASE 2.2(b):  $Z+(1,0) \notin F$.} Then, $Z$ must belong to $N[P] \setminus P$. This means that $e(k';l',l'+1) \in hb(F)$, which, as we, as we noted in the second paragraph of the proof, is not possible. See Figure 4.18 (c).  End of Case 2.2(b). $\square$

\endgroup 

\null

\noindent This completes the proof of Proposition 4.4. An immediate consequence of it is that the $A_1$-type following $L$ is the only $A_1$-type in $F$, so we can refer to it as \textit{the} $A_1$-type of $F$. 

\endgroup

\subsection{Turns}

\noindent In this section we show that every looping fat path $F$ must have a turn. We do this by showing that the shadow $h(F)$ of a looping fat path $F$ must have a (necessarily closed) turn and note that this would immediately imply that $F$ must have a turn. 

In Lemma 4.7, Lemma 4.8, and Corollaries 4.9 (a) and (b) below, we will often use the definition of the shadow of southern looping fat path, Proposition 4.4, and the fact that the $A_1$-type of $F$ is unique in $F$, and write (DsFP) whenever we appeal to them.

\null

\begingroup
\setlength{\intextsep}{0pt}
\setlength{\columnsep}{20pt}
\begin{wrapfigure}[]{l}{0cm}
\begin{adjustbox}{trim=0cm 0cm 0cm 0cm}

\begin{tikzpicture}[scale=1.25]

\begin{scope}[xshift=0cm, yshift=0cm]
{
\draw[gray,very thin, step=0.5cm, opacity=0.5] (0,0) grid (2.5,2);

\fill[blue!50!white, opacity=0.5] (0.5,0.5) rectangle (1,1.5);
\fill[blue!50!white, opacity=0.5] (1.5,0.5) rectangle (2,1.5);

\draw[blue, line width=0.5mm] (0.5,1.5)--++(0.5,0)--++(0,-1);
\draw[blue, line width=0.5mm] (1.5,1.5)--++(0,-1);

\draw[->, blue,line width=0.5mm] (1.5,1.5)--++(0.5,0);

\node[above] at (0.5,2) [scale=1]{\tiny{$-2$}};
\node[above] at (1,2) [scale=1]{\tiny{$-1$}};
\node[above] at (1.5,2) [scale=1]{\tiny{$k'$}};
\node[above] at (2,2) [scale=1]{\tiny{$+1$}};

\node[left] at (0,1.5) [scale=1]{\tiny{$\ell'$}};
\node[left] at (0,1) [scale=1]{\tiny{$-1$}};
\node[left] at (0,0.5) [scale=1]{\tiny{$-2$}};

\node at (0.75,1.25) [scale=0.8]{Y};
\node at (1.75,1.25) [scale=0.8]{X};

\node[below] at (1.25,0) [scale=0.75]{\small{Fig. 4.19.}};

}
\end{scope}

\end{tikzpicture}
\end{adjustbox}
\end{wrapfigure}

\noindent \textbf{Lemma 4.7.} Let $H$ be a Hamiltonian cycle of an $m \times n$ grid graph $G$ and let $h(F)$ be the shadow of a looping fat path of $G$. Then $h(F)$ has at least one turn $T_1$ such that both leaves of $T_1$ belong to $F$.

\null

\noindent \textit{Proof.} For definiteness assume that $F=G\langle N[P(X,Y)] \rangle $ is a southern looping fat path following $R(k'-1,l'+1)$ southward, with $X=R(k',l'-1)$ and $Y=R(k'-2,l'-1)$. By Lemma 4.3, $E(h(F))$ is a cycle. Orient  $E(h(F))$ as a directed trail $\overrightarrow{K}$ so that the first edge of $\overrightarrow{K}$ 
is  $(v(k',l'), v(k'+1,l'))$. With this orientation we can give a direction  - $N,S,E$ or $W$ - to edges in $\overrightarrow{K}$, defined as the position of the head of an edge relative to its tail. See Figure 4.19.

\endgroup

Our choice of direction for the first edge and Lemma 1.6 imply that $\text{Boxes}(\Phi(\overrightarrow{K}, \text{right})) \subset \text{Boxes}(F)$ so the boxes of $h(F)$ are on the right side of the oriented edges of $K$. We call this fact (RSK) for reference. We sweep the edges of $K$ in the direction of the orientation starting at $v(k',l')$. We observe that we must encounter at least one west edge $e_W$, since $Y$ is west of $X$. Let $e_W=e_0=e(k-1, k;l)$ be the first west edge encountered, let $e_1$ be the edge preceding $e_0$ in the sweep, let $e_j$ be the edge preceding the edge $e_{j-1}$ in the sweep and let $(v(k',l'), v(k'+1,l'))=e_s$.

\null

\noindent In this proof we will use the fact that $e_W$ is the first west edge encountered (1stW) several times. 

\null 

\begingroup
\setlength{\intextsep}{0pt}
\setlength{\columnsep}{10pt}
\begin{wrapfigure}[]{r}{0cm}
\begin{adjustbox}{trim=0cm 0.5cm 0cm 0cm}

\begin{tikzpicture}[scale=1.5]

\begin{scope}[xshift=0cm, yshift=0 cm]
{

\draw[gray,very thin, step=0.5cm, opacity=0.5] (0,0) grid (1,1.5);

\fill[blue!40!white, opacity=0.5] (0,0.5) rectangle (0.5,1.5);
\fill[blue!40!white, opacity=0.5] (0,0) rectangle (1,0.5);

\begin{scope}
[very thick,decoration={
    markings,
    mark=at position 0.6 with {\arrow{>}}}
    ]
    
    \draw[postaction={decorate}, blue, line width=0.5mm] (1,0)--++(-0.5,0);
    \draw[postaction={decorate}, blue, line width=0.5mm] (1,0.5)--++(0,-0.5);
    \draw[postaction={decorate}, blue, line width=0.5mm] (0.5,0.5)--++(0.5,0);
    
    \draw[postaction={decorate}, blue, line width=0.5mm] (0.5,1.5)--++(0,-0.5);
    \draw[postaction={decorate}, blue, line width=0.5mm] (0.5,1)--++(0,-0.5);

\end{scope}

\draw[blue, line width=0.5mm] (0,0.5)--++(0,1)--++(0.5,0);

\node[left] at (0,0) [scale=1]{\tiny{$\ell$}};
\node[above] at (1,1.5) [scale=1]{\tiny{$k$}};

\node at (0.65,0.15) [scale=1]{\small{$e_{_0}$}};
\node at (0.65,1.15) [scale=1]{\small{$e_{_4}$}};

\node at (0.25,1.25) [scale=1]{\small{X}};

}

\node[below] at (0.5,0) [scale=0.8]{\small{\begin{tabular}{c} Fig. 4.20 (a). Case  \\ 1.1.  $\Phi(e_4, \text{right})=X$.
 \end{tabular}}};;

\end{scope}

\begin{scope}[xshift=2cm, yshift=0 cm]
{

\draw[gray,very thin, step=0.5cm, opacity=0.5] (0,0) grid (1.5,1.5);

\fill[blue!40!white, opacity=0.5] (0,0.5) rectangle (0.5,1.5);
\fill[blue!40!white, opacity=0.5] (0,0) rectangle (1,0.5);

\begin{scope}
[very thick,decoration={
    markings,
    mark=at position 0.6 with {\arrow{>}}}
    ]
    
    \draw[postaction={decorate}, blue, line width=0.5mm] (1,0)--++(-0.5,0);
    \draw[postaction={decorate}, blue, line width=0.5mm] (1,0.5)--++(0,-0.5);
    \draw[postaction={decorate}, blue, line width=0.5mm] (0.5,0.5)--++(0.5,0);
    
    \draw[postaction={decorate}, blue, line width=0.5mm] (0.5,1.5)--++(0,-0.5);
    \draw[postaction={decorate}, blue, line width=0.5mm] (0.5,1)--++(0,-0.5);

\end{scope}

\draw[blue, line width=0.5mm] (0,0.5)--++(0,1)--++(0.5,0);

\draw[blue, line width=0.5mm] (1,0.5)--++(0,1)--++(0.5,0);

\draw[fill=blue] (1,0.5) circle [radius=0.05];

\node[left] at (0,0) [scale=1]{\tiny{$\ell$}};
\node[above] at (1,1.5) [scale=1]{\tiny{$k$}};

\node at (0.65,0.15) [scale=1]{\small{$e_{_0}$}};
\node at (0.65,1.15) [scale=1]{\small{$e_{_4}$}};

\node at (0.25,1.25) [scale=1]{\small{Y}};

\node at (1.25,1.25) [scale=1]{\small{X}};

}

\node[below] at (0.75,0) [scale=0.8]{\small{\begin{tabular}{c} Fig. 4.20 (b). Case  \\ 1.1.  $\Phi(e_4, \text{right})=Y$.
 \end{tabular}}};;

\end{scope}

\end{tikzpicture}
\end{adjustbox}
\end{wrapfigure}

\noindent By (DsFP),  $e_0$ was immediately preceded by a south edge or a north edge.

\null 

\noindent \textit{CASE 1: $e_1$ is southern.} We shall find a northeastern turn. By (DsFP) and (1stW), the preceding edge $e_2$ must be eastern; By (DsFP) and (1stW), $e_3$  has to be southern. By (1stW) $e_4$ cannot be western. Then $e_4$ is southern or $e_4$ is eastern.

\null 

\noindent \textit{CASE 1.1: $e_4$ is southern.} Then, by (DsFP) and (RSK), we have that $\Phi(e_4, \text{right})=X$ or $\Phi(e_4, \text{right})=Y$. But the former contradicts Proposition 4.4, and the latter implies that $\text{deg}_{h(F)}(v(k,l+1)=3$, contradicting Lemma 4.3. See Figure 4.20. Thus, $e_4$ must be eastern. End of Case 1.1. 

\null

\endgroup 

\noindent \textit{CASE 1.2: $e_4$ is eastern.} By (DsFP), $e_5$ is not eastern. Then $e_5$ is northern or $e_5$ is southern. 

\null

\noindent \textit{CASE 1.2(a): $e_5$ is northern.} Then there is a northeastern turn $T_1$ on the edges $e_0, ..., e_5$ with both leaves contained in $F$. See Figure 4.21 (a). End of Case 1.2 (a).

\null

\begingroup
\setlength{\intextsep}{0pt}
\setlength{\columnsep}{10pt}
\begin{wrapfigure}[]{l}{0cm}
\begin{adjustbox}{trim=0cm 0.25cm 0cm 0.5cm}

\begin{tikzpicture}[scale=1.5]

\begin{scope}[xshift=0cm, yshift=0 cm]
{

\draw[gray,very thin, step=0.5cm, opacity=0.5] (0,0) grid (1,1.5);

\fill[blue!40!white, opacity=0.5] (0,0.5) rectangle (0.5,1);
\fill[blue!40!white, opacity=0.5] (0,0) rectangle (1,0.5);

\begin{scope}
[very thick,decoration={
    markings,
    mark=at position 0.6 with {\arrow{>}}}
    ]
    
    \draw[postaction={decorate}, blue, line width=0.5mm] (1,0)--++(-0.5,0);
    \draw[postaction={decorate}, blue, line width=0.5mm] (1,0.5)--++(0,-0.5);
    \draw[postaction={decorate}, blue, line width=0.5mm] (0.5,0.5)--++(0.5,0);

    \draw[postaction={decorate}, blue, line width=0.5mm] (0,0.5)--++(0,0.5);
    \draw[postaction={decorate}, blue, line width=0.5mm] (0,1)--++(0.5,0);
    \draw[postaction={decorate}, blue, line width=0.5mm] (0.5,1)--++(0,-0.5);

\end{scope}

\node[left] at (0,0) [scale=1]{\tiny{$\ell$}};
\node[above] at (1,1.5) [scale=1]{\tiny{$k$}};

\node at (0.65,0.15) [scale=1]{\small{$e_{_0}$}};
\node at (0.15, 0.6) [scale=1]{\small{$e_{_5}$}};

}

\node[below] at (0.5,0) [scale=0.8]{\small{\begin{tabular}{c} Fig. 4.21(a). 
 \end{tabular}}};;

\end{scope}

\begin{scope}[xshift=1.75 cm]
{

\draw[gray,very thin, step=0.5cm, opacity=0.5] (0,0) grid (1.5,2);

\fill[blue!40!white, opacity=0.5] (0,0.5) rectangle (0.5,1.5);
\fill[blue!40!white, opacity=0.5] (0.5,0) rectangle (1,1);
\fill[blue!40!white, opacity=0.5] (1,0) rectangle (1.5,0.5);

\begin{scope}
[very thick,decoration={
    markings,
    mark=at position 0.6 with {\arrow{>}}}
    ]

    \draw[postaction={decorate}, blue, line width=0.5mm] (1.5,0)--++(-0.5,0);
    \draw[postaction={decorate}, blue, line width=0.5mm] (1.5,0.5)--++(0,-0.5);
    \draw[postaction={decorate}, blue, line width=0.5mm] (1,0.5)--++(0.5,0);
    
    \draw[postaction={decorate}, blue, line width=0.5mm] (0,1.5)--++(0.5,0);
    \draw[postaction={decorate}, blue, line width=0.5mm] (0.5,1.5)--++(0,-0.5);

\end{scope}

\draw[blue,line width=0.5mm] (0,0.5)--++(0,1);

{
\tikzset
  {
    myCircle/.style=
    {
      blue,
      path fading=fade out,
    }
  }

\foreach \x in {0,...,2}
\fill[myCircle,] (0.625+0.125*\x,0.875-0.125*\x) circle (0.075);

}

\node[right] at (1.5,0) [scale=1]{\tiny{$\ell$}};
\node[above] at (1.5,2) [scale=1]{\tiny{$k$}};
\node[above] at (0,2) [scale=1]{\tiny{$k'$}};
\node[right] at (1.5,1.5) [scale=1]{\tiny{$\ell'$}};

\node[right] at (0.5, 1.25) [scale=1]{\small{$e_{s-1}$}};

\node at (1.15,0.15) [scale=1]{\small{$e_{_0}$}};

\node at (0.25,1.25) [scale=1]{\small{X}};

}

\node[below] at (0.75,0) [scale=0.8]{\small{\begin{tabular}{c} Fig. 4.21(b). Case 1.2($b_1$). 
 \end{tabular}}};;

\end{scope}

\begin{scope}[xshift=4cm]
{

\draw[gray,very thin, step=0.5cm, opacity=0.5] (0,0) grid (1.5,2);

\fill[blue!40!white, opacity=0.5] (0,0.5) rectangle (0.5,1.5);
\fill[blue!40!white, opacity=0.5] (0.5,0) rectangle (1,1);
\fill[blue!40!white, opacity=0.5] (1,0) rectangle (1.5,0.5);

\begin{scope}
[very thick,decoration={
    markings,
    mark=at position 0.6 with {\arrow{>}}}
    ]
    
    \draw[postaction={decorate}, blue, line width=0.5mm] (1.5,0)--++(-0.5,0);
    \draw[postaction={decorate}, blue, line width=0.5mm] (1.5,0.5)--++(0,-0.5);
    \draw[postaction={decorate}, blue, line width=0.5mm] (1,0.5)--++(0.5,0);
    
    \draw[postaction={decorate}, blue, line width=0.5mm] (0,1)--++(0,0.5);
    \draw[postaction={decorate}, blue, line width=0.5mm] (0,1.5)--++(0.5,0);
    \draw[postaction={decorate}, blue, line width=0.5mm] (0.5,1.5)--++(0,-0.5);

\end{scope}

{
\tikzset
  {
    myCircle/.style=
    {
      blue,
      path fading=fade out,
    }
  }

\foreach \x in {0,...,2}
\fill[myCircle,] (0.625+0.125*\x,0.875-0.125*\x) circle (0.075);

}

\node[right] at (1.5,0) [scale=1]{\tiny{$\ell$}};
\node[above] at (1.5,2) [scale=1]{\tiny{$k$}};

\node[right] at (-0.05, 1.25) [scale=1]{\small{$e_{j_{_0}}$}};

\node at (1.15,0.15) [scale=1]{\small{$e_{_0}$}};

}

\node[below] at (0.75,0) [scale=0.8]{\small{\begin{tabular}{c} Fig. 4.21(c). Case 1.2($b_2$). 
 \end{tabular}}};;

\end{scope}

\end{tikzpicture}
\end{adjustbox}
\end{wrapfigure}

\noindent \textit{CASE 1.2 (b): $e_5$ is southern.} Let $Q(j)$ be the statement: ``$e_j$ is southern and $e_{j+1}$ is eastern''. Now, either $Q(j)$ is true for each $j \in \{1,3, ..., s-1\}$ (Case 1.2($b_1$)), or there is some $j_0 \in \{ 5,9, ..., s-1\}$  such that $Q(j)$ for each odd $j<j_0$, but $Q(j_0)$ is not true (Case 1.2($b_2$)).

\null 

\noindent \textit{CASE 1.2($b_1$).} Then we have a northeastern turn $T_1$ on the edges $e_0, ..., e_s, e(k'; l'-1, l')$ with both leaves contained in $F$. See Figure 4.21 (b). End of Case 1.2($b_1$).

\endgroup 

\null 

\noindent \textit{CASE 1.2($b_2$).} By (DsFP), $e_{j_0}$ is not eastern. If $e_{j_0}$ is southern, then we run into the same contradiction as in Case 1.1; and if $e_{j_0}$ is northern then we have a northeastern turn $T_1$ on the edges $e_0, ..., e_{j_0}$ with both leaves contained in $F$ (Figure 4.21 (c)). End of Case  1.2($b_2$). End of Case 1.2 (b). End of Case 1.2.

\endgroup 

\null 

\begingroup
\setlength{\intextsep}{0pt}
\setlength{\columnsep}{20pt}
\begin{wrapfigure}[]{r}{0cm}
\begin{adjustbox}{trim=0cm 0cm 0cm 0.25cm}

\begin{tikzpicture}[scale=1.5]

\begin{scope}[xshift=0cm]
{

\draw[gray,very thin, step=0.5cm, opacity=0.5] (0,0) grid (2,2);

\fill[blue!40!white, opacity=0.5] (0.5,0) rectangle (1,0.5);
\fill[blue!40!white, opacity=0.5] (0.5,0.5) rectangle (1.5,1);
\fill[blue!40!white, opacity=0.5] (1,1) rectangle (2,1.5);
\fill[blue!40!white, opacity=0.5] (0.5,1.5) rectangle (1.5,2);

\begin{scope}
[very thick,decoration={
    markings,
    mark=at position 0.6 with {\arrow{>}}}
    ]

    \draw[postaction={decorate}, blue, line width=0.5mm] (0.5,0)--++(0,0.5);
    \draw[postaction={decorate}, blue, line width=0.5mm] (0.5,0.5)--++(0,0.5);
    \draw[postaction={decorate}, blue, line width=0.5mm] (0.5,1)--++(0.5,0);
    \draw[postaction={decorate}, blue, line width=0.5mm] (1,1)--++(0,0.5);
    \draw[postaction={decorate}, blue, line width=0.5mm] (1,1.5)--++(-0.5,0);

\end{scope}

\draw[blue,line width=0.5mm] (1.5,2)--++(0,-0.5)--++(0.5,0);
\draw[blue,line width=0.5mm] (1.5,0.5)--++(0,0.5)--++(0.5,0);

\draw[blue, line width=0.5mm] (2,1)--++(0,0.5);

\draw[blue,line width=0.5mm] (0,2)--++(0,-0.5)--++(0.5,0)--++(0,0.5);

\draw[fill=blue] (0.5,1.5) circle [radius=0.05];


{
\draw[black, line width=0.15mm] (1.2,0.95)--++(0,0.1);
\draw[black, line width=0.15mm] (1.25,0.95)--++(0,0.1);
\draw[black, line width=0.15mm] (1.3,0.95)--++(0,0.1);

\draw[black, line width=0.15mm] (1.2,1.45)--++(0,0.1);
\draw[black, line width=0.15mm] (1.25,1.45)--++(0,0.1);
\draw[black, line width=0.15mm] (1.3,1.45)--++(0,0.1);

\draw[black, line width=0.15mm] (1.45,1.2)--++(0.1,0);
\draw[black, line width=0.15mm] (1.45,1.25)--++(0.1,0);
\draw[black, line width=0.15mm] (1.45,1.3)--++(0.1,0);
}

\node[left] at (0,1.5) [scale=1]{\tiny{$\ell$}};

\node[left] at (0,1) [scale=1]{\tiny{-1}};
\node[left] at (0,0.5) [scale=1]{\tiny{-2}};

\node[above] at (1,2) [scale=1]{\tiny{$k$}};
\node[above] at (1.5,2) [scale=1]{\tiny{$+1$}};

\node at (0.25,1.75) [scale=1]{\small{$L$}};
\node[above] at (0.75,1.5) [scale=1]{\small{$e_{_0}$}};
\node[right] at (1,1.25) [scale=1]{\small{$e_{_1}$}};
\node[above] at (0.75,1) [scale=1]{\small{$e_{_2}$}};
\node[left] at (0.5,0.75) [scale=1]{\small{$e_{_3}$}};

\node[right] at (0.5,0.25) [scale=1]{\small{$e_{_4}$}};

\node at (0.75,0.75) [scale=0.8]{\small{X}};

\node[below] at (1,0) [scale=0.8]{\small{\begin{tabular}{c} Fig. 4.22(a). Case 2.
 \end{tabular}}};;

}
\end{scope}

\begin{scope}[xshift=2.75cm]
{
\draw[gray,very thin, step=0.5cm, opacity=0.5] (0,0) grid (2.5,2.5);

\fill[blue!40!white, opacity=0.5] (0.5,0) rectangle (1,1);
\fill[blue!40!white, opacity=0.5] (1,0.5) rectangle (1.5,1);
\fill[blue!40!white, opacity=0.5] (1,1) rectangle (2,1.5);
\fill[blue!40!white, opacity=0.5] (1.5,1.5) rectangle (2.5,2);
\fill[blue!40!white, opacity=0.5] (1,2) rectangle (2,2.5);

\begin{scope}
[very thick,decoration={
    markings,
    mark=at position 0.6 with {\arrow{>}}}
    ]

    \draw[postaction={decorate}, blue, line width=0.5mm] (1,1)--++(0,0.5);
    \draw[postaction={decorate}, blue, line width=0.5mm] (1,1.5)--++(0.5,0);
    \draw[postaction={decorate}, blue, line width=0.5mm] (1.5,1.5)--++(0,0.5);
    \draw[postaction={decorate}, blue, line width=0.5mm] (1.5,2)--++(-0.5,0);

    \draw[postaction={decorate}, blue, line width=0.5mm] (0.5,1)--++(0.5,0);
    \draw[postaction={decorate}, blue, line width=0.5mm] (0.5,1.5)--++(0,-0.5);

\end{scope}

\draw[blue,line width=0.5mm] (1.5,0.5)--++(0,0.5)--++(0.5,0);
\draw[blue,line width=0.5mm] (2,2.5)--++(0,-0.5)--++(0.5,0);
\draw[blue,line width=0.5mm] (2,1)--++(0,0.5)--++(0.5,0);

\draw[blue,line width=0.5mm] (1,0.5)--++(0.5,0);
\draw[blue, line width=0.5mm] (2.5,1.5)--++(0,0.5);

\draw[blue, line width=0.5mm] (0,0.5)--++(0.5,0)--++(0,-0.5);
\draw[blue, line width=0.5mm] (1,0)--++(0,0.5);


{
\draw[black, line width=0.15mm] (1.7,1.45)--++(0,0.1);
\draw[black, line width=0.15mm] (1.75,1.45)--++(0,0.1);
\draw[black, line width=0.15mm] (1.8,1.45)--++(0,0.1);

\draw[black, line width=0.15mm] (1.7,1.95)--++(0,0.1);
\draw[black, line width=0.15mm] (1.75,1.95)--++(0,0.1);
\draw[black, line width=0.15mm] (1.8,1.95)--++(0,0.1);

\draw[black, line width=0.15mm] (1.95,1.7)--++(0.1,0);
\draw[black, line width=0.15mm] (1.95,1.75)--++(0.1,0);
\draw[black, line width=0.15mm] (1.95,1.8)--++(0.1,0);

\draw[black, line width=0.15mm] (0.45,0.7)--++(0.1,0);
\draw[black, line width=0.15mm] (0.45,0.75)--++(0.1,0);
\draw[black, line width=0.15mm] (0.45,0.8)--++(0.1,0);

\draw[black, line width=0.15mm] (0.95,0.7)--++(0.1,0);
\draw[black, line width=0.15mm] (0.95,0.75)--++(0.1,0);
\draw[black, line width=0.15mm] (0.95,0.8)--++(0.1,0);

\draw[black, line width=0.15mm] (0.7,0.45)--++(0,0.1);
\draw[black, line width=0.15mm] (0.75,0.45)--++(0,0.1);
\draw[black, line width=0.15mm] (0.8,0.45)--++(0,0.1);
}


{
\node[right] at (2.5,2) [scale=1]{\tiny{$\ell$}};
\node[right] at (2.5,1.5) [scale=1]{\tiny{-1}};
\node[right] at (2.5,1) [scale=1]{\tiny{-2}};
\node[right] at (2.5,0.5) [scale=1]{\tiny{-3}};
\node[right] at (2.5,0) [scale=1]{\tiny{-4}};

\node[above] at (0,2.5) [scale=1]{\tiny{$-3$}};
\node[above] at (0.5,2.5) [scale=1]{\tiny{$-2$}};
\node[above] at (1,2.5) [scale=1]{\tiny{$-1$}};
\node[above] at (1.5,2.5) [scale=1]{\tiny{$k$}};
\node[above] at (2,2.5) [scale=1]{\tiny{$+1$}};
\node[above] at (2.5,2.5) [scale=1]{\tiny{$+2$}};

\node[below] at (0.7,1) [scale=1]{\small{$e_{_4}$}};
\node[left] at (0.5,1.15) [scale=1]{\small{$e_{_5}$}};

\node[above] at (1.25,2) [scale=1]{\small{$e_{_0}$}};
\node[right] at (1.5,1.75) [scale=1]{\small{$e_{_1}$}};
\node[above] at (1.25,1.5) [scale=1]{\small{$e_{_2}$}};
\node[left] at (1,1.25) [scale=1]{\small{$e_{_3}$}};
}

\node[below] at (1.25,0) [scale=0.8]{\small{\begin{tabular}{c} Fig. 4.22(b). Case 2.1 \end{tabular}}};;
}

\end{scope}

\end{tikzpicture}
\end{adjustbox}
\end{wrapfigure}

\noindent \textit{CASE 2: $e_1$ is northern.} We shall find a southeastern turn with both leaves contained in $F$. By (1stW) and (DsFP), $e_2$ is eastern. By (DsFP), $e_3$ is northern. Note that (RSK) and Lemma 4.1 (d) imply that $\Phi(e_1, \text{right}) \in P$ and $e_1 \in H$. By Lemma 4.3, $e(k,k+1;l-1) \notin H$ and $e(k,k+1;l) \notin H$. If $e(k+2;l-1,l) \in H$, then $\Phi(e_1, \text{right})$ is switchable and in $P$, contradicting the definition of a fat path, so we may assume that $e(k+2;l-1,l) \notin H$. Then we must have that $S_{\downarrow}(k+1,l+1;k+2,l) \in H$,  $S_{\uparrow}(k+1,l-2;k+2,l-1) \in H$ and that $e(k+2;l-1,l) \in hb(F)$.

By (1stW), $e_4$ is not western. If $e_4$ is northern then (DsFP) and (RSK) imply that $\Phi(e_3, \text{right}) =X$. But then $L=R(k-1,l)$ and then $\deg_H(v(k-1,l))=3$, contradicting $H$ is Hamiltonian. See Figure 4.22 (a). Then $e_4$ must be eastern. By (DsFP), $e_5$ is not eastern. Then $e_5$ is southern or northern.

\endgroup 

\null 

\noindent \textit{CASE 2.1: $e_5$ is southern.} By (RSK) and Lemma 4.1 (d), we have that $\Phi(e_4, \text{right}) \in P$, and that $e_4 \in H$.  By Lemma 4.3, $e(k-2;l-3,l-2) \notin H$ and $e(k-1;l-3,l-2) \notin H$. If $e(k-2,k-1;l-3) \in H$, then $\Phi(e_4, \text{right})$ is switchable and in $P$, contradicting that $F$ is a fat path, so we may assume that $e(k-2,k-1;l-3) \notin H$. It follows that  $S_{\rightarrow}(k-3,l-3;k-2,l-4) \in H$ and that $S_{\uparrow}(k-1,l-4;k+2,l-1) \in H$. Then there is a southeastern turn on $e(k-2;l-4,l-3)$, $S_{\uparrow}(k-1,l-4;k+2,l-1) \in H$, $e(k+1,k+2;l)$ with both leaves in $F$. See Figure 4.22 (b). End of Case 2.1

\endgroup
 
\null

\begingroup
\setlength{\intextsep}{0pt}
\setlength{\columnsep}{20pt}
\begin{wrapfigure}[]{l}{0cm}
\begin{adjustbox}{trim=0cm 0.5cm 0cm 0cm}
\begin{tikzpicture}[scale=1.35]

\begin{scope}[xshift=0cm]
{
\draw[gray,very thin, step=0.5cm, opacity=0.5] (0,0) grid (3,3.5);

\fill[blue!40!white, opacity=0.5] (0.5,0.5) rectangle (1,1.5);
\fill[blue!40!white, opacity=0.5] (1,1) rectangle (1.5,2);
\fill[blue!40!white, opacity=0.5] (1.5,1.5) rectangle (2,2);
\fill[blue!40!white, opacity=0.5] (1.5,2) rectangle (2.5,2.5);
\fill[blue!40!white, opacity=0.5] (2,2.5) rectangle (3,3);
\fill[blue!40!white, opacity=0.5] (1.5,3) rectangle (2.5,3.5);

\draw[blue,line width=0.5mm] (0.5,0.5)--++(0,1);

\begin{scope}
[very thick,decoration={
    markings,
    mark=at position 0.6 with {\arrow{>}}}
    ]

    \draw[postaction={decorate}, blue, line width=0.5mm] (2,2.5)--++(0,0.5);
    \draw[postaction={decorate}, blue, line width=0.5mm] (2,3)--++(-0.5,0);
    \draw[postaction={decorate}, blue, line width=0.5mm] (1.5,2.5)--++(0.5,0);
    \draw[postaction={decorate}, blue, line width=0.5mm] (1,1.5)--++(0,0.5);
    \draw[postaction={decorate}, blue, line width=0.5mm] (0.5,1.5)--++(0.5,0);



\end{scope}

\draw[blue,line width=0.5mm] (2.5,3.5)--++(0,-0.5)--++(0.5,0);

\draw[blue,line width=0.5mm] (1,0.5)--++(0,0.5)--++(0.5,0)--++(0,0.5);

\draw[blue,line width=0.5mm] (2,2)--++(0.5,0)--++(0,0.5)--++(0.5,0);

\draw[blue,line width=0.5mm] (0.5,0.5)--++(0.5,0);

\draw[blue, line width=0.5mm] (3,2.5)--++(0,0.5);

{
\draw[black, line width=0.15mm] (2.45,2.7)--++(0.1,0);
\draw[black, line width=0.15mm] (2.45,2.75)--++(0.1,0);
\draw[black, line width=0.15mm] (2.45,2.8)--++(0.1,0);

\draw[black, line width=0.15mm] (2.2,2.45)--++(0,0.1);
\draw[black, line width=0.15mm] (2.25,2.45)--++(0,0.1);
\draw[black, line width=0.15mm] (2.3,2.45)--++(0,0.1);

\draw[black, line width=0.15mm] (2.2,2.95)--++(0,0.1);
\draw[black, line width=0.15mm] (2.25,2.95)--++(0,0.1);
\draw[black, line width=0.15mm] (2.3,2.95)--++(0,0.1);

\draw[black, line width=0.15mm] (0.2,0.45)--++(0,0.1);
\draw[black, line width=0.15mm] (0.25,0.45)--++(0,0.1);
\draw[black, line width=0.15mm] (0.3,0.45)--++(0,0.1);

\draw[black, line width=0.15mm] (0.45,0.2)--++(0.1,0);
\draw[black, line width=0.15mm] (0.45,0.25)--++(0.1,0);
\draw[black, line width=0.15mm] (0.45,0.3)--++(0.1,0);
}

{
\tikzset
  {
    myCircle/.style=
    {
      blue,
      path fading=fade out,
    }
  }

\foreach \x in {0,...,2}
\fill[myCircle,] (1.375+0.125*\x,1.875+0.125*\x) circle (0.075);
}

\node[right] at (3,2.5) [scale=1]{\tiny{-1}};
\node[right] at (3,3) [scale=1]{\tiny{$\ell$}};

\node[above] at (2,3.5) [scale=1]{\tiny{$k$}};
\node[above] at (2.5,3.5) [scale=1]{\tiny{+1}};
\node[above] at (3,3.5) [scale=1]{\tiny{+2}};

\node[above] at (0,3.5) [scale=1]{\tiny{$-1$}};
\node[above] at (0.5,3.5) [scale=1]{\tiny{$k'$}};
\node[above] at (1,3.5) [scale=1]{\tiny{+1}};

\node[right] at (3,1.5) [scale=1]{\tiny{$\ell'$}};
\node[right] at (3,1) [scale=1]{\tiny{-1}};
\node[right] at (3,0.5) [scale=1]{\tiny{-2}};
\node[right] at (3,0) [scale=1]{\tiny{-3}};

\node[above] at (1.65,3) [scale=1]{\small{$e_{_0}$}};
\node[right] at (1,1.75) [scale=1]{\small{$e_{s-1}$}};
\node at (0.75,1.25) [scale=1]{\small{X}};

\node[below] at (1.5,0) [scale=0.8]{\small{\begin{tabular}{c} Fig. 4.23. Case 2.2(a). \end{tabular}}};;
}

\end{scope}

\end{tikzpicture}
\end{adjustbox}
\end{wrapfigure}

\noindent \textit{CASE 2.2: $e_5$ is northern.} Let $Q(j)$ be the statement: ``$e_j$ is northern and $e_{j+1}$ is eastern''. Now, either $Q(j)$ is true for each $j \in \{1,3, ..., s-1\}$ (Case 2.2(a)), or there is some $j_0 \in \{ 5,7, ..., s-1\}$  such that $Q(j)$ for each odd $j<j_0$, but $Q(j_0)$ is not true (Case 2.2(b)).

\null

\noindent \textit{CASE 2.2(a).} Then (DsFP) and (RSK) imply that $\Phi(e_s, \text{right})=X$. This means that $e(k';l'-2,l'-1) \in H$, $e(k';l'-1,l') \in H$, and that $\Phi(e_{s-1}, \text{right}) \in h(F)$. By Proposition 4.4, $e(k';l'-3,l'-2) \notin H$. Note that if $e(k'-1,k';l'-2) \in H$, then $P(X,Y)$ is the $H$-path $X, X+(0,-1), X+(0,-2), X+(-1,-2), X+(-2,-2), X+(-2,-1), Y$. This contradicts our finding that $\Phi(e_{s-1}, \text{right}) \in h(F)$. Then it must be the case that $e(k'-1,k';l'-2) \notin H$. Then we must have $e(k', k'+1;l'-2) \in H$. It follows that $S_{\uparrow}(k'+1,l'-2;k+,l-1) \in H$. Then there is a southeastern turn $T_1$ on $e(k';l'-2,l'-1)$, $S_{\rightarrow}(k'+1,l'-2;k+2,l-1)$, $e(k+1,k+2;l)$ with both leaves contained in $F$. See Figure 4.23. End of Case 2.2(a). 

\endgroup 

\null

\noindent \textit{CASE 2.2(b).} By (DsFP), $e_{j_0}$ is not eastern. Suppose that $e_{j_0}$ is northern (in orange in Figure 4.25). The assumption that $Q(j_0)$ is false implies that $e_{j_0+1}$ is not eastern and (1stW) implies that $e_{j_0+1}$ is not western. It must be the case that $e_{j_0+1}$ is also northern. By (DsFP) and (RSK), we have that $\Phi(e_{j_0}, \text{right})=X$. But then $j_0-1=s$, contradicting the assumption that $j_0 \in \{5, ..., s-1\}$ (in orange in Figure 4.24). Thus $e_{j_0}$ cannot be northern. It follows that $e_{j_0}$ is southern. Let $e_{j_0}=e(k'';l'',l''+1)$.

\null 

\begingroup
\setlength{\intextsep}{0pt}
\setlength{\columnsep}{20pt}
\begin{wrapfigure}[]{r}{0cm}
\begin{adjustbox}{trim=0cm 0cm 0cm 0cm}

\begin{tikzpicture}[scale=1.35]

\begin{scope}[xshift=0cm]
{
\draw[gray,very thin, step=0.5cm, opacity=0.5] (0,0) grid (3,3.5);

\fill[blue!40!white, opacity=0.5] (0.5,0.5) rectangle (1,1.5);
\fill[blue!40!white, opacity=0.5] (1,1) rectangle (1.5,2);
\fill[blue!40!white, opacity=0.5] (1.5,1.5) rectangle (2,2);
\fill[blue!40!white, opacity=0.5] (1.5,2) rectangle (2.5,2.5);
\fill[blue!40!white, opacity=0.5] (2,2.5) rectangle (3,3);
\fill[blue!40!white, opacity=0.5] (1.5,3) rectangle (2.5,3.5);


\begin{scope}
[very thick,decoration={
    markings,
    mark=at position 0.6 with {\arrow{>}}}
    ]

    \draw[postaction={decorate}, blue, line width=0.5mm] (2,2.5)--++(0,0.5);
    \draw[postaction={decorate}, blue, line width=0.5mm] (2,3)--++(-0.5,0);
    \draw[postaction={decorate}, blue, line width=0.5mm] (1.5,2.5)--++(0.5,0);
    \draw[postaction={decorate}, blue, line width=0.5mm] (1,1.5)--++(0,0.5);
    \draw[postaction={decorate}, blue, line width=0.5mm] (0.5,1.5)--++(0.5,0);

    \draw[postaction={decorate}, blue, line width=0.5mm] (0.5,2)--++(0,-0.5);

    \draw[postaction={decorate}, orange!90!black, dotted, line width=0.5mm] (0.5,0.5)--++(0,0.5);

    \draw[postaction={decorate}, orange!90!black, dotted, line width=0.5mm] (0.5,1)--++(0,0.5);
\end{scope}

\draw[blue,line width=0.5mm] (2.5,3.5)--++(0,-0.5)--++(0.5,0);

\draw[blue,line width=0.5mm] (1,0.5)--++(0,0.5)--++(0.5,0)--++(0,0.5);

\draw[blue,line width=0.5mm] (2,2)--++(0.5,0)--++(0,0.5)--++(0.5,0);

\draw[blue,line width=0.5mm] (0,1)--++(0.5,0)--++(0,-0.5);

\draw[blue, line width=0.5mm] (3,2.5)--++(0,0.5);

{
\draw[black, line width=0.15mm] (2.45,2.7)--++(0.1,0);
\draw[black, line width=0.15mm] (2.45,2.75)--++(0.1,0);
\draw[black, line width=0.15mm] (2.45,2.8)--++(0.1,0);

\draw[black, line width=0.15mm] (2.2,2.45)--++(0,0.1);
\draw[black, line width=0.15mm] (2.25,2.45)--++(0,0.1);
\draw[black, line width=0.15mm] (2.3,2.45)--++(0,0.1);

\draw[black, line width=0.15mm] (2.2,2.95)--++(0,0.1);
\draw[black, line width=0.15mm] (2.25,2.95)--++(0,0.1);
\draw[black, line width=0.15mm] (2.3,2.95)--++(0,0.1);


\draw[black, line width=0.15mm] (0.95,1.2)--++(0.1,0);
\draw[black, line width=0.15mm] (0.95,1.25)--++(0.1,0);
\draw[black, line width=0.15mm] (0.95,1.3)--++(0.1,0);

\draw[black, line width=0.15mm] (0.7,0.95)--++(0,0.1);
\draw[black, line width=0.15mm] (0.75,0.95)--++(0,0.1);
\draw[black, line width=0.15mm] (0.8,0.95)--++(0,0.1);
}

{
\tikzset
  {
    myCircle/.style=
    {
      blue,
      path fading=fade out,
    }
  }

\foreach \x in {0,...,2}
\fill[myCircle,] (1.375+0.125*\x,1.875+0.125*\x) circle (0.075);
}

\node[right] at (3,2.5) [scale=1]{\tiny{-1}};
\node[right] at (3,3) [scale=1]{\tiny{$\ell$}};

\node[above] at (2,3.5) [scale=1]{\tiny{$k$}};
\node[above] at (2.5,3.5) [scale=1]{\tiny{+1}};
\node[above] at (3,3.5) [scale=1]{\tiny{+2}};

\node[above] at (0,3.5) [scale=1]{\tiny{$-1$}};
\node[above] at (0.5,3.5) [scale=1]{\tiny{$k''$}};
\node[above] at (1,3.5) [scale=1]{\tiny{+1}};

\node[right] at (3,1.5) [scale=1]{\tiny{$\ell''$}};
\node[right] at (3,1) [scale=1]{\tiny{-1}};
\node[right] at (3,0.5) [scale=1]{\tiny{-2}};
\node[right] at (3,0) [scale=1]{\tiny{-3}};

\node[above] at (1.65,3) [scale=1]{\small{$e_{_0}$}};
\node[left] at (0.5,1.75) [scale=1]{\small{{$e_{j_0}$}}};
\node[left] at (0.5,0.75) [scale=1]{\small{{\textcolor{orange}{$e_{j_0+1}$}}}};
\node[left] at (0.5,1.25) [scale=1]{\small{{\textcolor{orange}{$e_{j_0}$}}}};

\node at (0.75,1.25) [scale=1] {\small{{\textcolor{orange}{$X$}}}};

\node[below] at (1.5,0) [scale=0.8]{\small{\begin{tabular}{c} Fig. 4.24. Case 2.2(b). \end{tabular}}};;
}

\end{scope}

\end{tikzpicture}
\end{adjustbox}
\end{wrapfigure}

\noindent Using the same arguments as in Case 2.1, we find that $e(k'';l''-1,l'') \notin H$, $e(k''+1;l''-1,l'') \notin H$, $e(k'',k''+1;l''-1) \notin H$, and that $\Phi({j_0-1}, \text{right}) \in P$. It follows that $S_{\rightarrow}(k''-1,l''-1;k'',l''-2) \in H$ and that $S_{\uparrow}(k''+1,l''-2;k+2,l-1) \in H$. Then there is a southeastern turn $T_1$ on $e(k'';l''-2,l''-1)$, $S_{\uparrow}(k''+1,l''-2;k+2,l-1)$, $e(k+1,k+2;l)$.  with both leaves contained in $F$ (Figure 4.25). End of Case 2.2(b). End of Case 2.2. End of Case 2. $\square$ 

\null 

\noindent\textbf{Definition.} Let $G$ be an  $m\times n$ grid graph, let $H$ be a Hamiltonian cycle of $G$, and let $F$ be a looping fat path in $G$. We say that a turn $T$ of $h(F)$ is \textit{admissible} if:

(i) \  no leaf of $T$ is an end-box of $F$, and 

(ii) both leaves of $T$ belong to $F$.  

\null 

\endgroup

\noindent \textbf{Lemma 4.8.} Let $H$ be a Hamiltonian cycle of an $m\times n$ grid graph $G$ and let $h(F)$ be the shadow of a looping fat path of $G$. Then $h(F)$ has an admissible turn. 

\null 

\noindent \textit{Proof.} Let $F$, $X$, $Y$, $\overrightarrow{K}$ and $e_W, e_1, ... e_s$ be as in Lemma 4.7, including the assumption that $F$ is southern and (RSK).

\null 

\noindent \textit{CASE 1: $e_1$ is southern.} By Case 1 in Lemma 4.7, there is a northeastern turn $T_1$. 
We continue sweeping $\overrightarrow{K}$, beginning from $e_W$, until we find the first northern edge $e_N$ in the subtrail $\overrightarrow{K}(e_W,e_N)$ of $\overrightarrow{K}$, where $e_N=(\widehat{k}; \widehat{l}, \widehat{l}+1)= \widehat{e_0}$. 
We write (1stN) to refer to the fact that $e_N$ is the first northern edge encountered after $e_W$, whenever we appeal to it. Let $\widehat{e_1}$ be the edge preceding $\widehat{e_0}$ in the sweep, let $\widehat{e_j}$ be the edge preceding the edge $\widehat{e_{j-1}}$ in the sweep and let $\widehat{e_{t+1}}=e_W$. Then $\widehat{e_1}$ is western or $\widehat{e_1}$ is eastern. 

\null 

\begingroup 
\setlength{\intextsep}{0pt}
\setlength{\columnsep}{20pt}
\begin{wrapfigure}[]{r}{0cm}
\begin{adjustbox}{trim=0cm 0cm 0cm 0cm}
\begin{tikzpicture}[scale=1.5]
\begin{scope}[xshift=0cm]{
\draw[gray,very thin, step=0.5cm, opacity=0.5] (0,0) grid (3,2.5);

\fill[blue!50!white, opacity=0.5]  (0,0)--++(1.5,1.5)--++(1.5,-1.5);
\fill[orange!50!white, opacity=0.5](2.5,0)--++(-1,1)--++(1.5,1.5)--++(0,-2.5);


\draw[black, line width=0.25mm] (0,2)--++(0.5,0)--++(0,-1);
\draw[black, line width=0.25mm] (1.5,2)--++(-0.5,0)--++(0,-1);
\draw[orange, line width=0.5mm, opacity=1] (2.5,0)--++(-1,1)--++(1.5,1.5);
\draw[blue, line width=0.5mm, opacity=1] (0,0)--++(1.5,1.5)--++(1.5,-1.5);
\draw[red, dashed, line width=0.5mm, opacity=1] (0,0)--++(2.5,2.5);

{
\node[left] at (0,2) [scale=1]
{\tiny{0}};
\node[left] at (0,1.5) [scale=1]
{\tiny{-1}};
\node[left] at (0,1) [scale=1]
{\tiny{-2}};

\node[above] at (1,2.5) [scale=1]
{\tiny{0}};
\node[above] at (1.5, 2.5) [scale=1]
{\tiny{1}};
\node[above] at (2,2.5) [scale=1]
{\tiny{2}};

\node at (0.25,1.75) [scale=0.8] {\small{$Y$}};
\node at (1.25,1.75) [scale=0.8] {\small{$X$}};

\draw[fill=blue] (1.5,1) circle [radius=0.05];

\node at (1.5,0.25) [scale=0.8] {$U_2$};

\node at (2.5,1) [scale=0.8] {{$U_{1, \text{end}}$}};

\node[right] at (1.5,1) [scale=0.8] {\small{(a,b)}};
}

\node[below] at (1.5,0) [scale=0.8]{\small{\begin{tabular}{c} Fig. 4.25. Case 3.2(a) The line \\ $y-x+2=0$ in  red; $U_{1, \text{end}}$ shaded \\ orange, $U_2$ shaded blue. 
 \end{tabular}}};;

} \end{scope}

\end{tikzpicture}
\end{adjustbox}
\end{wrapfigure}

\noindent Before we consider each case, we will check that the subtrail $\overrightarrow{K}(\widehat{e_0}, \widehat{e_t})$ of $\overrightarrow{K}$ does not contain the right or left colinear edges of the $A_1$-type of $F$. To this end, we will translate $H$ by $(-k',-l')$ to simplify calculations. (DsFP) and (1stW) imply that for every eastern edge in the subtrail $\overrightarrow{K}(e_s, e_1)$ of $\overrightarrow{K}$ there is at most one northern or southern edge. Denote by ${v_\text{end}}$ the head of the edge $e_W$.
The assumption that $e_1$ is southern and the fact that a shortest turn has length two imply that ${v_\text{end}}$ is contained in the region $U_{1, \text{end}}$, determined by $x \geq 1 $ and $|y+2| \leq x-1$ (Eq.1). Let ${v_\text{end}}=v(a,b)$. It follows that $\overrightarrow{K}(\widehat{e_t}, \widehat{e_0})$ is contained in the region $U_2$ bounded by
$y \leq b+1$ and $|x-a| \leq b+1-y$ (Eq.2). See Figure 4.35. 

We will check that $U_2$ and the colinear edges of the $A_1$-type of $F$ lie on two different sides of the line $y=x-2$. By (Eq.1) we have that $b \leq a-3$ and by (Eq.2) we have that $y \leq x-a+b+1$. Let $(x,y) \in U_2$.  Then $y-x+2 \leq -a+b+3 \leq 0$, so $U_2$ lies below the line $y-x+2$. Plugging in the values of the coordinates of the vertices $A_1$-type $v(x_1,y_1)$, we see that they lie above the line $y=x-2$. This shows that $\overrightarrow{K}(\widehat{e_t}, \widehat{e_0})$ does not contain colinear edges. We will write (NCE) whenever we appeal to this fact. Note that (NCE) implies (i).

\endgroup

\null

\begingroup
\setlength{\intextsep}{0pt}
\setlength{\columnsep}{20pt}
\begin{wrapfigure}[]{l}{0cm}
\begin{adjustbox}{trim=0cm 0cm 0cm 0cm}

\begin{tikzpicture}[scale=1.5]

\begin{scope}[xshift=0cm]
{

\draw[gray,very thin, step=0.5cm, opacity=0.5] (0,0) grid (3,2.5);

\fill[blue!40!white, opacity=0.5] (0,0) rectangle (0.5,0.5);
\fill[blue!40!white, opacity=0.5] (0,0.5) rectangle (1,1);
\fill[blue!40!white, opacity=0.5] (0.5,1) rectangle (1.5,1.5);
\fill[blue!40!white, opacity=0.5] (1,1.5) rectangle (2,2);

\begin{scope}
[very thick,decoration={
    markings,
    mark=at position 0.6 with {\arrow{>}}}
    ]

    \draw[postaction={decorate}, blue, line width=0.5mm] (1.5,2)--++(0.5,0);
    \draw[postaction={decorate}, blue, line width=0.5mm] (2,2)--++(0,-0.5);
    \draw[postaction={decorate}, blue, line width=0.5mm] (2,1.5)--++(-0.5,0);

    \draw[postaction={decorate}, blue, line width=0.5mm] (1,1)--++(0,-0.5);
    \draw[postaction={decorate}, blue, line width=0.5mm] (1,0.5)--++(-0.5,0);
    \draw[postaction={decorate}, blue, line width=0.5mm] (0.5,0.5)--++(0,-0.5);
    \draw[postaction={decorate}, blue, line width=0.5mm] (0.5,0)--++(-0.5,0);
    \draw[postaction={decorate}, blue, line width=0.5mm] (0,0)--++(0,0.5);

    \draw[postaction={decorate}, green!50!black, dotted, line width=0.5mm] (2.5,2)--++(-0.5,0);

    \draw[postaction={decorate}, orange!90!black, dotted, line width=0.5mm] (1,1)--++(-0.5,0);

\end{scope}

{
\draw[green!50!black, line width=0.15mm] (2.45,1.7)--++(0.1,0);
\draw[green!50!black, line width=0.15mm] (2.45,1.75)--++(0.1,0);
\draw[green!50!black, line width=0.15mm] (2.45,1.8)--++(0.1,0);

\draw[green!50!black, line width=0.15mm] (2.45,2.2)--++(0.1,0);
\draw[green!50!black, line width=0.15mm] (2.45,2.25)--++(0.1,0);
\draw[green!50!black, line width=0.15mm] (2.45,2.3)--++(0.1,0);

\draw[green!50!black, line width=0.15mm] (2.7,1.95)--++(0,0.1);
\draw[green!50!black, line width=0.15mm] (2.75,1.95)--++(0,0.1);
\draw[green!50!black, line width=0.15mm] (2.8,1.95)--++(0,0.1);
}

{
\tikzset
  {
    myCircle/.style=
    {
      blue,
      path fading=fade out,
    }
  }

\foreach \x in {0,...,2}
\fill[myCircle,] (0.875+0.125*\x,1.375+0.125*\x) circle (0.075);

}

\node[left] at (0,0.25) [scale=1]{\small{$\widehat{e_{_0}}$}};

\node[above] at (2.25, 2) [scale=1]{\small{{\textcolor{green!50!black}{$\widehat{e_{j_0}}$}}}};

\node[above] at (1.75,2) [scale=1]{\small{$\widehat{e_{j_0}}$}};

\node[above] at (0.75,1) [scale=1]{\small{{\textcolor{orange!90!black}{$e_2$}}}};

\node[right] at (1,0.75) [scale=1]{\small{{\textcolor{orange!90!black}{$e_1$}}}};

\node[above] at (0.75,0.5) [scale=1]{\small{{\textcolor{orange!90!black}{$e_W$}}}};

\node[below] at (1.5,0) [scale=0.8]{\small{\begin{tabular}{c} Fig. 4.26. Case 1.1.
 \end{tabular}}};;

}
\end{scope}

\end{tikzpicture}
\end{adjustbox}
\end{wrapfigure}

\noindent \textit{CASE 1.1: $\widehat{e_1}$ is western.} Note that $\widehat{e_1} \neq e_W$, otherwise we get a cycle on $e_1, e_W, e_1, e_2$. By (DsFP), $\widehat{e_2}$ is not western and by (1stN), $\widehat{e_2}$ is not northern, so $\widehat{e_2}$ must be southern. By (NCE) $\widehat{e_3}$ cannot be southern. Then $\widehat{e_3}$ must be western. If $\widehat{e_3}=e_W$, then we have a southeastern turn $T_2$ on $\widehat{e_0}, ..., \widehat{e_3}, e_1, e_2$, satisfying (i) and, by (RSK), (ii) (in orange in Figure 4.27), so we may assume that $\widehat{e_3} \neq e_W$. By (DsFP), $\widehat{e_4}$ is not western and by (1stN), $\widehat{e_4}$ is not northern. Then $\widehat{e_4}$ is southern.

Let $Q(j)$ be the statement: ``$\widehat{e_j}$ is western and $\widehat{e_{j+1}}$ is southern''.

\endgroup 

\null 

\noindent Now, either $Q(j)$ is true for each $j \in \{1,3, ..., t+1\}$, or there is some $j_0 \in \{ 5,7, ..., t+1\}$  such that $Q(j)$ for each odd $j<j_0$, but $Q(j_0)$ is not true. If the former then we have a southeastern turn $T_2$ on $\widehat{e_0}, \widehat{e_1}, ..., \widehat{e_t}, e_W, e_1,e_2$ satisfying (i) and (ii) (blue in Figure 4.26), so assume the latter. By (NCE), $\widehat{e_{j_0}}$ is not southern. If $\widehat{e_{j_0}}$ is eastern then we have southeastern turn on $\widehat{e_0}, \widehat{e_1}, ..., \widehat{e_{j_0}}$  satisfying (i) and (ii) (blue in figure 4.26). Suppose then, that $\widehat{e_{j_0}}$ is western (green in figure 4.26). This is impossible: by (1stN). $\widehat{e_{j_0+1}}$ is not northern; since $Q(j_0)$ is false, $\widehat{e_{j_0+1}}$ is not southern; and by (DsFP), $\widehat{e_{j_0+1}}$ is not western. End of Case 1.1

\null

\begingroup
\setlength{\intextsep}{0pt}
\setlength{\columnsep}{20pt}
\begin{wrapfigure}[]{r}{0cm}
\begin{adjustbox}{trim=0cm 0cm 0cm 0.5cm}

\begin{tikzpicture}[scale=1.35]

\begin{scope}[xshift=0cm]
{

\draw[gray,very thin, step=0.5cm, opacity=0.5] (0,0) grid (3.5,3.5);

\fill[blue!40!white, opacity=0.5] (0.5,3) rectangle (1.5,3.5);
\fill[blue!40!white, opacity=0.5] (0.5,2.5) rectangle (1,3);
\fill[blue!40!white, opacity=0.5] (0.5,2) rectangle (1.5,2.5);
\fill[blue!40!white, opacity=0.5] (1,1.5) rectangle (2,2);
\fill[blue!40!white, opacity=0.5] (1.5,1) rectangle (2.5,1.5);
\fill[blue!40!white, opacity=0.5] (2,0.5) rectangle (3.5,1);
\fill[blue!40!white, opacity=0.5] (3,1) rectangle (3.5,1.5);

\fill[blue!40!white, opacity=0.5] (0,2.5)--++(0.5,0)--++(0,0.5)--++(-0.5,0);
\fill[blue!40!white, opacity=0.5] (2.5,0)--++(0.5,0)--++(0,0.5)--++(-0.5,0);

\begin{scope}
[very thick,decoration={
    markings,
    mark=at position 0.6 with {\arrow{>}}}
    ]

    \draw[postaction={decorate}, blue, line width=0.5mm] (1.5,3)--++(-0.5,0);
    \draw[postaction={decorate}, blue, line width=0.5mm] (1,3)--++(0,-0.5);
    \draw[postaction={decorate}, blue, line width=0.5mm] (1,2.5)--++(0.5,0);

    \draw[postaction={decorate}, blue, line width=0.5mm] (2,2)--++(0,-0.5);
    \draw[postaction={decorate}, blue, line width=0.5mm] (2,1.5)--++(0.5,0);
    \draw[postaction={decorate}, blue, line width=0.5mm] (2.5,1.5)--++(0,-0.5);
    \draw[postaction={decorate}, blue, line width=0.5mm] (2.5,1)--++(0.5,0);
    \draw[postaction={decorate}, blue, line width=0.5mm] (3,1)--++(0,0.5);

    \draw[postaction={decorate}, orange!90!black, dotted, line width=0.5mm] (0.5,3)--++(0.5,0);
    
\end{scope}

\draw[green!50!black, dotted, line width=0.5mm] (1.5,3)--++(-0.5,0);

\draw[blue, line width=0.5mm] (0,2.5)--++(0.5,0)--++(0,-0.5)--++(0.5,0)--++(0,-0.5)--++(0.5,0)--++(0,-0.5)--++(0.5,0)--++(0,-0.5)--++(0.5,0)--++(0,-0.5);
\draw[blue, line width=0.5mm] (0,3)--++(0.5,0);
\draw[blue, line width=0.5mm] (3,0)--++(0,0.5);
{
\draw[orange!90!black, line width=0.15mm] (0.45,2.7)--++(0.1,0);
\draw[orange!90!black, line width=0.15mm] (0.45,2.75)--++(0.1,0);
\draw[orange!90!black, line width=0.15mm] (0.45,2.8)--++(0.1,0);

\draw[orange!90!black, line width=0.15mm] (0.45,3.2)--++(0.1,0);
\draw[orange!90!black, line width=0.15mm] (0.45,3.25)--++(0.1,0);
\draw[orange!90!black, line width=0.15mm] (0.45,3.3)--++(0.1,0);

\draw[orange!90!black, line width=0.15mm] (0.2,2.95)--++(0,0.1);
\draw[orange!90!black, line width=0.15mm] (0.25,2.95)--++(0,0.1);
\draw[orange!90!black, line width=0.15mm] (0.3,2.95)--++(0,0.1);
}

{
\tikzset
  {
    myCircle/.style=
    {
      blue,
      path fading=fade out,
    }
  }

\foreach \x in {0,...,2}
\fill[myCircle,] (1.625-0.125*\x,1.875+0.125*\x) circle (0.075);

}

\node[left] at (0,1) [scale=0.8]{\small{$\widehat{\ell}$}};

\node[left] at (0,3) [scale=0.8]{\small{$\ell$}};
\node[above] at (1.5,3.5) [scale=0.8]{\small{k}};

\node[above] at (3,3.5) [scale=0.8]{\small{$\widehat{k}$}};

\node[right] at (3,1.25) [scale=1]{\small{$\widehat{e_{_0}}$}};

\node[above] at (0.75, 3) [scale=1]{\small{{\textcolor{orange!90!black}{$\widehat{e_{j_0}}$}}}};

\node[above] at (1.25,3) [scale=1]{\small{$\widehat{e_{t+1}}$}};

\node[below] at (1.35,3) [scale=1]{\small{ \textcolor{green!50!black}{$\widehat{e_{j_0}}$} }};

\node[below] at (1.75,0) [scale=0.8]{\small{\begin{tabular}{c} Fig. 4.27. Case 1.2 (a)  and (b) \end{tabular}}};

}
\end{scope}

\end{tikzpicture}
\end{adjustbox}
\end{wrapfigure}

\noindent \textit{CASE 1.2: $\widehat{e_1}$ is eastern.} By (DsFP), $\widehat{e_2}$ is not eastern and by (1stN), $\widehat{e_2}$ is not northern, so $\widehat{e_2}$ must be southern. By (NCE) $\widehat{e_3}$  is not southern and by (DsFP), $\widehat{e_3}$ is not western. Then $\widehat{e_3}$ must be eastern. (DsFP) and (1stN) imply that $\widehat{e_4}$ must be southern.

Let $Q(j)$ be the statement: ``$\widehat{e_j}$ is eastern and $\widehat{e_{j+1}}$ is southern''. Now, either $Q(j)$ is true for each $j \in \{1, 3, ..., t-1\}$, or there is some $j_0 \in \{ 5,7, ..., t-1\}$  such that $Q(j)$ is true for each odd $j<j_0$, but $Q(j_0)$ is not true. 
 
\null 

\noindent \textit{CASE 1.2 (a): $Q(j)$ is true for each $j \in \{1, 3, ..., t-1\}$.}  Then we have a southwestern turn on $\widehat{e_0}, \widehat{e_1}, ..., \widehat{e_{t-1}}, \widehat{e_t}, e_W$. Recall that $e_W=e(k-1, k;l)$. Observe that $R(k-2; l-1) \in P$, so $e(k-2;l-1,l) \notin H$. Similarly, $R(\widehat{k}-1,\widehat{l}-1) \in P$ and $e(\widehat{k}-1, \widehat{k};\widehat{l}-1) \notin H$. It follows 

\endgroup

\noindent that $R(k-3, l-1) \in F$, $R(\widehat{k}-1,\widehat{l}) \in F$ and that there is a southeastern turn $T_2$ on $e(k-3,k-2;l)$, $S_{\rightarrow}(k-3,l-1;\widehat{k}-1, \widehat{l}-2)$, $e(\widehat{k}; \widehat{l}-2;\widehat{l}-1)$ satisfying (i) and (ii) (blue in Figure 4.27). End of Case 1.2 (a)

\null 

\noindent \textit{CASE 1.2 (b): There is some $j_0 \in \{ 5,7, ..., t-1\}$  such that $Q(j)$ is true for each odd $j<j_0$, but $Q(j_0)$ is not true.} If $\widehat{e_{j_0}}$ is western then we have a southeastern turn on $\widehat{e_0}, \widehat{e_1}, ..., \widehat{e_{j_0}}$. Then as in Case 1.2 (a), there is a southeastern turn $T_2$ satisfying (i) and (ii) (in blue in Figure 4.27, with $e_{j_0}$ in green). 

By (DsFP), $\widehat{e_{j_0}}$ is not southern. Suppose then $\widehat{e_{j_0}}$ is eastern. This is impossible: by (1stN). $\widehat{e_{j_0+1}}$ is not northern; since $Q(j_0)$ is false, $\widehat{e_{j_0+1}}$ is not southern; and by (DsFP), $\widehat{e_{j_0+1}}$ is not eastern (in orange in Figure 4.27). End of Case 1.2 (b). End of Case 1.2.

\endgroup

\null 

\noindent \textit{CASE 2: $e_1$ is northern.}   By Case 2 in Lemma 4.7, $\overrightarrow{K}$, has a southeastern turn $T_1$. We continue sweeping $K$, beginning from $e_W$, until we find the first southern edge $e_S$ in the subtrail $\overrightarrow{K}(e_W,e_S)$ of $\overrightarrow{K}$.  where $e_S=(\widehat{k}; \widehat{l}, \widehat{l}+1)= \widehat{e_0}$. 
We write (1stS) to refer to the fact that $e_S$ is the first southern edge encountered after $e_W$, whenever we appeal to it. Let $\widehat{e_1}$ be the edge preceding $\widehat{e_0}$ in the sweep, let $\widehat{e_j}$ be the edge preceding the edge $\widehat{e_{j-1}}$ in the sweep and let $\widehat{e_{t+1}}=e_W$. Then $\widehat{e_1}$ is western or $\widehat{e_1}$ is eastern.

\null

\noindent \textit{CASE 2.1: $\widehat{e_1}$ is western.} Note that the assumption that $e_1$ is northern implies that $\widehat{e_1} \neq e_W$, otherwise there is a cycle $e_2, e_1, e_W, \widehat{e_0}$. By (DsFP), $\widehat{e_2}$ is not western and by 1stS, $\widehat{e_2}$ is not southern, so $\widehat{e_2}$ must be northern. By (DsFP)
$\widehat{e_3}$ is not northern or eastern. Then $\widehat{e_3}$ must be western. By 1stS $\widehat{e_4}$ is not southern, and by (DsFP), $\widehat{e_4}$ is not western. Then $\widehat{e_4}$ must be northern. Let $Q(j)$ be the statement: ``$\widehat{e_j}$ is western and $\widehat{e_{j+1}}$ is northern. Then either $Q(j)$ is true for each $j \in \{1,3, ..., t-1\}$ or there is some $j_0 \in \{ 5,7, ..., t-1\}$  such that $Q(j)$ for each odd $j<j_0$, but $Q(j_0)$ is not true.

\begingroup
\setlength{\intextsep}{0pt}
\setlength{\columnsep}{15pt}
\begin{wrapfigure}[]{l}{0cm}
\begin{adjustbox}{trim=0cm 0cm 0cm 0.25cm}
\begin{tikzpicture}[scale=1.5]

\begin{scope}[xshift=0cm]
{

\draw[gray,very thin, step=0.5cm, opacity=0.5] (0,0) grid (3.5,3.5);

\fill[blue!40!white, opacity=0.5] (1,2) rectangle (2,2.5);
\fill[blue!40!white, opacity=0.5] (1.5,1.5) rectangle (2.5,2);
\fill[blue!40!white, opacity=0.5] (2,1) rectangle (3,1.5);
\fill[blue!40!white, opacity=0.5] (2.5,0.5) rectangle (3,1);
\fill[blue!40!white, opacity=0.5] (2,0) rectangle (3,0.5);

\fill[blue!40!white, opacity=0.5] (0,2.5)--++(1.5,0)--++(0,0.5)--++(-1.5,0);
\fill[blue!40!white, opacity=0.5] (0.5,3)--++(0.5,0)--++(0,0.5)--++(-0.5,0);
\fill[blue!40!white, opacity=0.5] (3,0.5)--++(0.5,0)--++(0,0.5)--++(-0.5,0);

\begin{scope}
[very thick,decoration={
    markings,
    mark=at position 0.6 with {\arrow{>}}}
    ]

    \draw[postaction={decorate}, blue, line width=0.5mm] (2,0.5)--++(0.5,0);
    \draw[postaction={decorate}, blue, line width=0.5mm] (2.5,0.5)--++(0,0.5);
    \draw[postaction={decorate}, blue, line width=0.5mm] (2.5,1)--++(-0.5,0);

    \draw[postaction={decorate}, blue, line width=0.5mm] (1.5,1.5)--++(0,0.5);
    \draw[postaction={decorate}, blue, line width=0.5mm] (1.5,2)--++(-0.5,0);
    \draw[postaction={decorate}, blue, line width=0.5mm] (1,2)--++(0,0.5);
    \draw[postaction={decorate}, blue, line width=0.5mm] (1,2.5)--++(-0.5,0);
    \draw[postaction={decorate}, blue, line width=0.5mm] (0.5,2.5)--++(0,-0.5);

    \draw[postaction={decorate}, orange, dotted, line width=0.5mm] (3,0.5)--++(-0.5,0);

\end{scope}

\draw[blue, line width=0.5mm] (1,3.5)--++(0,-0.5)--++(0.5,0)--++(0,-0.5)--++(0.5,0)--++(0,-0.5)--++(0.5,0)--++(0,-0.5)--++(0.5,0)--++(0,-0.5)--++(0.5,0);

\draw[blue, line width=0.5mm] (0.5,3)--++(0,0.5);
\draw[blue, line width=0.5mm] (3,0.5)--++(0.5,0);

\draw[orange!90!black, dotted, line width=0.5mm] (2,0.5)--++(0.5,0)--++(0,0.5)--++(-0.5,0);

{
\draw[orange, line width=0.15mm] (2.95,0.8)--++(0.1,0);
\draw[orange, line width=0.15mm] (2.95,0.75)--++(0.1,0);
\draw[orange, line width=0.15mm] (2.95,0.7)--++(0.1,0);

\draw[orange, line width=0.15mm] (2.95,0.3)--++(0.1,0);
\draw[orange, line width=0.15mm] (2.95,0.25)--++(0.1,0);
\draw[orange, line width=0.15mm] (2.95,0.2)--++(0.1,0);

\draw[orange, line width=0.15mm] (3.2,0.55)--++(0,-0.1);
\draw[orange, line width=0.15mm] (3.25,0.55)--++(0,-0.1);
\draw[orange, line width=0.15mm] (3.3,0.55)--++(0,-0.1);
}

{
\tikzset
  {
    myCircle/.style=
    {
      blue,
      path fading=fade out,
    }
  }

\foreach \x in {0,...,2}
\fill[myCircle,] (1.875+0.125*\x,1.625-0.125*\x) circle (0.075);

}

\node[left] at (0.5,2.25) [scale=1]{\small{$\widehat{e_{_0}}$}};

\node[below] at (2.25,0.5) [scale=1]{\small{$e_2$}};
\node[right] at (2.5, 0.75) [scale=1]{\small{$e_1$}};
\node[above] at (2.25, 1) [scale=1]{\small{$ \widehat{e_{t+1}}$}};

\node[right] at (1.5,1.75) [scale=1]{\small{{$\widehat{e_4}$}}};

\node[right] at (3.5,2) [scale=0.8]{\small{$\widehat{\ell}$}};
\node[right] at (3.5,1) [scale=0.8]{\small{$\ell$}};

\node[above] at (2.5,3.5) [scale=0.8]{\small{k}};
\node[above] at (0.5,3.5) [scale=0.8]{\small{$\widehat{k}$}};

\node[below] at (2.75, 0.5) [scale=1]{\small{{\textcolor{orange}{$\widehat{e_{j_0}}$}}}};

\node[below] at (1.75,0) [scale=0.8]{\small{\begin{tabular}{c}  Fig. 4.28. Case 2.1 (a) and (b).
 \end{tabular}}};

}
\end{scope}

\end{tikzpicture}
\end{adjustbox}
\end{wrapfigure}

\null

\noindent \textit{CASE 2.1 (a): $Q(j)$ is true for each $j \in \{1,3, ..., t-1\}$.} Then there is a northeastern turn on $\widehat{e_0}, \widehat{e_1}, ..., \widehat{e_{t-1}},\widehat{e_t}, e_W, e_1, e_2$. As in Case 1.2 (a), we have that $R(\widehat{k},\widehat{l}+1) \in P$, $R(k,l-1) \in P$, $e(\widehat{k}, \widehat{k}+1; \widehat{l}+2) \notin H$ and $e(k+1;l-1,l) \notin H$. Then there is a northeastern turn $T_2$ on $e(\widehat{k};  \widehat{l}+2,\widehat{l}+3)$, $S_{\downarrow}(\widehat{k}+1, \widehat{l}+3; k+2,l)$, $e(k+1,k+2;l-1)$ satisfying (i) and (ii) (in blue in figure 4.28). End of Case 2.1(a).

\null 

\noindent \textit{CASE 2.1 (b): There is some $j_0 \in \{ 5,7, ..., t-1\}$  such that $Q(j)$ for each odd $j<j_0$, but $Q(j_0)$ is not true.} If $\widehat{e_{j_0}}$ is eastern then we have a northeastern turn on $\widehat{e_0}, \widehat{e_1}, ..., \widehat{e_{j_0}}$. Then, as in Case 2.1(a), there is a northeastern turn $T_2$ satisfying (i) and (ii). 

By (DsFP), $\widehat{e_{j_0}}$ is not southern. Suppose then, that $\widehat{e_{j_0}}$ is western. This is impossible: by 1stS. $\widehat{e_{j_0+1}}$ is not southern; since 

\endgroup 

\noindent $Q(j_0)$ is false, $\widehat{e_{j_0+1}}$ is not northern; and by (DsFP), $\widehat{e_{j_0+1}}$ is not western (in orange in Figure 4.28). End of Case 2.1. End of Case 2.1(b). End of Case 2.1.

\null 
\null

\begingroup
\setlength{\intextsep}{0pt}
\setlength{\columnsep}{20pt}
\begin{wrapfigure}[]{r}{0cm}
\begin{adjustbox}{trim=0cm 0cm 0cm 0cm}

\begin{tikzpicture}[scale=1.5]

\begin{scope}[xshift=0cm]
{

\draw[gray,very thin, step=0.5cm, opacity=0.5] (0,0) grid (3,2.5);

\fill[blue!40!white, opacity=0.5] (1,0.5) rectangle (2,1);
\fill[blue!40!white, opacity=0.5] (1.5,1) rectangle (2.5,1.5);
\fill[blue!40!white, opacity=0.5] (2,1.5) rectangle (3,2);
\fill[blue!40!white, opacity=0.5] (2.5,2) rectangle (3,2.5);

\begin{scope}
[very thick,decoration={
    markings,
    mark=at position 0.6 with {\arrow{>}}}
    ]

    \draw[postaction={decorate}, blue, line width=0.5mm] (1.5,0.5)--++(-0.5,0);
    \draw[postaction={decorate}, blue, line width=0.5mm] (1,0.5)--++(0,0.5);
    \draw[postaction={decorate}, blue, line width=0.5mm] (1,1)--++(0.5,0);

    \draw[postaction={decorate}, blue, line width=0.5mm] (2,1.5)--++(0,0.5);
    \draw[postaction={decorate}, blue, line width=0.5mm] (2,2)--++(0.5,0);
    \draw[postaction={decorate}, blue, line width=0.5mm] (2.5,2)--++(0,0.5);
    \draw[postaction={decorate}, blue, line width=0.5mm] (2.5,2.5)--++(0.5,0);
    \draw[postaction={decorate}, blue, line width=0.5mm] (3,2.5)--++(0,-0.5);

    \draw[postaction={decorate}, green!50!black, dotted, line width=0.5mm] (0.5,0.5)--++(0.5,0);

\end{scope}

\draw[orange!90!black, dotted, line width=0.5mm] (1.5,0.5)--++(-0.5,0)--++(0,0.5)--++(0.5,0);

{
\draw[green!50!black, line width=0.15mm] (0.55,0.8)--++(-0.1,0);
\draw[green!50!black, line width=0.15mm] (0.55,0.75)--++(-0.1,0);
\draw[green!50!black, line width=0.15mm] (0.55,0.7)--++(-0.1,0);

\draw[green!50!black, line width=0.15mm] (0.55,0.3)--++(-0.1,0);
\draw[green!50!black, line width=0.15mm] (0.55,0.25)--++(-0.1,0);
\draw[green!50!black, line width=0.15mm] (0.55,0.2)--++(-0.1,0);

\draw[green!50!black, line width=0.15mm] (0.3,0.55)--++(0,-0.1);
\draw[green!50!black, line width=0.15mm] (0.25,0.55)--++(0,-0.1);
\draw[green!50!black, line width=0.15mm] (0.2,0.55)--++(0,-0.1);
}

{
\tikzset
  {
    myCircle/.style=
    {
      blue,
      path fading=fade out,
    }
  }

\foreach \x in {0,...,2}
\fill[myCircle,] (1.625-0.125*\x,1.625-0.125*\x) circle (0.075);

}

\node[right] at (3,2.25) [scale=1]{\small{$\widehat{e_{_0}}$}};

\node[below] at (1.25,0.5) [scale=1]{\small{$\widehat{e_{j_0}}$}};

\node[left] at (2,1.75) [scale=1]{\small{{$\widehat{e_4}$}}};

\node[above] at (1.25,0.5) [scale=1]{\small{{\textcolor{orange!90!black}{$e_{_W}$}}}};



\node[below] at (0.75, 0.5) [scale=1]{\small{{\textcolor{green!50!black}{$\widehat{e_{j_0}}$}}}};

\node[below] at (1.5,0) [scale=0.8]{\small{\begin{tabular}{c} Fig. 4.29. Case 2.2.
 \end{tabular}}};

}
\end{scope}

\end{tikzpicture}

\end{adjustbox}
\end{wrapfigure}

\noindent \textit{CASE 2.2: $\widehat{e_1}$ is eastern.} By 1stS and (DsFP), $\widehat{e_2}$ is northern.  By (DsFP), $\widehat{e_3}$ is not western. An argument analogous to (NCE-1) in Case 1 can be used to show that $T_2$ and the $A_1$-type lie on two different sides of the line $y=2-x$. In this case, we have that the region $U_{1,\text{end}}$ containing $v_{\text{end}}$ is determined by $x \geq 1$ and $|y-2| \leq x-1$, and the region $U_2$ containing $T_2$, as defined in the next paragraph, is determined by $y\geq b-1$ and $|x-a| \leq y-b+1$. We will refer to this argument as (NCE-2). Note that by (NCE-2), $\widehat{e_3}$ is not northern. Then $\widehat{e_3}$ must be eastern. By (DsFP) and (1stS) $\widehat{e_4}$ is not southern or eastern. Then $\widehat{e_4}$ must be northern.

Let $Q(j)$ be the statement: ``$\widehat{e_j}$ is eastern and $\widehat{e_{j+1}}$ is northern''. Now, either $Q(j)$ is true for each $j \in \{1,3, ..., t-1\}$, or there is some $j_0 \in \{ 5,7, ..., t-1\}$  such that $Q(j)$ for each odd $j<j_0$, but $Q(j_0)$ is not true. If the former then we have a northwestern turn $T_2$ on $\widehat{e_0}, \widehat{e_1}, ..., \widehat{e_{t-1}},\widehat{e_t}, e_W$ satisfying (i) and (ii), (blue in figure 4.29 with $\widehat{e_{t-1}},\widehat{e_t}, e_W$ dotted orange) so assume the latter. If $\widehat{e_{j_0}}$ is western then again we have a northwestern turn $T_2$ on $\widehat{e_0}, \widehat{e_1}, ..., \widehat{e_{j_0}}$ satisfying (i) and (ii) (blue in figure 4.29).  By (NCE-2), $\widehat{e_{j_0}}$ is not southern. Then, suppose that $\widehat{e_{j_0}}$ is western (green in figure 4.33). This is impossible: by 1stS. $\widehat{e_{j_0+1}}$ is not southern; since $Q(j_0)$ is false, $\widehat{e_{j_0+1}}$ is not northern; and by (DsFP), $\widehat{e_{j_0+1}}$ is not eastern. End of Case 2.2.

\endgroup 

\null 

\noindent \textbf{Corollary 4.9} Let $G$ be an $m \times n$ grid graph, let $H$ be a Hamiltonian cycle of $G$, let $F=G\langle N[P(X,Y)] \rangle $ be a looping fat path in $G$, let $T$ be an admissible turn of $F$, let $L$ be a leaf of $T$, and let $L'$ be the $H$-neighbour of $L$ in $F$. Then:

(a) $d(T) \geq 3$, and 

(b)  $L \in N[P]\setminus P$ and $L' \in P$.

\null 

\noindent  \textit{Proof of (a).} We prove the contrapositive. For definiteness, assume that $T$ is northeastern with northern leaf $L_N=R(k,l-1)$. Suppose that $d(T) < 3$. Then $d(T)=2$ and the eastern leaf of $T$ must be $L_E=R(k+1,l-2)$. Note that  $L_N+(0,-1) \in F$, otherwise $L_E, ..., L_N, L_N+(0,-1)$ is an $H$-cycle.

Now, by Proposition 4.4, $e(k;l-2,l-1)$ and $e(k,k+1;l-2)$ cannot both belong to $H$. Then, either exactly one of $e(k;l-2,l-1)$ and $e(k,k+1;l-2)$ belongs to $H$, or neither does. 

\null 

\begingroup 
\setlength{\intextsep}{0pt}
\setlength{\columnsep}{20pt}
\begin{wrapfigure}[]{l}{0cm}
\begin{adjustbox}{trim=0cm 0.25cm 0cm 0.5cm}
\begin{tikzpicture}[scale=1.5]

\begin{scope}[xshift=0cm]{
\draw[gray,very thin, step=0.5cm, opacity=0.5] (0,0) grid (2,1);

\fill[blue!50!white, opacity=0.5] (1,0) rectangle (1.5,1);
\fill[blue!50!white, opacity=0.5] (1.5,0) rectangle (2,0.5);

\draw[blue, line width=0.5mm] (1,0)--++(0,1);
\draw[blue, line width=0.5mm] (1.5,0)--++(0.5,0);
\draw[blue, line width=0.5mm] (1.5,1)--++(0,-0.5)--++(0.5,0);

\draw[orange, line width=0.5mm] (1,1)--++(0.5,0);

\draw[orange, line width=0.5mm] (0,1)--++(0.5,0)--++(0,-1);

{
\draw[black, line width=0.15mm] (1.2,-0.05)--++(0,0.1);
\draw[black, line width=0.15mm] (1.25,-0.05)--++(0,0.1);
\draw[black, line width=0.15mm] (1.3,-0.05)--++(0,0.1);
}

{
\node[left] at (0,1) [scale=1]{\tiny{$\ell$}};
\node[left] at (0,0.5) [scale=1]{\tiny{$-1$}};
\node[left] at (0,0) [scale=1]{\tiny{$-2$}};

\node[above] at (1,1) [scale=1]{\tiny{$k$}};
\node[above] at (1.5,1) [scale=1]{\tiny{$+1$}};

\node at (1.25,0.75) [scale=0.8] {\small{$L_N$}};
\node at (1.75,0.25) [scale=0.8] {\small{$L_E$}};
}

\node[below] at (1,0) [scale=0.8]{\small{\begin{tabular}{c} Fig. 4.30. Case 1. \end{tabular}}};;

} \end{scope}

\begin{scope}[xshift=2.75cm]{
\draw[gray,very thin, step=0.5cm, opacity=0.5] (0,0) grid (1.5,1.5);

\fill[blue!50!white, opacity=0.5] (0.5,0.5) rectangle (1,1.5);
\fill[blue!50!white, opacity=0.5] (1,0.5) rectangle (1.5,1);

\draw[blue, line width=0.5mm] (0.5,1)--++(0,0.5);
\draw[blue, line width=0.5mm] (1,0.5)--++(0.5,0);
\draw[blue, line width=0.5mm] (1,1.5)--++(0,-0.5)--++(0.5,0);

\draw[orange, line width=0.5mm] (0,1)--++(0.5,0);
\draw[orange, line width=0.5mm] (0,0.5)--++(0.5,0)--++(0,-0.5);
\draw[orange, line width=0.5mm] (1,0)--++(0,0.5);

{
\draw[black, line width=0.15mm] (0.7,0.45)--++(0,0.1);
\draw[black, line width=0.15mm] (0.75,0.45)--++(0,0.1);
\draw[black, line width=0.15mm] (0.8,0.45)--++(0,0.1);

\draw[black, line width=0.15mm] (0.45,0.7)--++(0.1,0);
\draw[black, line width=0.15mm] (0.45,0.75)--++(0.1,0);
\draw[black, line width=0.15mm] (0.45,0.8)--++(0.1,0);
}

{
\node[left] at (0,1.5) [scale=1]{\tiny{$\ell$}};
\node[left] at (0,1) [scale=1]{\tiny{$-1$}};
\node[left] at (0,0.5) [scale=1]{\tiny{$-2$}};

\node[above] at (0.5,1.5) [scale=1]{\tiny{$k$}};
\node[above] at (1,1.5) [scale=1]{\tiny{$+1$}};

\node at (0.75,1.25) [scale=0.8] {\small{$L_N$}};
\node at (1.25,0.75) [scale=0.8] {\small{$L_E$}};
}

\node[below] at (0.75,0) [scale=0.8]{\small{\begin{tabular}{c} Fig. 4.31. Case 2. \end{tabular}}};;

} \end{scope}

\end{tikzpicture}
\end{adjustbox}
\end{wrapfigure}

\noindent \textit{CASE 1: exactly one of $e(k;l-2,l-1)$ and $e(k,k+1;l-2)$ belongs to $H$.} By symmetry, we may assume WLOG that  $e(k;l-2,l-1) \in H$ and $e(k,k+1;l-2) \notin H$. Note that the assumption that  $e(k,k+1;l-2) \notin H$ implies that $F$ is northern. It follows that $L_N$ is an end-box of $P(X,Y)$. See figure 4.30. End of Case 1.

\null 

\noindent \textit{CASE 2: neither $e(k;l-2,l-1)$ nor $e(k,k+1;l-2)$ belongs to $H$.} Then $e(k-1,k;l-1)$, $e(k+1;l-3,l-2)$ and \ $S_{\rightarrow}(k-1,l-2;k,l-3)$ belong to $H$. By Lemma 4.1(d), $L_N+(0,-1)$ must be long to $P(X,Y)$. It follows that at least one of $L_N$, $L_N+(0,-2)$, $L_E$ and $L_E+(-2,0)$ belongs to $P(X,Y)$ and is switchable, contradicting the assumption that $F$ is a looping fat path. See figure 4.31. End of Case 2. End of proof for (a).

\endgroup

\null 

\begingroup
\setlength{\intextsep}{0pt}
\setlength{\columnsep}{20pt}
\begin{wrapfigure}[]{r}{0cm}
\begin{adjustbox}{trim=0cm 0cm 0 0cm}
\begin{tikzpicture}[scale=1.5]

\begin{scope}[xshift=0cm, yshift=0cm]

\draw[gray,very thin, step=0.5cm, opacity=0.5] (0,0) grid (1,1);

\fill[blue!40!white, opacity=0.5] (0,0)--++(1,0)--++(0,0.5)--++(-1,0);
\fill[green!40!white, opacity=0.5] (0,0.5)--++(0.5,0)--++(0,0.5)--++(-0.5,0);

\draw[blue, line width=0.5mm] (0.5,1)--++(0,-0.5)--++(0.5,0)--++(0,-0.5);
\draw[blue, line width=0.5mm] (0,0.5)--++(0,0.5);
\draw[blue, line width=0.5mm] (0,0)--++(0.5,0);
\draw[blue, line width=0.5mm] (0,1)--++(0.5,0);

{
\draw[black, line width=0.15mm] (0.45,0.2)--++(0.1,0);
\draw[black, line width=0.15mm] (0.45,0.25)--++(0.1,0);
\draw[black, line width=0.15mm] (0.45,0.3)--++(0.1,0);

}

{
\node[right] at (1,0) [scale=1]{\tiny{-1}};
\node[right] at (1,0.5) [scale=1]{\tiny{$b$}};
\node[right] at (1,1) [scale=1]{\tiny{+1}};

\node[above] at (0,1) [scale=1]{\tiny{$a$}};
\node[above] at (0.5,1) [scale=1]{\tiny{+1}};
}
\node at (0.25,0.75) [scale=0.8]{\small{$L_{\textrm{N}}$}};

{
\node[below] at (0.5,0) [scale=0.8]{\small{\begin{tabular}{c} Fig. 4.32(a). \\   $e(a,a{+}1;b{+}1) {\in} H$. \end{tabular}}};;
}

\end{scope}

\begin{scope}[xshift=1.75cm, yshift=0cm]

\draw[gray,very thin, step=0.5cm, opacity=0.5] (0,0) grid (1,1);

\fill[blue!40!white, opacity=0.5] (0,0)--++(1,0)--++(0,0.5)--++(-1,0);
\fill[green!40!white, opacity=0.5] (0,0.5)--++(0.5,0)--++(0,0.5)--++(-0.5,0);

\draw[blue, line width=0.5mm] (0.5,1)--++(0,-0.5)--++(0.5,0)--++(0,-0.5);
\draw[blue, line width=0.5mm] (0,0.5)--++(0,0.5);
\draw[blue, line width=0.5mm] (0,0)--++(0.5,0);

{
\draw[black, line width=0.15mm] (0.45,0.2)--++(0.1,0);
\draw[black, line width=0.15mm] (0.45,0.25)--++(0.1,0);
\draw[black, line width=0.15mm] (0.45,0.3)--++(0.1,0);

\draw[black, line width=0.15mm] (0.2,0.95)--++(0,0.1);
\draw[black, line width=0.15mm] (0.25,0.95)--++(0,0.1);
\draw[black, line width=0.15mm] (0.3,0.95)--++(0,0.1);
}

{
\node[right] at (1,0) [scale=1]{\tiny{-1}};
\node[right] at (1,0.5) [scale=1]{\tiny{$b$}};
\node[right] at (1,1) [scale=1]{\tiny{+1}};

\node[above] at (0,1) [scale=1]{\tiny{$a$}};
\node[above] at (0.5,1) [scale=1]{\tiny{+1}};
}
\node at (0.25,0.75) [scale=0.8]{\small{$L_{\textrm{N}}$}};

{
\node[below] at (0.5,0) [scale=0.8]{\small{\begin{tabular}{c} Fig. 4.32(b). \\   $e(a,a{+}1;b{+}1) {\notin} H$. \end{tabular}}};;
}

\end{scope}

\end{tikzpicture}
\end{adjustbox}
\end{wrapfigure}

\noindent \textit{Proof of (b).} Let $T$ be an admissible turn. For definiteness, assume that $T$ is northeastern with northern leaf $L_N=R(a,b)$ By Corollary 4.9(a), $d(T) \geq 3$. Then we have that $e(a,a+1;b-1) \in H$, and that $e(a+1;b-1,b) \notin H$. By (RSK), $L_N+(0,-1)$ and $L_N+(1,-1)$ belong to $F$. This means that $L_N+(0,-1)=L_N'$ is the $H$-neighbour of $L_N$ in $F$. Now, either $e(a,a+1;b+1) \in H$ or $e(a,a+1;b+1) \notin H$. If $e(a,a+1;b+1) \in H$, then, since $L_N$ is not an end-box of $P$, $L_N \in N[P]\setminus P$. Then, by Lemma 4.1 (b), $L_N' \in P$. And if $e(a,a+1;b+1) \notin H$, then $L_N$ is switchable, so $L_N \in N[P] \setminus P$. Since $L_N' \in F$, by Lemma 4.1 (b), $L_N' \in P$. Either way we have that $L_N' \in P$ and  and $L_N \in N[P]\setminus P$. See Figure 4.32. End of proof for (b). $\square$

\null 
\null 

\subsection{Turn weakenings.}

\begingroup
\setlength{\intextsep}{0pt}
\setlength{\columnsep}{10pt}
\begin{wrapfigure}[]{r}{0cm}
\begin{adjustbox}{trim=0cm 0cm 0cm 0.75cm} 
\begin{tikzpicture}[scale=1.25]
\begin{scope}[xshift=0cm] 
{
\draw[gray,very thin, step=0.5cm, opacity=0.5] (0,0) grid (4,3.5);

\draw[green!50!black, fill=green!50!white, dotted, line width=0.65mm, opacity=0.5] (0.5,2.5)--++(0.5,0)--++(0,-0.5)--++(0.5,0)--++(0,-0.5)--++(0.5,0)--++(0,-0.5)--++(0.5,0)--++(0,-0.5)--++(1.5,0)--++(0,3)--++(-3.5,0)--++(0,-1);

\draw[blue, line width=0.5mm] (1,2.5)--++(0,-0.5)--++(0.5,0)--++(0,-0.5)--++(0.5,0)--++(0,-0.5)--++(0.5,0)--++(0,-0.5)--++(-0.5,0);
\draw[blue, line width=0.5mm] (0.5,2.5)--++(0,-0.5);

\draw[orange, line width=0.5mm] (1.5,3)--++(0,-0.5)--++(0.5,0)--++(0,-0.5)--++(0.5,0)--++(0,-0.5)--++(0.5,0)--++(0,-0.5)--++(0.5,0);
\draw[orange, line width=0.5mm] (3,0.5)--++(0.5,0);
\draw[orange, line width=0.5mm] (1,2.5)--++(0,0.5);

{

\node[left] at (0,2.5) [scale=1]{\tiny{$\ell$}};
\node[left] at (0,0.5) [scale=1]{\tiny{$\ell'$}};

\node[left] at (0,3.5) [scale=1]{\tiny{$n-1$}};

\node[above] at (0.5,3.5) [scale=1]
{\tiny{$k$}};
\node[above] at (2.5,3.5) [scale=1]
{\tiny{$k'$}};
\node[above] at (4,3.5) [scale=1]
{\tiny{$m-1$}};

\node at (0.75,2.25) [scale=1]
{\small{$L_N$}};
\node at (2.25,0.75) [scale=1]
{\small{$L_E$}};

\node[below] at (2,-0.1) [scale=0.75] 
{\small{\begin{tabular}{c}  Fig. 4.33.  A half-open turn $T$ in blue, its lengthening \\ $T'$ in orange, Sector($T$) shaded in green. \end{tabular}}};;

}

}

\end{scope}

\end{tikzpicture}
\end{adjustbox}
\end{wrapfigure}

\textbf{Definitions.} Let $H$ be a Hamiltonian cycle of an $m \times n$ grid graph $G$. Let $T$ be a northeastern turn of $H$ on $\{e(k;l-1,l), S_{\downarrow}(k+1,l;k',l'+1), e(k'-1,k';l')\}$ and let $L_N$ and $L_E$ be the leaves of $T$. Define a \textit{weakening} of $T$ to be a cascade $\mu_1, ..., \mu_s$, where $ \mu_s$ is the first valid nontrivial move in the cascade that has the form $L \mapsto L'$, or $L' \mapsto L$, with $L=L_N$ or $L=L_E$. We note that $\mu_1, ..., \mu_{s-1}$ must avoid $T$, as $T$ has no switchable boxes or leaves, other than $L_1$ and $L_2$. We call a weakening consisting of three or less moves a \textit{short weakening}. We call the subgraph $S_{\downarrow}(k+1,l;k',l'+1)$ \textit{the stairs-part of T} and denote it by $S(T)$. We say that $T$ has a \textit{lengthening} $T'$ if $T'$ is a northeastern turn of $H$ such that: 
 
 a) $d(T') \geq d(T)$ and
 
 b) $S(T') \supseteq S(T)+(1,1)$. 

\noindent Analogous definitions apply to southeastern, southwestern and northwestern turns. We note that if $T'$ is a lengthening of $T$, then $T'$ is unique.  Given a turn $T_0$ let $\mathcal{T}(T_0)=\mathcal{T}$ be a set of lengthenings such that:

1. The turn $T_0 \in \mathcal{T}$

2. The turn $T_j \in \mathcal{T}$ if and only if $T_j$ is a lengthening of the turn $T_{j-1}$.

\noindent Define the \index{sector of T@sector of $T$|textbf}\textit{sector of T} to be the induced subgraph of $G$ bounded by $e(k,k+1;l)$, $S_{\downarrow}(k+1,l;k',l'), e(k'-1,k';l')$, and the segments $[(k',l'), (m-1,l')]$, $[(m-1,l'), (m-1,n-1)]$, $[(m-1,n-1), (k,n-1)]$, $[(k,n-1), (k,l)]$, and denote it by Sector($T$). See Figure 4.33.  Analogous definitions apply to sectors of southeastern, southwestern and northwestern turns.

\null

\noindent \textbf{Lemma 4.10.} Let $H$ be a Hamiltonian cycle of an $m\times n$ grid graph $G$, and let $T$ be a turn in $H$ with $d(T)\geq 3$. Then:

I.  \ $T$ has a short weakening or $T$ has a lengthening.

II. If $T'$ is a lengthening of $T$ and $T'$ has a weakening of length at most $s$, then $T$ has a weakening of 

\hspace{0.4cm} length at most $s+1$, with  $s+1 \leq \min(m,n)$.

\null

\noindent We prove Lemma 4.10 after we use it to prove Proposition 4.11. 

\null

\noindent \textbf{Proposition 4.11.} Let $H$ be a Hamiltonian cycle of an $m\times n$ grid graph $G$, and let $T$ be a turn in $H$ with $d(T)\geq 3$. Then $T$ has a weakening of length at most $\min(m,n)$.

\null 

\noindent \textit{Proof.} Let $T=T_0$ be an admissible turn of $H$. 
If $T_0$ has a short weakening, then we're done, so we assume $T_0$ has no short weakening. By I in Lemma 4.10, $T_0$ has a lengthening $T_1$. So, $T_1\in \mathcal{T}$, where $\mathcal{T}=\mathcal{T}(T_0)$. Since $m,n < \infty$, we have that $|\mathcal{T}| < \infty$. Let $\mathcal{T}=\{T_0, T_1, ..., T_j\}$. 
Then $T_j$ has no lengthening; thus, by I of Lemma 4.10, it must have a short weakening. Then, by induction and II on Lemma 4.10, $T_0$ has a weakening. The bound follows immediately. $\square$

\null 

\noindent \textit{Proof of Lemma 4.10.} We first remark that none of the moves we use throughout this proof fit the description of the moves in Observation 3.4 (i) and (ii) in Section 3. We will use this fact repeatedly and implicitly.

Let $H$ be a Hamiltonian cycle of $G$ and let $T$ be a turn of $H$ with $d(T) \geq 3$. For definiteness, assume that $T$ is northeastern and that $T$ is on $\{e(k;l-1,l), S_{\downarrow}(k+1,l;k',l'+1), e(k'-1,k';l')\}$. Let $L_{\textrm{N}}$ be the northern leaf of $T$ and let $L_{\textrm{E}}$ be the eastern leaf of $T$. Since $d(T)\geq 3$, $m-1 \geq k+3$ and $0\leq l-3$. $L_{\textrm{N}}$ can be open or closed, so there are two cases to check.

\null 

\begingroup
\setlength{\intextsep}{10pt}
\setlength{\columnsep}{20pt}
\begin{wrapfigure}[]{l}{0cm}
\begin{adjustbox}{trim=0cm 0.5cm 0cm 0.5cm}
\begin{tikzpicture}[scale=1.5]

\begin{scope}[xshift=0cm, yshift=0cm]
\draw[gray,very thin, step=0.5cm, opacity=0.5] (0,0) grid (1.5,1);
\draw[yellow, line width =0.2mm] (0,1) -- (1.5,1);

{
\node[right] at (1.5,1) [scale=1]{\tiny{$\ell$}};
\foreach \x in {1,...,2}
\node[right] at (1.5,1 -0.5*\x) [scale=1]
{\tiny{-\x}};


\foreach \x in {1,...,2}
\node[above] at (+0.5*\x, 1) [scale=1]
{\tiny{+\x}};
\node[above] at (0, 1) [scale=1] {\tiny{k}};

}

{

\draw[blue, line width=0.5mm] (0,0.5)--++(0,0.5)--++(0.5,0)--++(0,-0.5)--++(0.5,0)--++(0,-0.5)--++(0.5,0);
}

\draw[] (1,1) circle [radius=0.05];

{

\node at  (0.25,0.75 ) [scale=0.8]{\small{$L_{\textrm{N}}$}};
\node[below] at  (0.75,0) [scale=0.75]{\small{Fig. 4.34(a). Case 1}};

}
\end{scope}

\begin{scope}[xshift=2.25cm, yshift=0cm]
\draw[gray,very thin, step=0.5cm, opacity=0.5] (0,0) grid (1.5,2);
\draw[yellow, line width=0.2mm] (0,2) -- (1.5,2);

{

\node[right] at (1.5,1.5) [scale=1]{\tiny{$\ell$}};
\foreach \x in {1,...,3}
\node[right] at (1.5,1.5 - 0.5*\x) [scale=1]{\tiny{-\x}};

\foreach \x in {1,...,3}
\node[above] at (0.5*\x, 2) [scale=1]{\tiny{+\x}};
\node[above] at (0, 2) [scale=1]{\tiny{k}};

}


\fill[blue!40!white, opacity=0.5] (0,1.5) rectangle (0.5,0.5);
\fill[blue!40!white, opacity=0.5] (0.5,1) rectangle (1,0);
\fill[blue!40!white, opacity=0.5] (1,0) rectangle (1.5,0.5);

\draw[blue, line width=0.5mm] (0,1)--++(0,0.5);

\draw[blue, line width=0.5mm] (0.5,1.5)--++(0,-0.5)--++(0.5,0)--++(0,-0.5)--++(0.5,0)--++(0,0.5);

\draw[blue, line width=0.5mm] (1,0)--++(-0.5,0)--++(0,0.5)--++(-0.5,0);

\draw[blue, line width=0.5mm] (0,1.5)--++(0.5,0);

\fill[blue!40!white, opacity=0.5] (1,1.5) rectangle (1.5,2);

\draw[blue, line width=0.5mm] (1,2)--++(0,-0.5)--++(0.5,0)--++(0,0.5);

\draw[blue, line width=0.5mm] (1,0)--++(0.5,0);


{
\node[below] at (0.75,0) [scale=0.75]{\small{Fig. 4.34(b). Case 1.1.}};
\node at (0.25,1.25) [scale=0.8]{\small{$L_{\textrm{N}}$}};
\node at (1.25,0.25) [scale=0.8]{\small{$L_{\textrm{E}}$}};
}

\end{scope}

\end{tikzpicture}
\end{adjustbox}
\end{wrapfigure}

\noindent \textit{Proof of CASE 1: $L_{\textrm{N}}$ is closed. Proof of I.}  First we note $n-1 \neq l$, otherwise $H$ misses $v(k+2,l)$. Then we must have $S_{\downarrow}(k+2,l+1;k+3,l) \in H$.  Now, $n-1 = l+1$, $n-1 = l+2$, or $n-1 \geq l+3$.

\null 

\noindent \textit{CASE 1.1: n-1=l+1.} By Lemma 1.14, $L_N+(0,1) \in \text{int}(H)$. This implies that $L_N+(2,1)$ is a small cookie of $H$, so $e(k+3;l,l+1) \in H$. Then $e(k+3;l-2,l-1) \in H$. It follows that $L_N+(2,-2)=L_E$. But then $L_E \mapsto L_E+(0,1)$ is a short weakening of $T$. End of Case 1.1.

\endgroup

\begingroup
\setlength{\intextsep}{0pt}
\setlength{\columnsep}{10pt}
\begin{wrapfigure}[]{r}{0cm}
\begin{adjustbox}{trim=0cm 0cm 0cm 0cm}
\begin{tikzpicture}[scale=1.5]

\begin{scope}[xshift=0cm, yshift=0cm]
\draw[gray,very thin, step=0.5cm, opacity=0.5] (0,0) grid (2,1.5);
\draw[yellow, line width=0.2mm] (0,1.5) -- (2,1.5);

{

\node[left] at (0,0.5) [scale=1]{\tiny{$\ell$}};
\node[left] at (0,1) [scale=1]{\tiny{+1}};
\node[left] at (0,1.5) [scale=1]{\tiny{+2}};

\node[above] at (0.5, 1.5) [scale=1]{\tiny{k}};
\node[above] at (1, 1.5) [scale=1]{\tiny{+1}};

}


\fill[blue!40!white, opacity=0.5] (0,0)--++(0.5,0)--++(0,1)--++(-0.5,0);
\fill[blue!40!white, opacity=0.5] (0.5,0.5)--++(0.5,0)--++(0,0.5)--++(-0.5,0);
\fill[blue!40!white, opacity=0.5] (1,0)--++(0.5,0)--++(0,1)--++(-0.5,0);
\fill[blue!40!white, opacity=0.5] (1.5,0)--++(0.5,0)--++(0,0.5)--++(-0.5,0);

\draw[blue, line width=0.5mm] (0.5,0)--++(0,0.5)--++(0.5,0)--++(0,-0.5)--++(0.5,0);

\draw[blue, line width=0.5mm] (0.5,1)--++(0.5,0);

\draw[blue, line width=0.5mm] (1.5,1)--++(0,-0.5)--++(0.5,0);

{
}


\node at (0.75,0.25) [scale=0.8]{\small{$L_{\textrm{N}}$}};

\node[below] at (1,0) [scale=0.75] {\small{\begin{tabular}{c}Fig. 4.35 (a). Case 1.2 (a): \\ $L_N+(0,2) \in \text{ext}(H)$. \end{tabular}}};;

\end{scope}

\begin{scope}[xshift=2.75cm, yshift=0cm]
\draw[gray,very thin, step=0.5cm, opacity=0.5] (0,0) grid (2,1.5);
\draw[yellow, line width=0.2mm] (0,1.5) -- (2,1.5);

{

\node[left] at (0,0.5) [scale=1]{\tiny{$\ell$}};
\node[left] at (0,1) [scale=1]{\tiny{+1}};
\node[left] at (0,1.5) [scale=1]{\tiny{+2}};

\node[above] at (0.5, 1.5) [scale=1]{\tiny{k}};
\node[above] at (1, 1.5) [scale=1]{\tiny{+1}};

}


\fill[blue!40!white, opacity=0.5] (0.5,1)--++(0.5,0)--++(0,0.5)--++(-0.5,0);
\fill[blue!40!white, opacity=0.5] (0.5,0)--++(0.5,0)--++(0,0.5)--++(-0.5,0);
\fill[blue!40!white, opacity=0.5] (1.5,0.5)--++(0.5,0)--++(0,0.5)--++(-0.5,0);

\draw[blue, line width=0.5mm] (0.5,1.5)--++(0,-0.5)--++(0.5,0)--++(0,0.5);
\draw[blue, line width=0.5mm] (0.5,0)--++(0,0.5)--++(0.5,0)--++(0,-0.5);
\draw[blue, line width=0.5mm] (2,0.5)--++(-0.5,0)--++(0,0.5)--++(0.5,0);

{
}


\node at (0.75,0.25) [scale=0.8]{\small{$L_{\textrm{N}}$}};

\node[below] at (1,0) [scale=0.75] {\small{\begin{tabular}{c}Fig. 4.35 (b). Case 1.2 (a): \\ $L_N+(0,2) \in \text{int}(H)$. \end{tabular}}};;

\end{scope}

\begin{scope}[xshift=0cm, yshift=-3cm]
\draw[gray,very thin, step=0.5cm, opacity=0.5] (0,0) grid (2,1.5);
\draw[yellow, line width=0.2mm] (0,1.5) -- (2,1.5);

{

\node[left] at (0,0.5) [scale=1]{\tiny{$\ell$}};
\node[left] at (0,1) [scale=1]{\tiny{+1}};
\node[left] at (0,1.5) [scale=1]{\tiny{+2}};

\node[above] at (0.5, 1.5) [scale=1]{\tiny{k}};
\node[above] at (1, 1.5) [scale=1]{\tiny{+1}};

}



\draw[blue, line width=0.5mm] (0.5,0)--++(0,0.5)--++(0.5,0)--++(0,-0.5)--++(0.5,0);
\draw[blue, line width=0.5mm] (1,1.5)--++(0,-0.5)--++(0.5,0)--++(0,-0.5)--++(0.5,0);
\draw[blue, line width=0.5mm] (0,1)--++(0.5,0)--++(0,0.5);
\draw[blue, line width=0.5mm] (0.5,1.5)--++(0.5,0);

\draw[] (1.5,1.5) circle [radius=0.05];

{
\draw[black, line width=0.15mm] (0.7,0.95)--++(0,0.1);
\draw[black, line width=0.15mm] (0.75,0.95)--++(0,0.1);
\draw[black, line width=0.15mm] (0.8,0.95)--++(0,0.1);
}


\node at (0.75,0.25) [scale=0.8]{\small{$L_{\textrm{N}}$}};

\node[below] at (1,0) [scale=0.75] {\small{\begin{tabular}{c}Fig. 4.36 (a). Case 1.2 (b): \\ $e(k,k+1;l+2) \in H$. \end{tabular}}};;

\end{scope}

\begin{scope}[xshift=2.75cm, yshift=-3cm]
\draw[gray,very thin, step=0.5cm, opacity=0.5] (0,0) grid (2,1.5);
\draw[yellow, line width=0.2mm] (0,1.5) -- (2,1.5);

{

\node[left] at (0,0.5) [scale=1]{\tiny{$\ell$}};
\node[left] at (0,1) [scale=1]{\tiny{+1}};
\node[left] at (0,1.5) [scale=1]{\tiny{+2}};

\node[above] at (0.5, 1.5) [scale=1]{\tiny{k}};
\node[above] at (1, 1.5) [scale=1]{\tiny{+1}};

}


\fill[blue!40!white, opacity=0.5] (0,0)--++(0.5,0)--++(0,1)--++(-0.5,0);
\fill[blue!40!white, opacity=0.5] (0.5,0.5)--++(0.5,0)--++(0,1)--++(-0.5,0);
\fill[blue!40!white, opacity=0.5] (1,0)--++(0.5,0)--++(0,1)--++(-0.5,0);
\fill[blue!40!white, opacity=0.5] (1.5,0)--++(0.5,0)--++(0,0.5)--++(-0.5,0);

\draw[blue, line width=0.5mm] (0.5,0)--++(0,0.5)--++(0.5,0)--++(0,-0.5)--++(0.5,0);
\draw[blue, line width=0.5mm] (1,1.5)--++(0,-0.5)--++(0.5,0)--++(0,-0.5)--++(0.5,0);
\draw[blue, line width=0.5mm] (0,1)--++(0.5,0)--++(0,0.5);
\draw[blue, line width=0.5mm] (1,1.5)--++(0.5,0);

{
\draw[black, line width=0.15mm] (0.7,0.95)--++(0,0.1);
\draw[black, line width=0.15mm] (0.75,0.95)--++(0,0.1);
\draw[black, line width=0.15mm] (0.8,0.95)--++(0,0.1);

\draw[black, line width=0.15mm] (0.7,1.45)--++(0,0.1);
\draw[black, line width=0.15mm] (0.75,1.45)--++(0,0.1);
\draw[black, line width=0.15mm] (0.8,1.45)--++(0,0.1);
}


\node at (0.75,0.25) [scale=0.8]{\small{$L_{\textrm{N}}$}};

\node[below] at (1,0) [scale=0.75] {\small{\begin{tabular}{c}Fig. 4.36 (b). Case 1.2 (b): \\ $e(k,k+1;l+2) \notin H$. \end{tabular}}};;

\end{scope}

\end{tikzpicture}
\end{adjustbox}
\end{wrapfigure}

\noindent \textit{CASE 1.2: n-1= l+2.} Either $e(k,k+1;l+1) \in H$ or $e(k,k+1;l+1) \notin H$. 

\null 

\noindent \textit{CASE 1.2(a): $e(k,k+1;l+1) \in H$.} Either $L_N+(0,2) \in \text{int}(H)$ or $L_N+(0,2) \in \text{ext}(H)$. If $L_N+(0,2) \in \text{int}(H)$, then $L_N+(0,1) \mapsto L_N$ is a short weakening of $T$. Suppose then, that $L_N+(0,2) \in \text{ext}(H)$. This implies that $L_N+(0,2)$ is a small cookie. Then $L_N+(0,1) \mapsto L_N+(0,2)$, $L_N \mapsto L_N+(1,1)$ is a short weakening. End of Case 1.2 (a).

\null 

\noindent \textit{CASE 1.2(b): $e(k,k+1;l+1) \notin H$.} Then $S_{\rightarrow}(k-1,l+1;k,l+2) \in H$ and $S_{\downarrow}(k+1,l+2;k+2,l+1) \in H$. Note that if $e(k,k+1;l+2) \in H$, then $H$ misses $v(k+2,l+2)$, so we may assume that $e(k,k+1;l+2) \notin H$. It follows that $e(k+1,k+2;l+2) \in H$. Then $L_N+(0,2) \mapsto L_N+(1,2)$, $L_N+(0,1) \mapsto L_N$ is a short weakening. End of Case 1.2(b). End of Case 1.2.

\endgroup 

\null 

\begingroup
\setlength{\intextsep}{0pt}
\setlength{\columnsep}{15pt}
\begin{wrapfigure}[]{l}{0cm}
\begin{adjustbox}{trim=0cm 0.5cm 0cm 0.5cm}
\begin{tikzpicture}[scale=1.5]

\draw[gray,very thin, step=0.5cm, opacity=0.5] (0,0) grid (3,3);

{
\foreach \x in {1,...,2}
\node[left] at (0,+0.5*\x+2) [scale=1]
{\tiny{\x}};
\node[left] at (0,2) [scale=1]{\tiny{$\ell$}};

\node[left] at (0,1) [scale=1]{\tiny{+1}};
\node[left] at (0,0.5) [scale=1]{\tiny{$\ell'$}};

\node[above] at (0.5, 3) [scale=1] {\tiny{k}};
\node[above] at (1, 3) [scale=1] {\tiny{+1}};

\foreach \x in {1,...,2}
\node[above]  at (+0.5*\x+2, 3) [scale=1]
{\tiny{\x$'$}};
\node[above]  at (2, 3) [scale=1] {\tiny{k$'$}};
}

{

\draw[blue, line width=0.5mm] (0.5,1.5)--++(0,0.5)--++(0.5,0)--++(0,-0.5);

\draw[blue, line width=0.5mm] (1.5,1)--++(0.5,0)--++(0,-0.5)--++(-0.5,0);

}

{

\draw[orange, line width=0.5mm](0,2.5)--++(0.5,0)--++(0,0.5);
\draw[orange, line width=0.5mm](2.5,0)--++(0,0.5)--++(0.5,0);

\draw[orange, line width=0.5mm](1,3)--++(0,-0.5)--++(0.5,0)--++(0,-0.5);
\draw[orange, line width=0.5mm](3,1)--++(-0.5,0)--++(0,0.5)--++(-0.5,0);
}

{
\tikzset
  {
    myCircle/.style=    {orange}
  }
 
\foreach \x in {0,...,2}
\fill[myCircle, orange,] (1.625+0.125*\x,1.875-0.125*\x) circle (0.05);

\foreach \x in {0,...,2}
\fill[myCircle, blue,] (1.125+0.125*\x,1.375-0.125*\x) circle (0.05);

}


{
\draw[black, line width=0.15mm] (0.7,2.45)--++(0,0.1);
\draw[black, line width=0.15mm] (0.75,2.45)--++(0,0.1);
\draw[black, line width=0.15mm] (0.8,2.45)--++(0,0.1);

\draw[black, line width=0.15mm] (2.45,0.7)--++(0.1,0);
\draw[black, line width=0.15mm] (2.45,0.75)--++(0.1,0);
\draw[black, line width=0.15mm] (2.45,0.8)--++(0.1,0);

}

{

\node at  (0.75,1.75 ) [scale=0.8]{\small{$L_{\textrm{N}}$}};

\node at  (0.75,2.75 ) [scale=0.8]{\small{$ \hat{L}_{\textrm{N}}$}};

\node at  (1.75,0.75 ) [scale=0.8]{\small{$L_{\textrm{E}}$}};

\node[below] at  (1.5,0 ) [scale=0.75]{\small{Fig. 4.37. Case 1.3(a).}};

}

\end{tikzpicture}
\end{adjustbox}
\end{wrapfigure}

\noindent \textit{CASE 1.3: $n-1 \geq l+3$.} By Case 1.2, we may assume that $e(k,k+1;l+1) \notin H$, $S_{\rightarrow}(k-1,l+1;k,l+2) \in H$ and that $S_{\downarrow}(k+1,l+2;k+2,l+1) \in H$. Now $L_{\textrm{E}}$ is either open or closed.

\null 

\noindent \textit{CASE 1.3(a): $L_{\textrm{E}}$ is closed.} By previous cases and symmetry we may assume that $m-1 \geq k'+3$. Using symmetry once more, we may assume that $e(k'+1;l',l'+1) \notin H$. Then the turn $\hat{T}$  on $\{e(k; l+1,l+2), S_{\downarrow}(k+1,l+2;k'+2,l'+1), e(k'+1,k'+2;l')\}$ is in $H$ and it is a lengthening of $T$. End of proof of I for Case 1.3(a).

\null 

\noindent \textit{Proof of II for Case 1.3 (a).} WLOG assume that the last move $\mu_s$ of a weakening $\mu_1,..., \mu_s$ of $\hat{T}$ is $Z \mapsto \hat{L}_{\textrm{N}}$, where $\hat{L}_{\textrm{N}}$ is the northern leaf of $\hat{T}$.  Then $\mu_1,..., \mu_s, L_N+(0,1) \mapsto L_N$, is a weakening of $T$. 

\noindent It remains to check that $s+1 \leq \min(m,n)$. Since the $j^{\text{th}}$ lengthening $T_j$ in $\mathcal{T}(T)$ is $j$ units north and east of $T$, and $d(T) \geq 3$, there can be at most $\min(m,n)-3$ such lengthenings. Since a turn with no lengthening has a short weakening, $s+1 \leq 3+ \min(m,n)-3 =\min(m,n)$. End of II for Case 1.3(a). End of Case 1.3(a).

\endgroup 

\begingroup
\setlength{\intextsep}{0pt}
\setlength{\columnsep}{10pt}
\begin{wrapfigure}[]{r}{0cm}
\begin{adjustbox}{trim=0cm 0.25cm 0cm 0.5cm}
\begin{tikzpicture}[scale=1.5]

\draw[gray,very thin, step=0.5cm, opacity=0.5] (0,0) grid (2.5,2.5);

{
\foreach \x in {1,...,2}
\node[left] at (0,0.5*\x+1.5) [scale=1]{\tiny{\x}};
\node[left] at (0,1.5) [scale=1]{\tiny{$\ell$}};

\node[left] at (0,0.5) [scale=1]{\tiny{+1}};
\node[left] at (0,0) [scale=1]{\tiny{$\ell'$}};

\node[above]  at (0.5, 2.5) [scale=1]{\tiny{k}};
\node[above]  at (1, 2.5) [scale=1]{\tiny{+1}};

\foreach \x in {1,...,1}
\node[above]  at (0.5*\x+2, 2.5) [scale=1]{\tiny{\x$'$}};
\node[above]  at (2, 2.5) [scale=1]{\tiny{k$'$}};
}

{
\draw[blue, line width=0.5mm] (0.5,1)--++(0,0.5)--++(0.5,0)--++(0,-0.5);

\draw[blue, line width=0.5mm] (1.5,0.5)--++(0.5,0);
\draw[blue, line width=0.5mm] (1.5,0)--++(0.5,0);
}

{
\draw[orange, line width=0.5mm](0,2)--++(0.5,0)--++(0,0.5);

\draw[orange, line width=0.5mm](1,2.5)--++(0,-0.5)--++(0.5,0)--++(0,-0.5);
\draw[orange, line width=0.5mm](2,1)--++(0.5,0);
\draw[orange, line width=0.5mm](2,0.5)--++(0.5,0);
}

{
\tikzset
  {
    myCircle/.style=    {orange}
  }

\foreach \x in {0,...,2}
\fill[myCircle, orange,] (1.625+0.125*\x,1.375-0.125*\x) circle (0.05);

\foreach \x in {0,...,2}
\fill[myCircle, blue,] (1.125+0.125*\x,0.875-0.125*\x) circle (0.05);
}

{
\draw[black, line width=0.15mm] (0.7,1.95)--++(0,0.1);
\draw[black, line width=0.15mm] (0.75,1.95)--++(0,0.1);
\draw[black, line width=0.15mm] (0.8,1.95)--++(0,0.1);

\draw[black, line width=0.15mm] (1.95,0.2)--++(0.1,0);
\draw[black, line width=0.15mm] (1.95,0.25)--++(0.1,0);
\draw[black, line width=0.15mm] (1.95,0.3)--++(0.1,0);

\draw[black, line width=0.15mm] (1.95,0.7)--++(0.1,0);
\draw[black, line width=0.15mm] (1.95,0.75)--++(0.1,0);
\draw[black, line width=0.15mm] (1.95,0.8)--++(0.1,0);
}

{
\node at (0.75,1.25) [scale=0.8]{\small{$L_{\textrm{N}}$}};

\node at (0.75,2.25) [scale=0.8]{\small{$\hat{L}_{\textrm{N}}$}};

\node at (1.75,0.25) [scale=0.8]{\small{$L_{\textrm{E}}$}};

\node[below] at (1.25,0) [scale=0.75]{\small{Fig. 4.38. Case 1.3(b).}};
}

\end{tikzpicture}
\end{adjustbox}
\end{wrapfigure}

\null

\noindent \textit{CASE 1.3(b): $L_{\textrm{E}}$ is open.} Either $m-1 < k'+2$ or $m-1 \geq k'+2$. It will follow from Case 2 that if a turn has on open leaf adjacent to the boundary or at distance one away from the boundary, then we can find a weakening outright. Therefore, we may assume that $m-1 \geq k'+2$. 

If $e(k';l'+1,l'+2) \in H$, then there is a weakening 
$L_E \mapsto L_E+(0,1)$, so we may assume that $e(k';l'+1,l'+2) \notin H$. Then the turn $\hat{T}$  on $\{e(k; l+1,l+2), S_{\downarrow}(k+1,l+2;k'+1,l'+2), e(k',k'+1;l'+1)\}$ is in $H$ and it is a lengthening of $T$. End of proof of I for Case 1.3.

\null 

\noindent \textit{Proof of II for Case 1.3(b).}Let $\hat{L_N}$ and $\hat{L_E}$ be the northern and eastern leaves of $\hat{T}$ respectively, and let $\mu_1,..., \mu_s$ be a weakening of $\hat{T}$. If $\mu_s$ is the move $X \mapsto \hat{L}_N$, then, as in the Case 1.3(a),  $\mu_1,..., \mu_s, L_N+(0,1) \mapsto L_N$, is a weakening of $T$. Suppose then that $\mu_s$ is the move $Z' \mapsto \hat{L}_E$. Then $\mu_1,..., \mu_s, L_E \mapsto L_E+(0,1)$, is a weakening of $T$. The argument that $s+1 \leq \min(m,n)$ is the same as the one in Case 1.3(a), so we omit it. End of proof of II for Case 1.3(b) End of Case 1.3(b). End of Case 1.3. End of Case 1.

\endgroup

\null 

\noindent \textit{CASE 2: $L_{\textrm{N}}$ is open. Proof of I.} If $n-1=l$ then we must have $e(k+1, k+2;l) \in H$. Then $L_{\textrm{N}} \mapsto L_{\textrm{N}}+(1,0)$ is a weakening. Therefore, we may assume that $n-1 >l$.

\null 

\noindent \textit{CASE 2.1: $n-1=l+1$.} Either $e(k+1, k+2;l) \in H$ or  $e(k+1, k+2;l) \notin H$.

\begingroup
\setlength{\intextsep}{0pt}
\setlength{\columnsep}{10pt}
\begin{wrapfigure}[]{l}{0cm}
\begin{adjustbox}{trim=0cm 0cm 0cm 0.5cm}
\begin{tikzpicture}[scale=1.5]

\begin{scope}[xshift=0cm, yshift=0cm]
\draw[gray,very thin, step=0.5cm, opacity=0.5] (0,0) grid (1,1);
\draw[yellow, line width =0.2mm] (0,1) -- (1,1);

{

\node[right] at (1,1) [scale=1]{\tiny{$\ell$}};

\node[above] at (0, 1) [scale=1] {\tiny{k}};
\node[above] at (0.5, 1) [scale=1] {\tiny{+1}};

}

{
\fill[blue!40!white, opacity=0.5] (0,0)  rectangle (1,0.5);
\fill[blue!40!white, opacity=0.5] (0,0.5)  rectangle (0.5,1);

\draw[blue, line width=0.5mm] (0,0.5)--++(0,0.5);
\draw[blue, line width=0.5mm] (1,1)--++(-0.5,0)--++(0,-0.5)--++(0.5,0);

}

{

\draw[black, line width=0.15mm] (0.25,0.95)--++(0,0.1);
\draw[black, line width=0.15mm] (0.2,0.95)--++(0,0.1);
\draw[black, line width=0.15mm] (0.3,0.95)--++(0,0.1);

}

{
\node at  (0.25,0.75 ) [scale=0.8]{\small{$L_{\textrm{N}}$}};

\node[below] at (0.5,0) [scale=0.75]{\small{\begin{tabular}{c} Fig. 4.39 (a).  \\ Case  2. \\$n-1=l.$  \end{tabular}}};

}

\end{scope}

\begin{scope}[xshift=1.5cm, yshift=0cm]
\draw[gray,very thin, step=0.5cm, opacity=0.5] (0,0) grid (1.5,1.5);
\draw[yellow, line width =0.2mm] (0,1.5) -- (1.5,1.5);

{

\node[right] at (1.5,1) [scale=1]{\tiny{$\ell$}};
\node[right] at (1.5,0.5) [scale=1]{\tiny{$-1$}};

\node[above] at (0.5, 1.5) [scale=1] {\tiny{k}};
\node[above] at (1, 1.5) [scale=1] {\tiny{+1}};

}

{
\fill[blue!40!white, opacity=0.5] (0,0.5)  rectangle (0.5,1);
\fill[blue!40!white, opacity=0.5] (1,0.5)  rectangle (1.5,1);

\draw[blue, line width=0.5mm] (0,1)--++(0.5,0)--++(0,-0.5)--++(-0.5,0);
\draw[blue, line width=0.5mm] (1.5,1)--++(-0.5,0)--++(0,-0.5)--++(0.5,0);

}

{

\draw[black, line width=0.15mm] (0.75,0.95)--++(0,0.1);
\draw[black, line width=0.15mm] (0.7,0.95)--++(0,0.1);
\draw[black, line width=0.15mm] (0.8,0.95)--++(0,0.1);

}

{
\node at  (0.75,0.75 ) [scale=0.8]{\small{$L_{\textrm{N}}$}};

\node[below] at (0.75,0) [scale=0.75]{\small{\begin{tabular}{c} Fig. 4.39 (b). Case  \\ 2.1($a_1$). $L_N+(-1,0)$ is \\ a small cookie.  \end{tabular}}};

}

\end{scope}

\begin{scope}[xshift=3.5cm, yshift=0cm]
\draw[gray,very thin, step=0.5cm, opacity=0.5] (0,0) grid (1.5,1.5);
\draw[yellow, line width =0.2mm] (0,1.5) -- (1.5,1.5);

{

\node[right] at (1.5,1) [scale=1]{\tiny{$\ell$}};
\node[right] at (1.5,0.5) [scale=1]{\tiny{$-1$}};

\node[above] at (0.5, 1.5) [scale=1] {\tiny{k}};
\node[above] at (1, 1.5) [scale=1] {\tiny{+1}};

}

{
\fill[blue!40!white, opacity=0.5] (0,0.5)  rectangle (0.5,1);
\fill[blue!40!white, opacity=0.5] (1,0.5)  rectangle (1.5,1);

\draw[blue, line width=0.5mm] (0.5,1)--++(0,-0.5);
\draw[blue, line width=0.5mm] (1.5,1)--++(-0.5,0)--++(0,-0.5)--++(0.5,0);

}

{

\draw[black, line width=0.15mm] (0.75,0.95)--++(0,0.1);
\draw[black, line width=0.15mm] (0.7,0.95)--++(0,0.1);
\draw[black, line width=0.15mm] (0.8,0.95)--++(0,0.1);

}

{
\node at  (0.75,0.75 ) [scale=0.8]{\small{$L_{\textrm{N}}$}};

\node[below] at (0.75,0) [scale=0.75]{\small{\begin{tabular}{c} Fig. 4.39 (c). Case  \\ 2.1($a_1$). $L_N+(-1,0)$ is \\ not a small cookie.  \end{tabular}}};

}

\end{scope}

\end{tikzpicture}
\end{adjustbox}
\end{wrapfigure}

\null 

\noindent \textit{CASE 2.1(a): $e(k+1, k+2;l) \in H$.} Then $e(k,k+1;l+1) \in H$ and $e(k+1,k+2;l+1) \in H$. Then by Lemma 1.14, $L_N+(0,1) \in \textrm{int}H$, and so $L_N \in \textrm{int}H$ and $L_N+(1,0) \in \textrm{ext}H$. Either $k=0$, or $k>0$.

\null 

\noindent \textit{CASE 2.1($a_1$): $k>0$.} Then $L_N+(-1,0) \in \textrm{ext}H$, and $L_N+(-1,0)$ is either  is a small cookie of $H$ or it is not. If the former, then $L_N \mapsto L_N+(-1,0)$ is a short weakening; and if the latter then $L_N \mapsto L_N+(1,0)$ is a short weakening. See figures 4.39 (b) and (c). End of Case 2.1($a_1$).

\endgroup 

\begingroup
\setlength{\intextsep}{0pt}
\setlength{\columnsep}{10pt}
\begin{wrapfigure}[]{r}{0cm}
\begin{adjustbox}{trim=0cm 0.25cm 0cm 0.5cm}
\begin{tikzpicture}[scale=1.5]

\begin{scope}[xshift=0cm, yshift=0cm]
\draw[gray,very thin, step=0.5cm, opacity=0.5] (0,0) grid (2,2.5);
\draw[yellow, line width =0.2mm] (0,2.5) -- (2,2.5);

{

\node[left] at (0,2) [scale=1]{\tiny{$\ell$}};
\node[left] at (0,1.5) [scale=1]{\tiny{$-1$}};
\node[left] at (0,1) [scale=1]{\tiny{$-2$}};
\node[left] at (0,0.5) [scale=1]{\tiny{$-3$}};
\node[left] at (0,0) [scale=1]{\tiny{$-4$}};

\node[above] at (0, 2.5) [scale=1] {\tiny{0}};
\node[above] at (0.5, 2.5) [scale=1] {\tiny{1}};
\node[above] at (1, 2.5) [scale=1] {\tiny{2}};
\node[above] at (1.5, 2.5) [scale=1] {\tiny{3}};
\node[above] at (2, 2.5) [scale=1] {\tiny{4}};
}

{
\fill[blue!40!white, opacity=0.5] (0,0.5)--++(0.5,0)--++(0,0.5)--++(-0.5,0);
\fill[blue!40!white, opacity=0.5] (1,0)--++(1,0)--++(0,0.5)--++(-1,0);
\fill[blue!40!white, opacity=0.5] (0.5,1.5)--++(1,0)--++(0,0.5)--++(-1,0);
\fill[blue!40!white, opacity=0.5] (1,1)--++(1,0)--++(0,0.5)--++(-1,0);

\draw[blue, line width=0.5mm] (1.5,2.5)--++(0.5,0);

\draw[blue, line width=0.5mm] (0,0)--++(0,0.5)--++(0.5,0)--++(0,0.5)--++(-0.5,0)--++(0,1.5)--++(1,0);

\draw[blue, line width=0.5mm] (1,2)--++(-0.5,0)--++(0,-0.5)--++(0.5,0)--++(0,-0.5)--++(0.5,0);

\draw[blue, line width=0.5mm] (1,0)--++(0,0.5)--++(0.5,0);

\draw[blue, line width=0.5mm] (1.5,2)--++(0,-0.5)--++(0.5,0);

}

{

\draw[black, line width=0.15mm] (0.25,1.95)--++(0,0.1);
\draw[black, line width=0.15mm] (0.2,1.95)--++(0,0.1);
\draw[black, line width=0.15mm] (0.3,1.95)--++(0,0.1);

\draw[black, line width=0.15mm] (1.75,1.95)--++(0,0.1);
\draw[black, line width=0.15mm] (1.7,1.95)--++(0,0.1);
\draw[black, line width=0.15mm] (1.8,1.95)--++(0,0.1);

\draw[black, line width=0.15mm] (1.45,0.2)--++(0.1,0);
\draw[black, line width=0.15mm] (1.45,0.25)--++(0.1,0);
\draw[black, line width=0.15mm] (1.45,0.3)--++(0.1,0);

\draw[black, line width=0.15mm] (1.45,1.2)--++(0.1,0);
\draw[black, line width=0.15mm] (1.45,1.25)--++(0.1,0);
\draw[black, line width=0.15mm] (1.45,1.3)--++(0.1,0);
}

{
\node at  (1.25,0.75) [scale=0.8]{\small{$L_{\textrm{E}}$}};
\node at  (0.25,1.75 ) [scale=0.8]{\small{$L_{\textrm{N}}$}};

\node[below] at (1,0) [scale=0.75]{\small{\begin{tabular}{c} Fig. 4.40. Case  2.1$(a_2)$ \end{tabular}}};

}

\end{scope}

\end{tikzpicture}
\end{adjustbox}
\end{wrapfigure}

\null

\noindent \textit{CASE 2.1($a_2$): $k=0$.} Then $e(0;l,l+1) \in H$, $e(0;l-2,l-1) \in H$, and $S_{\rightarrow}(0,l-2;1,l-3) \in H$. This implies that $0\leq l-4$ and that $S_{\leftarrow}(1,l-3;0,l-4) \in H$. Then we must have $S_{\uparrow}(2,l-4;3,l-3) \in H$ as well, that $L_E=L_N+(2,-2)$ and that $L_E+(0,-1) \in \textrm{ext}H$. See Figure 4.40. Note that if $L_E+(0,-1)$ is a small cookie, then $L_E \mapsto L_E+(0,-1)$ is a short weakening, so we may assume that $L_E+(0,-1)$ is not a small cookie. Note that this implies that $0\leq l-5$. 

If $e(3;l-4,l-3) \in H$, then again $L_E \mapsto L_E+(0,-1)$ is a short weakening. Similarly, if $e(3;l-2,l-1) \in H$, then $L_E \mapsto L_E+(0,1)$ is a short weakening. Therefore we only need to check the case where $e(3;l-4,l-3) \notin H$ and $e(3;l-2,l-1) \notin H$. Then $L_E+(1,-1) \in \textrm{ext}H$, and by the assumption that $H$ is Hamiltonian, $m-1 \geq 4$. Then we have $S_{\downarrow}(3,l;4,l-1) \in H$.

\begingroup
\setlength{\intextsep}{0pt}
\setlength{\columnsep}{20pt}
\begin{wrapfigure}[]{l}{0cm}
\begin{adjustbox}{trim=0.25cm 0.25cm 0cm 0cm}
\begin{tikzpicture}[scale=1.5]

\begin{scope}[xshift=0cm, yshift=0cm]
\draw[gray,very thin, step=0.5cm, opacity=0.5] (0,0) grid (2.5,2.5);
\draw[yellow, line width =0.2mm] (0,2.5) -- (2,2.5);

{

\node[left] at (0,2) [scale=1]{\tiny{$\ell$}};
\node[left] at (0,1.5) [scale=1]{\tiny{$-1$}};
\node[left] at (0,1) [scale=1]{\tiny{$-2$}};
\node[left] at (0,0.5) [scale=1]{\tiny{$-3$}};
\node[left] at (0,0) [scale=1]{\tiny{$-4$}};

\node[above] at (0, 2.5) [scale=1] {\tiny{$0$}};
\node[above] at (0.5, 2.5) [scale=1] {\tiny{1}};
\node[above] at (1, 2.5) [scale=1] {\tiny{2}};
\node[above] at (1.5, 2.5) [scale=1] {\tiny{3}};
\node[above] at (2, 2.5) [scale=1] {\tiny{4}};
}

{
\fill[blue!40!white, opacity=0.5] (2,2)--++(0.5,0)--++(0,0.5)--++(-0.5,0);
\fill[blue!40!white, opacity=0.5] (0,0.5)--++(0.5,0)--++(0,0.5)--++(-0.5,0);
\fill[blue!40!white, opacity=0.5] (1,0)--++(1,0)--++(0,0.5)--++(-1,0);
\fill[blue!40!white, opacity=0.5] (0.5,1.5)--++(1,0)--++(0,0.5)--++(-1,0);
\fill[blue!40!white, opacity=0.5] (1,1)--++(1,0)--++(0,0.5)--++(-1,0);

\fill[blue!40!white, opacity=0.5] (1.5,0.5)--++(1,0)--++(0,0.5)--++(-1,0);

\draw[blue, line width=0.5mm] (1.5,2.5)--++(0.5,0);

\draw[blue, line width=0.5mm] (0,0)--++(0,0.5)--++(0.5,0)--++(0,0.5)--++(-0.5,0)--++(0,1.5)--++(1,0);

\draw[blue, line width=0.5mm] (1,2)--++(-0.5,0)--++(0,-0.5)--++(0.5,0)--++(0,-0.5)--++(0.5,0);

\draw[blue, line width=0.5mm] (1,0)--++(0,0.5)--++(0.5,0);

\draw[blue, line width=0.5mm] (1.5,2)--++(0,-0.5)--++(0.5,0);

\draw[blue, line width=0.5mm] (1.5,0.5)--++(0,0.5);
\draw[blue, line width=0.5mm] (2,0)--++(0,0.5)--++(0.5,0);
\draw[blue, line width=0.5mm] (2,1.5)--++(0,-0.5)--++(0.5,0)--++(0,0.5);
\draw[blue, line width=0.5mm] (2,2.5)--++(0,-0.5)--++(0.5,0)--++(0,0.5);

}

{

\draw[black, line width=0.15mm] (0.25,1.95)--++(0,0.1);
\draw[black, line width=0.15mm] (0.2,1.95)--++(0,0.1);
\draw[black, line width=0.15mm] (0.3,1.95)--++(0,0.1);

\draw[black, line width=0.15mm] (1.75,1.95)--++(0,0.1);
\draw[black, line width=0.15mm] (1.7,1.95)--++(0,0.1);
\draw[black, line width=0.15mm] (1.8,1.95)--++(0,0.1);

\draw[black, line width=0.15mm] (1.45,0.2)--++(0.1,0);
\draw[black, line width=0.15mm] (1.45,0.25)--++(0.1,0);
\draw[black, line width=0.15mm] (1.45,0.3)--++(0.1,0);

\draw[black, line width=0.15mm] (1.45,1.2)--++(0.1,0);
\draw[black, line width=0.15mm] (1.45,1.25)--++(0.1,0);
\draw[black, line width=0.15mm] (1.45,1.3)--++(0.1,0);

\draw[black, line width=0.15mm] (1.95,0.7)--++(0.1,0);
\draw[black, line width=0.15mm] (1.95,0.75)--++(0.1,0);
\draw[black, line width=0.15mm] (1.95,0.8)--++(0.1,0);
}

{
\node at  (1.25,0.75) [scale=0.8]{\small{$L_{\textrm{E}}$}};
\node at  (0.25,1.75 ) [scale=0.8]{\small{$L_{\textrm{N}}$}};

\node[below] at (1.25,0) [scale=0.75]{\small{\begin{tabular}{c} Fig. 4.41. Case  2.1$(a_2).(i)$ \end{tabular}}};

}

\end{scope}

\end{tikzpicture}
\end{adjustbox}
\end{wrapfigure}

\noindent Note that either $e(2,3;l+1) \in H$ and $e(2,3;l) \in H$, or $e(2;l,l+1) \in H$ and $e(3;l,l+1) \in H$. Either way, we must have $e(3,4;l) \notin H$ and $e(3,4;l+1) \in H$. Now, either $e(3;l-3,l-2) \in H$ or $e(3;l-3,l-2) \notin H$.

\null

\noindent \textit{CASE 2.1($a_2$).(i): $e(3;l-3,l-2) \in H$.} Then $L_E+(1,0) \in \textrm{ext}H$. By Lemma 1.14, this implies that $m-1 \geq 5$. If $e(4;l-3,l-2) \in H$, then $L_E+(1,0) \mapsto L_E$ is a short weakening, so we may assume that $e(4;l-3,l-2) \notin H$. Then $S_{\downarrow}(4,l-1;5,l-2) \in H$ and $S_{\uparrow}(4,l-4;5,l-3) \in H$. We must also have that $S_{\downarrow}(4,l+1;5,l)\in H$ and $e(5;l,l+1)\in H$. Then $e(5;l-2,l-1)\in H$ as well. See Figure 4.41. Then $L_E+(2,0) \mapsto L_E+(2,1)$, $L_E+(1,0) \mapsto L_E$ is a short weakening of $T$. End of Case 2.1($a_2$).(i).

\endgroup 

\null

\begingroup
\setlength{\intextsep}{00pt}
\setlength{\columnsep}{10pt}
\begin{wrapfigure}[]{r}{0cm}
\begin{adjustbox}{trim=0cm 0.5cm 0cm 0.5cm}
\begin{tikzpicture}[scale=1.5]

\begin{scope}[xshift=0cm, yshift=0cm]
\draw[gray,very thin, step=0.5cm, opacity=0.5] (0,0) grid (2.5,2.5);
\draw[yellow, line width =0.2mm] (0,2.5) -- (2,2.5);

{

\node[left] at (0,2) [scale=1]{\tiny{$\ell$}};
\node[left] at (0,1.5) [scale=1]{\tiny{$-1$}};
\node[left] at (0,1) [scale=1]{\tiny{$-2$}};
\node[left] at (0,0.5) [scale=1]{\tiny{$-3$}};
\node[left] at (0,0) [scale=1]{\tiny{$-4$}};

\node[above] at (0, 2.5) [scale=1] {\tiny{$0$}};
\node[above] at (0.5, 2.5) [scale=1] {\tiny{1}};
\node[above] at (1, 2.5) [scale=1] {\tiny{2}};
\node[above] at (1.5, 2.5) [scale=1] {\tiny{3}};
\node[above] at (2, 2.5) [scale=1] {\tiny{4}};
}

{
\fill[blue!40!white, opacity=0.5] (2,1)--++(0.5,0)--++(0,0.5)--++(-0.5,0);
\fill[blue!40!white, opacity=0.5] (2,2)--++(0.5,0)--++(0,0.5)--++(-0.5,0);
\fill[blue!40!white, opacity=0.5] (0,0.5)--++(0.5,0)--++(0,0.5)--++(-0.5,0);
\fill[blue!40!white, opacity=0.5] (1,0)--++(1,0)--++(0,0.5)--++(-1,0);
\fill[blue!40!white, opacity=0.5] (0.5,1.5)--++(1,0)--++(0,0.5)--++(-1,0);
\fill[blue!40!white, opacity=0.5] (1,1)--++(1,0)--++(0,0.5)--++(-1,0);

\draw[blue, line width=0.5mm] (1.5,2.5)--++(0.5,0);

\draw[blue, line width=0.5mm] (0,0)--++(0,0.5)--++(0.5,0)--++(0,0.5)--++(-0.5,0)--++(0,1.5)--++(1,0);

\draw[blue, line width=0.5mm] (1,2)--++(-0.5,0)--++(0,-0.5)--++(0.5,0)--++(0,-0.5)--++(1,0);

\draw[blue, line width=0.5mm] (1,0)--++(0,0.5)--++(1,0);

\draw[blue, line width=0.5mm] (1.5,2)--++(0,-0.5)--++(0.5,0);

\draw[blue, line width=0.5mm] (1,2.5)--++(0.5,0);
\draw[blue, line width=0.5mm] (1,2)--++(0.5,0);
\draw[blue, line width=0.5mm] (2,2.5)--++(0,-0.5)--++(0.5,0)--++(0,0.5);

\draw[blue, line width=0.5mm] (2,1.5)--++(0.5,0);
}

{

\draw[black, line width=0.15mm] (0.25,1.95)--++(0,0.1);
\draw[black, line width=0.15mm] (0.2,1.95)--++(0,0.1);
\draw[black, line width=0.15mm] (0.3,1.95)--++(0,0.1);

\draw[black, line width=0.15mm] (1.75,1.95)--++(0,0.1);
\draw[black, line width=0.15mm] (1.7,1.95)--++(0,0.1);
\draw[black, line width=0.15mm] (1.8,1.95)--++(0,0.1);

\draw[black, line width=0.15mm] (1.45,0.2)--++(0.1,0);
\draw[black, line width=0.15mm] (1.45,0.25)--++(0.1,0);
\draw[black, line width=0.15mm] (1.45,0.3)--++(0.1,0);

\draw[black, line width=0.15mm] (1.45,1.2)--++(0.1,0);
\draw[black, line width=0.15mm] (1.45,1.25)--++(0.1,0);
\draw[black, line width=0.15mm] (1.45,1.3)--++(0.1,0);

\draw[black, line width=0.15mm] (1.95,1.7)--++(0.1,0);
\draw[black, line width=0.15mm] (1.95,1.75)--++(0.1,0);
\draw[black, line width=0.15mm] (1.95,1.8)--++(0.1,0);

\draw[black, line width=0.15mm] (1.95,1.2)--++(0.1,0);
\draw[black, line width=0.15mm] (1.95,1.25)--++(0.1,0);
\draw[black, line width=0.15mm] (1.95,1.3)--++(0.1,0);
}

{
\node at  (1.25,0.75) [scale=0.8]{\small{$L_{\textrm{E}}$}};
\node at  (0.25,1.75 ) [scale=0.8]{\small{$L_{\textrm{N}}$}};

\node[below] at (1.25,0) [scale=0.75]{\small{\begin{tabular}{c} Fig. 4.42(a). Case  2.2$.a_2(ii)_1$. \end{tabular}}};

}

\end{scope}

\begin{scope}[xshift=3.25cm, yshift=0cm]
\draw[gray,very thin, step=0.5cm, opacity=0.5] (0,0) grid (2.5,2.5);
\draw[yellow, line width =0.2mm] (0,2.5) -- (2,2.5);

{

\node[left] at (0,2) [scale=1]{\tiny{$\ell$}};
\node[left] at (0,1.5) [scale=1]{\tiny{$-1$}};
\node[left] at (0,1) [scale=1]{\tiny{$-2$}};
\node[left] at (0,0.5) [scale=1]{\tiny{$-3$}};
\node[left] at (0,0) [scale=1]{\tiny{$-4$}};

\node[above] at (0, 2.5) [scale=1] {\tiny{$0$}};
\node[above] at (0.5, 2.5) [scale=1] {\tiny{1}};
\node[above] at (1, 2.5) [scale=1] {\tiny{2}};
\node[above] at (1.5, 2.5) [scale=1] {\tiny{3}};
\node[above] at (2, 2.5) [scale=1] {\tiny{4}};
}

{
\fill[blue!40!white, opacity=0.5] (1,2)--++(0.5,0)--++(0,0.5)--++(-0.5,0);
\fill[blue!40!white, opacity=0.5] (2,1)--++(0.5,0)--++(0,0.5)--++(-0.5,0);
\fill[blue!40!white, opacity=0.5] (2,2)--++(0.5,0)--++(0,0.5)--++(-0.5,0);
\fill[blue!40!white, opacity=0.5] (0,0.5)--++(0.5,0)--++(0,0.5)--++(-0.5,0);
\fill[blue!40!white, opacity=0.5] (1,0)--++(1,0)--++(0,0.5)--++(-1,0);
\fill[blue!40!white, opacity=0.5] (0.5,1.5)--++(1,0)--++(0,0.5)--++(-1,0);
\fill[blue!40!white, opacity=0.5] (1,1)--++(1,0)--++(0,0.5)--++(-1,0);

\draw[blue, line width=0.5mm] (1.5,2.5)--++(0.5,0);

\draw[blue, line width=0.5mm] (0,0)--++(0,0.5)--++(0.5,0)--++(0,0.5)--++(-0.5,0)--++(0,1.5)--++(1,0);

\draw[blue, line width=0.5mm] (1,2)--++(-0.5,0)--++(0,-0.5)--++(0.5,0)--++(0,-0.5)--++(1,0);

\draw[blue, line width=0.5mm] (1,0)--++(0,0.5)--++(1,0);

\draw[blue, line width=0.5mm] (1.5,2)--++(0,-0.5)--++(0.5,0);

\draw[blue, line width=0.5mm] (1,2)--++(0,0.5);
\draw[blue, line width=0.5mm] (1.5,2)--++(0,0.5);
\draw[blue, line width=0.5mm] (2,2.5)--++(0,-0.5)--++(0.5,0)--++(0,0.5);

\draw[blue, line width=0.5mm] (2,1.5)--++(0.5,0);
}

{

\draw[black, line width=0.15mm] (0.25,1.95)--++(0,0.1);
\draw[black, line width=0.15mm] (0.2,1.95)--++(0,0.1);
\draw[black, line width=0.15mm] (0.3,1.95)--++(0,0.1);

\draw[black, line width=0.15mm] (1.75,1.95)--++(0,0.1);
\draw[black, line width=0.15mm] (1.7,1.95)--++(0,0.1);
\draw[black, line width=0.15mm] (1.8,1.95)--++(0,0.1);

\draw[black, line width=0.15mm] (1.45,0.2)--++(0.1,0);
\draw[black, line width=0.15mm] (1.45,0.25)--++(0.1,0);
\draw[black, line width=0.15mm] (1.45,0.3)--++(0.1,0);

\draw[black, line width=0.15mm] (1.45,1.2)--++(0.1,0);
\draw[black, line width=0.15mm] (1.45,1.25)--++(0.1,0);
\draw[black, line width=0.15mm] (1.45,1.3)--++(0.1,0);

\draw[black, line width=0.15mm] (1.95,1.7)--++(0.1,0);
\draw[black, line width=0.15mm] (1.95,1.75)--++(0.1,0);
\draw[black, line width=0.15mm] (1.95,1.8)--++(0.1,0);

\draw[black, line width=0.15mm] (1.95,1.2)--++(0.1,0);
\draw[black, line width=0.15mm] (1.95,1.25)--++(0.1,0);
\draw[black, line width=0.15mm] (1.95,1.3)--++(0.1,0);
}

{
\node at  (1.25,0.75) [scale=0.8]{\small{$L_{\textrm{E}}$}};
\node at  (0.25,1.75 ) [scale=0.8]{\small{$L_{\textrm{N}}$}};

\node[below] at (1.25,0) [scale=0.75]{\small{\begin{tabular}{c} Fig. 4.42(b). Case  2.2$.a_2(ii)_2$. \end{tabular}}};

}

\end{scope}

\end{tikzpicture}
\end{adjustbox}
\end{wrapfigure}

\noindent \textit{CASE 2.1($a_2$).(ii): $e(3;l-3,l-2) \notin H$.} Then $e(3,4;l-3) \in H$ and $e(3,4;l-2) \in H$. Now, either $e(2,3;l+1) \in H$ and $e(2,3;l) \in H$, or $e(2;l,l+1) \in H$ and $e(3;l,l+1) \in H$.

\null 

\noindent \textit{CASE 2.1($a_2$).$(ii)_1$:  $e(2,3;l+1) \in H$ and $e(2,3;l) \in H$.} Then $e(4;l-2,l-1) \notin H$. Note that if $e(4;l-1,l)\in H$, then $L_E+(1,1) \mapsto L_E+(1,2)$, $L_E \mapsto L_E+(0,1)$ is a short weakening of $T$, so we may assume that. $e(4;l-1,l)\notin H$. Then $e(4,5;l-1) \in H$, $S_{\downarrow}(4,l+1;5,l) \in H$, and $e(5;l,l+1) \in H$. Then $L_E+(2,2) \mapsto L_E+(2,3)$, $L_E+(1,1) \mapsto L_E+(1,2)$, $L_E \mapsto L_E+(0,1)$ is a short weakening of $T$. End of Case 2.1($a_2$).$(ii)_1$.

\null 

\noindent \textit{CASE 2.1($a_2$).$(ii)_2$:  $e(2;l,l+1) \in H$ and $e(3;l,l+1) \in H$.} Now, if $e(4;l-2,l-1) \notin H$, than we can use the same argument and find the same cascades as in Case 2.1($a_2$).$(ii)_1$ (Fig, 4.22(b); and if $e(4;l-2,l-1) \in H$, then we must have that $e(5;l-2,l-1) \in H$ as well. Then $L_E+(2,1) \mapsto L_E+(1,1)$, Then $L_E \mapsto L_E+(0,1)$ is a short weakening of $T$. End of Case 2.1($a_2$).$(ii)_2$. End of Case 2.1($a_2$).$(ii)$.
End of Case 2.1($a_2$). End of Case 2.1(a).

\endgroup

\null 

\noindent \textit{CASE 2.1(b): $e(k+1, k+2;l) \notin H$.} Then we must have $e(k+1;l,l+1) \in H$ and $S_{\downarrow}(k+2,l+1;k+3,l) \in H$. Now, either $e(k+1,k+2;l+1) \notin H$ or $e(k+1,k+2;l+1) \in H$.

\begingroup
\setlength{\intextsep}{0pt}
\setlength{\columnsep}{10pt}
\begin{wrapfigure}[]{r}{0cm}
\begin{adjustbox}{trim=0cm 0.5cm 0cm 0.75cm}
\begin{tikzpicture}[scale=1.5]

\begin{scope}[xshift=0cm, yshift=0cm]
\draw[gray,very thin, step=0.5cm, opacity=0.5] (0,0) grid (1.5,1);
\draw[yellow, line width =0.2mm] (0,1) -- (1,1);

{

\node[right] at (1.5,0.5) [scale=1]{\tiny{$\ell$}};


\foreach \x in {1,...,3}
\node[above] at (+0.5*\x, 1) [scale=1]
{\tiny{+\x}};
\node[above] at (0, 1) [scale=1] {\tiny{k}};

}

{

\fill[blue!40!white, opacity=0.5] (0.5,0)  rectangle (1,1);

\draw[blue, line width=0.5mm] (0,0)--++(0,0.5);
\draw[blue, line width=0.5mm] (1,0)--++(-0.5,0)--++(0,1);
\draw[blue, line width=0.5mm] (1,0.5)--++(0,0.5)--++(0.5,0);

\draw[blue, line width=0.5mm] (1,0.5)--++(0.5,0);

}

{

\draw[black, line width=0.15mm] (0.25,0.45)--++(0,0.1);
\draw[black, line width=0.15mm] (0.2,0.45)--++(0,0.1);
\draw[black, line width=0.15mm] (0.3,0.45)--++(0,0.1);

\draw[black, line width=0.15mm] (0.75,0.45)--++(0,0.1);
\draw[black, line width=0.15mm] (0.7,0.45)--++(0,0.1);
\draw[black, line width=0.15mm] (0.8,0.45)--++(0,0.1);

\draw[black, line width=0.15mm] (0.75,0.95)--++(0,0.1);
\draw[black, line width=0.15mm] (0.7,0.95)--++(0,0.1);
\draw[black, line width=0.15mm] (0.8,0.95)--++(0,0.1);

}

{
\node at  (0.25,0.25 ) [scale=0.8]{\small{$L_{\textrm{N}}$}};

\node[below] at (0.7,0) [scale=0.75] {\small{\begin{tabular}{c}Fig. 4.43(a). Case 2.1($b_1$). \end{tabular}}};;

}

\end{scope}

\begin{scope}[xshift=2.25cm, yshift=0cm]
\draw[gray,very thin, step=0.5cm, opacity=0.5] (0,0) grid (1.5,2);
\draw[yellow, line width=0.2mm] (0,2) -- (1.5,2);

{

\node[right] at (1.5,1.5) [scale=1]{\tiny{$\ell$}};
\node[right] at (1.5,1) [scale=1]{\tiny{-1}};
\node[right] at (1.5,0.5) [scale=1]{\tiny{-2}};
\node[right] at (1.5,0) [scale=1]{\tiny{-3}};

\foreach \x in {1,...,3}
\node[above] at (0.5*\x, 2) [scale=1]{\tiny{+\x}};
\node[above] at (0, 2) [scale=1]{\tiny{k}};

}

{
\fill[blue!40!white, opacity=0.5] (1,1.5) rectangle (1.5,2);

\fill[blue!40!white, opacity=0.5] (0,1) rectangle (0.5,2);
\fill[blue!40!white, opacity=0.5] (0,0.5) rectangle (1,1);
\fill[blue!40!white, opacity=0.5] (0.5,0) rectangle (1.5,0.5);

\draw[blue, line width=0.5mm] (0,1)--++(0,0.5);
\draw[blue, line width=0.5mm] (1,1)--++(-0.5,0)--++(0,1)--++(0.5,0);

\draw[blue, line width=0.5mm] (1,1)--++(0,-0.5)--++(0.5,0)--++(0,0.5);
\draw[blue, line width=0.5mm] (1,0)--++(0.5,0);

\draw[blue, line width=0.5mm] (1,2)--++(0,-0.5)--++(0.5,0)--++(0,0.5);
}

{
\draw[black, line width=0.15mm] (0.25,1.45)--++(0,0.1);
\draw[black, line width=0.15mm] (0.2,1.45)--++(0,0.1);
\draw[black, line width=0.15mm] (0.3,1.45)--++(0,0.1);

\draw[black, line width=0.15mm] (0.25,1.95)--++(0,0.1);
\draw[black, line width=0.15mm] (0.2,1.95)--++(0,0.1);
\draw[black, line width=0.15mm] (0.3,1.95)--++(0,0.1);

}

{
\node at (0.25,1.25) [scale=0.8]{\small{$L_{\textrm{N}}$}};

\node at (1.25,0.25) [scale=0.8]{\small{$L_{\textrm{E}}$}};

\node[below] at (0.75,0) [scale=0.75]{\small{\begin{tabular}{c} Fig. 4.43(b).  Case 2.1($b_2$).  \end{tabular}}};

}

\end{scope}

\end{tikzpicture}
\end{adjustbox}
\end{wrapfigure}

\null 

\noindent \textit{CASE 2.1($b_1$): $e(k+1, k+2;l+1) \notin H$.} Then we must have $e(k+2,k+3;l+1) \in H$ and that $L_N+(1,1) \in \textrm{ext}H$ is the neck of the large cookie. Then $k>0$, or $k=0$.

If $k>0$, then after $L_N+(1,1) \mapsto L_N+(2,1)$, we are back to Case 2.1($a_1$). And if $k=0$ then we are effectively in the same scenario as in Case 2.1($a_2$), with the additional, inconsequential assumption that $L_N+(1,1) \in \textrm{ext}H$ is the neck of the large cookie. End of Case 2.1($b_1$).

\null 

\noindent \textit{CASE 2.1($b_2$): $e(k+1, k+2;l+1) \in H$.} Then $e(k,k+1;l+1) \notin H$ and $S_{\downarrow}(k+2,l+1;k+3,l) \in H$. Then we must have that $e(k+3;l,l+1)\in H$ as well. This implies that $e(k+3;l-2,l-1)\in H$. It follows that $L_N+(2,-2)=L_E$, and that $L_E$ is open. Then $L_E \mapsto L_E+(0,1)$ is a short weakening. End of Case 2.1($b_2$). End of Case 2.1(b) End of Case 2.1.

\endgroup 

\null

\begingroup
\setlength{\intextsep}{0pt}
\setlength{\columnsep}{20pt}
\begin{wrapfigure}[]{l}{0cm}
\begin{adjustbox}{trim=0cm 0cm 0cm 0.5cm}
\begin{tikzpicture}[scale=1.35]

\draw[gray,very thin, step=0.5cm, opacity=0.5] (0,0) grid (3,2.5);

{
\foreach \x in {1,...,2}
\node[left] at (0,0.5*\x+1.5) [scale=1]{\tiny{\x}};
\node[left] at (0,1.5) [scale=1]{\tiny{$\ell$}};

\node[left] at (0,0.5) [scale=1]{\tiny{+1}};
\node[left] at (0,0) [scale=1]{\tiny{$\ell'$}};

\node[above] at (0.5, 2.5) [scale=1]{\tiny{k}};
\node[above] at (1, 2.5) [scale=1]{\tiny{+1}};

\foreach \x in {1,...,1}
\node[above] at (0.5*\x+2, 2.5) [scale=1]{\tiny{\x$'$}};
\node[above] at (2, 2.5) [scale=1]{\tiny{k$'$}};
}

{
\draw[blue, line width=0.5mm] (0.5,1)--++(0,0.5);

\draw[blue, line width=0.5mm] (1,1)--++(0,0.5);

\draw[blue, line width=0.5mm] (1.5,0.5)--++(0.5,0);
\draw[blue, line width=0.5mm] (1.5,0)--++(0.5,0);
}

{
\draw[orange, line width=0.5mm](1,1.5)--++(0,0.5);
\draw[orange, line width=0.5mm](1.5,1.5)--++(0,0.5);

\draw[orange, line width=0.5mm](2,1)--++(0.5,0);
\draw[orange, line width=0.5mm](2,0.5)--++(0.5,0);
}

{
\tikzset
  {
    myCircle/.style=    {orange}
  }

\foreach \x in {0,...,2}
\fill[myCircle, orange,] (1.625+0.125*\x,1.375-0.125*\x) circle (0.05);

\foreach \x in {0,...,2}
\fill[myCircle, blue,] (1.125+0.125*\x,0.875-0.125*\x) circle (0.05);
}

{
\draw[black, line width=0.15mm] (1.2,1.45)--++(0,0.1);
\draw[black, line width=0.15mm] (1.25,1.45)--++(0,0.1);
\draw[black, line width=0.15mm] (1.3,1.45)--++(0,0.1);

\draw[black, line width=0.15mm] (1.95,0.7)--++(0.1,0);
\draw[black, line width=0.15mm] (1.95,0.75)--++(0.1,0);
\draw[black, line width=0.15mm] (1.95,0.8)--++(0.1,0);

}

{
\node at (0.75,1.25) [scale=0.8]{\small{$L_{\textrm{N}}$}};

\node at (1.25,1.75) [scale=0.8]{\small{$\hat{L}_{\textrm{N}}$}};

\node at (1.75,0.25) [scale=0.8]{\small{$L_{\textrm{E}}$}};

\node[below] at (1.5,0) [scale=0.75]{\small{Fig. 4.44.  Case 2.2.}};
}

\end{tikzpicture}
\end{adjustbox}
\end{wrapfigure}

\noindent \textit{CASE 2.2: $n-1 \geq l+2$.} By previous cases we may assume that $m-1 \geq k'+2$, $e(k+1, k+2;l) \notin H$, $e(k+1;l,l+1) \in H$ and $S_{\downarrow}(k+2,l+1;k+3,l) \in H$. If $L_{\textrm{E}}$ is closed, then we're done by Case 1, so we may assume that $L_{\textrm{E}}$ is open. By Case 2.1, we may assume that $e(k';l'+1,l'+2) \notin H$. Then the turn $\hat{T}$  on $\{e(k+1; l,l+1), S_{\downarrow}(k+2,l+1; k'+1,l'+2), e(k',k'+1;l'+1)\}$ is in $H$ and it is a lengthening of $T$. End of proof of I for Case 2.2.

The proof of II for Case 2.2 is the same as the proof of II for Case 1.3 (a). End of proof for Case 2. $\square$

\endgroup

\null 

\noindent \textbf{Observation 4.12.} All turn weakenings found in Lemma 4.10 are contained in Sector(T). $\square$

\null

\noindent \textbf{Lemma 4.13.} Let $G$ be an $m \times n$ grid graph, and let $H$ be a Hamiltonian cycle of $G$. Let $F$ be a looping fat path of $G$, anchored at some outermost small cookie $C$. Then $F$ has an admissible turn $T$ such that Sector$(T)$ and the $j$-stack of $A_0$'s following $C$ are disjoint.

\null

\begingroup
\setlength{\intextsep}{0pt}
\setlength{\columnsep}{20pt}
\begin{wrapfigure}[]{r}{0cm}
\begin{adjustbox}{trim=0cm 0cm 0cm 0cm} 
\begin{tikzpicture}[scale=1.25]

\begin{scope}[xshift=0cm] 
{

\draw[gray,very thin, step=0.5cm, opacity=0.5] (0,0) grid (2.5,3.5);

\fill[blue!50!white, opacity=0.5] (1,1.5)--++(0,1)--++(0.5,0)--++(0,-0.5)--++(0.5,0)--++(0,-0.5)--++(0.5,0)--++(0,-0.5)--++(-0.5,0)--++(0,-0.5)--++(-0.5,0)--++(0,-0.5)--++(-0.5,0)--++(0,1)--++(0.5,0)--++(0,0.5)--++(-0.5,0);

\draw[blue, line width=0.5mm] (0.5,3.5)--++(0,-0.5)--++(0.5,0)--++(0,0.5);

\draw[blue, line width=0.5mm] (0,2.5)--++(0.5,0)--++(0,-1);

\draw[blue, line width=0.5mm] (1,1.5)--++(0,1)--++(0.5,0)--++(0,-0.5)--++(0.5,0)--++(0,-0.5);

\draw[orange, line width=0.5mm] (2,1.5)--++(0.5,0)--++(0,-0.5)--++(-0.5,0)--++(0,-0.5)--++(-0.5,0)--++(0,-0.5)--++(-0.5,0)--++(0,0.5);

{

\node at (1.25,2.25) [scale=1]
{\small{$X$}};
\node at (0.75,3.25) [scale=1]
{\small{$L$}};

\node at (1.75,0.75) [scale=1]
{\small{$T_2$}};

\node[below] at (1.25,-0.1) [scale=0.75] 
{\small{\begin{tabular}{c}  Fig. 4.45(a). An illustration of \\ Case 1 with $T_2$  southeastern.  \end{tabular}}};;

}

}

\end{scope}

\begin{scope}[xshift=3.25cm] 
{

\draw[gray,very thin, step=0.5cm, opacity=0.5] (0,0) grid (3,4.5);

\fill[blue!50!white, opacity=0.5](1,2.5)--++(0,1)--++(0.5,0)--++(0,-0.5)--++(0.5,0)--++(0,-0.5)--++(0.5,0)--++(0,-0.5)--++(-0.5,0)--++(0,-0.5)--++(0.5,0)--++(0,-0.5)--++(0.5,0)--++(0,-1)--++(-0.5,0)--++(0,0.5)--++(-0.5,0)--++(0,0.5)--++(-0.5,0)--++(0,0.5)--++(-0.5,0)--++(0,0.5)--++(0.5,0)--++(0,0.5);

\draw[blue, line width=0.5mm] 
(0.5,4.5)--++(0,-0.5)--++(0.5,0)--++(0,0.5);

\draw[blue, line width=0.5mm] 
(0,3.5)--++(0.5,0)--++(0,-1);

\draw[blue, line width=0.5mm] 
(1,2.5)--++(0,1)--++(0.5,0)--++(0,-0.5)--++(0.5,0)--++(0,-0.5)--++(0.5,0)--++(0,-0.5)--++(-0.5,0)--++(0,-0.5)--++(0.5,0)--++(0,-0.5)--++(0.5,0)--++(0,0.5);

\draw[orange, line width=0.5mm](1.5,2)--++(-0.5,0)--++(0,-0.5)--++(0.5,0)--++(0,-0.5)--++(0.5,0)--++(0,-0.5)--++(0.5,0)--++(0,-0.5)--++(0.5,0)--++(0,0.5);

\node at (1.25,3.25) [scale=1]
{\small{$X$}};
\node at (0.75,4.25) [scale=1]
{\small{$L$}};
\node at (1.75,1.25) [scale=1]
{\small{$T_2$}};

\node[below] at (1.5,-0.1) [scale=0.75] 
{\small{\begin{tabular}{c}  Fig. 4.45(b). An illustration of \\ Case 1 with $T_2$  southwestern.  \end{tabular}}};;

}

\end{scope}

\end{tikzpicture}
\end{adjustbox}
\end{wrapfigure}

\noindent  \textit{Proof.} For definiteness, assume that $C$ is a small northern cookie, with $L=R(k'-1,l'+1)$. Let $F=G\langle N[P(X,Y)] \rangle$ be the looping $H$ path following $L$, as in Lemma 4.7, with $X=R(k',l'-1)$. Let $\overrightarrow{K}$ and $e_W, e_1, ... e_s$ be as in Lemma 4.7 as well. By Lemma 4.7, $F$ has a turn $T_1$. By proof of Lemma 4.7, either $T_1$ is a northeastern turn with $X$ as its northern leaf, or it is not. 

\null 

\noindent  \textit{CASE 1: $T_1$ is a northeastern turn with $X$ as its northern leaf.} Then $e_1$ is southern. By proof of Lemma 4.8, $T_2$ is either southeastern or southwestern. In either case, we note that $\text{Sector}(T_2)$ is south of the stack of $A_0$'s. By Lemma 4.8, $T_2$ is admissible. End of Case 1. 

\endgroup 

\null

\noindent  \textit{CASE 2: $T_1$ is not a northeastern turn with $X$ as its northern leaf.} By proof of Lemma 4.7, $T_1$ is a northeastern turn or $T_1$ is a southeastern turn.

\null 

\begingroup
\setlength{\intextsep}{0pt}
\setlength{\columnsep}{20pt}
\begin{wrapfigure}[]{l}{0cm}
\begin{adjustbox}{trim=0cm 0.75cm 0cm 0.5cm} 
\begin{tikzpicture}[scale=1.25]

\begin{scope}[xshift=0cm] 
{

\draw[gray,very thin, step=0.5cm, opacity=0.5] (0,0) grid (3.5,2.5);

\fill[blue!50!white, opacity=0.5] (1,1)--++(0.5,0)--++(0,0.5)--++(-0.5,0);
\fill[blue!50!white, opacity=0.5] (2,1)--++(0.5,0)--++(0,0.5)--++(-0.5,0);
\fill[blue!50!white, opacity=0.5] (1,0.5)--++(2,0)--++(0,0.5)--++(-2,0);
\fill[blue!50!white, opacity=0.5] (2.5,0)--++(1,0)--++(0,0.5)--++(-1,0);

\draw[blue, line width=0.5mm] (0.5,2.5)--++(0,-0.5)--++(0.5,0)--++(0,0.5);

\draw[blue, line width=0.5mm] (0,1.5)--++(0.5,0)--++(0,-1);

\draw[blue, line width=0.5mm] (1,0.5)--++(0,1)--++(0.5,0)--++(0,-0.5)--++(0.5,0);

\draw[orange, line width=0.5mm] (2,1)--++(0,0.5)--++(0.5,0)--++(0,-0.5)--++(0.5,0)--++(0,-0.5)--++(0.5,0)--++(0,-0.5)--++(-0.5,0);

{
\node at (1.25,1.25) [scale=1]
{\small{$X$}};
\node at (0.75,2.25) [scale=1]
{\small{$L$}};
\node at (2.75,0.75) [scale=1]
{\small{$T_1$}};

\node[below] at (1.75,-0.1) [scale=0.75] 
{\small{\begin{tabular}{c}  Fig. 4.46(a). An illustration of Case 2.1.  \end{tabular}}};;

}

}

\end{scope}

\begin{scope}[xshift=4.25cm] 
{

\draw[gray,very thin, step=0.5cm, opacity=0.5] (0,0) grid (3,3);

\fill[blue!50!white, opacity=0.5] (1,1)--++(0,0.5)--++(0.5,0)--++(0,0.5)--++(0.5,0)--++(0,0.5)--++(1,0)--++(0,-0.5)--++(-0.5,0)--++(0,-0.5)--++(-0.5,0)--++(0,-0.5)--++(-0.5,0)--++(0,-0.5)--++(-0.5,0)--++(0,0.5);

\draw[blue, line width=0.5mm] (0.5,2.5)--++(0,-0.5)--++(0.5,0)--++(0,0.5);

\draw[blue, line width=0.5mm] (0,1.5)--++(0.5,0)--++(0,-1);

\draw[blue, line width=0.5mm] (1,1)--++(0,0.5)--++(0.5,0)--++(0,0.5)--++(0.5,0)--++(0,0.5)--++(-0.5,0);

\draw[orange, line width=0.5mm] (2.5,2.5)--++(0.5,0)--++(0,-0.5)--++(-0.5,0)--++(0,-0.5)--++(-0.5,0)--++(0,-0.5)--++(-0.5,0)--++(0,-0.5)--++(-0.5,0)--++(0,0.5);

\draw[red, dashed, line width=0.5mm, opacity=1] (0,0.5)--++(2.5,2.5);

\node[left] at (0,1.5) [scale=1]{\tiny{$0$}};
\node[above] at (1,3) [scale=1] {\tiny{0}};

{
\node at (1.25,1.25) [scale=1]
{\small{$X$}};
\node at (0.75,2.25) [scale=1]
{\small{$L$}};
\node at (2.25,1.75) [scale=1]
{\small{$T_1$}};

\node[below] at (1.5,-0.1) [scale=0.75] 
{\small{\begin{tabular}{c}  Fig. 4.46(b). An illustration of  \\ Case 2.2.  The line $x=y$ in red. \end{tabular}}};;

}

}

\end{scope}

\end{tikzpicture}
\end{adjustbox}
\end{wrapfigure}

\noindent  \textit{CASE 2.1: $T_1$ is a northeastern turn.} It follows from the proof of Lemma 4.7 that $\text{Sector}(T_1)$ is east of the stack of $A_0$'s, that both leaves of $T_1$ are in $F$, and that neither leaf of $T_1$ is an end-box of $P(X,Y)$. End of Case 2.1.

\null

\noindent  \textit{CASE 2.2: $T_1$ is a southeastern turn.}  It follows from the proof of Lemma 4.7 that $\text{Sector}(T_1)$ is southeast of the stack of $A_0$'s. More precisely, we can check that $\text{Sector}(T_1)$ is below the line $y = x + (l' - k')$, and the j-stack of $A_0$'s is above the line $y = x + (l' - k')$. By proof of Lemma 4.7, we have that both leaves of $T_1$ are in $F$, and that no leaf of $T_1$ is an end-box of $P(X,Y)$. End of Case 2.2. End of Case 2. $\square$

\null 

\noindent Now we are ready to give a proof of Lemma 3.13.

\null 

\noindent \textbf{Lemma 3.13.} Let $G$ be an $m \times n$ grid graph, and let $H$ be a Hamiltonian cycle of $G$. Let $C$ be a small cookie of $G$.  Assume that $G$ has only one large cookie, and that there is a $j$-stack of $A_0$ starting at the $A_0$-type containing $C$. Let $L$ be the leaf in the top ($j^{\text{th}}$) $A_0$ of the stack, and assume that $L$ is followed by an $A_1$-type. Let $X$ and $Y$ be the boxes adjacent to the middle-box of the $A_1$-type that are not its $H$-neighbours. If $P(X,Y)$ has no switchable boxes, then either:

(i)  there is a cascade of length at most $\min(m,n)$, which avoids the stack of $A_0$'s, and after 

\hspace{0.4 cm} which $P(X,Y)$ gains a switchable box, or 

(ii) there is a cascade of length at most $\min(m,n)+1$, that collects $L$ and avoids the stack of $A_0$'s.

\null

\noindent \textit{Proof.} Suppose that $P(X,Y)$ has no switchable boxes. Then $P(X,Y)$ is contained in a looping fat path $F=G\langle N[P(X,Y)] \rangle$. By Lemma 4.13, $F$ has an admissible turn $T$ such that $\text{Sector}(T)$ and the $j$-stack of $A_0$'s are disjoint. Then, by Corollary 4.9(a), $d(T) \geq 3$. By Proposition 4.11, $T$ has a weakening $\mu_1, ..., \mu_s$. By Observation 4.12, $\mu_1, ..., \mu_s$ is contained in $\text{Sector}(T)$, and thus it avoids the j-stack of $A_0$'s.

\null 

\begingroup
\setlength{\intextsep}{10pt}
\setlength{\columnsep}{20pt}

\begin{adjustbox}{trim=0cm 0cm 0cm 0cm}
\begin{tikzpicture}[scale=1.45]

\begin{scope}[xshift=0cm, yshift=0cm]
\draw[gray,very thin, step=0.5cm, opacity=0.5] (0,0) grid (1.5,1.5);

\draw [->,black, very thick] (1.7,0.75)--(2.3,0.75);
\node[above] at  (2,0.75) [scale=0.8]{\small{$\mu_s$}};

{

\node[right] at (1.5,1.5) [scale=1]{\tiny{+1}};
\node[right] at (1.5,1) [scale=1]{\tiny{b}};
\node[right] at (1.5,0.5) [scale=1]{\tiny{-1}};
\node[right] at (1.5,0) [scale=1]{\tiny{-2}};

\node[above] at (0, 1.5) [scale=1]{\tiny{a}};
\node[above] at (0.5, 1.5) [scale=1]{\tiny{+1}};
\node[above] at (1, 1.5) [scale=1]{\tiny{+2}};
}


\fill[green!40!white, opacity=0.5] (0,1) rectangle (0.5,1.5);
\fill[blue!40!white, opacity=0.5] (0,0.5) rectangle (1,1);

\draw[blue, line width=0.5mm] (0,1)--++(0,0.5);
\draw[blue, line width=0.5mm] (0.5,1.5)--++(0,-0.5)--++(0.5,0)--++(0,-0.5)--++(0.5,0);
\draw[blue, line width=0.5mm] (0,0.5)--++(0.5,0)--++(0,-0.5);


{
\draw[black, line width=0.15mm] (0.7,0.45)--++(0,0.1);
\draw[black, line width=0.15mm] (0.75,0.45)--++(0,0.1);
\draw[black, line width=0.15mm] (0.8,0.45)--++(0,0.1);

\draw[black, line width=0.15mm] (0.45,0.7)--++(0.1,0);
\draw[black, line width=0.15mm] (0.45,0.75)--++(0.1,0);
\draw[black, line width=0.15mm] (0.45,0.8)--++(0.1,0);

\draw[black, line width=0.15mm] (-0.05,0.7)--++(0.1,0);
\draw[black, line width=0.15mm] (-0.05,0.75)--++(0.1,0);
\draw[black, line width=0.15mm] (-0.05,0.8)--++(0.1,0);
}

\node at (0.25,1.25) [scale=0.8]{\small{$L_{\textrm{N}}$}};

\node[below] at (0.75,0) [scale=0.75] 
{\small{\begin{tabular}{c} Fig. 4.47(a). $L_N{+}(0,{-}1)$ \\ is not an end-box. \end{tabular}}};;

\end{scope}

\begin{scope}[xshift=2.5cm, yshift=0cm]
\draw[gray,very thin, step=0.5cm, opacity=0.5] (0,0) grid (1.5,1.5);

{

\node[right] at (1.5,1.5) [scale=1]{\tiny{+1}};
\node[right] at (1.5,1) [scale=1]{\tiny{b}};
\node[right] at (1.5,0.5) [scale=1]{\tiny{-1}};
\node[right] at (1.5,0) [scale=1]{\tiny{-2}};

\node[above] at (0, 1.5) [scale=1]{\tiny{a}};
\node[above] at (0.5, 1.5) [scale=1]{\tiny{+1}};
\node[above] at (1, 1.5) [scale=1]{\tiny{+2}};
}


\fill[green!40!white, opacity=0.5] (0,1) rectangle (0.5,1.5);
\fill[blue!40!white, opacity=0.5] (0,0.5) rectangle (1,1);

\draw[blue, line width=0.5mm] (0,1)--++(0.5,0);
\draw[blue, line width=0.5mm] (0,1.5)--++(0.5,0);
\draw[blue, line width=0.5mm] (0.5,1)--++(0.5,0)--++(0,-0.5)--++(0.5,0);
\draw[blue, line width=0.5mm] (0,0.5)--++(0.5,0)--++(0,-0.5);

{
\draw[black, line width=0.15mm] (0.7,0.45)--++(0,0.1);
\draw[black, line width=0.15mm] (0.75,0.45)--++(0,0.1);
\draw[black, line width=0.15mm] (0.8,0.45)--++(0,0.1);

\draw[black, line width=0.15mm] (0.45,0.7)--++(0.1,0);
\draw[black, line width=0.15mm] (0.45,0.75)--++(0.1,0);
\draw[black, line width=0.15mm] (0.45,0.8)--++(0.1,0);

\draw[black, line width=0.15mm] (-0.05,0.7)--++(0.1,0);
\draw[black, line width=0.15mm] (-0.05,0.75)--++(0.1,0);
\draw[black, line width=0.15mm] (-0.05,0.8)--++(0.1,0);

}

\node at (0.25,1.25) [scale=0.8]{\small{$L_{\textrm{N}}$}};

\node[below] at (0.75,0) [scale=0.75] 
{\small{\begin{tabular}{c} Fig. 4.47(b). After $\mu_s$. \end{tabular}}};;

\end{scope}

\begin{scope}[xshift=5cm, yshift=0cm]
\draw[gray,very thin, step=0.5cm, opacity=0.5] (0,0) grid (2,2);

\draw [->,black, very thick] (2.2,0.75)--(2.8,0.75);
\node[above] at  (2.5,0.75) [scale=0.8]{\small{$\mu_s$}};

\node[right] at (2,2) [scale=1]{\tiny{+1}};
\node[right] at (2,1.5) [scale=1]{\tiny{b}};
\node[right] at (2,1) [scale=1]{\tiny{-1}};
\node[right] at (2,0.5) [scale=1]{\tiny{-2}};

\node[above] at (0.5, 2) [scale=1]{\tiny{-1}};
\node[above] at (1, 2) [scale=1]{\tiny{a}};
\node[above] at (1.5, 2) [scale=1]{\tiny{+1}};

\fill[green!40!white, opacity=0.5] (1,1.5) rectangle (1.5,2.0);
\fill[blue!40!white, opacity=0.5] (1,1.0) rectangle (2,1.5);

\draw[blue, line width=0.5mm] (1,1.5)--++(0,0.5);
\draw[blue, line width=0.5mm] (1.5,2.0)--++(0,-0.5)--++(0.5,0)--++(0,-0.5);
\draw[blue, line width=0.5mm] (1,1)--++(0.5,0)--++(0,-0.5);
\draw[blue, line width=0.5mm] (1,1)--++(0,0.5);

\draw[blue, line width=0.5mm] (0,1)--++(0.5,0)--++(0,1);

\draw[blue, line width=0.5mm] (0.5,0)--++(0,0.5)--++(0.5,0)--++(0,-0.5);

\draw[black, line width=0.15mm] (1.7,0.95)--++(0,0.1);
\draw[black, line width=0.15mm] (1.75,0.95)--++(0,0.1);
\draw[black, line width=0.15mm] (1.8,0.95)--++(0,0.1);

\draw[black, line width=0.15mm] (1.45,1.2)--++(0.1,0);
\draw[black, line width=0.15mm] (1.45,1.25)--++(0.1,0);
\draw[black, line width=0.15mm] (1.45,1.3)--++(0.1,0);

\node at (1.25,1.75) [scale=0.8]{\small{$L_{\textrm{N}}$}};
\node at (1.25,1.25) [scale=0.8]{\small{$Y$}};
\node at (0.25,1.25) [scale=0.8]{\small{$X$}};
\node at (0.75,0.25) [scale=0.8]{\small{$L$}};

\node[below] at (1,0) [scale=0.75] 
{\small{\begin{tabular}{c} Fig. 4.48(a). $L_N{+}(0,{-}1)$ \\ is an end-box. \end{tabular}}};

\end{scope}

\begin{scope}[xshift=8cm, yshift=0cm]
\draw[gray,very thin, step=0.5cm, opacity=0.5] (0,0) grid (2,2);

\node[right] at (2,2) [scale=1]{\tiny{+1}};
\node[right] at (2,1.5) [scale=1]{\tiny{b}};
\node[right] at (2,1) [scale=1]{\tiny{-1}};
\node[right] at (2,0.5) [scale=1]{\tiny{-2}};

\node[above] at (0.5, 2) [scale=1]{\tiny{-1}};
\node[above] at (1, 2) [scale=1]{\tiny{a}};
\node[above] at (1.5, 2) [scale=1]{\tiny{+1}};

\fill[green!40!white, opacity=0.5] (1,1.5) rectangle (1.5,2.0);
\fill[blue!40!white, opacity=0.5] (1,1.0) rectangle (2,1.5);

\draw[blue, line width=0.5mm] (1,1.5)--++(0.5,0);
\draw[blue, line width=0.5mm] (1,2)--++(0.5,0);

\draw[blue, line width=0.5mm] (1.5,1.5)--++(0.5,0)--++(0,-0.5);
\draw[blue, line width=0.5mm] (1,1)--++(0.5,0)--++(0,-0.5);
\draw[blue, line width=0.5mm] (1,1)--++(0,0.5);

\draw[blue, line width=0.5mm] (0,1)--++(0.5,0)--++(0,1);

\draw[blue, line width=0.5mm] (0.5,0)--++(0,0.5)--++(0.5,0)--++(0,-0.5);

\draw[black, line width=0.15mm] (1.7,0.95)--++(0,0.1);
\draw[black, line width=0.15mm] (1.75,0.95)--++(0,0.1);
\draw[black, line width=0.15mm] (1.8,0.95)--++(0,0.1);

\draw[black, line width=0.15mm] (1.45,1.2)--++(0.1,0);
\draw[black, line width=0.15mm] (1.45,1.25)--++(0.1,0);
\draw[black, line width=0.15mm] (1.45,1.3)--++(0.1,0);

\node at (1.25,1.75) [scale=0.8]{\small{$L_{\textrm{N}}$}};
\node at (1.25,1.25) [scale=0.8]{\small{$Y$}};
\node at (0.25,1.25) [scale=0.8]{\small{$X$}};
\node at (0.75,0.25) [scale=0.8]{\small{$L$}};

\node[below] at (1,0) [scale=0.75] 
{\small{\begin{tabular}{c} Fig. 4.48(b). After $\mu_s$. \end{tabular}}};

\end{scope}

\end{tikzpicture}
\end{adjustbox}

\null

\noindent For definiteness, assume that $T$ is northeastern with northern leaf $L_N=R(a,b)$ and that $\mu_s$ is the move $L_N \mapsto L_N'$ or the move $L_N' \mapsto L_N$. By Corollary 4.9(b), $L_N \in N[P] \setminus P$ and $L_N+(0,-1) \in P$. Note that we must have $S_{\rightarrow}(a,b-1;a+1,b-2) \in H$, $e(a,a+1;b) \notin H$ and $e(a+1;b-1,b)\notin H$. Now, either $L_N+(0,-1)$ is an end-box of $P$ or it is not. If the latter then, $e(a;b-1,b \notin H)$. Then, after $\mu_s$, $L_N+(0,-1) \in P(X,Y)$ is switchable. The fact that $s \leq \min(m,n)$, follows immediately from Proposition 4.11. Thus, in this case, (i) holds. 

Suppose then, that $L_N+(0,-1)$ is an end-box of $P$, say $Y$. This implies that $e(a;b-1,b)\in H$, that $F$ is northern,  and that $L=L_N+(-1,-3)$. 
Then $\mu_s$ can be followed by $Y+(-1,0) \mapsto Y$, $L+(0,1) \mapsto L$, which collects $L$. To check the length of the cascade, first we note that $b+1 \geq 4$, and that the eastern leaf of $T$ has x-coordinate at least $5$. It follows that $|\mathcal{T}(T)|\leq \min(m,n)-4$. It follows from the proof of Lemma 4.10 that $T$ has a weakening of length at most $\min(m,n)-4+3$. Since we need two additional moves after $\mu_s$ to collect $L$, the length of the cascade is at most $\min(m,n)-4+3+2=\min(m,n)+1$ moves. Thus, in this case, (ii) holds. $\square$

\endgroup

\null 

{



}

\null

\section*{Acknowledgments}

\noindent This paper is based on part of the author’s PhD dissertation at York University, 
written under the supervision of Professor Neal Madras. 

\null

\noindent The author is grateful to his advisor, Neal Madras, for his patience and guidance through several rewrites, and for his thoughtful and insightful feedback on each of them. The author also thanks Nathan Clisby for creating an interactive online tool for visualizing and experimenting with Hamiltonian paths and cycles\cite{clisbyHamiltonian}.

\null 

\noindent This research was supported in part by a Discovery Grant from NSERC Canada to my advisor Neal Madras.

\null 

\section*{Declarations}

\noindent\textbf{Funding} \\
The author received no external funding for this work.\\[6pt]

\noindent\textbf{Competing interests} \\
The author has no relevant financial or non-financial interests to disclose.\\[6pt]

\noindent\textbf{Availability of data and materials} \\
Data sharing not applicable to this article as no datasets were generated or analyzed during the current study.\\[6pt]

\noindent\textbf{Authors' contributions} \\
The author is the sole contributor to all aspects of this work.

\null

\bibliography{references}

\null 
\null 
\null

\null

\null 
\null 
\null 
\null

\end{document}